\PassOptionsToClass {plain,biber}{nowfnt}
\documentclass[CIT,biber]{nowfnt}

\usepackage{amsmath}
\usepackage{amsthm,amsfonts}

\usepackage{xspace,exscale,relsize}
\usepackage{amssymb}
\usepackage{subcaption}
\usepackage{booktabs}
\usepackage{multirow,soul}

\usepackage{algorithm,algorithmic}

\usepackage{macro}
\graphicspath{{figures/}}
\usetikzlibrary{arrows}
\usepackage{pgfplots}
\usepgfplotslibrary{colormaps,fillbetween}

\newcommand{\stepa}[1]{\overset{\rm (a)}{#1}}
\newcommand{\stepb}[1]{\overset{\rm (b)}{#1}}

\issuesetup
{
}

\maintitleauthorlist
{
    Yihong Wu \\
    Department of Statistics and Data Science, Yale University\\
    \url{yihong.wu@yale.edu} \\    
    \and
    Pengkun Yang \\
    Center for Statistical Science, Tsinghua University \\
    \url{yangpengkun@tsinghua.edu.cn}\\
}

\author[1]{Wu,Yihong}
\author[2]{Yang,Pengkun}

\affil[1]{Department of Statistics and Data Science, Yale University, New Haven, CT; \url{yihong.wu@yale.edu}}
\affil[2]{Department of Electrical Engineering, Princeton University, Princeton, NJ; \url{pengkuny@princeton.edu}}

\title{Polynomial methods in statistical inference: theory and practice}

\begin{document}
\makeabstracttitle

\begin{abstract}
%
%
%

This survey provides an exposition of a suite of techniques based on the theory of polynomials, collectively referred to as polynomial methods, which have recently been applied to address several challenging problems in statistical inference successfully. Topics including polynomial approximation, polynomial interpolation and majorization, moment space and positive polynomials, orthogonal polynomials and Gaussian quadrature are discussed, with their major probabilistic and statistical applications in property estimation on large domains and learning mixture models. These techniques provide useful tools not only for the design of highly practical algorithms with provable optimality, but also for establishing the fundamental limits of the inference problems through the method of moment matching. The effectiveness of the polynomial method is demonstrated in concrete problems such as entropy and support size estimation, distinct elements problem, and learning Gaussian mixture models.

\end{abstract}




\chapter{Introduction}
\label{chap:intro}

Modern data-analytic applications frequently involve complex and high-dimensional statistical models. 
For example, applications such as natural language processing, genetics, and neuroscience deal with datasets naturally viewed as being sampled from \emph{probability distributions over a large domain}.
A number of real-world signal processing and machine learning tasks rest upon data-driven procedures for 
estimating \emph{distributional properties} (functionals of the data-generating distribution), including \emph{entropy} for understanding the neural coding \cite{spikes-book,amigo2004estimating,Aktulga07,GKB06,SKSB98,NBS04,knudson2013spike}, \emph{mutual information} for image registration in fMRI \cite{pluim2003mutual,kybic2004high,tsai1999analysis,tsai2004mutual} and learning graphical models \cite{CL68,jiao2016beyond}, etc. 
For these tasks, the key challenge is to accurately estimate the property even when the domain size far exceeds the sample size and the distribution itself is impossible to learn.


Another prominent example of complex statistical models deals with \emph{mixture models}, which are useful to model the effects of latent variables and form the basis of many clustering algorithms. 
The simplest mixture models is perhaps the Gaussian mixture model, introduced by Pearson in 1894 to model the presence of hidden subpopulations within an overall population. 
Despite the seemingly innocuous nature of the Gaussian mixture models, many difficult challenges arise, such as the vanishing Fisher information leading to nonparametric rates, the nonexistence of maximum likelihood estimator in location-scale mixtures, etc. 
For this reason, it proves to be a fertile ground for innovations in statistical methodologies, including the method of moments \cite{Pearson1894}, the 
Expectation-Maximization (EM) algorithm \cite{DLR1977}, the Generalized Method of Moments \cite{Hansen1982}, etc.
Despite the vast literature and recent breakthroughs, many problems as basic as optimal estimation rates remain open in finite mixture models.

Recently, several challenging problems in property estimation and mixture models have been successfully resolved using methods based on the theory of polynomials, in particular, polynomial approximation, interpolation, as well as moments and positive polynomials. 
They provide useful tools not only for the design of algorithms that are both statistically optimal and computationally efficient, but also in establishing the fundamental limits of the inference problems. 
This survey aims to provide an exposition of these techniques, which are collectively referred to as the \emph{polynomial method}, 
as well as their application in statistical inference.

\section{Background on polynomial methods}


The theory of polynomials is a rich subject in mathematics of both algebraic and analytic flavor. It forms the foundation of and has diverse applications in many subjects including optimization, combinatorics, 
coding theory, control theory, digital signal processing, game theory, statistics and machine learning, etc, leading to many deep theoretical results and highly practical algorithms. 
In this survey, we mainly focus on 
polynomial approximation, interpolation, and positive polynomials that will be introduced below.

\paragraph{Polynomial approximation and interpolation.}
One of the most well-understood subjects in approximation theory, polynomial approximation 
aims at approximating a given complicated function, in either a local or global sense, using algebraic or trigonometric polynomials of a certain degree. 
For instance, the Taylor expansion characterizes the local behavior of a smooth function and provide the foundation for optimization techniques such as gradient descent and the Newton-Raphson method \cite{Mitchell1997} and kernel-based methods in statistical inference \cite{Tsybakov09,hardle2012wavelets}; trigonometric polynomials represent functions in the frequency domain through Fourier analysis, which are the theoretical underpinnings for digital signal processing and wireless transmission \cite{Oppenheim1999,TV2005}. 
A closely related topic is polynomial interpolation, which can be viewed as achieving zero approximation error on a discrete set of points.

In property estimation, the functional to be estimated can be highly nonsmooth and classical methods requires a large sample size in order to be accurate. 
In such settings, polynomial approximation and interpolation provide a useful primitive for constructing better estimates by first approximating the original functional by 
a polynomial and then estimate the polynomial approximant. 
Besides the approximation error which is the primary concern in approximation theory, other properties of the polynomial approximant such as the magnitude of its coefficients are also crucial for bounding the statistical error.

\paragraph{Moments and positive polynomials.}
The theory of moments plays a key role in the developments of analysis, probability, statistics, and optimization. We refer the readers to the classics \cite{ST1943,KS1953,KN1977} and the more recent monographs \cite{Lasserre2009,Schmudgen17} for a detailed treatment. In statistical inference, the method of moments was originally introduced by Pearson \cite{Pearson1894} for mixture models, which constructs estimates by solving polynomial equations. 
Due to its conceptual simplicity and flexibility, especially in models without the complete specification of the joint distribution of data, method of moments and its extensions have been widely applied in practice, for instance, to analyze economic and financial data \cite{Hall2005}. 
In probability and optimization literature, the classical moment problem refers to determining whether a probability distribution is determined by all of its moments. 
Solution to the moment problem requires understanding the moment space, which is the convex set formed by moments of probability distributions. 
The moment space satisfies many geometric properties (such as the Cauchy-Schwarz and H\"older inequalities) and a complete description can be phrased in terms of positive polynomials, which are further related to sums of squares and semidefinite programming.
Together with techniques based on polynomial interpolation, this structural information can be leveraged to 
design moment-based methods for learning mixture models that are statistically optimal, robust to model misspecification, and highly practical.


\section{Polynomial methods for designing estimators}
\label{sec:intro-estimator}
We will apply the above polynomial methods to the tasks of estimating distributional properties and learning mixture models 
with the goal of constructing estimators with good statistical performance.

\paragraph{Estimating distributional properties on large domains. }
Given samples drawn from an unknown distribution $P$ on a large domain, the goal is to estimate a specific property of that distribution, such as various information measures including the Shannon entropy, R\'enyi entropy, and the support size. This falls under the category of \emph{functional estimation} \cite{Rao2014}, where we are not interested in directly estimating the high-dimensional parameter (the data-generating distribution $P$) per se, but rather a function thereof. Estimating a distributional functional has been intensively studied in nonparametric statistics, including estimating a scalar function of a regression function or density such as linear functionals \cite{Stone80,DL91}, quadratic functionals \cite{laurent1996efficient,CL05}, $L_q$ norm \cite{LNS99}, etc. 

To estimate a functional, perhaps the most natural idea is the ``plug-in'' approach, namely, first estimate the parameter and then substitute into the function. As frequently observed in the functional estimation literature, the plug-in estimator can suffer from severe bias (see \cite{Efron82,Berkson80} and the references therein). Indeed, although the plug-in estimate is typically asymptotically efficient and minimax (cf., \eg, \cite[Sections 8.7 and 8.9]{VdV00}) for fixed domain size, it can be highly suboptimal in high dimensions, where, due to the large alphabet and resource constraints, we are constantly contending with the difficulty of \emph{undersampling} in applications such as
\begin{itemize}
    \item Natural language processing: The vast vocabulary size of natural languages, compounded by the frequent use of bigrams and trigrams in practice \cite{MS1999}, leads to an effective alphabet size far exceeding the sample size.  A well-known example from corpus linguistics is that about half of the words in the Shakespearean canon only appeared once \cite{ET76};
    \item Neuroscience: in 		
		analyzing neural spike trains, natural stimuli generate neural responses of high timing precision resulting in a massive space of meaningful responses \cite{Berry13051997,mainen1995reliability,SLSKB97};

    \item Network traffic analysis: many customers or website users are only seen a small number of times \cite{benevenuto2009characterizing}.
\end{itemize}

Statistical inference on large domains
has a rich history in information theory, statistics and computer science, with early contributions dating back to Fisher, Good and Turing,  Efron and Thisted, etc \cite{FCW43,Good1953,ET76,TE87} 
and recent renewed interests on compression, prediction, classification and estimation on large alphabets \cite{OSZ04,BS09,KWTV13,WVK11,VV13}; however, none of the aforementioned results allows a general understanding of the fundamental limits of estimating information quantities of large distributions. 
While there exists a vast literature on information-theoretic approaches to the statistical inference of high-dimensional parameters 
\cite{LeCam73,IKbook,pinsker.minimax,Birge83,Yu97,YB99}, a systematic theory 
for estimating their low-dimensional functionals remains
severely under-developed, especially in the \emph{sublinear regime} where the sample size is far less than the domain size so that the underlying distribution is impossible to learn but certain low-dimensional features can nevertheless be estimated accurately.

In this survey, we will investigate a few prototypical problems in estimating distributional properties such as the Shannon entropy and the support size.
These properties can be easily estimated if the sample size far exceeds the support size of the underlying distribution, but how can it be done if the observations are relatively scarce, especially in the \emph{sublinear regime} where the sample size is far less than the domain size? It turns out the theory of polynomial approximation provides a principled approach to construct an optimal estimator.
To illustrate this program let us consider the problem of estimating a function $f(p)$ based on $n$ independent observations drawn from Bernoulli distribution with mean $p$, or equivalently, the sufficient statistic $N\sim \Binom(n,p)$. This simple setting forms the basis of designing estimators for distributional properties in Chapters~\ref{chap:approx} -- \ref{chap:unseen}.
Given any estimator $\hat f(N)$, its mean is given by
\[
\Expect[\hat f(N)]=\sum_{j=0}^n f(j)\binom{n}{j}p^j (1-p)^{n-j},
\]
which is a degree-$n$ polynomial in $p$.
Consequently, unless the function $f$ is a polynomial, there exists no unbiased estimator for $f(p)$.
Conversely, given any degree-$n$ polynomial $\tilde f$, we can always construct an unbiased estimator for $\tilde f(p)$ by combining the unbiased estimator of each monomial (see, \eg, \prettyref{eq:unbiased-iid} in \prettyref{sec:functional}). 
These observations suggest that, for the purpose of reducing the bias, we should first find a polynomial $\tilde f$ of degree at most $n$ such that the approximation error $|f(p)-\tilde f(p)|$ is small for every possible values of $p$, and then construct an unbiased estimator $\hat f(N)$ for $\tilde f(p)$.
Fixing $L \leq n$, the best degree-$L$ polynomial $\tilde f$  that minimizes the worst-case approximation error can be found by solving the following optimization problem:
\begin{equation}
\label{eq:poly-approx-general}
\inf_{\lambda_0,\dots,\lambda_L}\sup_p \abs{f(p)-\sum_{i=0}^n\lambda_i p^i};
\end{equation}
this is known as the \emph{best uniform polynomial approximation} problem which will be discussed at length in \prettyref{sec:uniform}. 
Although the approximation error decays with the degree, typically we cannot choose it to be as large as $n$ since the estimation error of monomials grows rapidly with the degree.
Therefore, the degree $L$ must be chosen appropriately (often logarithmic in the sample size $n$) so as to balance the approximation error and the estimation error (the bias-variance tradeoff). 
This method was pioneered by Lepski, Nemirovski, and Spokoiny \cite{LNS99} for nonparametric regression and further developed in Cai and Low \cite{CL11} for the Gaussian sequence model.
We will elaborate on the high-level ideas in \prettyref{chap:approx} and illustrate the effectiveness of this approach in Chapters~\ref{chap:entropy} and \ref{chap:unseen} for specific problems.

\paragraph{Learning Gaussian mixtures. }
Sampling from a mixture model can be viewed as being a two-step process: 
first draw a latent parameter $\theta\sim \nu$; then draw an observation $X\sim P_\theta$. 
The marginal distribution of each sample is 
\begin{equation}
\label{eq:gm-mixture}
\pi_\nu=\int P_\theta\diff\nu(\theta).
\end{equation}
We refer to $\nu$ as the mixing distribution and $\pi_\nu$ as the mixture distribution.
A finite mixture model has a discrete mixing distribution of finite support and a mixture distribution of the form $\sum_i w_i P_{\theta_i}$.
The key question in mixture model is the following: 
If we are only given unlabeled data from the mixture model, 
can we reconstruct the parameters in each component accurately and efficiently? 
Furthermore, in the regime where it is impossible to learn the labels with small misclassification rate, is it still possible to learn the mixing distribution and the mixture distribution accurately?

In the special case that each $P_\theta$ is a Gaussian distribution, this is the problem of learning Gaussian mixtures, a classical problem in statistics dating back to the work of Pearson \cite{Pearson1894}.
In addition, methods for learning Gaussian mixtures are widely used as part of the core machine learning toolkit, such as the popular {\sf scikit-learn} package in {\sf Python} \cite{scikit-learn}, Google's {\sf Tensorflow} \cite{abadi2016tensorflow}, and Spark's {\sf MLlib} \cite{meng2016mllib}; however, few provable guarantees are available. It is only recently proved in \cite{KMV2010,MV2010} that a mixture of constant number of components can be learned in polynomial time using a polynomial number of observations. 
The optimal rate for learning finite Gaussian location mixtures is recently determined in \cite{HK2015,WY18,DWYZ20} and for location-scale mixture only for the special case of two components \cite{HP15}. Is there a systematic way to obtain the sharp error rates and how to efficiently and optimally learn a Gaussian mixture? We will investigate the moment methods for the optimal estimation of Gaussian mixtures, where we learn a discrete mixing distribution by learning its moments. 
The key observation is that as opposed to the vanilla method of solving moment equations, 
the moment estimates should be first denoised based on the geometry of the moment space, and the denoising step can be efficiently carried out through convex optimization (semidefinite programming). The learned moments can be then converted to a discrete distribution by the efficient algorithm of Gaussian quadrature. 
This approach will be presented in Chapters~\ref{chap:framework}--\ref{chap:gm}.

\section{Polynomial methods for determining theoretical limits}
\label{sec:intro-lb}
Another focus of this survey is to investigate the fundamental limits of statistical inference, that is, the optimal estimation error among all estimators regardless of computational costs. While the use of polynomial methods on the constructive side is admittedly natural, the fact that it also arises in the optimal lower bound is perhaps surprising. 

To give a precise definition of the fundamental limits, we begin with an account of the general framework for statistical inference. We assume that the sample $X_1,\ldots,X_n$ are independently generated from an unknown distribution $P$ that belongs to a collection of distributions $\calP$. The goal is to estimate a certain property $T(P)$ of the distribution $P$.

In this survey we consider the following two types of problems:
\begin{itemize}
	\item \emph{Estimating distributional properties: }
	$T(P)$ is a functional of the unknown discrete distribution $P=(p_1,p_2,\ldots)$, such as the Shannon entropy 
	\begin{equation}
	H(P)=\sum_i p_i\log \frac{1}{p_i}
	\label{eq:entropy}
	\end{equation}
	 and the support size 
	\begin{equation}
	S(P)=\sum_i \indc{p_i>0}
	\label{eq:supportsize}
	\end{equation}
	
	\item \emph{Learning Gaussian mixtures: } 
	$P$ is a Gaussian mixture and $T(P)$ represents the parameters, including the mean, variance, and the mixing weights, of each Gaussian component. Equivalently, $T(P)$ can be viewed as the mixing distribution of the mixture model (see \prettyref{chap:framework}).  
\end{itemize}
Given a loss function $\ell(\hat T, T(P))$ that measures the accuracy of an estimator $\hat T$, the decision-theoretic fundamental limit is defined as the \emph{minimax risk}
\begin{equation}
    \label{eq:R-minimax}
    R^*_n\triangleq \inf_{\hat{T}}\sup_{P \in \calP}\Expect_P[\ell( \hat{T},T(P))],     
\end{equation}
where the infimum is taken over all estimators $\hat T$ measurable with respect to $X_1,\dots,X_n$ drawn independently from $P$. Examples of the loss function include the quadratic loss $\ell(x,y)=\Norm{x-y}_2^2$ and the zero-one loss $\ell(x,y)=\indc{\Norm{x-y}_2>\epsilon}$ for a desired accuracy $\epsilon$. For the zero-one loss, we also consider the \emph{sample complexity}: 
\begin{definition}
	\label{def:sample-complexity}
	For a desired accuracy $\epsilon$ and confidence $1-\delta$, the sample complexity is the minimal sample size $n$ such that there exists an estimator $\hat T$ based on $n$ independent and identically distributed (\iid) observations drawn from a distribution $P$ such that $\Prob[\ell( \hat{T},T(P)) < \epsilon ] \geq 1-\delta$ for any $P\in\calP$. 
\end{definition}
In this survey, our primary goal is to characterize the minimax risk \prettyref{eq:R-minimax} within universal constant factors, which is known as the \emph{minimax rate}; we will also consider the sample complexity in \prettyref{def:sample-complexity}. This task entails an upper bound achieved by certain estimators, preferably a computationally efficient one, and a matching minimax lower bound that applies to all estimators.

A general program for obtaining lower bounds is based on a reduction of estimation to testing (Le Cam's method); cf.~\prettyref{sec:d-mm}. If there are two distributions $P$ and $Q$ that cannot be reliably distinguished based on a given number of independent observations, while $T(P)$ and $T(Q)$ are different, then any estimate suffers a maximum risk at least proportional to the distance between $T(P)$ and $T(Q)$. Furthermore, sometimes one needs to consider a pair of randomized distributions in which case one needs to construct two distributions (priors) on the space of  distributions (also known as fuzzy hypothesis testing in \cite{Tsybakov09}). Here the polynomial method enters the scene again: statistical closeness between two distributions can be bounded by comparing their moments.
More precisely, the strategy is to choose two priors with matching moments up to a certain degree, which ensures the induced distributions of data are impossible to test. The minimax lower bound is then given by the maximal separation in the expected functional values subject to the moment matching condition. 
For example, it pertains to the optimal value of the following type of moment matching problem:
\begin{equation}
\label{eq:poly-dual}
\begin{aligned}
\sup \quad & \Expect_\nu[f(X)]-\Expect_{\nu'}[f(X)],\\
\textrm{s.t.}\quad & \Expect_\nu[X^j]=\Expect_{\nu'}[X^j],\quad j=0,\dots,L,\\
& \nu,\nu'\textrm{ are supported on }[a,b],
\end{aligned}
\end{equation}
where the supremum is over all pairs of distributions, and the function $f$, the degree $L$, and the interval $[a,b]$ are problem specific. 
We will discuss how to choose those parameters, construct a pair of least favorable priors from the optimal solution, and then derive the minimax lower bound in Chapters \ref{chap:entropy} and \ref{chap:unseen}. 
It turns out this optimization problem is the \emph{dual} problem of the best polynomial approximation that arises in the design of polynomial-based estimator in \prettyref{sec:intro-estimator}. 
In the introduction, let us first look into the relation to polynomial method.
Below we formally derive the duality, and we leave the discussion on strong duality and the correspondence between primal and dual solutions to \prettyref{sec:dual-best}.
By introducing the Lagrangian multipliers $\lambda_1,\dots,\lambda_L$, we optimize the Lagrangian function by
\begin{align*}
&\phantom{{}=}\sup_{\nu,\nu'}\Expect_\nu[f(X)]-\Expect_{\nu'}[f(X)]-\sum_{j=1}^L\lambda_i(\Expect_\nu[X^j]-\Expect_{\nu'}[X^j])\\
&=\sup_{\nu,\nu'}\Expect_\nu\qth{f(X)-\sum_{j=1}^L\lambda_iX^i}-\Expect_{\nu'}\qth{f(X)-\sum_{j=1}^L\lambda_iX^i}\\
&=\sup_{x\in[a,b]}\pth{f(x)-\sum_{j=1}^L\lambda_ix^i}-\min_{x\in[a,b]}\pth{f(x)-\sum_{j=1}^L\lambda_ix^i}.
\end{align*}
We can introduce another variable $\lambda_0$ that does not impact the optimal value and formulate the dual problem as 
\begin{align}
&\phantom{{}=}\inf_{\lambda_0,\dots,\lambda_L}\sup_{x\in[a,b]}\pth{f(x)-\sum_{j=0}^L\lambda_ix^i}-\min_{x\in[a,b]}\pth{f(x)-\sum_{j=0}^L\lambda_ix^i}\nonumber\\
&=2\inf_{\lambda_0,\dots,\lambda_L}\sup_{x\in[a,b]}\abs{f(x)-\sum_{j=0}^L\lambda_ix^i}.\label{eq:poly-primal}
\end{align}
This last formulation is precisely the best polynomial approximation problem \prettyref{eq:poly-approx-general}.
For this reason, estimators constructed using the method of polynomial approximation frequently comes with a matching lower bound that certifies their statistical optimality. 
The connection is precisely the duality between polynomial approximation and moment matching.

The method of moment matching can be similarly carried out for learning mixture models. Typically, there is a minimal number $L$ of moments that identifies a finite mixture model, which depends on the order (the number of components) of the mixture model. A statistical lower bound can then be obtained by constructing a pair of distributions with matching $L-1$ moments. This naturally matches the performance of the ``most economical'' moment-based estimators that learns the mixture distribution using the minimal number of moments.
We will discuss this approach in \prettyref{chap:gm}.

\section{Organization}
\begin{figure}[ht]%
\centering
\tikzset{
    implies/.style={double,double equal sign distance,-implies},
    equiv/.style={implies-implies,double equal sign distance},
}
\usetikzlibrary{shapes.misc,automata,positioning,arrows}

\begin{tikzpicture}[node distance=2cm,scale = 0.8]
	\node at (6, 6.6) (a1) [align=center, rectangle, draw,thick]{polynomial approximation};
	\node at (0, 3) (a2) [align=center, rectangle, draw,thick]{polynomial \\ interpolation};
	\node at (6, -0.5) (a3) [align=center, rectangle, draw,thick]{moment space \&\\ positive polynomials};

	\node at (6, 4.5) (b1) [align=center, rounded rectangle, draw, thick]{\large Property Estimation};
	\node at (6, 3) (b2) [align=center, rectangle, draw, thick]{moment matching};
	\node at (6, 1.5) (b3) [align=center, rounded rectangle, draw, thick]{\large Learning Mixtures};
	\draw[->,>=latex,thick] (b2.north) -- (b1.south);
	\draw[->,>=latex,thick] (b2.south) -- (b3.north);

	\draw[->,>=latex,thick] (a1.south) -- (b1.north);
	\draw[->,>=latex,thick] (a2.east) to [bend left] (b1.west);

	\draw[->,>=latex,thick] (a2.east) to [bend right] (b3.west);
	\draw[->,>=latex,thick] (a3.north) -- (b3.south);

	\node at (12, 3) (c) [align=center, rectangle, draw,thick]{orthogonal ~\\ polynomials};
	\draw[->,>=latex,thick] (c.west) to [bend right] (b1.east);
	\draw[->,>=latex,thick] (c.west) to [bend left] (b3.east);

\end{tikzpicture}%

\caption{Statistical applications of polynomial methods.\label{fig:org}}
\end{figure}
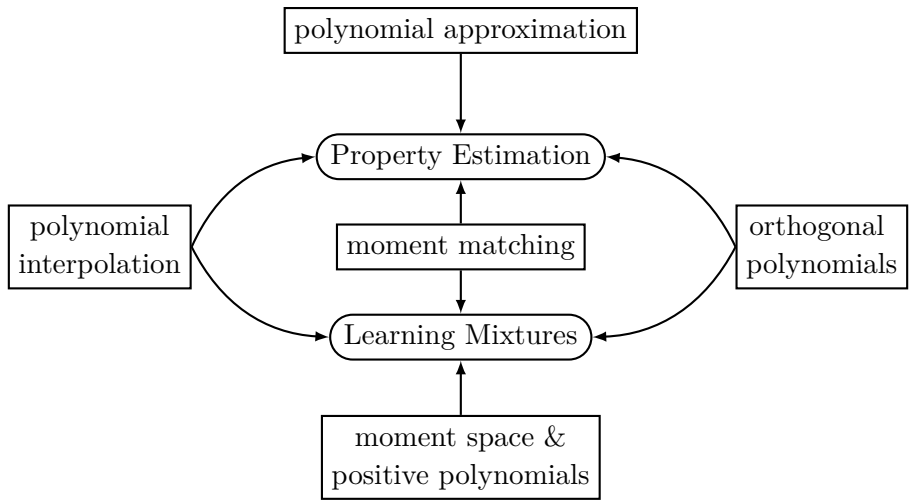

In this survey, we present several tools from the theory of polynomials and their applications in statistical problems. \prettyref{chap:background} provides a brief introduction to the necessary background in the theory of polynomials, including polynomial approximation, interpolation and majorization, theory of moments and positive polynomials, orthogonal polynomials, and Gaussian quadrature. \prettyref{fig:org} describes how these techniques are used in specific statistical applications.

The first statistical application is in the topic of property estimation. \prettyref{chap:approx} introduces some common framework and techniques, including 
 Poisson sampling, 
approximation-theoretic construction of statistical estimators, and minimax lower bounds based on moment matching. 
We then apply these techniques to two representative problems: The problem of entropy estimation is studied in details in \prettyref{chap:entropy}; In \prettyref{chap:unseen}, we study the estimation of the unseen, including estimating the support size and the distinct elements problem.

The second statistical application is learning Gaussian mixture models using moment methods. 
A general framework for mixture models and various moment comparison theorems are developed \prettyref{chap:framework}, which form the underpinnings of our statistical theory. 
Most of these results do not depend on properties of Gaussians and are applicable to general mixture models. 
\prettyref{chap:gm} describes algorithms for Gaussian mixture models and their statistical guarantees, complemented by matching lower bounds.



\section{Notations}
\label{sec:notation}
For $ k\in \naturals $, let $ [k]\triangleq \sth{1,\dots,k} $.
We use standard big-$O$ notations, e.g., for any positive sequences $\{a_n\}$ and $\{b_n\}$, $a_n=O(b_n)$ or $a_n \lesssim b_n$ if $a_n \leq C b_n$ for some absolute constant $C>0$, $a_n=o(b_n)$ or $a_n \ll b_n$ or if $\lim a_n/b_n = 0$.
We write $o_{\delta}(1)$  as $\delta\to0$ to indicate convergence that is uniform in all other parameters.
The notations $a\wedge b$ and $a\vee b$ stand for $\min\{a,b\}$ and $\max\{a,b\}$, respectively.
For a probability measure $\pi$ on the real line, let $F_{\pi}$ denote its cumulative distribution function (CDF), with $F_\pi(t) \triangleq \pi((-\infty,t])$.
A distribution $\pi$ is called $\sigma$-subgaussian if $\Expect_{\pi}[e^{tX}]\le \exp(t^2\sigma^2/2)$ for all $t\in\reals$.
For matrices $A \succeq B$ stands for $A-B$ being positive semidefinite.
The Euclidean ball centered at $x\in\reals^d$ of radius $r$ is denoted by $B(x,r)$.

Denote by $\Binom(n,p)$ the binomial distribution with $n$ Bernoulli trials and success probability $p$. 
For $P=(p_1,\ldots,p_k)$, denote by $\Multinom(n,P)$ the multinomial distribution with $n$ trials where each trial has outcome $i$ with probability $p_i$.
Denote by $N(\mu,\sigma^2)$ the normal distribution with mean $\mu$ and variance $\sigma^2$ and let $\phi(x)\triangleq \frac{1}{\sqrt{2\pi}}e^{-x^2/2}$ denote the standard normal density. Denote by $\Poi(\mu)$ the Poisson distribution with mean $\mu$.

We recall the definition of the following $f$-divergences (cf.~\cite[Chap.~2]{Tsybakov09} for details). For probability distributions $P$ and $Q$, 
the Kullback-Leibler (KL) divergence is $D(P\|Q) \triangleq \int dP \log \frac{dP}{dQ}$ if $P\ll Q$ and $\infty$ otherwise; 
the $\chi^2$-divergence is defined as $\chi^2(P\|Q) \triangleq \int dP (\frac{dP}{dQ}-1)^2$ if $P\ll Q$ and $\infty$ otherwise;
the squared Hellinger distance is $H^2(P,Q) \triangleq \int (\sqrt{\frac{dP}{d\mu}}-\sqrt{\frac{dQ}{d\mu}})^2 d\mu$ and 
the total variation distance is $\TV(P,Q) \triangleq \int |\frac{dP}{d\mu}-\frac{dQ}{d\mu}| d\mu$, for any dominating measure $\mu$ such that $P\ll \mu$ and $Q\ll \mu$.



\chapter{Background}
\label{chap:background}

In this chapter we introduce the necessary background on the theory of polynomials that are useful for statistical inference. 
We focus on three types of topics: \emph{polynomial approximation}, \emph{polynomial interpolation}, and \emph{moments and positive polynomials}. The major probabilistic and statistical applications of the first topic are in estimating properties of distributions (Chapters \ref{chap:approx}--\ref{chap:unseen}), while the second and third will be applied to learning mixture models (Chapters \ref{chap:framework}--\ref{chap:gm}).
For a comprehensive survey on the theory of polynomials see the monographs by Prasolov \cite{Prasolov2009} and Timan \cite{timan63}. We focus on algebraic (ordinary) polynomials in one variable
and briefly discuss trigonometric polynomials.
See \cite{Reimer2012} for extensions to multivariate polynomials.

Polynomials of one variable $x$ are functions of the form
\[
    p_n(x)=a_0+a_1x+a_2x^2+\dots+a_nx^n,
\]
where $n\in\integers_+$, $a_0,a_1,\dots,a_n$ are arbitrary real or complex coefficients. The degree of a polynomial is the highest power in $x$ with a nonzero coefficient. The set $\calP$ of all polynomials is a vector space with countably infinite dimension; if one restricts to polynomials of degree at most $n$, then it is a vector space of $n+1$ dimensions, denoted by $\calP_n$.

The canonical basis for the space of polynomials consists of \emph{monomials} $\{1,x,x^2,\ldots\}$. Any set of $n+1$ polynomials $\{p_0,p_1,\dots,p_n\}$ such that each $p_m$ has degree $m$ can serve as a basis for the polynomials space $\calP_n$, and every polynomial of degree at most $n$ can be uniquely represented by a linear combination of these polynomials via a change of basis.

Trigonometric polynomials are functions in $\theta$ of the form 
\[
    p_n(\theta)=\sum_{k=0}^n(a_k\cos k\theta+b_k\sin k\theta),
\]
with coefficients $a_k$ and $b_k$. The degree of a trigonometric polynomial is the largest $k$ such that $a_k$ and $b_k$ are not both zero. The functions $\cos k\theta$ and $\frac{\sin(k+1)\theta}{\sin\theta}$ are ordinary polynomials in $\cos\theta$, named Chebyshev polynomials of the first and second kind, respectively \cite{timan63}:
\begin{equation}
    \label{eq:cheby}
    \cos k\theta=T_k(\cos \theta),\quad \frac{\sin (k+1)\theta}{\sin\theta}=U_k(\cos \theta).
\end{equation}


\section{Uniform approximation}
\label{sec:uniform}
Approximation theory studies how well functions can be approximated by simpler ones.
In statistical applications, simpler functions are often easier to estimate.
Polynomials are among the most well-studied approximants. 
As mentioned in Section \ref{sec:intro-estimator} and \ref{sec:intro-lb}, polynomial approximation serve as a proxy for constructing estimators for complicated functionals, and the best uniform approximation error is connected to the minimax lower bound. 

In this section we provide a brief summary of some of the key results from the theory of polynomial approximation. 
We start by recalling a fundamental theorem on the denseness of polynomials:
\begin{theorem}[Weierstrass]
\label{thm:Weierstrass}
    Given a function $f$ that is continuous on the interval $[a,b]$, and any $\epsilon>0$, there exists a polynomial $p$ such that
    \[
        |f(x)-p(x)|<\epsilon,\quad \forall x\in[a,b].
    \]
    If $f$ is continuous and has the period $2\pi$, then there exists a trigonometric polynomial $q$ such that
    \[
        |f(x)-q(x)|<\epsilon,\quad \forall x.
    \]
\end{theorem}

This theorem has been proved in many different ways, and can be generalized to the approximation of multivariate continuous functions in a closed bounded region. For more information on this theorem, we refer to \cite[Chapter 1]{timan63}. In the first case of the theorem, an elegant constructive proof is via \emph{Bernstein polynomials} that approximate a continuous function $f$ on $[0,1]$:
\begin{equation}
    \label{eq:Bernstein}
    B_n(x)=\sum_{k=0}^nf(k/n)\binom{n}{k}x^k(1-x)^{n-k}.
\end{equation}
As explained next, Bernstein polynomials have a clear probabilistic interpretation in terms of ``coin flips''. 
The approximation of a function $f$ using Bernstein polynomials can be characterized in terms of its \emph{modulus of continuity} 
\begin{equation}
    \label{eq:mod-cont}
    \omega_f(\delta)=\sup\{f(x)-f(y):|x-y|\le \delta\}.
\end{equation}

\begin{theorem}[T. Popoviciu]
    \label{thm:bernstein}
    Given a continuous function $f$ on $[0,1]$,
    \[
        \sup_{0\le x\le 1}|f(x)-B_n(x)|\le \frac{5}{4}\omega_f(n^{-\frac{1}{2}}).
    \]
\end{theorem}
\begin{proof}
Note that the Bernstein polynomial \prettyref{eq:Bernstein} can be written as 
    \begin{equation}
		    B_n(x)=\Expect[f(N/n)],
		\label{eq:Bernstein-prob}
		\end{equation}
    where $N\sim\Binom(n,x)$. In other words, Bernstein polynomials are the mean of a ``plug-in'' estimator for $f(x)$ based on $n$ \iid coin flips with bias $x$. 
		For any $\delta>0$,
    \begin{gather*}
        |f(x)-B_n(x)|\le \Expect|f(x)-f(N/n)|\indc{|x-\frac{N}{n}|\le \delta}+\\\Expect|f(x)-f(N/n)|\indc{|x-\frac{N}{n}|> \delta}.
    \end{gather*}
    To prove an upper bound of the right-hand side, we note that $|f(x)-f(y)|\le 1+\Floor{\frac{|x-y|}{\delta}}\omega_f(\delta)$. Then we have 
    \[
        |f(x)-B_n(x)| \le \omega_f(\delta)+\frac{\omega_f(\delta)}{\delta}\Expect\abs{x-\frac{N}{n}}\indc{|x-\frac{N}{n}|> \delta}.
    \]
		The second term can be bounded using the variance of the binomial distribution as
    \[
        \Expect\abs{x-\frac{N}{n}}\indc{|x-\frac{N}{n}|> \delta}
        < \frac{1}{\delta}\Expect\abs{x-\frac{N}{n}}^2=\frac{x(1-x)}{n\delta}
        \le \frac{1}{4n\delta}.
    \]
    The desired statement follows by choosing $\delta = n^{-1/2}$.
\end{proof}

In general, approximation by Bernstein polynomials is not optimal. In fact, for Bernstein polynomials the rate $\omega_f(n^{-\frac{1}{2}})$ in \prettyref{thm:bernstein} is tight.\footnote{This can be shown by considering $f(x)=|x-\frac{1}{2}|$, for which $\omega_f(\delta) \asymp \delta$; on the other hand, from \prettyref{eq:Bernstein-prob}, it is clear that $B_n(\frac{1}{2}) = \frac{1}{n} \Expect[|N-\frac{n}{2}|] = \Theta(\frac{1}{\sqrt{n}})$ by the central limit theorem.}
In \prettyref{sec:plug-in}, the suboptimality of the Bernstein polynomials also explains the severe bias of the plug-in estimator.
A major result in the approximation theory, Jackson's theorem states that any continuous function can be uniformly approximated with error $\omega_f(n^{-1})$. 
\begin{theorem}[Jackson]
\label{thm:jackson}
    Given a continuous function $f$ on $[0,1]$, there exists a polynomial $P_n$ of degree at most $n$ such that
    \[
        \sup_{0\le x\le 1}|f(x)-P_n(x)|\le 3\omega_f(n^{-1}).
    \]
\end{theorem}
Note that it is clear from \prettyref{eq:Bernstein-prob} that for positive function $f$, Bernstein polynomials are also positive. In contrast, Jackson's construction uses trigonometric polynomials which have both positive and negative parts.
This is analogous to the well-known phenomenon in nonparametric statistics that nonnegative kernels are insufficient to leverage higher-order smoothness \cite{Tsybakov09}. 
For more details on \prettyref{thm:jackson} we refer to \cite{timan63}. Generalizations and extensions, called Jackson-type theorems, provide approximation guarantees in terms of various notions of modulus of continuity. 
See \cite[Chapter V]{timan63} and \cite[Chapter 7]{DL93} for more constructive approximations.

\subsection{Best uniform approximation}
\label{sec:best}
The study of best uniform approximation was initiated by Chebyshev. 
Let $f$ be a continuous function on an interval $[a,b]$. Consider its best uniform approximation by degree-$n$ polynomials and denote the best approximation error by 
\begin{equation}
    \label{eq:best}
    E_n(f,[a,b]) =\inf_{P\in\mathcal{P}_n}\sup_{x\in[a,b]}|f(x)-P(x)|.
\end{equation}
Jackson-type theorems (such as \prettyref{thm:jackson}) provide upper bounds on the best approximation error \prettyref{eq:best} in terms of various moduli of continuity. 
Conversely, it is also possible to use moduli of continuity to lower bound the best approximation error, although this is typically carried out indirectly. 
This type of impossibility results are needed for proving statistical lower bound (see \prettyref{sec:d-mm} and \prettyref{lmm:sep}).
However, there is no converse theorem in terms of $\omega_f$ in \prettyref{eq:mod-cont} due to the special behavior of the best approximating polynomial near the boundary of the approximation interval, as firstly observed by Nikolsky \cite{nikolsky1946mean}; see further discussions in \cite[Sec.~6.1.7]{timan63} and \cite[Chapter 8]{DL93}.
A modulus of continuity with refined measurements of the smoothness near the boundary is needed to establish impossibility results. 
Below we present one such result due to K. Ivanov \cite{Ivanov1983} (see also \cite[Section 3.4]{petrushev2011rational}); see \cite{DK2012} for similar results in terms of the Ditzian-Totik moduli of smoothness.
\begin{theorem}
    \label{thm:Ivanov}
    Define 
    \begin{align*}
    & \Delta_n(x)\triangleq\frac{1}{n}\sqrt{1-x^2}+\frac{1}{n^2},\\
    & \tau_1(f,\Delta_n)\triangleq\sup\{|f(x)-f(y)|: x,y\in[-1,1],|x-y|\le \Delta_n(x)\}.
    \end{align*}
    Then, there exist absolute constants $c_1$ and $c_2$ such that
    \begin{align}
    & E_n(f,[-1,1])\le c_1\tau_1(f,\Delta_n), \label{eq:ivanov1}\\
    & \tau_1(f,\Delta_n) \le \frac{c_2}{n}\sum_{s=0}^nE_s(f,[-1,1]). \label{eq:ivanov2}
    \end{align}
\end{theorem}
Note that the converse result \prettyref{eq:ivanov2} is in terms of the best approximation error averaged over all degrees. 
To produce a concrete lower bound on an individual approximation error with degree say $L$, one needs to use \prettyref{eq:ivanov2} in conjunction with the upper bound \prettyref{eq:ivanov1}. Indeed, by the monotonicity of $E_s$ in the degree $s$, we have
\[
\frac{n-L}{n} E_L(f,[-1,1]) \geq \frac{1}{c_2} \tau_1(f,\Delta_n) - \frac{c_1}{n} \sum_{s=0}^L \tau_1(f,\Delta_s)
\]
and optimize over $n \geq L$; for a concrete example see \cite[Appendix F]{WY14}.

It is known that the infimum in \prettyref{eq:best} is achieved by a unique polynomial (see, \eg, \cite[Chapter 3]{DL93}), with the following remarkable characterization:
\begin{theorem}[Chebyshev alternation theorem]
    \label{thm:cheby}
    A polynomial $P_n\in\calP_n$ is the best uniform approximation of a continuous function $f$ on $[a,b]$ by $\calP_n$ if and only if there exists $n+2$ points $x_j$, $a\le x_0<\dots<x_{n+1}\le b$ such that $f(x_j)-P_n(x_j)=\pm \sup_{x\in[a,b]}|f(x)-P(x)|$ with successive changes of sign, \ie, $f(x_{j+1})-P_n(x_{j+1})=-(f(x_{j})-P_n(x_{j}))$ for $j=0,\dots,n$.
\end{theorem}
In addition to the ordinary polynomials, the above characterization holds for any real \emph{Haar system}
such as the trigonometric polynomials. See \cite[Section 3.3 -- 3.5]{DL93} for a proof of this theorem and more information.

For certain special cases the exact value of the best approximation error and the explicit formula of the best polynomial approximant are known; see \cite[Section 2.11]{timan63} for examples with explicit solutions. We shall give one example due to Chebyshev, which will be used in \prettyref{chap:unseen} for the statistical problem of estimating the unseen. 
\begin{theorem}
    \label{thm:cheby-best}
    For $n\in\naturals$, the degree-$n$ monic polynomial (i.e., with leading coefficient equal to one) that deviates the least from zero over $[-1,1]$ is $\frac{1}{2^{n-1}}T_n(x)$, where $T_n$ is the Chebyshev polynomial of the first kind given by \prettyref{eq:cheby}. Furthermore, the value of its deviation is $\sup_{x\in[-1,1]} |T_n(x)|=\frac{1}{2^{n-1}}$.
\end{theorem}
\begin{proof}
    Observe that the problem is equivalent to finding the best polynomial of degree $n-1$ to approximate the monomial $x^n$ over $[-1,1]$:
    \[
        \inf_{a_0,\dots,a_{n-1}}\sup_{x\in[-1,1]}|x^n-a_{n-1}x^{n-1}-\dots-a_1x_1-a_0|.
    \]
    The polynomial $\frac{1}{2^{n-1}}T_n(x)$ is monic with maximum magnitude $\frac{1}{2^{n-1}}$. Furthermore, the Chebyshev polynomial $T_n$ successively attains $1$ or $-1$ at $\cos(k\pi/n)$ for $k=0,\dots,n$. The optimality of $\frac{1}{2^{n-1}}T_n(x)$ follows from \prettyref{thm:cheby}.
\end{proof}

In general, for a given function $f$, there is no known close-form formula for its best polynomial approximation; nevertheless, many fast algorithms have been developed.
Note that the optimization problem in \prettyref{eq:best} can be rewritten as a linear program (LP) with $n+2$ decision variables and infinitely many constraints: 
\begin{equation}
    \label{eq:primal-approx}
    \begin{aligned}
        E_n(f,[a,b])=\min\quad & t\\
        \textrm{s.t.}\quad & a_0+a_1x+\dots+a_nx^n-t\le f(x),\quad x\in[a,b],\\
        & a_0+a_1x+\dots+a_nx^n+t\ge f(x),\quad x\in[a,b].
    \end{aligned}
\end{equation}
Thanks to Chebyshev's alternation theorem (\prettyref{thm:cheby}), instead of enforcing the constraints for all $x \in [a,b]$, it suffices to do so for the alternating points (maxima of the approximation error) corresponding to the optimal polynomial. This motivates an iteration scheme called the \emph{Remez algorithm} (\prettyref{algo:remez}), which successively updates the polynomial by solving a linear system and the constraint sets by the local maxima of the approximation error. 
\begin{algorithm}[ht]
    \caption{Remez algorithm.}
    \label{algo:remez}
    \begin{algorithmic}[1]
        \REQUIRE a continuous function $f$, an interval $[a,b]$, a degree $n$.
        \ENSURE a polynomial $P$ of degree at most $n$.
        \STATE{Initialize $n+2$ points $a\le x_0<x_1<\dots<x_{n+1}\le b$. }
        \REPEAT
            \STATE{Solve the system of linear equations
            \[
                f(x_j)-Q_n(x_j)=(-1)^j\delta,\quad j=0,\dots,n+1,
            \]
            where $Q_n(x)=\sum_{i=0}^n a_ix^i$, with respect to unknowns $\delta,a_0,\dots,a_{n}$.
            }
            \STATE{Find $\xi$ and $d$ such that 
            \[
                |f(\xi)-Q_n(\xi)|=\max_{x\in[a,b]}|f(x)-Q_n(x)|=d.
            \]
            }
            \STATE{Update the sequence $x_0<\dots<x_{n+1}$ by replacing one $x_j$ by $\xi$ so that $f-Q_n$ successively changes sign.}
        \UNTIL{stopping criterion is satisfied.}
        \STATE{Report $Q_n$.}
    \end{algorithmic}
\end{algorithm}
See \cite{petrushev2011rational} for the proof of correctness and convergence rates of the Remez algorithm.

\section{Duality and moment matching}
\label{sec:dual-best}
We have shown in \prettyref{eq:poly-dual} -- \prettyref{eq:poly-primal} that the dual program of the (infinite-dimensional) LP \prettyref{eq:primal-approx} is the following moment matching problem
\begin{equation}
    \label{eq:dual-mm}
    \begin{aligned}
 2 E_n(f,[a,b]) =       \sup \quad & \Expect_\nu[f(X)]-\Expect_{\nu'}[f(X)],\\
        \textrm{s.t.}\quad & \Expect_\nu[X^j]=\Expect_{\nu'}[X^j],\quad j=0,\dots,n,\\
        & \nu,\nu'\textrm{ are supported on }[a,b],
    \end{aligned}
\end{equation}
where the supremum is over pairs of distributions $\nu$ and $\nu'$.
The strong duality between \prettyref{eq:primal-approx} and \prettyref{eq:dual-mm} can be verified using the general theory of convex optimization (see \cite[pp.~48--50]{Rockafellar1974}) or by Chebyshev's alternating theorem. In the primal problem, as a consequence of Chebyshev's characterization in \prettyref{thm:cheby}, 
there exist $n+2$ points where the constraints are binding for the optimal solution. 
Consequently, in the dual problem \prettyref{eq:dual-mm}, the optimal $\nu_*$ and $\nu'_*$ are supported on those points by complementary slackness. The dual solution can be obtained accordingly from the primal solution: 
\begin{theorem}
    \label{thm:mm-best}
    Denote by $P^*$ the best polynomial that achieves $E_n(f,[a,b])$ in \prettyref{eq:best}. Suppose $P^*\ne f$ and the maximum deviation of $P^*$ from $f$ is attained at $a\le x_0<\dots<x_{n+1}\le b$ such that $f(x_i)-P^*(x_i)$ successively changes sign. 
    The dual optimal solution of \prettyref{eq:dual-mm} is 
    \begin{gather*}
        \nu_*(x_i)=\frac{2w_i}{w_0+w_1+\dots+w_{n+1}},\quad f(x_i)>P^*(x_i),\\
        \nu'_*(x_i)=\frac{2w_i}{w_0+w_1+\dots+w_{n+1}},\quad f(x_i)<P^*(x_i),
    \end{gather*}
    where $w_i=(\prod_{j\ne i}|x_i-x_j|)^{-1}$.
\end{theorem}
\begin{proof}
    Note that $\nu_*$ is supported on either $\{x_0,x_2,\dots\}$ or $\{x_1,x_3,\dots\}$ and $\nu'_*$ is supported on the rest. Denote $\epsilon = E_n(f,[a,b])=\sup_{x \in[a,b]} |f(x)-P^*(x)|$. Then $f-P^*$ is almost surely $\epsilon$ and $-\epsilon$ under $\nu_*$ and $\nu'_*$, respectively.

    We first verify the feasibility. 
    Note that $\prod_{j\ne i}(x_i-x_j)$ has alternating signs for $i=0,\dots,n+1$.
    Hence the moment matching constraints in \prettyref{eq:dual-mm} is equivalent to $\sum_{i=0}^{n+1} \frac{x_i^m}{\prod_{j\ne i}(x_i-x_j)}=0$ for $m=0,\dots,n$.
    For each $m\in\{0,1,\dots,n\}$, consider the polynomial $P(x)=\sum_{i=0}^{n+1} x_i^m\frac{\prod_{j\ne i}(x-x_j)}{\prod_{j\ne i}(x_i-x_j)}$ of degree at most $n+1$. Then $P(x)$ coincides with $x^m$ on $n+2$ distinct points $x_0,\dots,x_{n+1}$. Hence $P(x)\equiv x^m$ and $\sum_{i=0}^{n+1} \frac{x_i^m}{\prod_{j\ne i}(x_i-x_j)}=0$.
    In particular, the special case of $m=0$ shows that $\sum_i \nu_*(x_i)=\sum_i \nu'_*(x_i)$. Since $\sum_i \nu_*(x_i)+\sum_i \nu'_*(x_i)=2$ by construction, this verifies that $\nu_*$ and $\nu'_*$ are well-defined probability distributions.

    For optimality it suffices to show a zero duality gap: 
    \[
        \Expect_{\nu_*}[f]-\Expect_{\nu'_*}[f]
        =\Expect_{\nu_*} [f-P^*]-\Expect_{\nu'_*} [f-P^*]
        =2\epsilon.              
    \]
    The first equality is due to the moment matching constraints. 
\end{proof}

\section{Polynomial interpolation: Lagrange and Newton form}
\label{sec:interpol}
Interpolation is a method of estimating the value of a function within the range of a discrete set of data points. 
Given $(x_i,f_i)$ for $i=0,\dots,n$, the interpolation problem amounts to finding a simple function  $P$ such that
\begin{equation}
\label{eq:interp}
P(x_i)=f_i,\quad i=0,\dots,n.
\end{equation}
Examples of the simple function $P$ include ordinary polynomials and trigonometric polynomials. Interpolation also offers a useful primitive for approximating a given function by  interpolating it on a sagaciously chosen set of points. For a comprehensive survey on related topics, see \cite{Davis1975,rivlin2003introduction}.

The main result for polynomial interpolation in one dimension is the following:
\begin{theorem}
\label{thm:interp}
    Given distinct data points $(x_i,f_i)$ for $i=0,\dots,n$, there exists a unique interpolating polynomial $P$ of degree at most $n$ such that
    \begin{equation}
      P(x_i)=f_i,\quad i=0,\dots,n.
      \end{equation}
\end{theorem}
\begin{proof}
The existence is given by Lagrange or Newton formula discussed next.
For uniqueness, given two interpolating polynomials $P$ and $P'$ of degree at most $n$, the polynomial $Q=P-P'$ is of degree at most $n$ satisfying $Q(x_i)=0$ for $i=0,\dots,n$. Thus $Q\equiv 0$. 
\end{proof}

Interpolating polynomials are the main tool to construct estimator for the distinct elements problem in \prettyref{sec:distinct} and to prove moment comparison theorems in \prettyref{chap:framework}. For these applications it is critical to have a good control over the coefficients of the interpolating polynomial. 
To this end, we analyze the explicit formula of the interpolant: the \emph{Lagrange formula} and the \emph{Newton formula}.
\begin{itemize}
\item The \emph{Lagrange formula} for the interpolating polynomial $P$ is explicitly constructed in terms of the Lagrange basis:
\begin{equation}
    L_i(x) \triangleq \prod_{j\ne i}\frac{x-x_j}{x_i-x_j}=
    \begin{cases}
        1,& x=x_i,\\
        0,& x=x_j,j\ne i.
    \end{cases}
\label{eq:lagrange-basis}
\end{equation}
By linearity, we obtain the interpolation polynomial satisfying \prettyref{eq:interp}
\begin{equation}
    \label{eq:Lagrange}
    P(x)=\sum_{i=0}^nf_iL_i(x).
\end{equation}

\item The \emph{Newton formula} for the interpolating polynomial is of the form
\begin{gather}
    P(x)=a_0+a_1(x-x_0)+a_2(x-x_0)(x-x_1)+\dots\nonumber\\
    +a_n(x-x_0)\cdots(x-x_{n-1}).
    \label{eq:interpolation-Newton}
\end{gather}
The coefficients of the Newton form \prettyref{eq:interpolation-Newton} can be successively calculated by
\begin{align*}
  &f_0=P(x_0)=a_0,\\
  &f_1=P(x_1)=a_0+a_1(x_1-x_0),\\
  &\dots
\end{align*}
\end{itemize}

\begin{figure}[ht]
    \centering
    \begin{subfigure}[b]{.45\textwidth}
        \includegraphics[width=\linewidth]{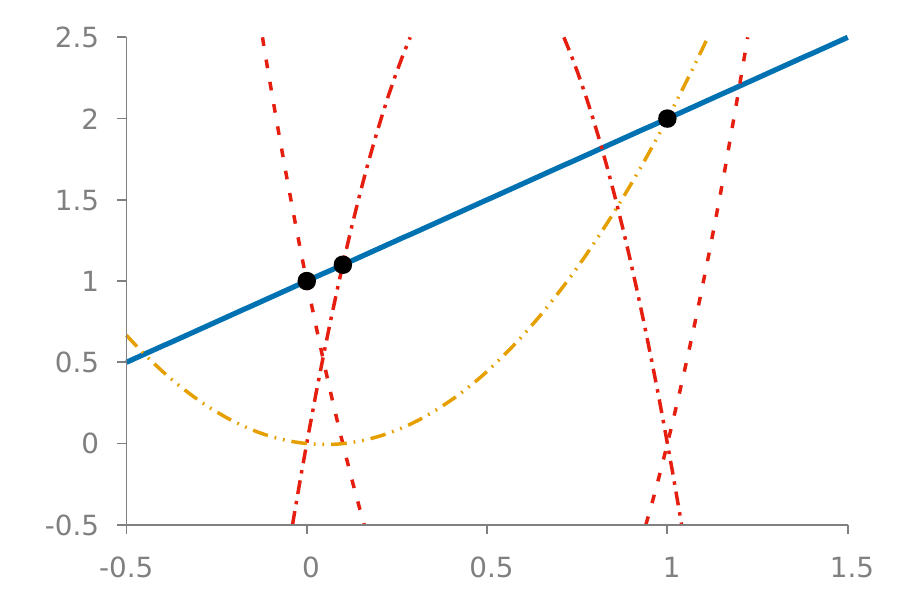}
        \caption{Lagrange formula}
        \label{fig:Lagrange}
    \end{subfigure}
    \hfill
    \begin{subfigure}[b]{.45\textwidth}
        \includegraphics[width=\linewidth]{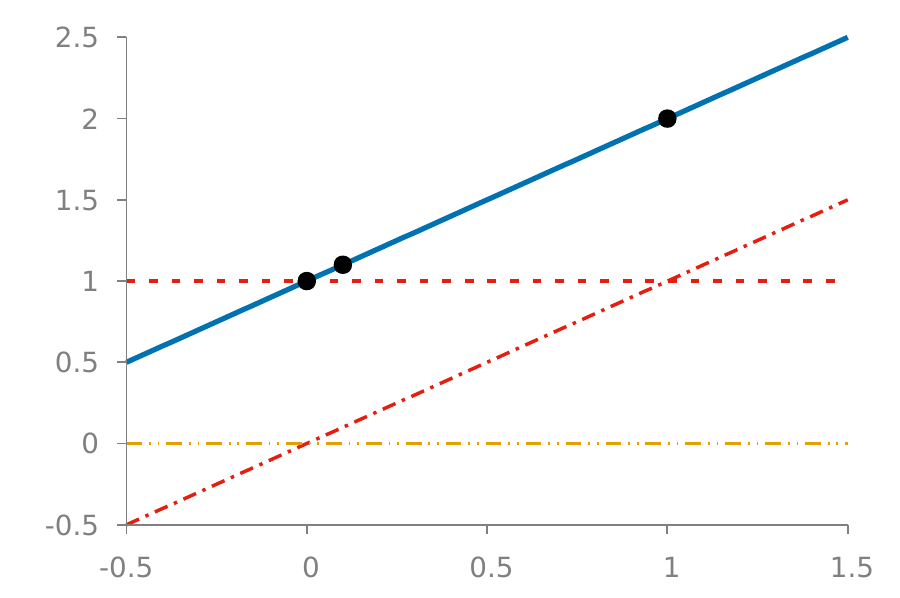}
        \caption{Newton formula}
        \label{fig:Newton} 
    \end{subfigure}
    \caption{Interpolation on three data points $(0,1)$, $(0.1,1.1)$, and $(1,2)$, shown in black dots. (\subref{fig:Lagrange}) Three terms in Lagrange formula are shown in dashed lines, summing up to the interpolating polynomial in the solid line. (\subref{fig:Newton}) Three terms in Newton formula are shown in dashed lines, and the same interpolating polynomial in solid.}
    \label{fig:formulas}
\end{figure}

In numerical analysis, the Newton form of polynomial interpolation is usually introduced for computational considerations so that, unlike the Lagrange form, one does not need to recompute all coefficients when an extra node is introduced \cite{stoer.2002}. 
For our statistical applications of learning mixture models in Chapters \ref{chap:framework} and \ref{chap:gm}, the Newton form turns out to be crucial, which offers better bound on the coefficients of the interpolating polynomials because it takes into account the cancellation between each terms in the polynomial.
Indeed, in the Lagrange form \prettyref{eq:Lagrange}, if two nodes are very close, then each term can be arbitrarily large, even if $f$ itself is a smooth function. In contrast, each term of \prettyref{eq:interpolation-Newton} is stable when $f$ is smooth since the coefficients are closely related to derivatives. The following example illustrates this point:
\begin{example}[Lagrange versus Newton form]
    \label{ex:Lagrange-Newton}
    Given three points $x_1=0, x_2=\epsilon, x_3=1$ with $f(x_1)=1, f(x_2)=1+\epsilon, f(x_3)=2$, the interpolating polynomial is $P(x)=x+1$. The next equation gives the interpolating polynomial in Lagrange's and Newton's form respectively. The interpolations are illustrated in \prettyref{fig:formulas}.
    \begin{align*}
      \text{Lagrange: }&P(x)=\frac{(x-\epsilon)(x-1)}{\epsilon}+(1+\epsilon)\frac{x(x-1)}{\epsilon(\epsilon-1)}+2\frac{x(x-\epsilon)}{1-\epsilon};\\
      \text{Newton: } &P(x)=1+x+0.
    \end{align*}
\end{example}

In general, the coefficients $a_k$ of Newton formula \prettyref{eq:interpolation-Newton} coincide with the divided differences $f_{0\dots k}$ that are recursively defined as 
\begin{equation}
    \label{eq:div-recursion}
    f_{i_0i_1\dots i_k}=\frac{f_{i_1\dots i_k}-f_{i_0\dots i_{k-1}}}{x_{i_k}-x_{i_0}}.
\end{equation}
The above recursion can be calculated with the help of \emph{Neville's diagram} as shown in \prettyref{fig:neville0} (cf.~\cite[Section 2.1.2]{stoer.2002}):
\begin{figure}[ht]%
\centering
\begin{center}
    \begin{tikzpicture}[every node/.style={draw,shape=circle,fill,inner sep=1.5pt,align=left},scale=0.6,transform shape]
        \def\num{3}
        \def\shift{0.4}
        \def\lshift{1.5}
        \foreach \n in {0,...,\num} {
            \foreach \k in {0,...,\n} {
                \node at (-2*\n,-2*\k+\n){};
            }
        }

        \pgfmathtruncatemacro\numd{\num-1}
        \foreach \n in {1,...,\numd} {
            \foreach \k in {1,...,\n} {
                \draw (-2*\n-2,-2*\k+2+\n-1) -- (-2*\n,-2*\k+2+\n) -- (-2*\n-2,-2*\k+2+\n+1);
            }
        }
        
        \foreach \k in {1,...,\num} {
            \pgfmathtruncatemacro\kd{\k-1}
            \node[draw=none, fill=none] at (-2*\num-\lshift,-2*\k+2+\num){$x_{\kd}$};
            \node[draw=none, fill=none] at (-2*\num,-2*\k+2+\num+\shift){$f_{\kd}$};
            \node[draw=none, fill=none] at (-2*\k,-\k+2-1){$\vdots$};
        }
        
        \foreach \k in {1,...,\numd} {
            \pgfmathtruncatemacro\kd{\k-1}
            \node[draw=none, fill=none] at (-2*\num+2,-2*\k+2+\num-1+\shift){$f_{\kd\k}$};
            \draw (-2*\k,-\k) -- (-2*\k-2,-\k-1);
        }

        \node[draw=none, fill=none] at (-2,1+\shift){$f_{012}$};
        \node[draw=none, fill=none] at (-0,0+\shift){$f_{0\dots n}$};
        \node[draw=none, fill=none] at (-6-\lshift,-3){$x_{n}$};
        \node[draw=none, fill=none] at (-6-\lshift,-2){$\vdots$};
        \node[draw=none, fill=none] at (-6,-3-\shift){$f_{n}$};
        \draw[dotted] (-2*.6,-1*.6) -- (-2*.4,-1*.4);
        \draw[dotted] (-2*.6,1*.6) -- (-2*.4,1*.4);
        \node[draw=none, fill=none] at (-6,4){$k=0$};
        \node[draw=none, fill=none] at (-4,4){$1$};
        \node[draw=none, fill=none] at (-2,4){$2$};
        \node[draw=none, fill=none] at (-1,4){$\dots$};
        \node[draw=none, fill=none] at (-0,4){$n$};
        \draw (-8,3.7) -- (1,3.7);
        \draw (-6.8,4.5) -- (-6.8,-4);
    \end{tikzpicture}
\end{center}%
\caption{Neville's diagram for computing the interpolation polynomial in the Newton form. 
}%
\label{fig:neville0}%
\end{figure}
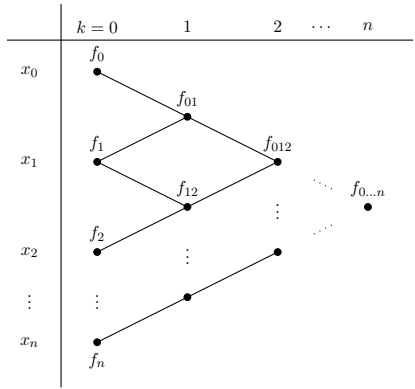
In Neville's diagram, the $k\Th$ order divided differences are computed in the $k\Th$ column, and are determined by the previous column and the interpolation nodes $x_0,\dots,x_n$. The coefficients in \prettyref{eq:interpolation-Newton} are found in the top diagonal.

If the data points correspond to values $f_i=f(x_i)$ of a given function $f$ on a set of distinct points (commonly referred to as \emph{nodes}) $\{x_0,\dots,x_{n}\}$, the divided difference $f_{i_0\dots i_k}$ can be viewed as a multivariate function of the nodes $x_{i_0},\dots,x_{i_k}$ and is denoted by 
\begin{equation}
\label{eq:div-diff-f}
f[x_{i_0},\dots,x_{i_k}]=f_{i_0\dots i_k}.
\end{equation}
If $f$ is $k$-times differentiable, then its $k\Th$ order divided difference admits the following integral representation (known as the \emph{Genocchi-Hermite formula}):
\begin{gather}
    f[x_0,\dots,x_k]=\int_0^1\int_0^{s_1}\dots\int_0^{s_{k-1}}f^{(k)}((1-s_1)x_0+\dots\nonumber\\+(s_{k-1}-s_k)x_{k-1}+s_kx_k)\diff s_k\dots\diff s_1,
    \label{eq:div-GH}
\end{gather}
which, in the special case of $k=1$, reduces to the fundamental theorem of calculus: $f[x_0,x_1]= \int_0^1 f'((1-s)x_0+sx_1) ds$. 
For details on this formula and other representations of the divided differences, see \cite{deBoor2005}.

Next we study the interpolation error. 
The remainder in the polynomial interpolation can be conveniently expressed in terms of the divided differences
\begin{equation}
    \label{eq:remainder}
    R(x)\triangleq f(x)-P(x)=f[x_0,\dots,x_n,x]\prod_{i=0}^n(x-x_i).
\end{equation}
If the function $f$ is $(n+1)$-times differentiable, then the remainder term can be represented using \prettyref{eq:div-GH} by
\begin{equation}
    \label{eq:remainder-xi}
    R(x)=\frac{f^{(n+1)}(\xi)}{(n+1)!}\prod_{i=0}^n(x-x_i),
\end{equation}
for some $\xi$ in the convex hull of $\{x_0,\dots,x_n,x\}$.

By using the interpolant as a (suboptimal) approximant, the remainder offers a convenient upper bound for the \emph{pointwise} approximation error, which is particularly useful when an explicit nonasymptotic bound is needed; for a statistical application see \prettyref{thm:tv-bound}.
The formula \prettyref{eq:remainder-xi} can be applied to analyze the approximation error of the interpolation polynomial for a given set of nodes. Superior to equidistant nodes, the \emph{Chebyshev nodes} consist of the zeros of Chebyshev polynomials, which, on the interval $[a,b]$, are given by
\[
    x_i=\frac{b+a}{2}+\frac{b-a}{2}\cos\pth{\frac{2k+1}{2n+2}\pi},\quad k=0,\ldots,n,
\]
The interpolating polynomial with respect to the Chebyshev nodes satisfies the following error bound (see \cite[Eq. (4.7.28)]{Atkinson89})
\begin{equation}
    \label{eq:interp-cheby}
    |R(x)|\le \frac{\max_{x\in[a,b]}|f^{(n+1)}(x)|}{2^n(n+1)!}\pth{\frac{b-a}{2}}^{n+1},\quad x\in[a,b],
\end{equation}
In fact, it is known that for any continuous function, the uniform approximation error of the interpolating polynomial corresponding to the Chebyshev nodes is within a logarithmic factor of the best approximation error, while that of the equidistant nodes can be off by an exponential factor \cite[Chapter 4]{rivlin2003introduction}.


\section{Hermite interpolation and majorizing polynomials}
\label{sec:hermite}
Polynomial interpolation can be generalized to interpolate not only the values of the function but also its derivatives; this is known as \emph{Hermite interpolation}.
For certain functions, e.g.~step functions, Hermite interpolation can be used to construct majorizing (resp.~minorizing) polynomials that are pointwise larger (resp.~smaller) than the given function.
The majorizing and minorizing polynomials naturally produce sandwich bounds for proving moment comparison theorems in \prettyref{chap:framework}.

We start with the counterpart to \prettyref{thm:interp} for Hermite interpolation:
\begin{theorem}
    Given distinct real numbers $x_0<x_1<\dots<x_n$, and values $f_i^{(k)}$ for $i=0,\dots,n$ and $k=0,\dots,m_i$, there exists a unique polynomial of degree at most $N=n+\sum_im_i$ such that
    \[
        P^{(k)}(x_i)=f_i^{(k)},\quad i=0,\dots,n,~k=0,\dots,m_i.
    \]
\end{theorem}
\begin{proof}
The existence is given by the generalized Lagrange or Newton formula introduced next.
    Given two interpolating polynomials $P$ and $P'$ of degree at most $N$, the polynomial $Q=P-P'$ is of degree at most $N$ and satisfies $Q^{(k)}(x_i)=0$ for $i=0,\dots,n$ and $k=0,\dots,m_i$. Therefore, each $x_i$ is a root of $Q$ of multiplicities $m_i+1$. Since $\sum_i(m_i+1)>N$, $Q\equiv 0$, and the uniqueness follows. 
\end{proof}

Analogous to the Lagrange formula \prettyref{eq:Lagrange}, the interpolating polynomial can be explicitly constructed with the help of the generalized Lagrange polynomials $L_{i,k}$ satisfying 
\[
    L_{i,k}^{(k')}(x_{i'})=
    \begin{cases}
        1,&i=i',k=k',\\
        0,&\text{otherwise}.
    \end{cases}
\]
For an explicit formula of the generalized Lagrange polynomials, see \cite[pp.~52--53]{stoer.2002}. The Hermite interpolating polynomial can then be expressed as a linear combination
\[
    P(x)=\sum_{i,k}f_i^{(k)}L_{i,k}(x).
\]

The Newton formula \prettyref{eq:interpolation-Newton} can also be extended by using generalized divided differences when repeated nodes are present:
\begin{equation}
    \label{eq:general-div}
    f[x_0,\dots,x_k]=\frac{f^{(k)}(x_0)}{k!},\quad x_0=x_1=\dots=x_k.
\end{equation}
To this end, we define an expanded sequence by repeating each $x_i$ for $k_i$ times:
\begin{equation}
    \label{eq:new-sequence}
    \underbrace{x_0=\ldots =x_0}_{k_0}<\underbrace{x_1=\ldots =x_1}_{k_1}<\ldots<\underbrace{x_m=\ldots=x_m}_{k_m}.
\end{equation}
The Hermite interpolating polynomial is obtained by \prettyref{eq:interpolation-Newton} using this new sequence and generalized divided differences, which can also be calculated from the Neville's diagram by replacing differences by derivatives whenever encountering repeated nodes. When the data points are sampled from a given function $f$, the remainder formulas \prettyref{eq:remainder} and \prettyref{eq:remainder-xi} hold verbatim.


Next we give an explicit example using Hermite interpolation to construct majorizing polynomials for step functions, which will be used to prove moment comparison theorems in \prettyref{chap:framework}.
Similar constructions can be carried out for functions other than step functions; cf.~\cite[Theorem 5.4]{Freud}.

\begin{example}[Hermite interpolation as polynomial majorization]
    \label{ex:Hermite}
    Let $f(x)=\indc{x\le 0}$. We want to find a polynomial majorization $P\ge f$ such that $P(x)=f(x)$ on $x=\pm 1$. To this end we interpolate $f$ on $\{-1,0,1\}$ with the following constraints:
    \begin{center}
        \begin{tabular}{ c | c  c  c }
            \toprule
            $x$ & $-1$ & $0$ & $1$ \\
            \hline
            $P(x)$  & 1 & 1 & 0 \\
            $P'(x)$ & 0 & \text{any} & 0 \\
            \bottomrule
        \end{tabular}
    \end{center}
    The resulting interpolating polynomial $P$ has degree four and majorizes $f$ \cite[p.~65]{Akhiezer1965}. To see this, we note that $P'(\xi)=0$ for some $\xi\in (-1,0)$ by Rolle's theorem. Since $P'(-1)=P'(1)=0$, $P$ has no other stationary point than $-1,\xi,1$, and thus decreases monotonically in $(\xi,1)$. Hence, $-1,1$ are the only local minimum points of $P$, and thus $P\ge f$ everywhere. The polynomial $P$ is shown in \prettyref{fig:step-approx}.
Similarly, one can construct the minorizing polynomial (see \prettyref{fig:major-minor} in \prettyref{sec:compare-wass}).

    To explicitly compute the majorizing polynomial, we first find the expanded sequence $-1,-1,0,1,1$ per \prettyref{eq:new-sequence}. Applying Newton formula \prettyref{eq:interpolation-Newton} with generalized divided differences from the Neville's diagram \prettyref{fig:Neville}, we obtain $P(x)=1-\frac{1}{4}x(x+1)^2+\frac{1}{2}x(x+1)^2(x-1)$.

    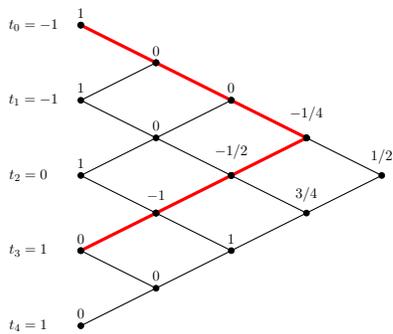
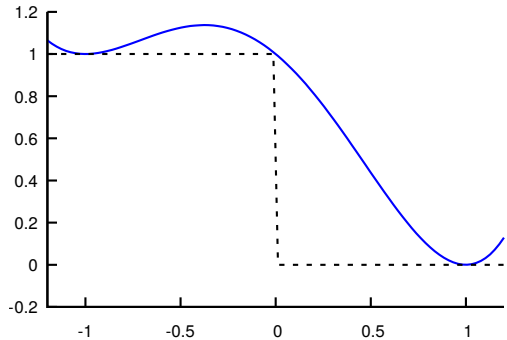
\begin{figure}[ht]
        \centering
        \begin{subfigure}[b]{.48\textwidth}
            \begin{tikzpicture}[every node/.style={draw,shape=circle,fill,inner sep=1.5pt,align=left},scale=0.5,transform shape]
    \def\labels{{1,1,1,0,0},{0,0,$-1$,0},{0,$-1/2$,1},{$-1/4$,3/4},{1/2}}
    \foreach[count=\n] \x in \labels{
        \foreach[count=\k] \y in \x{
            \node[label={[yshift=0cm]\y}] at (2*\n,-2*\k-\n) {};
            \ifthenelse{\n > 1}{
                \draw ({2*(\n-1)},{-2*\k-(\n-1)}) -- (2*\n,-2*\k-\n) -- ({2*(\n-1)},{-2*(\k+1)-(\n-1)});
            }{}
        }
    }

    \def\labels{$t_0=-1$,$t_1=-1$,$t_2=0$,$t_3=1$,$t_4=1$}
    \foreach[count=\k] \y in \labels{
        \node[draw=none, fill=none, anchor=west] at (0,-2*\k-1) {\y};
    }

    \foreach \n in {1,2,3} {
        \draw[red,very thick]
        (2*\n,-2-\n) node[black,thin] {} -- ({2*(\n+1)},{-2-(\n+1)}) node[black,thin] {}
        (2*\n,{-2*(5-\n)-\n}) node[black,thin] {} -- ({2*(\n+1)},{-2*(4-\n)-(\n+1)}) node[black,thin] {};
    }



\end{tikzpicture}
            \caption{Neville's diagram.}
            \label{fig:Neville}
        \end{subfigure}
        ~
        \begin{subfigure}[b]{.48\textwidth}
        \includegraphics[width=\textwidth]{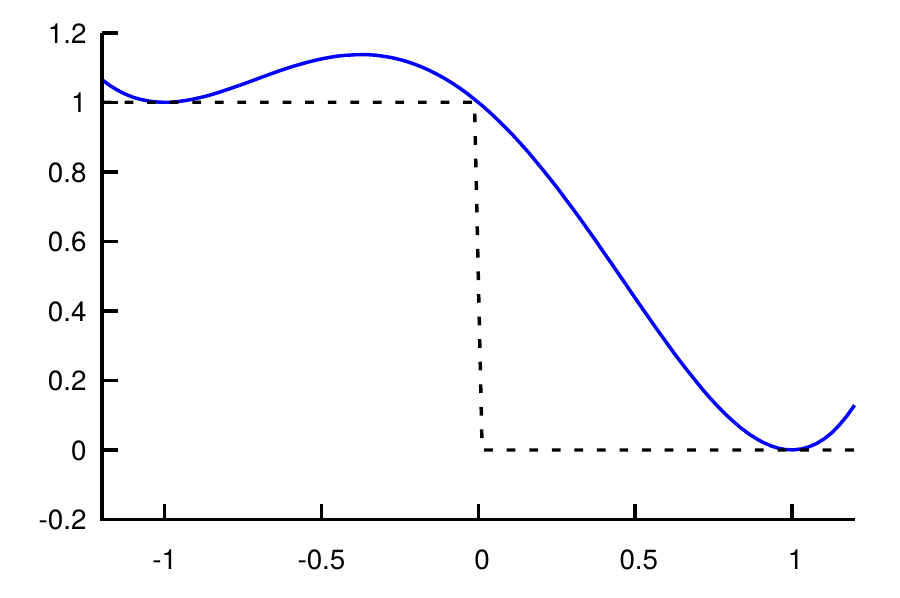}
        \caption{Hermite interpolation.}
        \label{fig:step-approx}
    \end{subfigure}
    \caption{Neville's diagram and Hermite interpolation. In (\subref{fig:Neville}), values are recursively calculated from left to right. For example, the red thick line shows that $f[-1,-1,0,1]$ is calculated by $\frac{-1/2-0}{1-(-1)}=-1/4$.}
    \end{figure}
\end{example}

\section{Moments and positive polynomials}
\label{sec:moment}
As explained in \prettyref{sec:dual-best}, moments of probability distributions arise in the dual program of best uniform approximation. 
In this section, we introduce some basic theory of moments and discuss the geometry of the moment space. 
Optimization over the moment space will prove central for both the theory and algorithm of the \emph{method of moments}, which is the topic of \prettyref{chap:gm}.

We start by introducing some notations.
The moment of a probability distribution $\mu$ is defined as 
\begin{equation}
    \label{eq:def-moment}
    m_k(\mu)=\int x^k\diff \mu(x),
\end{equation}
with $m_0(\mu) = 1$.
The $n\Th$ moment vector of a distribution $\mu$ is an $n$-tuple 
\begin{equation}
\bfm_n(\mu)=(m_1(\mu),\dots,m_n(\mu)).     
\label{eq:moment-vec}
\end{equation}
This is sometimes referred to as the truncated moment sequence which is the first $n$ terms of the full moment sequence.
The $n\Th$ moment space on $K\subseteq \reals$ is defined as
\begin{equation}
\calM_n(K)=\{\bfm_n(\mu):\mu \text{ is supported on }K\},
\label{eq:moment-space}
\end{equation}
which is a convex set since it is the convex hull of the moment curve $\{(x,x^2,\dots,x^r):x\in K\}$. 
The \emph{moment matrix} of order $n$ is a Hankel matrix of size $(n+1)\times (n+1)$ given by 
\begin{equation}
\label{eq:moment-mat}
    \bfM_n(\mu)= 
        \Expect_\mu[\bfX \bfX^\top]=
    \begin{bmatrix}
        m_0 & m_1 & \cdots & m_n \\
        m_1 & m_2 & \cdots & m_{n+1} \\
        \vdots  & \vdots  & \ddots & \vdots  \\
        m_n & m_{n+1} & \cdots & m_{2n} 
    \end{bmatrix},
\end{equation}
where $\bfX=(1,X,X^2,\ldots,X^n)^\top$; apparently moment matrices are positive semidefinite (PSD).

The moment space satisfies many geometric constraints such as the Cauchy-Schwarz (e.g.~$m_1^2 \leq m_2$) and H\"older inequalities (e.g.~$m_2^3 \leq m_3^2$).
The main results we will prove in this section is \prettyref{thm:moment-psd} that shows the finite dimensional PSD characterization of the moment space on a compact interval $K=[a,b]$, namely $\calM_n([a,b])$.
In \prettyref{chap:gm}, it is an essential step to build a fast algorithm for learning Gaussian mixtures.  
To state the result we abbreviate the Hankel matrix with entries $m_i,m_{i+1},\dots,m_j$ by 
\[
    \bfM_{i,j}=
    \begin{bmatrix}
        m_i & m_{i+1} & \cdots & m_{\frac{i+j}{2}} \\
        m_{i+1} & m_{i+2} & \cdots & m_{\frac{i+j}{2}+1} \\
        \vdots  & \vdots  & \ddots & \vdots  \\
        m_{\frac{i+j}{2}} & m_{\frac{i+j}{2}+1} & \cdots & m_{j} 
    \end{bmatrix},\quad i+j~\text{is even}.
\]
The moment matrix in \prettyref{eq:moment-mat} is $\bfM_n=\bfM_{0,2n}$.

\begin{figure}[ht]
    \centering
    \usetikzlibrary{arrows}

\begin{tikzpicture}[>=latex',scale=3.5,font=\scriptsize]  
    \pgfmathsetmacro{\m}{1.2}
    \draw[color=gray!50,fill=gray!50,domain=0:1] plot ({\x},{\x*\x});
    \draw[->](0,0)--(0,\m) node[above]{$m_2$};
    \draw[->](0,0)--(\m,0) node[right]{$m_1$};
    \draw (1,0) node[below] {$1$}--(1,1)--(0,1);  
\end{tikzpicture}
    \caption{The moment space $\calM_2([0,1])$ corresponds to the shaded region.}
    \label{fig:moment-space}
\end{figure}
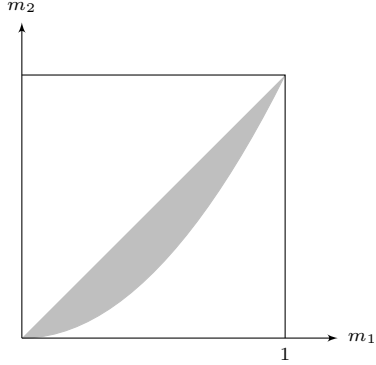

\begin{theorem}
    \label{thm:moment-psd}
    A vector $\bfm_n=(m_1,\dots,m_n)$ is in the moment space $\calM_n([a,b])$ if and only if
    \begin{equation}
        \label{eq:moment-psd}
        \begin{cases}
            ~\bfM_{0,n}\succeq 0,\quad (a+b)\bfM_{1,n-1}\succeq ab\bfM_{n-2}+\bfM_{2,n}, & \text{ $n$ even},\\
            ~b\bfM_{0,n-1}\succeq \bfM_{1,n} \succeq a\bfM_{0,n-1}, & \text{ $n$ odd}.
        \end{cases}
    \end{equation}
\end{theorem}
\begin{example}[Moment spaces on {$[0,1]$}]
    $\calM_2([0,1])$ is simply described by $m_1\ge m_2\ge 0$ and $m_2\ge m_1^2 $ as shown in \prettyref{fig:moment-space}. $\calM_3([0,1])$ is described by
    \[
        \begin{bmatrix}
            1   & m_1   \\
            m_1 & m_2  
        \end{bmatrix}
        \succeq 
        \begin{bmatrix}
            m_1 & m_2   \\
            m_2 & m_3  
        \end{bmatrix}
        \succeq 0.
    \]
    Using Sylvester's criterion (see \cite[Theorem 7.2.5]{horn-2nd}), they are equivalent to 
    \begin{align*}
        &0\le m_1\le 1, \quad m_2\ge m_3\ge 0,\\
        &m_1m_3\ge m_2^2, \quad (1-m_1)(m_2-m_3)\ge (m_1-m_2)^2,
    \end{align*}
        which can be further simplified to  $m_1^2\le m_2 \le m_1$ and $\frac{m_2^2}{m_1}\le m_3\le m_2-\frac{(m_1-m_2)^2}{1-m_1}$.
    The necessity of the above inequalities are apparent: the first two follow from the range $[0,1]$, and the last two follow from the Cauchy-Schwarz inequality. It turns out that they are also sufficient.
\end{example}

In the remaining of this section we prove \prettyref{thm:moment-psd} that reveals the intimate connection between the geometry of the moment space and the theory of positive polynomials. 
Note that a vector $(m_1,\ldots,m_n)$ can be viewed as values of a linear functional $L$ acting on monomials in $\calP_n$, such that $L(x\mapsto x^j)=m_j$ and $L(x\mapsto 1)=1$. It is a valid moment vector if there exists a representation probability measure $\mu$ such that $L(p)=\int p\diff\mu$ for every polynomial $p\in\calP_n$. Apparently, if the vector is valid, then for every positive polynomial $p\ge 0$ we have $L(p)\ge 0$. Next theorem shows that the converse also holds (see \cite[Theorem 2.14]{Rudin2006}).
\begin{theorem}[Riesz-Haviland]
    \label{thm:represent-truncated}
    Let $K\subseteq \reals$ be compact. If $L$ is a linear functional on $\calP_n$ such that $L(p)\ge 0$ for every $p\ge 0$ on $K$, then there exists a representing measure $\mu$ for $L$, \ie, $L(p)=\int p\diff\mu$ for every $p\in\calP_n$.
\end{theorem} 



The above theorems can be generalized to multiple dimensions (see \cite[pp.~17--18]{Schmudgen17} for proofs of these results); however, in general, an efficient (so that one can efficiently optimize over) characterization of positive polynomials is not known in multiple dimensions. Fortunately, for one dimension, positive polynomials can be described by sum of squares, leading to a PSD characterization of the moment space:
\begin{proposition}[{\cite[Propositions 3.1--3.3]{Schmudgen17}}]
    \label{prop:positive-sos}
		Denote by $\calS_n^2$ the set of finite sum of squares of polynomials in $\calP_n$. Then
    \begin{itemize}
        \item $p\ge 0$ on $\reals$, $\deg(p)=2n$  $\Leftrightarrow$ $p(x)=f(x)^2+g(x)^2$, $f,g\in\calP_n$.
        \item $p\ge 0$ on $[0,\infty)$, $\deg(p)=2n$  $\Leftrightarrow$ $p(x)=f(x)+xg(x)$, $f\in\calS_n^2, g\in\calS_{n-1}^2$.
        \item $p\ge 0$ on $[0,\infty)$, $\deg(p)=2n+1$  $\Leftrightarrow$ $p(x)=f(x)+xg(x)$, $f,g\in\calS_n^2$.
        \item $p\ge 0$ on $[a,b]$, $\deg(p)=2n$  $\Leftrightarrow$ $p(x)=f(x)+(b-x)(x-a)g(x)$, $f\in\calS_n^2, g\in\calS_{n-1}^2$.
        \item $p\ge 0$ on $[a,b]$, $\deg(p)=2n+1$  $\Leftrightarrow$ $p(x)=(b-x)f(x)+(x-a)g(x)$, $f,g\in\calS_n^2$.
    \end{itemize}
\end{proposition}

Using the above results, next we derive the characterization of the moment space $\calM_n([a,b])$ in \prettyref{thm:moment-psd} that was obtained in \cite[Theorem 3.1]{ST1943}. Other cases can be obtained analogously (see \cite[Part II--III]{Schmudgen17} or \cite[Chapter 3]{Lasserre2009}). 
\begin{proof}[Proof of \prettyref{thm:moment-psd}]
    If $n$ is even, by \prettyref{thm:represent-truncated} and \prettyref{prop:positive-sos}, $\bfm_n\in\calM_n([a,b])$ if and only if $L(p^2)\ge 0$ for every $p\in\calP_n$ and $L((b-x)(a-x)q^2(x))\ge 0$ for every $q\in \calP_{n-1}$. These are equivalent to $\bfM_{0,n}\succeq 0$ and $(a+b)\bfM_{1,n-1}\succeq ab\bfM_{0,n-2}+\bfM_{2,n}$, respectively. 

    If $n$ is odd, then $\bfm_n\in\calM_n([a,b])$ if and only if $L((x-a)p^2(x))\ge 0$ and $L((b-x)p^2(x))\ge 0$ for every $p\in\calP_n$. These are equivalent to $b\bfM_{0,n-1}\succeq \bfM_{1,n} \succeq a\bfM_{0,n-1}$. 
\end{proof}
\begin{remark}
    Alternatively, the characterization of the moment space in \prettyref{thm:moment-psd} can be obtained from the recursive properties of Hankel matrices; cf.~\cite{CF1991}.
\end{remark}

Moment matrices of discrete distributions satisfy more structural properties that are useful for learning finite mixture models. For instance, if $\mu$ is a \emph{$k$-atomic distribution} (supported on $k$ atoms), then its moment matrix of any order has rank at most $k$ -- see from \prettyref{eq:moment-mat}-- and is determined completely by ${\bf m}_{2k-1}(\mu)$. The number of atoms can be characterized using the determinants of moment matrices (see \cite[p.~362]{Uspensky37} or \cite[Theorem 2A]{Lindsay1989}):   
\begin{theorem}
    \label{thm:supp-detM}
    A sequence $ m_1,\dots,m_{2r} $ is the moments of a distribution with exactly $ r $ atoms if and only if 
    $ \det(\bfM_{r-1})>0 $          and $ \det(\bfM_r)=0 $.
\end{theorem}

\section{Orthogonal polynomials}
\label{sec:orthogo}


The theory of orthogonal polynomials is another classical topic with many applications, and we refer to the reader to the monographs \cite{orthogonal.poly,Gautschi2004,Ismail2005}.
In this section we recall a few definitions and constructions that will be used later for statistical problems.
\begin{definition}
    A set of functions $\{f_1,\ldots,f_n\}$ is orthogonal under the positive measure $\mu$ if 
    \[
        \int f_if_j\diff \mu=0,\quad i\ne j.
    \]
    It is orthonormal if in addition $\int f_i^2 \diff \mu=1$ for each $i$.
\end{definition}
Given a set of linear independent functions, an orthonormal set can be obtained by the Gram-Schmidt process. In the next subsection we will review some classical orthogonal polynomials under commonly used measures. 

\subsection{Classical orthogonal polynomials}

\paragraph{Chebyshev polynomials}
Recall Chebyshev polynomials (of the first kind) $T_n$ of degree $n$ defined in \prettyref{eq:cheby}.
They are orthogonal with respect to the weight function $(1-x^2)^{-1/2}$:
\begin{gather*}
    \int_{-1}^1T_n(x)T_m(x)(1-x^2)^{-1/2}\diff x=\int_0^\pi \cos(n\theta)\cos(m\theta)\diff \theta\\=
    \begin{cases}
        0,&n\ne m,\\
        \pi,&n=m=0,\\
        \pi/2,&n=m\ne 0.
    \end{cases}
\end{gather*}
An explicit formula of the Chebyshev polynomials is given by: 
\[
    T_n(x)=\frac{n}{2}\sum_{k=0}^{\Floor{n/2}}(-1)^k\frac{(n-k-1)!}{k!(n-2k)!}(2x)^{n-2k}.
\]
The approximation-theoretic properties of Chebyshev polynomials (such as \prettyref{thm:cheby-best}) will be used in \prettyref{sec:supp} where 
the optimal estimator of the unseen is constructed based on approximating step functions by Chebyshev polynomials. 

\paragraph{Hermite polynomials}
Hermite polynomials, denoted by $H_n$, are orthogonal under the standard normal distribution, \ie, for $Z\sim N(0,1)$, 
\begin{equation}
    \label{eq:Hermite-ortho}
    \Expect[H_n(Z)H_m(Z)]=\int H_n(x)H_m(x)\phi(x)\diff x=
    \begin{cases}
        n!&n=m,\\
        0&n\ne m,
    \end{cases}
\end{equation}
where $\phi(x)= \frac{1}{\sqrt{2\pi}}e^{-x^2/2}$ is the standard normal density. Hermite polynomials have the following formula
\begin{equation}
    \label{eq:Hermite1}
    H_n(x)=\Expect(x+{\bf i}Z)^n=n!\sum_{j=0}^{\floor{n/2}} \frac{(-1/2)^n }{n!(n-2j)!}x^{n-2j},
\end{equation}
where ${\bf i}=\sqrt{-1}$. Hermite polynomials are the (unique) unbiased estimate for monomials of the normal mean:
\begin{equation}
    \label{eq:Hermite-unbiased}
    \Expect[H_n(\mu+Z)]=\mu^n,\quad Z\sim N(0,1).
\end{equation}
The exponential generating function of Hermite polynomials is \cite[22.9.17]{AS64}
\begin{equation}
    \label{eq:EGF-Hermite}
    \sum_{j\ge 0}H_j(x)\frac{u^j}{j!}=\frac{\phi(x-u)}{\phi(x)}=e^{-\frac{u^2}{2}+xu}.
\end{equation}
Those properties will feature prominently in \prettyref{chap:gm} for learning Gaussian mixtures.

\paragraph{Laguerre polynomials}
The Laguerre polynomials are orthogonal under the exponential distribution (\ie, with respect to the weight function $e^{-x}$), given by the following formula:
\begin{equation}
    \label{eq:Laguerre}
    \calL_n(x)=\sum_{k=0}^n\binom{n}{k}\frac{(-x)^k}{k!}.
\end{equation}
The generalized Laguerre polynomials $\calL_n^{(k)}$ are orthogonal with respect to the weight function $x^{-k}e^x$, which can be obtained from the Rodrigues representation:
\begin{equation}
    \label{eq:Rodrigues-Lag}
    \calL_n^{(k)}(x)=\frac{x^{-k}e^x}{n!}\fracdk{}{x}{n}(e^{-x}x^{n+k})=(-1)^k\fracdk{}{k}{x}\calL_{n+k}(x), \quad k\in\naturals .
\end{equation}
Then the simple Laguerre polynomial in \prettyref{eq:Laguerre} corresponds to $ \calL_n=\calL_n^{(0)} $.
The orthogonality relation is given by
\[
    \int_0^\infty x^ke^{-x}\calL_n^{(k)}(x)\calL_m^{(k)}(x)=
    \begin{cases}
        \frac{\Gamma(n+k+1)}{n!},&n=m,\\
        0,&n\ne m.
    \end{cases}
\]
These properties will be applied to statistical lower bounds in \prettyref{sec:d-mm}, together with the following upper bound of Laguerre polynomials \cite[22.14.13]{AS64}
\begin{equation}
    \label{eq:AS64-Lag}
    |\calL_n^{(k)}(x)|\le \binom{n+k}{n}e^{x/2},\quad x\ge 0,~k\in \naturals .
\end{equation}


\paragraph{Discrete Chebyshev polynomials}
The discrete Chebyshev polynomials, denoted by $\{t_0,\ldots,t_{n-1}\}$, are orthogonal with respect to the counting measure over the discrete set $ \sth{0,1,\dots,n-1} $ with the following formula \cite[Sec.~2.8]{orthogonal.poly}: for $m=0,1,\dots,n-1$,
\begin{equation}
    t_m(x)
    \triangleq \frac{1}{m!}\Delta^mp_m(x)
    =\frac{1}{m!} \sum_{j=0}^{m}(-1)^j\binom{m}{j}p_m(x+m-j),
    \label{eq:tn-full}
\end{equation} 
where 
\begin{equation}
    p_m(x)\triangleq (x)_m (x-n)_m
    \label{eq:def-pm}
\end{equation}
and $ \Delta^m $ denotes the $ m\Th $ order forward difference. The orthogonality is given by (cf.~\cite[Sec.~2.8.2, 2.8.3]{orthogonal.poly}):
\begin{equation}
    \label{eq:ortho-discrete-cheby}
    \sum_{j=0}^{n-1}t_m(j)t_\ell(j)=
    \begin{cases}
        0,& m\ne \ell,\\
        \frac{(n+m)_{2m+1}}{2m+1},&m=\ell.
    \end{cases}
\end{equation}

The discrete Chebyshev polynomials will be used in \prettyref{sec:distinct-linear} to construct solutions for the distinct elements problem, a special case of the support size estimation problem.

\section{Gaussian quadrature}
\label{sec:GQ}

Gaussian quadrature finds a discrete approximation for a given distribution in the sense of moments, and plays a crucial role in the efficient execution of our denoised method of moments in \prettyref{chap:gm}. 
The theory of Gaussian quadrature is also an important application of interpolating polynomials, orthogonal polynomials, and moment matrices.
Given a probability measure $\mu$ supported on $K\subseteq \reals$, a $k$-point Gaussian quadrature is a $k$-atomic distribution $\mu_k=\sum_{i=1}^kw_i\delta_{x_i}$, also supported on $K$, such that, for any polynomial $P$ of degree at most $2k-1$,
\begin{equation}
    \label{eq:quadrature-goal}
    \Expect_\mu P = \Expect_{\mu_k}P=\sum_{i=1}^{k}w_iP(x_i).
\end{equation}
It is known that Gaussian quadrature always exists and is uniquely determined by ${\bf m}_{2k-1}(\mu)$ (cf.~e.g.~\cite[Section 3.6]{stoer.2002}), which, in turn, shows that any valid moment vector of order $2k-1$ can be realized by a unique $k$-atomic distribution.
A basic algorithm to compute Gaussian quadrature is \prettyref{algo:quadrature} \cite{GW1969} and many variants with improved computational efficiency and numerical stability have been proposed; cf.~\cite[Chapter 3]{Gautschi2004}.
\begin{algorithm}[ht]
    \caption{Quadrature rule}
    \label{algo:quadrature}
    \begin{algorithmic}[1]
        \REQUIRE a vector of $ 2k-1 $ moments $ (m_1,\dots,m_{2k-1}) $.
        \ENSURE nodes $ x=(x_1,\dots,x_k) $ and weights $ w=(w_1,\dots,w_k) $.
        \STATE Define the following degree-$k$ polynomial $ \Phi $ \label{line:P-ortho}
        \begin{equation}
            \label{eq:Phi-ortho}
            \Phi(x)=\det\begin{bmatrix}
                1 & m_1 & \cdots & m_k \\
                \vdots  & \vdots  & \ddots & \vdots  \\
                m_{k-1} & m_{k} & \cdots & m_{2k-1} \\
                1 & x & \cdots & x^{k} 
            \end{bmatrix}.
        \end{equation}
        \STATE Let the nodes $ (x_1,\dots,x_k) $ be the roots of the polynomial $ \Phi $.
        \STATE Let the weights $ w=(w_1,\dots,w_k)  $ be
        \begin{equation*}
            w=\begin{bmatrix}
                1 & 1 & \cdots & 1 \\
                x_1 & x_2 & \cdots & x_{k} \\
                \vdots  & \vdots  & \ddots & \vdots  \\
                x_1^{k-1} & x_2^{k-1} & \cdots & x_k^{k-1} 
            \end{bmatrix}^{-1}
            \begin{bmatrix}
                1 \\
                m_1 \\
                \vdots  \\
                m_{k-1} 
            \end{bmatrix}.
        \end{equation*}
    \end{algorithmic}
\end{algorithm}
The next result shows the correctness of \prettyref{algo:quadrature}:
\begin{theorem}
\label{thm:GQ-correct}
Let $(m_1,\dots,m_{2k-1})$ be the moments of a distribution $\mu$ supported on at least $k$ atoms. 
Then \prettyref{algo:quadrature} with the input $(m_1,\dots,m_{2k-1})$ returns a $k$-atomic distribution $\mu_k$ such that $m_j(\mu_k)=m_j(\mu)$ for $j=1,\dots,2k-1$.
\end{theorem}

It is instructive to prove the correctness of \prettyref{algo:quadrature}, in order to familiarize the readers with the properties of moment matrices and interpolating polynomials.
We first note that, if $(m_1,\dots,m_{2k-1})$ are the first $2k-1$ moments of a $k$-atomic distribution, then its atoms coincide with the zeros of the polynomial $\Phi$ defined in \prettyref{eq:Phi-ortho}, and thus $\mu$ is exactly recovered by \prettyref{algo:quadrature}.
This is shown in the next lemma, for which we provide two proofs:
\begin{lemma}
\label{lmm:roots-atoms}
        Let $\pi$ be a distribution supported on $k$ distinct atoms. 
        Then these atoms are precisely the roots of the polynomial $\Phi$ in \prettyref{eq:Phi-ortho}, where $m_i=m_i(\pi)$.
\end{lemma}
\begin{proof}[Proof 1]
    Denote the distinct atoms of $\pi$ by $x_1,\dots,x_k$.
    Denote the matrix in \prettyref{eq:Phi-ortho} by $N_k(x)$ and $\Phi(x) = \det N_k(x)$. 
    Note that its first $k$ rows coincide with those of the moment matrix $\bfM_k$. 
    It follows from \prettyref{thm:supp-detM} that the first $k$ rows are linearly independent.
    Furthermore, by definition, the first $k$ rows are in the span of $v_1,\dots,v_k$, where $v_i \triangleq (1,x_i,x_i^2,\dots,x_i^k)$ and $(v_1,\dots,v_k)$ are linearly independent (Vandermonde matrix). Therefore each $v_i$ is in span of the first $k$ rows, and thus $\Phi(x_i)=0$.
\end{proof}
\begin{proof}[Proof 2]
    Denote the distinct atoms of $\pi$ by $x_1,\dots,x_k$ and let the random variable $X\sim \pi$. 
    Then $m_i=\Expect[X^i]$. 
    The goal is to show that $\Phi(X)=0$ almost surely.
    Under the distribution $\pi$, a function $f$ is determined by $(f(x_1),\dots,f(x_k))\in\reals^k$.
    We claim that $\Expect[\Phi(X)X^r]=0$ for $0 \leq r\le k-1$. 
    Then $\Phi$ must be almost surely zero since   $1,x,\dots,x^{k-1}$ are $k$ linearly independent functions.
    To justify the claim, note that
    for all $0\leq r \leq k-1$,
    \begin{equation}
        \label{eq:Phi-orthogonal}
        \Expect[\Phi(X)X^r] = \det\begin{bmatrix}
            1 & m_1 & \cdots & m_k \\
            \vdots  & \vdots  & \ddots & \vdots  \\
            m_{k-1} & m_{k} & \cdots & m_{2k-1} \\
            m_r & m_{r+1} & \cdots & m_{r+k} 
        \end{bmatrix}=0
    \end{equation}
        where the first equality follows from expanding the determinant with respect to the last row and taking expectations, and the second follows from the existence of two identical rows when $r\le k-1$.
\end{proof}

\begin{remark}
The assumption in \prettyref{thm:GQ-correct} of $\mu$ having at least $k$ atoms is not superfluous. Suppose $\mu$ is supported on $k'<k$ atoms. Then the polynomial \eqref{eq:Phi-ortho} is $\Phi\equiv 0$ since the first $k$ rows are linearly dependent. Nevertheless, the actual support size $k'$ can be determined by the rank of the moment matrix by \prettyref{thm:supp-detM}, and then $\mu$ can be exactly recovered using \prettyref{algo:quadrature} by inputing the first $2k'-1$ moments.	
\end{remark}

Finally, we prove \prettyref{thm:GQ-correct}:
\begin{proof}[Proof of \prettyref{thm:GQ-correct}]
By \prettyref{eq:Phi-orthogonal}, $\Phi$ is orthogonal to all polynomial $P\in \calP_{k-1}$ under $\mu$. 
Then it follows from orthogonality that $\Phi$ has $k$ distinct real roots, denoted by $x_1,\ldots,x_k$.
To see this, first note that since $\Phi$ has real coefficients, its non-real roots must occur in conjugate pairs. 
Furthermore, the leading term of $\Phi$ is $ \det(\bfM_{r-1}) x^k$, where $\det(\bfM_{k-1})>0$ since $\mu$ has at least $k$ atoms (\prettyref{thm:supp-detM}). Thus,  $\Phi$ can be factorized into the form $\Phi(x)=\prod_{i=1}^{k'} (x-x_i) r(x)$, where $x_i$'s are distinct and real and $r(x)\ge 0$. Suppose ${k'}<k$, then, by orthogonality, $0=\Expect[\Phi(X)\prod_{i=1}^{k'}(X-x_i)]=\Expect[\prod_{i=1}^{k'}(X-x_i)^2r(X)]$, contradicting the fact that $X\sim \mu$ has at least $k$ atoms.

Next, for any polynomial $P$ of degree $2k-1$, we have
\begin{equation}
P(x)=\Phi(x)Q(x)+R(x),
\label{eq:PQR}
\end{equation}
where $Q,R$ are polynomials of degree at most $k-1$. 
Since $\Phi(x_i)=0$, the polynomial $R$ can be expressed by the Lagrangian interpolation formula \prettyref{eq:Lagrange}
\[
    R(x)=\sum_{i=1}^{k}R(x_i)\frac{\prod_{j\ne i}(x-x_j)}{\prod_{j\ne i}(x_i-x_j)}
    =\sum_{i=1}^{k} P(x_i)\frac{\prod_{j\ne i}(x-x_j)}{\prod_{j\ne i}(x_i-x_j)}.
\]
By orthogonality, taking expectations on both sides of \prettyref{eq:PQR} yields: for any polynomial $P$ of degree at most $2k-1$, 
\begin{equation}
\Expect[P(X)]=\sum_{i=1}^{k}w_iP(x_i),\quad w_i \triangleq \frac{\Expect\prod_{j\ne i}(X-x_j)}{\prod_{j\ne i}(x_i-x_j)}.
\label{eq:quadrule}
\end{equation}
where $X\sim \mu$.
This shows that $\sum_{i=1}^k w_i \delta_{x_k}$ defines a $k$-point Gaussian quadrature for $\mu$, provided that we can show the weights $w_i$'s define a valid probability distribution. To show that $w_i\ge 0$, recall the Lagrange basis $L_i(x)=\frac{\prod_{j\ne i}(x-x_j)}{\prod_{j\ne i}(x_i-x_j)}$ defined in \prettyref{eq:lagrange-basis}, which is a polynomial of degree $k-1$ and satisfies $L_i(x_j) = \indc{i=j}$. Then $w_i=\Expect[L_i(X)]$. Since $\deg(L_i^2) = 2k-2$, by the quadrature rule \prettyref{eq:quadrule},
\[
    0\le \Expect[L_i^2(X)]=\sum_{j=1}^{k}w_jL_i^2(x_j)=w_i.
\]
Finally, $\sum_{j=1}^{k} w_j = 1$ follows from \prettyref{eq:quadrule} by taking $P \equiv 1$.
\end{proof}

\begin{remark}[Gaussian quadrature and orthogonal polynomials]
\label{rmk:GQ-OP}	
	As shown by \prettyref{eq:Phi-orthogonal}, if $\mu$ has at least $k$ atoms, then the degree-$k$ polynomial $\Phi$ defined in \prettyref{eq:Phi-ortho} is orthogonal to all polynomials of lower degrees. Thus to solve for the Gaussian quadrature in \prettyref{algo:quadrature}, instead of evaluating the determinant in \prettyref{eq:Phi-ortho}, one can also apply the Gram-Schmidt procedure to the monomials to find the orthogonal polynomials under $\mu$, which, thanks to \prettyref{lmm:roots-atoms}, are guaranteed to have all simple real roots.
\end{remark}

\subsection{Gaussian quadrature of standard normal}
\label{sec:GQ-normal}
In this subsection we present a few properties of the Gaussian quadrature of the standard normal distribution that will be used for learning Gaussian mixture models in \prettyref{chap:gm}. 
Let $g_k$ be the $k$-point Gaussian quadrature of $N(0,1)$. By \prettyref{rmk:GQ-OP}, the atoms of $g_k$ are precisely the roots of the Hermite polynomial $H_k$ defined in \prettyref{eq:Hermite1}.
\begin{lemma}
    \label{lmm:quadrature-moments-error}
     For $j\ge 2k$, we have $m_j(g_k)\le m_j(N(0,1))$ when $j$ is even, and $m_j(g_k)= m_j(N(0,1))=0$ otherwise. In particular, $g_k$ is $1$-subgaussian.
\end{lemma}
\begin{proof}
    By the uniqueness of the Gaussian quadrature and the symmetry of the Gaussian distribution, $g_k$ is a symmetric distribution.  
    Let $\nu=N(0,1)$. If $j$ is odd, $m_j(g_k)= m_j(\nu)=0$ by symmetry. 
        If $j\ge 2k$ and $j$ is even, the conclusion follows from the integral representation of the error term of Gaussian quadrature (see, \eg,~\cite[Theorem 3.6.24]{stoer.2002}):
    \[
        m_j(\nu)- m_j(g_k) =\frac{f^{(2k)}(\xi)}{(2k)!}\int \pi_k^2(x)\diff \nu(x),
    \]
    for some $\xi\in\reals$; here $f(x)=x^{j}$, $\{x_1,\ldots,x_k\}$ is the support of $g_k$, and $\pi_k(x) \triangleq \prod_i (x-x_i)$. 
    Thus $m_j(g_k) \le m_j(\nu)$ for all even $j$ and hence $g_k$ is $1$-subgaussian \cite[Lemma 2]{subgaussian}.
\end{proof}
\begin{lemma}
    \label{lmm:quadrature-l1}
    Let $g_k$ be the $k$-point Gaussian quadrature of $N(0,1)$. Then
    \[
        \Expect_{g_k}|X|\ge (4k+2)^{-1/2},\quad k\ge 2.
    \]
\end{lemma}
\begin{proof}
    Let $G_k\sim g_k$. Note that $|G_k|\le \sqrt{4k+2}$ using the bound on the zeros of Hermite polynomials \cite[p.~129]{orthogonal.poly}. The conclusion follows from $1=\Expect[G_k^2]\le \Expect|G_k|\sqrt{4k+2}$.
\end{proof}
\begin{lemma}
    \label{lmm:quadrature-moments-Hermite}
    Let $g_k$ be the $k$-point Gaussian quadrature of $N(0,1)$.
    Then $\Expect_{g_k}[H_{j}]=0$ for $j=1,\dots, 2k-1$, and $\Expect_{g_k}[H_{2k}]=-k!$, where $H_j$ is the Hermite polynomial of degree $j$ (see \prettyref{eq:Hermite1}).
\end{lemma}
\begin{proof}
    Let $Z\sim N(0,1)$ and $G_k\sim g_k$. By orthogonality of Hermite polynomials \prettyref{eq:Hermite-ortho} we have $\Expect[H_j(Z)]=0$ for all $j\ge 1$ and thus $\Expect[H_j(G_k)]=0$ for $j=1,\dots,2k-1$. Expand $H_k^2(x)$ as
    \[
        H_k^2(x)=H_{2k}(x)+a_{2k-1}H_{2k-1}(x)+\dots+a_1H_1(x)+a_0.
    \]
    Since $G_k$ is supported on the zeros of $H_k$, we have $0=\Expect[H_k^2(G_k)]=\Expect[H_{2k}(G_k)]+a_0$. The conclusion follows from $k!=\Expect[H_k^2(Z)]=a_0$ (see \prettyref{eq:Hermite-ortho}). 
\end{proof}


\chapter{Polynomial approximation methods}
\label{chap:approx}
Property estimation is a common task in statistical inference. 
Given data from an unknown distribution, frequently the quantity of interest is a certain property of the data-generating distribution rather than the distribution itself. 
To estimate a function of a distribution, one natural idea is a two-step approach, known as the \emph{plug-in} estimate: first estimate the distribution and then substitute it into the function. 
However, this estimator is often highly biased when there is not enough data to fully recover the complicated distribution \cite{Efron82,Berkson80}.

It is natural to expect that estimating a functional is simpler (in the sense of lower sample complexity) than learning the entire distribution. As such, it is possible to accurately
estimate a functional directly even when the distribution itself is impossible to estimate. Polynomial approximation provides a powerful tool for this task.
In this chapter, we outline the recipe for implementing the polynomial approximation methods and introduce the common techniques for estimating properties of probability distributions.
Specific topics including estimating the Shannon entropy \prettyref{eq:entropy} and the support size \prettyref{eq:supportsize} will be detailed in the next two chapters. 
Throughout these chapters, the design of optimal estimator and the proof of a matching minimax lower bound both rely on the apparatus of \emph{best polynomial approximation} previously discussed in Sections~\ref{sec:uniform} and \ref{sec:dual-best}: 
\begin{itemize}
    \item For the upper bound (\prettyref{sec:functional}), we find a polynomial that approximates the property of interest and then use the unbiased estimator of the polynomial approximant.
    The bias of the resulting estimator is at most the approximation error.
    The optimal rate is obtained by carefully choosing the degree of approximation to balance the approximation error (bias) and the stochastic error (variance);
    \item For the lower bound (\prettyref{sec:d-mm}), the least favorable pair of priors can be constructed from the moment matching problem which is the dual of best polynomial approximation.
\end{itemize}


\section{Multinomial (i.i.d.) and Poisson sampling model}
\label{sec:poi}


We start by introducing a general setup for property estimation problems and relevant sampling models. 
The goal is to estimate some property $T(P)$ of the unknown distribution $P$ over an alphabet of cardinality $k$ using \iid\ observations $X_1,\dots,X_n\sim P$. 
Without loss of generality, we shall assume that the alphabet is $[k]$. 
To investigate the decision-theoretic fundamental limit \prettyref{eq:R-minimax}, we consider the minimax quadratic risk:
\begin{equation}
     \label{eq:Rkn}
    R^*(k,n) \triangleq \inf_{\hat{T}}\sup_{P \in \calM_k}\Expect( \hat{T}-T(P) )^2,    
\end{equation}
where $\hat T$ is an estimator measurable with respect to $n$ \iid observations from $P$, and $\calM_k$ denotes the set of probability distributions on $[k]$.

To perform statistical inference on the unknown distribution $P$ or any functional thereof, a sufficient statistic is the histogram $ N\triangleq (N_1,\ldots,N_k)$, where
\begin{equation}
    \label{eq:histogram}
    N_j=\sum_{i=1}^n \indc{X_i=j}
\end{equation}
records the number of occurrences of $j \in [k]$ in the sample. Then $ N\sim \Multinom(n,P) $. 
If $T(P)$ is a permutation-invariant functional of the distribution, a further sufficient statistic for estimating $T(P)$ is the histogram of the histogram $N$:
\begin{equation}
    \label{eq:fp}
    \Phi_i=\sum_{j=1}^k\indc{N_j=i},
\end{equation}
also known as histogram order statistics \cite{Paninski03}, profile \cite{OSZ04}, or fingerprint \cite{VV10}, which is the number of symbols that appear exactly $i$ times in the sample.

The \iid~sampling model is also named multinomial sampling model after the distribution of the sufficient statistic $N$. Multinomial distributions are frequently difficult to work with because of the dependency in its coordinates. To remove the dependency, a commonly used technique is the so-called \emph{Poisson sampling} where we relax the sample size $n$ from being deterministic to a Poisson random variable $n'$ with mean $n$. Under this model, we first draw the sample size $ n'\sim \Poi(n) $, then draw $ n' $ \iid\ observations from the distribution $P$. The main benefit is that now the sufficient statistics $ N_i\inddistr \Poi(np_i) $ are independent, which can significantly simplify the analysis. For more sampling models (such as sampling without replacement) and their relations, see \cite[Appendix A]{WY2016sample}.

Analogous to the minimax risk \prettyref{eq:Rkn} under \iid sampling, we define its counterpart for the Poisson sampling model:
\begin{equation}
    \label{eq:Rknt}
    \tilde{R}^*(k,n) \triangleq \inf_{\hat{T}}\sup_{P \in \calM_k}\Expect( \hat{T}-T(P) )^2,    
\end{equation}
where $\hat T$ is an estimator measurable with respect to $ N_i\inddistr \Poi(np_i) $ for $ i=1, \dots, k $. In view of the exponential tail of Poisson distributions, the Poissonized sample size is concentrated near its mean $ n $ with high probability, which guarantees that the minimax risk under Poisson sampling is provably close to that with fixed sample size. This is made precise by the following result:
\begin{theorem}
    \label{thm:poisson-sampling}
    Let $T$ be a bounded functional such that $|T(P)| \leq A$ for all $P\in\calM_k$.
	For any $\alpha>0$ and $0<\beta<1$, 
    \begin{equation}
        \label{eq:RRt}
        \tilde{R}^*(k,(1+\alpha)n)-A^2e^{-n\alpha^2/4}\le R^*(k,n)\le \frac{\tilde{R}^*(k,(1-\beta)n)}{1-\exp(-n\beta^2/2)}.
    \end{equation}
\end{theorem}

As an illustrative application of \prettyref{thm:poisson-sampling}, upon setting $\alpha=1$ and $\beta=1/2$, we have for all $n\ge 2$,
\[
\tilde{R}^*(k,2n)-A^2e^{-n/4}\le R^*(k,n)\le 5\tilde{R}^*(k,n/2). 
\]
This crude bound allows us to compare the risks under the multinomial and Poisson sampling models and show that their sample complexities are within constant factors. 
By choosing the parameters $\alpha$ and $\beta$ more carefully, one can obtain a more refined comparison of sample complexities within a factor of $1+o(1)$, which is useful for analyzing the optimal constant, for example, in \prettyref{thm:sample}.

\begin{proof}[Proof of \prettyref{thm:poisson-sampling}]
    We first prove the right inequality of \prettyref{eq:RRt}. We use the Bayesian risk as a lower bound of the minimax risk. The risk under the Poisson sampling can be expressed as
    \begin{equation*}
        \tilde{R}^*(k,(1-\beta)n)
        =\inf_{\{\hat{T}_m\}}\sup_{P\in\calM_k}\Expect[\ell(\hat{T}_{n'},T(P))],
    \end{equation*}
    where $\{\hat T_m\}$ is a sequence of estimators, $ n'\sim\Poi((1-\beta)n) $ and $ \ell(x,y)\triangleq (x-y)^2 $ is the loss function. The Bayesian risk is a lower bound of the minimax risk:
    \begin{equation}
        \tilde{R}^*(k,(1-\beta)n)
        \ge \sup_\pi\inf_{\{\hat{T}_m\}}\Expect[\ell(\hat{T}_{n'},T(P))],
        \label{eq:bayesian-ref1}
    \end{equation}
    where $\pi$ is a prior over the parameter space $\calM_k$. For any sequence of estimators $ \{\hat{T}_m\} $,
    \begin{equation*}
        \Expect[\ell(\hat{T}_{n'},T)]
        =\sum_{m\ge 0}\Expect[\ell(\hat{T}_{m},T)]\Prob[n'=m]
        \ge \sum_{m=0}^{n}\Expect[\ell(\hat{T}_{m},T)]\Prob[n'=m].
    \end{equation*}
    Taking the infimum of both sides, we obtain
    \begin{align*}
        \inf_{\{\hat{T}_m\}}\Expect[\ell(\hat{T}_{n'},T)]
        &\ge \inf_{\{\hat{T}_m\}}\sum_{m=0}^{n}\Expect[\ell(\hat{T}_{m},T)]\Prob[n'=m]\\
        &= \sum_{m=0}^{n}\inf_{\hat{T}_m}\Expect[\ell(\hat{T}_{m},T)]\Prob[n'=m].
    \end{align*}
    Note that for any fixed prior $ \pi $, the function $ m\mapsto\inf_{\hat T_m}\Expect[\ell(\hat{T}_{m},T)] $ is decreasing. Therefore
    \begin{align}
      \inf_{\{\hat{T}_m\}}\Expect[\ell(\hat{T}_{n'},T)]
      &\ge \inf_{\hat{T}_n}\Expect[\ell(\hat{T}_{n},T)]\Prob[n'\le n] \nonumber \\
      &  \ge \inf_{\hat{T}_n}\Expect[\ell(\hat{T}_{n},T)](1-\exp(n(\beta+\log(1-\beta))))\nonumber\\
      &\ge \inf_{\hat{T}_n}\Expect[\ell(\hat{T}_{n},T)](1-\exp(-n\beta^2/2)),\label{eq:bayesian-ref2}
    \end{align}
    where we used the Chernoff bound for Poisson distributions (see, \eg, \cite[Theorem 5.4]{MU06}) and the fact that $ \log(1-x)\le -x-x^2/2 $ for $x\in(0,1)$. Taking the supremum over $ \pi $ on both sides of \prettyref{eq:bayesian-ref2}, the conclusion follows from \prettyref{eq:bayesian-ref1} and the minimax theorem (cf.~\eg \cite[Theorem 46.5]{Strasser85}).

    Next we prove the left inequality of \prettyref{eq:RRt}. Recall that $0 \leq R^*(k,m) \leq R^*(k,0)$ and $m \mapsto R^*(k,m)$ is decreasing. Therefore,
    \begin{align*}
        \tilde{R}^*(k,(1+\alpha)n) 
        &\leq \sum_{m > n}R^*(k,m)\Prob[n'=m]+\sum_{0\le m \le n}R^*(k,m)\Prob[n'=m]\\
        &\le R^*(k,n)+R^*(k,0)\Prob[n'\le n]\\
        &\le R^*(k,n)+R^*(k,0)\exp(-n(\alpha-\log(1+\alpha)))\\
        &\le R^*(k,n)+A^2\exp(-n\alpha^2/4),
    \end{align*}
    where $n'\sim \Poi((1+\alpha)n)$ and we used the Chernoff bound and the fact that $\log(1+x)\le x-x^2/4$ for $0<x<1$.
\end{proof}

\section{Property estimation via polynomial approximation}
\label{sec:functional}



To construct a good estimator, the main idea is to trade bias with variance using polynomial approximation. 
The technique of polynomial approximation has been previously used for estimating non-smooth functions ($L_q$-norms) in Gaussian models \cite{INK87,LNS99,CL11} and more recently for estimating information quantities (entropy and power sums) on large discrete alphabets \cite{WY14,JVHW15}. The design principle is to  approximate the non-smooth function on a given interval using algebraic or trigonometric polynomials for which unbiased estimators exist; the degree is chosen to balance the bias (approximation error) and the variance (stochastic error).

Under the \iid{} sampling model, it is shown in \prettyref{sec:intro-estimator} 
that to estimate a functional $T(P)$ using a sample of size $n$, an unbiased estimator exists if and only if $T(P)$ is a polynomial in $P$ of degree at most $n$. 
Similarly, under the Poisson sampling model, $T(P)$ admits an unbiased estimator if and only if $T$ is real analytic. 
Consequently, there exists no unbiased estimator 
for the entropy \prettyref{eq:entropy} or the support size \prettyref{eq:supportsize}, with or without Poissonized sampling. Therefore, a natural idea is to approximate the functional by polynomials which can be estimated unbiasedly.
To be more specific, a variety of problems of theoretical and practical importance entails estimating separable functionals
\begin{equation}
T(P) = \sum_{i=1}^k f(p_i)
\label{eq:TP-sep}
\end{equation}
for some univariate function $f$ such as $f(p)=p\log \frac{1}{p}$ for the entropy \prettyref{eq:entropy} and $f(p) = \indc{p>0}$ for the support size \prettyref{eq:supportsize}. A ``meta procedure'' for constructing an estimator of $T(P)$ is to approximate $f$ by a polynomial $\tilde f$, then apply the unbiased estimator for 
\[
\tilde T(P) = \sum_{i=1}^k \tilde f(p_i), 
\]
which can be obtained using unbiased estimators of monomials given as follows:
\begin{itemize}
    \item Multinomial sampling model: 
    let $(x)_m\triangleq \frac{x!}{(x-m)!}$ denotes the falling factorial. Then,
    \begin{equation}
        \label{eq:unbiased-iid}
        \Expect[(N)_m/(n)_m]=p^n,\quad N\sim \Binom(n,p);
    \end{equation}
    \item Poisson sampling model: 
    \begin{equation}
        \label{eq:unbiased-poi}
        \Expect[(N)_m/n^m]=p^n,\quad N\sim\Poi(np).
    \end{equation}
\end{itemize}
The bias of the resulting estimator is at most the approximation error 
\begin{equation}
\label{eq:TP-ub-uniform}
|T(P)-\tilde T(P)|
\le k \cdot \sup_{p}|f(p)-\tilde f(p)|,
\end{equation}
and the variance is the error of estimating $\tilde f$ using \prettyref{eq:unbiased-iid} or \prettyref{eq:unbiased-poi}, which is determined by the polynomial degree and coefficients.

Two important parameters need to be chosen in order to carry out the above program:
\begin{itemize}
	\item Polynomial degree: 
	The bias of the estimator can be upper bounded by the approximation error, which decays as the polynomial degree increases. On the other hand, both the coefficients of the polynomial approximant and the variance of the monomial estimator 
grows (typically exponentially) with the degree as well. Thus the choice of the degree aims to strike a good bias-variance balance. 
	
	\item Approximation interval: In order to reduce the bias, typically one needs to approximate the function $f$ on a small interval near its singularity as opposed to the entire unit interval. 
    Therefore typically one uses \prettyref{eq:TP-ub-uniform} to control the bias when $p$ is small, and for larger $p$ either relies on a different bound or resort to another estimator.
\end{itemize}
The above strategy is executed for entropy estimation in \prettyref{chap:entropy}. For the support size and related problems in \prettyref{chap:unseen} the construction can be simplified 
by directly optimizing the so-called linear estimators.

Finally, we mention that, depending on the nature of the problem, the polynomial approximation problem in question can be either over a continuous interval, which is an infinite-dimensional LP as described in \prettyref{sec:uniform}, or over a discrete set of points, which is a finite-dimensional LP. 
We will encounter the latter case in discrete problems such as the distinct elements problem in \prettyref{sec:distinct} where the probabilities are known to take discrete values in the ball-urn model (see \prettyref{prop:rate-w}).

\section{Lower bounds from moment matching}
\label{sec:d-mm}

While the use of best polynomial approximation on the constructive side is admittedly natural, the fact that it also arises in the optimal lower bound is perhaps surprising. As carried out in \cite{LNS99,CL11}, the strategy is to choose two priors with matching moments up to a certain degree, which ensures the impossibility to test. The minimax lower bound is then given by the maximal separation in the expected functional values subject to the moment matching condition. 
As explained in \prettyref{sec:dual-best}, this problem is the \emph{dual} of best polynomial approximation; cf.~\prettyref{eq:dual-mm}.
In this section, we first introduce a general strategy for minimax lower bounds, and specialize it to mixture models using moment matching techniques.


A general approach for obtaining lower bounds is based on a reduction from estimation to testing. Consider the estimation of some functional $T_P=T(P)$ 
based on observations sampled from the distribution $P$, which is known to belong to some class of distributions $\calM$. For an estimator $\hat T$, suppose the loss function is a metric\footnote{
        Similar lower bounds hold if $\rho$ is not a distance but satisfies the triangle inequality within a constant factor such as the quadratic loss (cf.~e.g.~\cite{Yu97}).
} 
$\rho(\hat T,T_P)$.
If the two hypotheses 
\[
    H_0: X\sim P,\quad H_1: X\sim Q,
\] 
cannot be reliably distinguished, then any estimator suffers a loss at least proportional to the separation of their functional values $\rho(T_P,T_{Q})$. This leads to the following lower bound known as Le Cam's \emph{two-point method}. 
\begin{theorem}[{\cite[Lemma 1]{Yu97}}]
    \label{thm:lb-two-pts}
    For any two distributions $Q,Q'\in\calM$, 
    \[
            \inf_{\hat T}\sup_{P\in \calM}\Expect \rho(T_P,\hat T)\ge \frac{1}{2}\rho(T_{Q},T_{Q'})(1-\TV(Q,Q')).
    \]
\end{theorem}

\prettyref{thm:lb-two-pts} can be generalized by introducing two composite hypotheses (also known as fuzzy hypotheses in \cite{Tsybakov09}):
\[
    H_0: P\in \calM_0,\quad H_1: P\in\calM_1,
\]
where $\calM_0,\calM_1\subseteq\calM$, such that $\rho(T_P,T_{Q})\ge d$ for any $P\in\calM_0$ and $Q\in\calM_1$. 
Similarly, if no test can distinguish the above two hypotheses reliably, then any estimate suffers a maximum risk at least proportional to $d$. 
For a parametric family $\calM=\{P_\theta: \theta\in \Theta\}$ and a mixing distribution $\nu$ on $\Theta$, denote the mixture distribution by 
\begin{equation}
    \label{eq:pi-mixture-param}
    \pi_{\nu}=\int P_\theta \diff \nu(\theta),
\end{equation}
 Similar to \prettyref{thm:lb-two-pts}, we obtain the following ``two-prior'' minimax lower bound: 
\begin{theorem}
    \label{thm:lb-two-priors}
    For any two distributions $\nu_0$ and $\nu_1$ supported on $\calM_0$ and $\calM_1$, respectively,
    \[
        \inf_{\hat T}\sup_{P\in \calM}\Expect \rho(T_P,\hat T)\ge \frac{1}{2}\rho(\calM_0,\calM_1)(1-\TV(\pi_{\nu_0},\pi_{\nu_1})),
    \]
    where $\rho(\calM_0,\calM_1)\triangleq \inf\{ \rho(T_P,T_Q): P \in \calM_0, Q \in \calM_1\}$.
\end{theorem}

There are two main ingredients in Le Cam's method: (1) functional value separation; (2) indistinguishability, \ie, small total variation distance.
It turns out these two goals can be simultaneously accomplished by the \emph{dual} of uniform approximation \prettyref{eq:dual-mm}, which enables us to construct two (discrete) distributions $P$ and $Q$ supported on a closed interval $[a,b]$ such that 
\begin{equation}
    \label{eq:f-mm-best}
    \Expect_\nu[f] -\Expect_{\nu'} [f] = 2E_L(f,[a,b]),
\end{equation}
and that $\nu$ and $\nu'$ match their first $L$ moments:
\begin{equation}
    \label{eq:m-match}
    \Expect_\nu [X^j]=\Expect_{\nu'} [X^j],\quad j=0,\dots,L.
\end{equation}
For many parametric families, the statistical distance between two mixtures of the form \prettyref{eq:pi-mixture-param} can be bounded by the moment matching condition \prettyref{eq:m-match} (see Theorems~\ref{thm:chi2-gm} and \ref{thm:tv-bound} below). 
The duality relationship \prettyref{eq:f-mm-best} is essentially the reason why methods based on polynomial approximation comes naturally with a matching minimax lower bound certifying their statistical optimality.

 Let us mention the duality between statistical lower bound and upper bound in fact holds more generally beyond the paradigm of polynomial method. 
This duality view is formalized and operationalized in \cite{JN09,PW18,PSW17-colt}, leading to more general and sometimes stronger results than those obtained here from polynomial approximation and moment matching here. 

Next we elaborate on the moment-based bound on statistical distance for Gaussian mixtures and Poisson mixtures. 
The resulting statistical lower bounds using \prettyref{eq:f-mm-best} and \prettyref{eq:m-match} for specific problems will be elaborated in the subsequent chapters.

\begin{figure}[ht]
    \centering
    \begin{subfigure}[b]{.48\textwidth}
        \caption{Mixing distributions}
        \includegraphics[width=\linewidth]{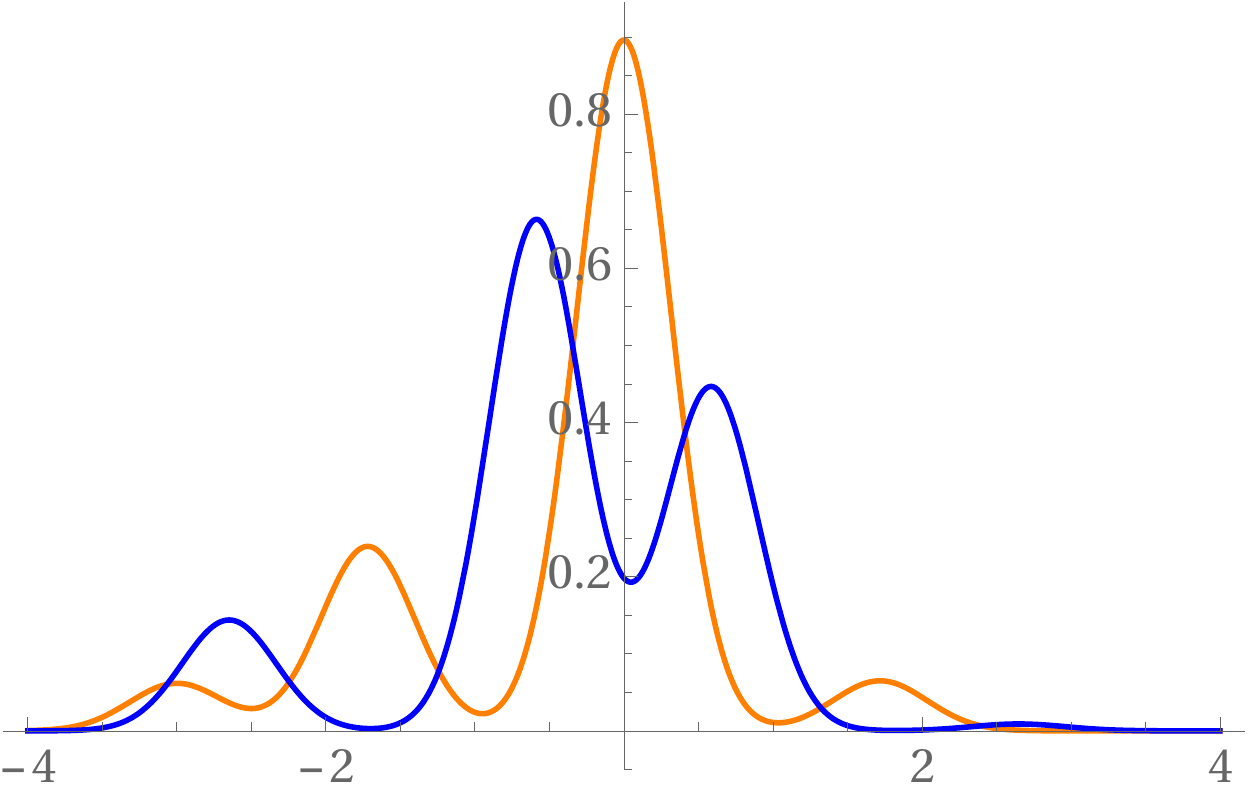}
        \label{fig:tv-latent}
    \end{subfigure}
    \begin{subfigure}[b]{.48\textwidth}
        \caption{Mixture distributions}
        \includegraphics[width=\linewidth]{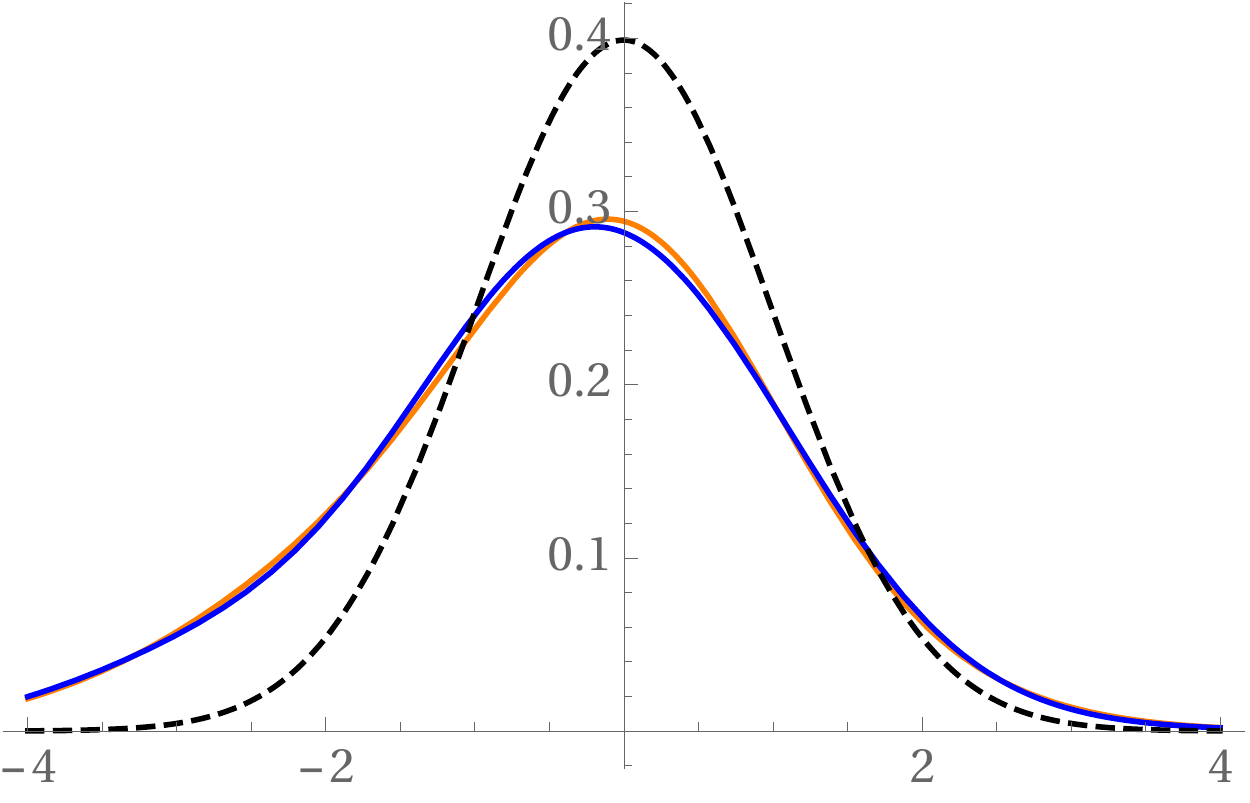}
        \label{fig:tv-mixture}
    \end{subfigure}
    \caption{Moment matching implies statistical closeness. In (\subref{fig:tv-latent}), two different mixing distributions have identical first six moments; in (\subref{fig:tv-mixture}), the mixing distributions are convolved with the standard normal (the black dashed line) and become almost indistinguishable. 
    }
    \label{fig:tv-latent-mixture}
\end{figure}

\paragraph{Gaussian mixtures. }
A Gaussian location mixture with mixing distribution $\nu$ is of the convolution form 
\[
    \pi_{\nu}=\int N(\theta,1) \diff \nu(\theta)=\nu*N(0,1).
\]
The next theorem gives an upper bound on the $\chi^2$-divergence between two Gaussian mixtures in terms of matching moments of the priors; 
see \prettyref{fig:tv-latent-mixture} for an illustration. 
Similar results have been previously obtained, for instance, by orthogonal expansion \cite{WV2010,CL11}, by Taylor expansion \cite{HP15,WY14}, and by best polynomial approximation \cite{WY15}. A more general moment comparison result in given in \prettyref{lmm:chi2-moments}.

\begin{theorem}
    \label{thm:chi2-gm}
    Suppose $m_j(\nu)=m_j(\nu')$ for $j=1,\dots,L $.
    \begin{enumerate}
        \item If $\nu$ and $\nu'$ are $\epsilon$-subgaussian for $\epsilon<1$, then
        \begin{equation}
            \label{eq:chi2-subg}
            \chi^2(\nu*N(0,1)\|\nu'*N(0,1))\le \frac{16}{\sqrt{L}}\frac{\epsilon^{2L+2}}{1-\epsilon^2}.
        \end{equation}
        \item If $\nu$ and $\nu'$ are supported on $[-\epsilon,\epsilon]$, then
        \begin{equation}
            \label{eq:chi2-bdd}
            \chi^2(\nu*N(0,1)\|\nu'*N(0,1))\le 4 e^{\epsilon^2/2} \pth{\frac{e\epsilon^2}{L+1}}^{L+1}.
        \end{equation}
    \end{enumerate}
\end{theorem}
\begin{proof}
    The proof is based on the orthogonal expansion of Gaussian mixture density via Hermite polynomials (cf.~\prettyref{sec:hermite}). 
		Since $\nu$ and $\nu'$ have the same mean, by the shift-invariance of the $\chi^2$-divergence, we assume without loss of generality that both distributions are centered.
    Let $\phi(x)=\frac{1}{\sqrt{2\pi}}e^{-\frac{x^2}{2}}$ denote the density of the standard normal distribution.
    The densities of the mixtures $\nu*N(0,1)$ and $\nu'*N(0,1)$ are 
    \begin{align}
        f(x)&=\int \phi(x-u)\diff \nu(u)=\phi(x)\sum_{j\ge 0}H_j(x)\frac{m_j(\nu)}{j!}, \label{eq:GM-orthogonal}\\
        f'(x)&=\int \phi(x-u)\diff \nu'(u)= \phi(x)\sum_{j\ge 0}H_j(x)\frac{m_j(\nu')}{j!},\label{eq:GM-orthogonal2}
    \end{align}
    respectively, where $\phi$ denotes the density of $N(0,1)$, and we used the fact in \prettyref{eq:EGF-Hermite} that $\phi(x-u)=\phi(x)\sum_{j\ge 0}H_j(x)\frac{u^j}{j!}$. Since $x\mapsto e^x$ is convex, applying Jensen's inequality yields that
    \[
        f'(x)=\phi(x)\Expect[\exp(U'x-U'^2/2)]
        \ge \phi(x)\exp(-\sigma^2/2),
    \]
    where $U'\sim\nu'$ and $\sigma$ denote the variance of $\nu'$.
    Consequently, applying the moment matching condition yields that
    \begin{align}
        &\phantom{{}={}}\chi^2(f\|f') = \int \frac{(f(x)-f'(x))^2}{f'(x)}\diff x\nonumber\\
        &\leq e^{\frac{\sigma^2}{2}} \expect{\pth{\sum_{j\ge L+1}H_j(Z)\frac{m_j(\nu)-m_j(\nu')}{j!}}^2}\nonumber\\       
        &=  e^{\frac{\sigma^2}{2}} \sum_{j\ge L+1}\frac{(m_j(\nu)-m_j(\nu'))^2}{j!},\label{eq:chi2-matching}
    \end{align}
    where $Z\sim N(0,1)$ and the last step follows from the orthogonality property of Hermite polynomials in \prettyref{eq:Hermite-ortho}.

    If $\nu$ and $\nu'$ are $\epsilon$-subgaussian, then $\sigma\le \epsilon$, and $|m_j(\nu)|, |m_j(\nu')|\le 2(\epsilon\sqrt{j/e})^j$ \cite{subgaussian}. Applying \eqref{eq:chi2-matching} yields that 
    \[
        \chi^2(f\|f')
        \le e^{\epsilon^2/2}\sum_{j\ge L+1}\frac{16\epsilon^{2j}}{\sqrt{2\pi j}}
        \le \frac{16}{\sqrt{L}}\sum_{j\ge L+1}\epsilon^{2j},
    \]
    where we used Stirling's approximation $j!>\sqrt{2\pi j}(j/e)^j$, and $e^{\epsilon^2/2}\le \sqrt{e}\le \sqrt{2\pi}$.
		
		If $\nu$ and $\nu'$ are supported on $[-\epsilon,\epsilon]$, then we have $\sigma \leq \epsilon$ and $|m_j(\nu)|, |m_j(\nu')|\le \epsilon^j$. 
		From \prettyref{eq:chi2-matching} we get 
		\begin{align*}
		\chi^2(f\|f')
		\leq 4 e^{\frac{\epsilon^2}{2}} \sum_{j\ge L+1} \frac{\epsilon^{2j}}{j!} 
		= 4 e^{-\frac{\epsilon^2}{2}} \prob{\Poi(\epsilon^2) \geq L+1} 
		\leq 4 e^{\frac{\epsilon^2}{2}} \pth{ \frac{e \epsilon^2}{L+1} }^{L+1},
		\end{align*}
		where the last inequality follows from the Chernoff bound for Poisson distribution (see e.g.~\cite[Thm 4.4]{MU06}).		
\end{proof}

\begin{remark}
    \prettyref{thm:chi2-gm} is tight in the subgaussian case. Specifically, when L is odd, there exists a pair of $\epsilon$-subgaussian distributions $\nu$ and $\nu'$ such that ${\bf m}_L(\nu)={\bf m}_L(\nu')$, while $\chi^2(\nu*N(0,1)\|\nu'*N(0,1))\ge \Omega_L(\epsilon^{2L+2})$. They can be constructed using Gaussian quadrature introduced in \prettyref{sec:GQ}. 
    To this end, let $L=2k-1$ and we set $\nu=N(0,\epsilon^2)$; 
    let $Z_k$ be distributed according to the $k$-point Gaussian quadrature of the standard normal distribution in \prettyref{sec:GQ-normal}, and $\nu'$ denote the law of $\epsilon Z_k$.
    Then ${\bf m}_{2k-1}(\nu)={\bf m}_{2k-1}(\nu')$, and $\nu'$ is also $\epsilon$-subgaussian (see \prettyref{lmm:quadrature-moments-error}). It is shown in \cite[(54)]{WV2010} that
    \[
        \chi^2(\nu'*N(0,1)\|\nu*N(0,1))=\sum_{j \ge 2k}\frac{1}{j!}\pth{\frac{\epsilon^2}{1+\epsilon^2}}^j|\Expect[H_j(Z_k)]|^2.
    \]
    where 
    $H_j$'s are the Hermite polynomials defined in \prettyref{eq:Hermite1}. 
	By \prettyref{lmm:quadrature-moments-Hermite}, we have $\Expect[H_{2k}(Z_k)]=-k!$ and thus for any $\epsilon<1$,
    \[
        \chi^2(\nu'*N(0,1)\|\nu*N(0,1))\ge \frac{(k!)^2}{(2k)!}\pth{\frac{\epsilon^2}{1+\epsilon^2}}^{2k}\ge (\Omega(\epsilon))^{4k}.
    \]
\end{remark}

\begin{remark}[Gaussian scale mixtures]
Similar results can be obtained for Gaussian scale mixtures of the form
\begin{equation}
    \label{eq:Gauss-scale}
    \sum_{i=1}^{k}w_iN(0,\sigma_i^2).
\end{equation}
In the case where $\sigma_i$'s are close to a fixed value, moment comparison of Gaussian scale mixtures can be established by applying \prettyref{eq:chi2-subg} with $\nu$ and $\nu'$ themselves being Gaussian scale mixtures with variances $O(\epsilon)$. See \prettyref{eq:gaussian-scale-lb} in \prettyref{sec:discuss} for details and applications.
\end{remark}


\paragraph{Poisson mixtures.}
Now we consider a Poisson mixture
\[
    \pi_{\nu} =\int \Poi(\lambda)\diff \nu(\lambda),
\]
with mixing distribution $\nu$ on $\reals_+$. This defines a discrete distribution over $\integers_+$. 
The following result bounds the total variation ($L_1$-distance between probability mass functions) between Poisson mixtures based on moment matching. 
As opposed to the proofs of \prettyref{thm:chi2-gm} (and \prettyref{lmm:chi2-moments}) which are based on orthogonal expansion, the proof here is by polynomial approximation.

\begin{theorem}[Poisson mixtures]
    \label{thm:tv-bound}
    Let $\nu$ and $\nu'$ be two distributions supported on $ [0,\Lambda] $.
    If $m_j(\nu)=m_j(\nu')$ for $j=1,\dots,L $, then
    \begin{equation}
        \TV(\pi_\nu,\pi_{\nu'})
        \le \frac{(\Lambda/2)^{L+1}}{(L+1)!}\pth{2+2^{\Lambda/2-L}+2^{\Lambda/(2\log 2)-L}}.
        \label{eq:tv-bound}
    \end{equation}
    In particular, $ \TV(\pi_\nu,\pi_{\nu'})\le (\frac{e\Lambda}{2L})^{L} $. Moreover, if $ L>\frac{e}{2}\Lambda $, then
    \[
        \TV(\pi_\nu,\pi_{\nu'})\le \frac{(\Lambda/2)^{L+1}}{(L+1)!}(2+o(1)), \quad \Lambda \diverge.
    \]
\end{theorem}
\begin{proof}
    Let $f_j(x)\triangleq \frac{e^{-x}x^j}{j!}$.
    For any polynomial $ P_{L,j} $ of degree $ L $ that uniformly approximates $ f_j $ over the interval $ [0,\Lambda] $ with approximation error $ e_j\triangleq\max_{x\in[0,\Lambda]}|f_j(x)-P_{L,j}(x)| $, we have, by moment matching, $\Expect_\nu [P_{L,j}] = \Expect_{\nu'} [P_{L,j}]$. Therefore
    \begin{align}
      &\phantom{{}={}}\TV(\pi_\nu,\pi_{\nu'})       
      =  \frac{1}{2}\sum_{j=0}^{\infty}|\Expect_\nu[f_j]-\Expect_{\nu'} [f_j]|\nonumber\\
      &\le \frac{1}{2}\sum_{j=0}^{\infty}|\Expect_{\nu} [f_j-P_{L,j}]|+|\Expect_{\nu'} [f_j-P_{L,j}]| 
      \le \sum_{j=0}^{\infty} e_j.\label{eq:tv-ub-approx}
    \end{align}
    Next we will bound the approximation error by taking $P_{j,L}$ to be Chebyshev interpolation polynomial \prettyref{eq:interp-cheby}.

    If $ j\le L+1 $, applying \prettyref{eq:Rodrigues-Lag} and the general Leibniz rule for derivatives yields that
    \begin{equation}
        f_j^{(L+1)}(x)
        =\fracdk{}{x}{j'}(\calL_j(x)e^{-x})
        =(-1)^{j'}e^{-x}\sum_{m=0}^{j'\wedge j}\binom{j'}{m}\calL_{j-m}^{(m)}(x).
        \label{eq:fj-diff}
    \end{equation}
    where $\calL_j$ is the degree-$j$ Laguerre polynomial, and $j'\triangleq |L+1-j|$. 
    Applying \prettyref{eq:AS64-Lag} yields that
    \begin{equation*}
        \abs{f_j^{(L+1)}(x)}
        \le e^{-x}\sum_{m=0}^{j'\wedge j}\binom{j'}{m}\binom{j}{j-m}e^{x/2}
        = e^{-x/2}\binom{L+1}{j}.
    \end{equation*}
    Therefore $ \max_{x\in [0,\Lambda]}|f_j^{(L+1)}(x)|\le \binom{L+1}{j} $ when $ j\le L+1 $.\footnote{This is in fact an equality. In view of \prettyref{eq:fj-diff} and the fact that $L_{j-m}^{(m)}(0)=\binom{j}{j-m}$ \cite[22.3]{AS64}, we have $ |f_j^{(L+1)}(0)|=\sum_{m}\binom{L+1-j}{m}\binom{j}{j-m}=\binom{L+1}{j} $.}
    Then, we apply \prettyref{eq:interp-cheby} and obtain that
    \begin{equation}
        \sum_{j=0}^{L+1} e_j
        \le \sum_{j=0}^{L+1}\frac{\binom{L+1}{j}(\Lambda/2)^{L+1}}{2^L(L+1)!}
        = \frac{2(\Lambda/2)^{L+1}}{(L+1)!}.
        \label{eq:approx-bd-jsmall}
    \end{equation}

    If $ j\ge L+2 $, again applying \prettyref{eq:Rodrigues-Lag} and \prettyref{eq:AS64-Lag} yields that
    \begin{align*}
        f_j^{(L+1)}(x)
        &=\frac{(L+1)!}{j!}x^{j'}e^{-x}\calL_{L+1}^{j'}(x)\\
        &\le \frac{(L+1)!}{j!}x^{j'}e^{-x}\binom{j}{L+1}e^{x/2}
        = \frac{1}{j'!}e^{-x/2}x^{j'}.
    \end{align*}
    The maximum of $e^{-x/2}x^{j'}$ for $x\in[0,\Lambda]$ occurs at $ x=(2j')\wedge \Lambda $.
    Therefore, for any $x\in [0,\Lambda]$,
    \begin{equation*}
        |f_j^{(L+1)}(x)|\le
        \begin{cases}
            \frac{1}{j'!}(2j'/e)^{j'},& L+1\le j\le L+1+\Lambda/2,\\
            \frac{1}{j'!}e^{-\Lambda/2}\Lambda^{j'},& j\ge L+1+\Lambda/2.
        \end{cases}
    \end{equation*}
    Then, applying \prettyref{eq:interp-cheby} and Stirling's approximation $ n!>\sqrt{2\pi n}(\frac{n}{e})^{n} $, we have
    \begin{align}
        \sum_{\substack{j\ge L+2\\ j<L+1+\Lambda/2}}e_j
        & \le \frac{(\Lambda/2)^{L+1}}{2^L(L+1)!}\sum_{1\le j<\frac{\Lambda}{2}}\frac{2^{j}}{\sqrt{2\pi j}}\le \frac{(\Lambda/2)^{L+1}2^{\Lambda/2}}{2^L(L+1)!},\label{eq:approx-bd-jmid}\\
        \sum_{j\ge L+1+\Lambda/2}e_j
        & \le \frac{(\Lambda/2)^{L+1}e^{-\Lambda/2}}{2^{L}(L+1)!}\sum_{j\ge \Lambda/2}\frac{\Lambda^{j}}{j!}\le \frac{(\Lambda/2)^{L+1}e^{\Lambda/2}}{2^{L}(L+1)!}.\label{eq:approx-bd-jbig}
    \end{align}

    Assembling the three ranges of summations in \prettyref{eq:approx-bd-jsmall}--\prettyref{eq:approx-bd-jbig} in the total variation bound \prettyref{eq:tv-ub-approx}, we obtain
    \begin{equation*}
        \TV(\pi_\nu,\pi_{\nu'})
        \le \frac{(\Lambda/2)^{L+1}}{(L+1)!}\pth{2+2^{\Lambda/2-L}+2^{\Lambda/(2\log 2)-L}}.
    \end{equation*}
    Finally, applying Stirling's approximation $ n!>\sqrt{2\pi n}(\frac{n}{e})^{n} $, we conclude that $ \TV(\pi_\nu,\pi_{\nu'})\le (\frac{e\Lambda}{2L})^{L} $.
    If $ L>\frac{e}{2}\Lambda>\frac{\Lambda}{2\log 2}>\frac{\Lambda}{2} $, then $ 2^{\Lambda/2-L}+2^{\Lambda/(2\log 2)-L}=o(1) $.    
\end{proof}


\chapter{Entropy estimation}
\label{chap:entropy}

In this chapter, we apply the method of polynomial approximation to the problem of entropy estimation.
The Shannon entropy \cite{shannon} of a discrete distribution $P$ is defined as 
\[
    H(P)=\sum_i p_i\log\frac{1}{p_i}.
\]
Entropy estimation has found numerous applications across various fields, such as psychology \cite{Attneave1959}, neuroscience \cite{spikes-book}, physics \cite{Vinck12}, telecommunication \cite{PW96}, biomedical research \cite{porta2001entropy}, etc. Furthermore, it serves as the building block for estimating other information measures expressible in terms of entropy, such as mutual information and directed information, which are instrumental in machine learning applications such as learning graphical models \cite{CL68,QKC13,jiao2013di,Bresler15}. However, accurate estimation of the Shannon entropy is a difficult task for large alphabet due to many symbols being unobserved.

In \prettyref{sec:plug-in}, we first discuss the maximum likelihood estimate, which is also known as the \emph{empirical entropy}. 
As introduced in \prettyref{sec:functional}, this is the plug-in estimator in functional estimation, for which we substitute the estimated distribution into the functional. Despite being a natural approach, this method suffers from a large bias when the sample size is small and can be highly suboptimal. 

Next, in \prettyref{sec:h-esti}, we construct an improved estimator by applying the polynomial approximation method in \prettyref{chap:approx}. To establish its statistical optimality, we consider the minimax quadratic risk:
\begin{equation}
    R_H^*(k,n)\triangleq \inf_{\hat{H}}\sup_{P \in \calM_k}\Expect_P[( \hat{H}-H(P) )^2],     
\end{equation}
where $\calM_k$ denotes the set of probability distributions on $[k]$, and $\hat H$ is an estimator measurable with respect to $n$ independent observations sampled from $P$. 
%
%
%
The main result is the following characterization of the \emph{minimax rate}, \ie, approximation of the minimax risk $R_H^*(k,n)$ within universal constant factors:
\begin{theorem}
    \label{thm:entropy}
    If $n \gtrsim \frac{k}{\log k}$, then
    \begin{equation}
        \label{eq:h-main}
        R_H^*(k,n) \asymp \pth{\frac{k }{n \log k}}^2 + \frac{\log^2 k}{n},
    \end{equation}
		which can be achieved by a linear-time estimator.
    If $n \lesssim \frac{k}{\log k}$, there exists no consistent estimators, \ie, $R_H^*(k,n) \gtrsim 1$.
\end{theorem}

The upper bound is achieved by the polynomial approximation method in \prettyref{sec:h-esti}, and the lower bound is obtained by a moment matching problem in \prettyref{sec:h-lb}.
To interpret the minimax rate \prettyref{eq:h-main}, we note that the second term corresponds to the classical ``parametric'' term proportional to $\frac{1}{n}$, which is governed by the variance and the central limit theorem (CLT). The first term corresponds to the squared bias, which is the dominating term for large domains. Note that $R_H^*(k,n) \asymp (\frac{k }{n \log k})^2$ if and only if $n \lesssim \frac{k^2}{\log^4 k}$, where the bias dominates. As a consequence, \prettyref{thm:entropy} implies that to estimate the entropy within $\epsilon$ bits with probability, say 0.9, the minimal sample size is given by
\begin{equation}
    \label{eq:H-sample-complexity}
    n \asymp \frac{\log^2 k}{\epsilon^2} \vee \frac{k}{\epsilon \log k}.    
\end{equation}

\begin{figure}[ht]
    \centering
    \begin{tikzpicture}[font=\tiny]

        \shade[top color = white, bottom color=black!40!white] (0,0) rectangle (2,0.8);
        \shade[top color = white, bottom color=red!40!white] (2,0) rectangle (7,0.8);
        \shade[top color = white, bottom color=green!40!white] (7,0) rectangle (11,0.8);
        \shade[bottom color = white, top color=black!40!white] (0,0) rectangle (3,-0.8);
        \shade[bottom color = white, top color=red!40!white] (3,0) rectangle (9,-0.8);
        \shade[bottom color = white, top color=green!40!white] (9,0) rectangle (11,-0.8);
        
        \draw [->](0,0) -- (11,0);
        \draw (11.1,0.2) node {\footnotesize $ n $};

        \draw (5.5,0.9) node {\footnotesize MSE of the optimal estimator};
        \draw (2,5pt)-- (2,0pt);
        \draw (2,0.5) node {$ \frac{k}{\log k} $};
        \draw (7,5pt)-- (7,0pt);
        \draw (7,0.5) node {$ \frac{k^2}{\log^4 k} $};
        \draw (1,0.3) node {Inconsistent};
        \draw (4.5,0.3) node {Bias dominates};
        \draw (9,0.3) node {Variance dominates};
        
        \draw (5.5,-0.9) node {\footnotesize MSE of the empirical entropy};
        \draw (3,0pt)-- (3,-5pt);
        \draw (3,-0.5) node {$ k $};
        \draw (9,0pt)-- (9,-5pt);
        \draw (9,-0.5) node {$ \frac{k^2}{\log^2 k} $};
        \draw (1.5,-0.3) node {Inconsistent};
        \draw (6,-0.3) node {Bias dominates};
        \draw (10.2,-0.3) node {Variance dominates};

\end{tikzpicture}
    \caption{Classification and comparison of optimal entropy estimator and the empirical entropy. }
    \label{fig:H-compare}
\end{figure}
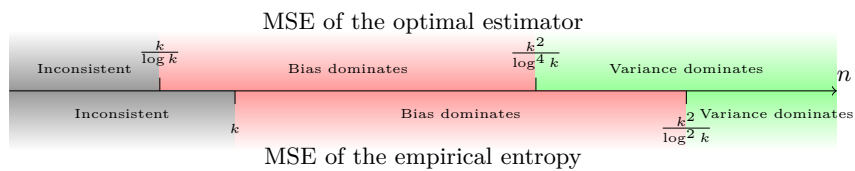

In comparison with \prettyref{thm:entropy}, the worst-case quadratic risk of the empirical entropy is analyzed in \prettyref{sec:plug-in} and is given by
\begin{theorem}
    \label{thm:emprical-entropy}
    Let $\hat P_n$ denote the empirical distribution of $n$ \iid~samples. If $n \gtrsim k$, then
    \begin{equation}
        \label{eq:Rplug-rate}
        \sup_{P:S(P)\le k}\Expect(H-H(\hat P_n))^2 \asymp \pth{\frac{k }{n}}^2 + \frac{\log^2 k}{n},
    \end{equation}
    If $n \lesssim k$, then $H(\hat P_n)$ is inconsistent, \ie, the left-hand side of \prettyref{eq:Rplug-rate} is $\Omega(1)$.
\end{theorem}
Again, the first and second terms in the risk again correspond to the squared bias and variance respectively. Comparing \prettyref{eq:h-main} and \prettyref{eq:Rplug-rate}, we reach the following verdict on the plug-in estimator: Empirical entropy is rate-optimal, \ie, achieving the minimax rate, if and only if we are in the ``data-rich'' regime of $n = \Omega(\frac{k^2}{\log^2 k})$. In the ``data-starved'' regime of $n = o\big(\frac{k^2}{\log^2 k}\big)$, empirical entropy is strictly rate-suboptimal. The comparison between the optimal estimator and the empirical entropy is illustrated in \prettyref{fig:H-compare}. 

\section{Empirical entropy and Bernstein polynomials}
\label{sec:plug-in}
Given $n$ independent observations $X_1,\dots,X_n$ from a discrete distribution $P=(p_1,\dots,p_k)$, the maximum likelihood estimate of the distribution is the empirical distribution
\[
    \hat P_n= (\hat p_1,\dots,\hat p_k)
\]
with $\hat p_i=N_i/n$, where $N_i$ records the number of occurrences of the $i\Th$ symbol. We refer to $N\triangleq (N_1,\dots,N_k)$ as the \emph{histogram}; both $N$ and $\hat P_n$ are sufficient statistics.
Then the empirical entropy is 
\begin{equation}
    \label{eq:H-emp}
    H(\hat P_n)=\sum_i \hat p_i\log \frac{1}{\hat p_i}.
\end{equation}
In this section, we will analyze the bias and the variance of $H(\hat P_n)$ in Propositions~\ref{prop:h-plug-bias} and \ref{prop:h-plug-var}, respectively.

\subsection{Bias of the empirical entropy}
\begin{figure}[ht]
    \centering
    \begin{subfigure}[b]{.48\textwidth}
        \includegraphics[width=\linewidth]{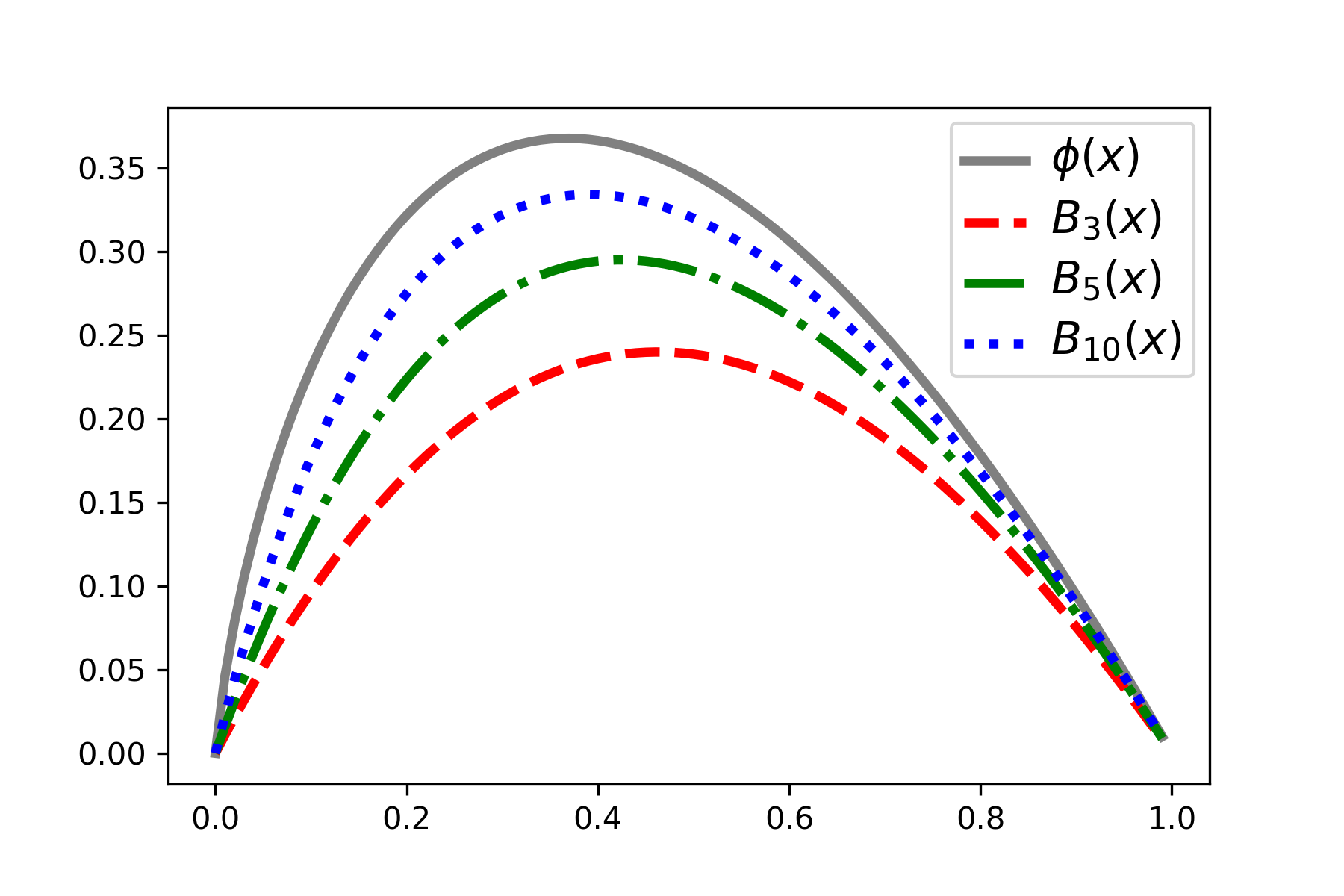}
        \caption{Bernstein polynomials}
        \label{fig:Bn}
    \end{subfigure}
    \hfill
    \begin{subfigure}[b]{.48\textwidth}
        \includegraphics[width=\linewidth]{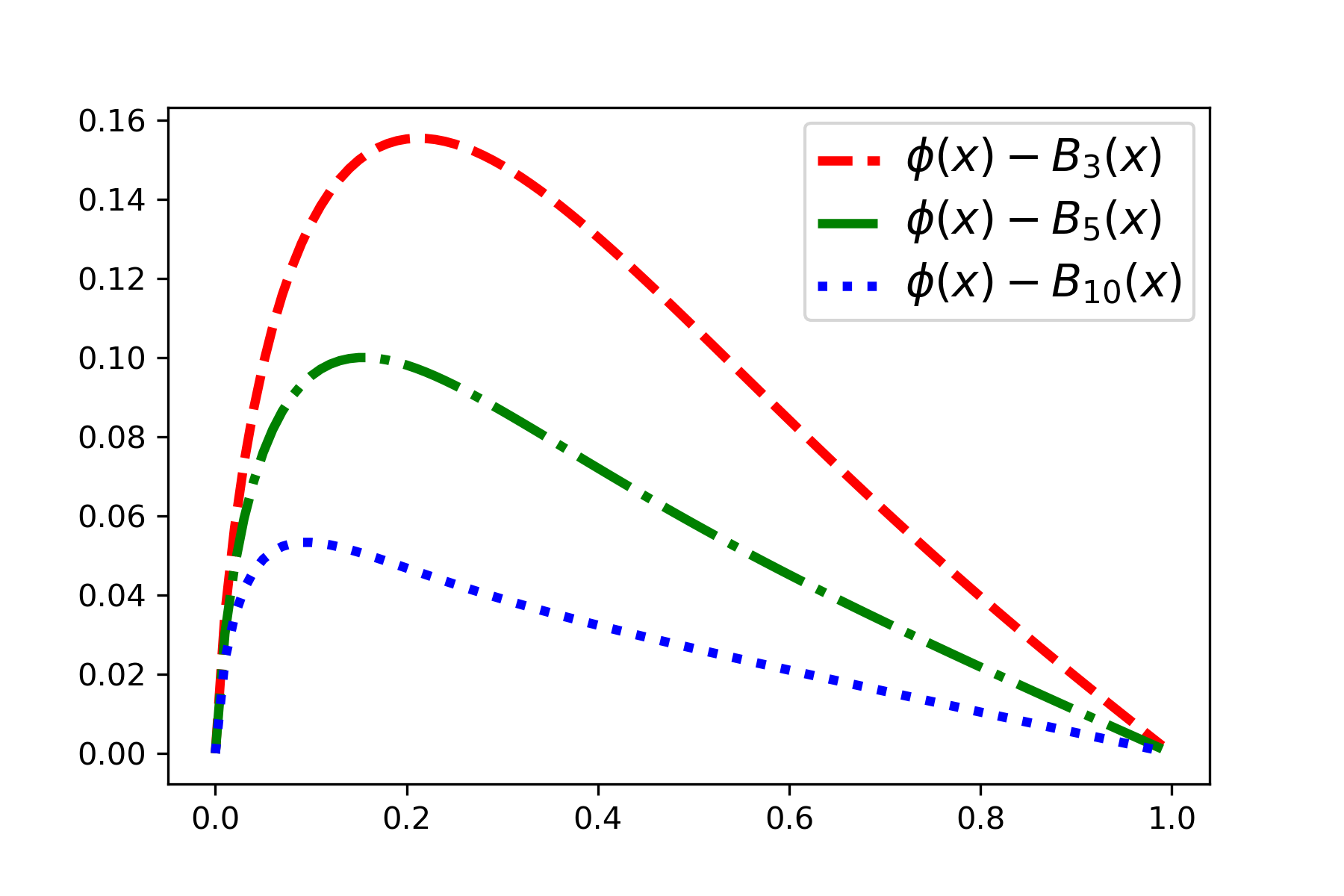}
        \caption{Approximation errors}
        \label{fig:Bn-error} 
    \end{subfigure}
    \caption{Illustration of Bernstein polynomial approximation of $\phi$ of degree $3,5$ and $10$. (\subref{fig:Bn}) shows the actual Bernstein polynomials, and (\subref{fig:Bn-error}) shows the errors of approximation.}
    \label{fig:Bernstein}
\end{figure}

It is straightforward to connect the bias of empirical entropy to the Bernstein polynomial approximation \prettyref{eq:Bernstein} by
\begin{align}
  \Expect[H(\hat P_n)]-H(P)
  &=\sum_i \pth{\sum_{j=0}^{n}\phi(j/n)\binom{n}{j}p_i^j(1-p_i)^{n-j} -\phi(p_i)}\nonumber\\
  &=\sum_i (B_n(p_i)-\phi(p_i)),\label{eq:bias-Bn}
\end{align}
where $\phi(x)=x\log\frac{1}{x}$ and $B_n$ is the Bernstein polynomial of degree $n$ to approximate $\phi$ using the formula \prettyref{eq:Bernstein}. See an illustration of Bernstein approximation in \prettyref{fig:Bernstein}.
We shall next derive several results on the bias of the empirical entropy using the Bernstein approximation. 

\begin{lemma}
    If $f$ is convex on $[0,1]$, then the corresponding Bernstein polynomial approximation \prettyref{eq:Bernstein} satisfies the following inequalities:
    \[
        B_n(x) \ge B_{n+1}(x) \ge f(x).
    \]
    The inequalities are strict if $f$ is strictly convex.
\end{lemma}
\begin{proof}
    Applying the formula of Bernstein polynomials \prettyref{eq:Bernstein}, we can calculate that (see \cite[pp. 309--310]{DL93})
    \[
        B_n(x)-B_{n+1}(x)=\frac{x(1-x)}{n(n+1)}\sum_{k=0}^{n-1}f\qth{\frac{k}{n},\frac{k+1}{n+1},\frac{k+1}{n}}\binom{n-1}{k}x^k(1-x)^{n-k},
    \]
    where $f[\frac{k}{n},\frac{k+1}{n+1},\frac{k+1}{n}]$ is the divided difference; cf.~\prettyref{eq:div-diff-f}. 
		By  \prettyref{eq:div-GH}, this divided difference is nonnegative when $f$ is convex.
\end{proof}
Note that $\phi$ is strictly concave on $[0,1]$. In this case we have 
\begin{equation}
    \label{eq:Bn-less}
    B_n(x)<B_{n+1}(x)<\phi(x), \quad 0<x<1.
\end{equation}
See \prettyref{fig:Bernstein} for an illustration. We conclude from \prettyref{eq:bias-Bn} that the empirical entropy is always \emph{underbiased}, and the magnitude of the bias is strictly decreasing as the sample size increases \cite{Paninski03}.

Bernstein approximation has the following asymptotic formula: 
\begin{lemma}
    Fix $x\in[0,1]$. If $f$ is bounded, differentiable in a neighborhood of $x$, and $f''(x)$ exists, then
    \begin{equation}
        \label{eq:Bn-asymp}
        \lim_{n\diverge}n(B_n(x)-f(x))=\frac{x(1-x)}{2}f''(x).
    \end{equation}
\end{lemma}
\begin{proof}
    By Taylor's expansion, for $t\in[0,1]$
    \[
        f(t)=f(x)+f'(x)(t-x)+\frac{f''(x)}{2}(t-x)^2+r(t-x)(t-x)^2,
    \]
    where in the remainder $r(y)$ is bounded and vanishes with $y$. Note that $B_n(x)=\Expect f(\hat p)$ where $\hat p=N/n$ and $N\sim \Binom(n,x)$. Then 
    \[
        B_n(x)-f(x)=\frac{x(1-x)}{2n}f''(x)+\Expect[r(\hat p-x)(\hat p-x)^2]. 
    \]
    Next we show the last term is $o(1/n)$. 
    Let $M(\delta)\triangleq \sup_{|x|\le \delta}|r(x)|$ that is bounded and vanishes with $\delta$.
    Then
    \begin{align*}
    \Expect[r(\hat p-x)(\hat p-x)^2]
    &\le M(\delta)\Expect[(\hat p-x)^2]+M(1) \Prob[|\hat p-x|>\delta]\\
    &\le M(\delta)\frac{1}{n}+2M(1)\exp(-n\delta^2), 
    \end{align*}
    where the last inequality used the Chernoff bound. 
    The conclusion follows by letting $\delta=n^{-1/4}$.
\end{proof}
In entropy estimation, $\phi''(x)=-1/x$. By using \prettyref{eq:bias-Bn}, for a \emph{fixed} distribution $P$, the asymptotic bias of the empirical entropy as $n$ diverges is given by
\begin{equation}
    \label{eq:hatH-asymp-bias}
    \Expect[H(\hat P_n)]-H(P)=\sum_i \frac{p_i-1}{2n}(1+o(1))
    =\frac{1-S(P)}{2n}(1+o(1)),
\end{equation}
where $S(P)$ denotes the support size of $P$. This asymptotic formula motivates the well-known bias reduction to the empirical entropy, named the Miller-Madow estimator \cite{Miller55}:
\begin{equation}
    \label{eq:hatH-MM}
    \hat H_{\rm MM}=\hat H_{\rm plug}+\frac{\hat S-1}{2n},
\end{equation}
where $\hat S$ is the number of observed distinct symbols. For higher order asymptotic expansions of the bias, as well as various types of bias reduction, see \cite{Harris75}. The formula \prettyref{eq:hatH-asymp-bias} also holds when the fixed distribution assumption is relaxed to $n \min_i p_i\diverge$ \cite[Theorem 5]{Paninski03}.

\begin{figure}[ht]
    \centering
    \begin{subfigure}[b]{.48\textwidth}
        \includegraphics[width=\linewidth]{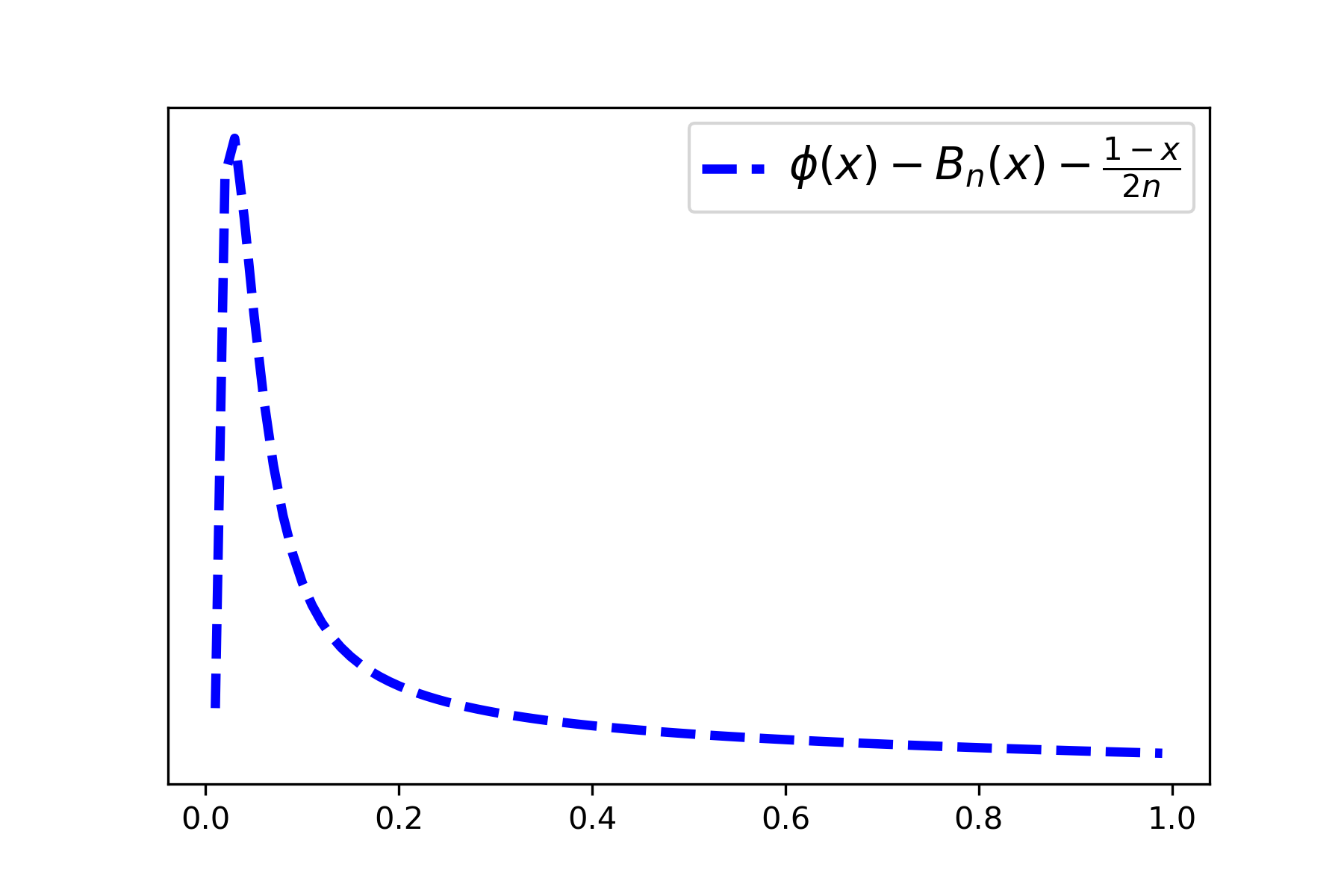}
        \caption{Approximation error}
        \label{fig:Bn-MM-diff}
    \end{subfigure}
    \hfill
    \begin{subfigure}[b]{.48\textwidth}
        \includegraphics[width=\linewidth]{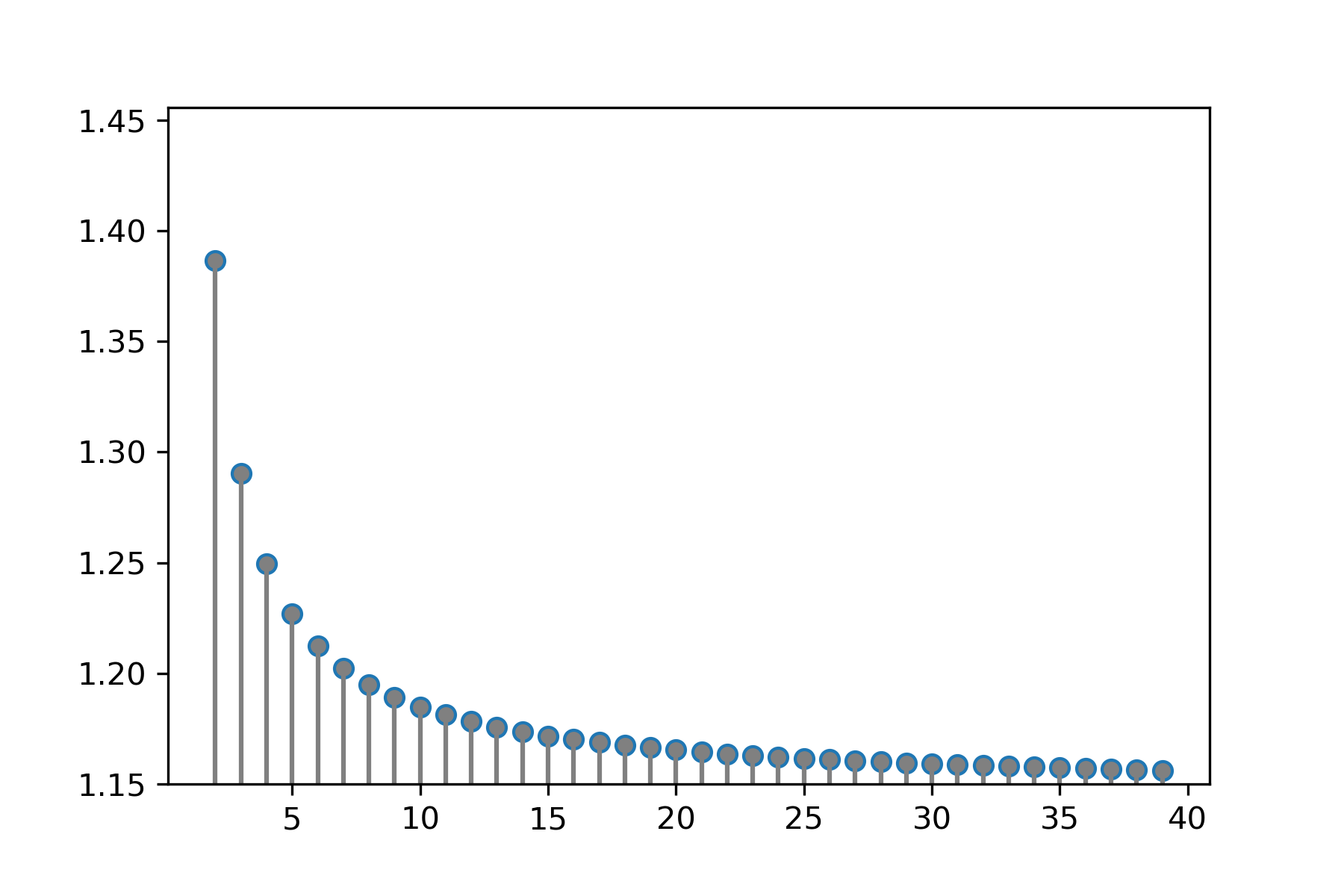}
        \caption{Non-uniform convergence}
        \label{fig:Bn-1/n} 
    \end{subfigure}
    \caption{Illustration of the non-uniform convergence of \prettyref{eq:Bn-asymp}. Panel (\subref{fig:Bn-MM-diff}) shows $\phi(x)-(B_n(x)+\frac{1-x}{2n})$ for $n=50$. Panel (\subref{fig:Bn-1/n}) further shows $\frac{2n}{1-x}(\phi(x)-B_n(x))$ evaluated at $x=1/n$ for different $n$. The sequence of values is not converging to one.  }
    \label{fig:Bn-asymp}
\end{figure}

However, the asymptotic estimate \prettyref{eq:Bn-asymp} is not uniform over $[0,1]$ (see \prettyref{eq:Bn-lb}, and also an illustration in \prettyref{fig:Bn-asymp}). 
Furthermore, when $S(P)$ is comparable or far exceeds the sample size, the asymptotic estimate of the bias in \prettyref{eq:hatH-asymp-bias} is no longer true. In fact, the bias can be expressed via \prettyref{eq:H-emp} as the expected KL divergence from the empirical to the true distribution. This shows that the empirical entropy is always underbiased. 
Applying \prettyref{eq:H-emp} yields that
\begin{equation}
    \label{eq:bias-KL}
    H(P)-\Expect[H(\hat P_n)]=\Expect[D(\hat P_n\| P)],
\end{equation}
where $D(\cdot\|\cdot)$ denotes the KL divergence (cf.~\prettyref{sec:notation}). This representation offers another way to see that the empirical entropy is underbiased.
Furthermore, from \prettyref{eq:bias-KL} we obtain the following upper bound of the bias \cite[Proposition 1]{Paninski03}: 
\begin{proposition}
\label{prop:h-plug-bias}
    \[
        0\le H(P)-\Expect[H(\hat P_n)]\le \log\pth{1+\frac{S(P)-1}{n}}.
    \]
\end{proposition}
\begin{proof}
    The KL divergence can be upper bounded by the $\chi^2$-divergence as follows \cite{probability.metrics}:
    \[
        D(\hat P_n\| P)\le \log(1+\chi^2(\hat P_n\| P)).
    \]
    By Jensen's inequality we have 
    \[
        \Expect[D(\hat P_n\| P)]\le \log(1+\Expect[\chi^2(\hat P_n\| P)]).
    \]
    The expectation in the right-hand side of the above inequality is
    \[
        \Expect[\chi^2(\hat P_n\| P)]
        = \sum_i \frac{\Expect(\hat p_i-p_i)^2}{p_i}
        = \sum_i \frac{1-p_i}{n}
        = \frac{S(P)-1}{n}.
    \]
    Applying \prettyref{eq:bias-KL} yields the conclusion.
\end{proof}

We next discuss the tightness of the previous bias analysis of the empirical entropy using the Bernstein polynomial \prettyref{eq:bias-Bn} again. We first state a lower bound on the Bernstein approximation obtained in \cite[Theorem 5]{BS2004}:
\begin{lemma}
    For $x\ge 15/n$,
    \begin{align}
      |B_n(x)-\phi(x)|
      &\ge \frac{1-x}{2n}+\frac{1}{12n^2x}-\frac{x}{12n^2}-\frac{1}{2n^3x^2}\nonumber\\
      &\ge \frac{1-x}{2n}+\frac{1}{20n^2x}-\frac{x}{12n^2}.\label{eq:Bn-lb}
    \end{align}
\end{lemma}
Consequently, using \prettyref{eq:bias-Bn} and \prettyref{eq:Bn-less}, for a distribution $P$ with $p_i\ge 15/n$, we have\footnote{
    For a fixed distribution, as $n$ diverges, it is obtained in \cite[(14)]{Harris75} that $$H(P)-H(\hat P_n)=\frac{S(P)-1}{2n}+\frac{1}{12n^2}\pth{\sum_i\frac{1}{p_i}-1}+O(n^{-3}).$$
}
\[
    |H(\hat P_n)-H(P)|\ge \frac{S(P)-1}{2n}+\frac{1}{20n^2}\pth{\sum_i\frac{1}{p_i}}-\frac{1}{12n^2}.
\]
From the above lower bound and the monotonicity in \prettyref{eq:Bn-less}, for the uniform distribution over $k$ elements, the bias of the empirical entropy is at least $\Omega(\min\{\frac{k}{n},1\})$.

\subsection{Variance of the empirical entropy}
Now we evaluate the variance of the empirical entropy. Note that empirical entropy is a linear estimate
\begin{equation}
    \label{eq:hatH-linear}
    H(\hat P)=\sum_{i}g(N_i)=\sum_j \Phi_j g(j),
\end{equation}
where $g(j)=\phi(j/n)$ and $\Phi_j$ is the $j$th fingerprint defined in \prettyref{eq:fp}, which
denotes the number of elements that appeared exactly $j$ times.
A variance upper bound can be obtained by the Efron-Stein-Steele inequality \cite{steele86}:
\[
    \var H(\hat P_n)\le \frac{n}{2}\Expect(\Delta g(\tilde N_{X_1})-\Delta g(\tilde N_{X_1'}))^2,
\]
where $X_1'$ is another independent sample from $P$, $\tilde N_i$ counts the occurrences of symbol $i$ in $X_2,\dots,X_n$, and $\Delta g(j)$ denotes the difference $g(j+1)-g(j)$. Applying the triangle inequality yields that 
\begin{equation}
    \label{eq:var-ub1}
    \var H(\hat P_n)\le n\Expect(\Delta g(\tilde N_{X_1}))^2.
\end{equation}
Another way of writing the above upper bound is 
\begin{equation}
    \label{eq:var-ub2}
    \var H(\hat P_n)\le n\sum_i p_i \Expect(\Delta g(\tilde N_i))^2
    = \sum_i \Expect N_i(\Delta g(N_i-1))^2,
\end{equation}
where $N_i\sim \Binom(n,p_i)$ and $g(j)=0$ for $j<0$. We have the following result on the variance of empirical entropy:
\begin{proposition}
\label{prop:h-plug-var}
    \[
        \var H(\hat P_n) \le \frac{\log^2 (\min\{n, eS(P)\})}{n}.
    \]
\end{proposition}
\begin{proof}
    Let $g(j)=\phi(j/n)$. The difference $\Delta g(j)$ can be uniformly upper bounded by $\frac{\log n}{n}$ in magnitude for every $j=0,\dots,n-1$, and thus by \prettyref{eq:var-ub1} we obtain that 
    \[
        \var H(\hat P_n)\le \frac{\log^2n}{n}.
    \]
    The derivative of $\phi$ over $[\frac{j}{n},\frac{j+1}{n}]$ is at most $\max\{|\log\frac{e j}{n}|,|\log \frac{e(j+1)}{n}|\}$ in magnitude. This yields a refined upper bound for $j= 1,\dots,n-1$:
    \[
        |\Delta g(j)| \le \frac{\max\{|\log(e j/n)|,1\}}{n}.
    \]
    Combining with the uniform upper bound $\frac{\log n}{n}$, we get 
    \[
        |\Delta g(j)| \le \frac{1}{n}\log \frac{e n}{j+1},\quad j=0,\dots n-1.
    \]
    Applying \prettyref{eq:var-ub2} yields that 
    \[
        \var H(\hat P_n)\le \frac{1}{n^2}\sum_i \Expect N_i \log^2 (en/N_i).
    \]
    Note that $x\mapsto x \log^2(e/x)$ is concave on $[0,1]$, and $x\mapsto \log^2(ex)$ is concave on $[1,\infty)$. We obtain that
    \[
        \sum_i \Expect [\log^2 (en/N_i)N_i/n ]\le \sum_i p_i \log^2(e/p_i)
        \le \log^2(e S(P)),
    \]
    according to Jensen's inequality.
\end{proof}

We will show in \prettyref{sec:h-lb} that, when the distribution is supported on $k$ elements, the MSE of any estimate using $n$ independent observations is $\Omega(\frac{\log^2 k}{n})$ in the worst case (see \prettyref{prop:lb1}). This lower bound also applies to the empirical entropy. 
Combined with the results in this section, we have proved the characterization of the worst-case MSE of the empirical entropy in \prettyref{thm:emprical-entropy}.

\section{Optimal entropy estimation on large domains}
\label{sec:h-esti}
The analysis in the last section shows that the empirical entropy is asymptotically optimal for distributions on a fixed alphabet as $n\diverge$. Specifically, \prettyref{eq:Rplug-rate} shows that the mean squared error of the empirical entropy is $O(\frac{\log^2 k}{n})$ when $n\ge \frac{k^2}{\log^2 k}$, which is optimal. However, the empirical entropy suffers from large bias when $n=O(k)$. In this section, we construct a minimax rate-optimal estimator that uses polynomial approximation to trade bias with variance optimally. 
We first present the estimator and its statistical guarantee, and then evaluate the estimator through numerical experiments.

\subsection{Optimal estimator via best polynomial approximation}
The major difficulty of entropy estimation lies in the bias due to insufficient sample size. Recall that the entropy is given by $H(P)=\sum \phi(p_i)$, where $\phi(x) = x \log \frac{1}{x}$. 
It is easy to see that 
the expectation of any estimator $T: [k]^n \to \reals_+$ is a polynomial of the underlying distribution $P$ and, consequently, no unbiased estimator for the entropy exists (see \cite[Proposition 8]{Paninski03}). This observation motivates us to approximate $\phi$ by a polynomial of degree $L$, say $g_L$, for which we pay a price in bias as the approximation error but yield the benefit of zero bias in estimating $g_L(p_i)$; cf.~\eqref{eq:unbiased-iid}.
While the approximation error clearly decays with the degree $L$, it is not unexpected that the variance of the unbiased estimator for $g_L(p_i)$ increases (exponentially) with $L$ as well as with the probability mass $p_i$. Therefore we only apply the polynomial approximation scheme to small $p_i$ and directly use the plug-in estimator for large $p_i$, since the ``signal-to-noise ratio'' is sufficiently large.

Next we describe the estimator in details. In view of the relationship between the risks with fixed and Poisson sample size in \prettyref{thm:poisson-sampling} and the fact that $0\le H\le \log k$, we shall assume the Poisson sampling model to simplify the analysis, where we first draw $ n'\sim \Poi(2n) $ and then draw a sample $X=(X_1,\ldots,X_{n'})$ of size $ n' $ independently from $P$. We split the sample roughly equally and use the first half for deciding whether to use the polynomial estimator or the plug-in estimator and the second half for the actual estimation. Specifically, for each $ X_i $ we draw an independent fair coin $ B_i\iiddistr \Bern\pth{\frac{1}{2}} $. We split the sample $ X $ according to the value of $ B_i$'s into two sets and compute the histogram of each subsample. That is, define $ N=(N_1,\dots,N_k) $ and $ N'=(N_1',\dots,N_k') $ by
\begin{align*}
    N_i=\sum_{j=1}^{n'}\indc{X_j=i}\indc{B_j=0},\quad  N_i'=\sum_{j=1}^{n'}\indc{X_j=i}\indc{B_j=1}.
\end{align*}
Then $ N$ and $N'$ are independent, where $ N_i,N_i'\iiddistr\Poi\pth{np_i} $.

Let $c_0,c_1,c_2>0$ be constants to be specified. Let $L= \floor{c_0\log k}$. Denote the best polynomial of degree $ L $ to uniformly approximate $ x\log\frac{1}{x} $ on $ \qth{0,1} $ by 
\begin{equation}
    \label{eq:pL}
    p_L(x)=\sum_{m=0}^{L}a_mx^m .
\end{equation} 
Through a change of variables, we see that the best polynomial of degree $ L $ to approximate $ x\log\frac{1}{x} $ on $ [0, \beta]$, where $\beta=\frac{c_1\log k}{n}$, is 
\[
    P_L(x)\triangleq \sum_{m=0}^L a_m\beta^{1-m}x^m-x\log\beta.
\]
Recall \prettyref{eq:unbiased-poi} an unbiased estimator for the monomials of the Poisson mean: $ \Expect[(X)_m]=\lambda^m $ where $ X\sim\Poi(\lambda) $. Consequently, the polynomial of degree $L$,
\begin{equation}
    \label{eq:gL}
    g_L(N_i)
    \triangleq \frac{1}{n}\pth{\sum_{m=0}^L \frac{a_m}{\pth{c_1\log k}^{m-1}} (N_i)_{m}
    - N_i\log\beta},
\end{equation}    
is an unbiased estimator for $ P_L(p_i) $.

Define a preliminary estimator of entropy $ H(P)=\sum_{i=1}^{k}\phi(p_i) $ by
\begin{equation}
    \label{eq:tH}
    \tilde{H}\triangleq\sum_{i=1}^{k}\pth{g_L(N_i)\indc{N_i'\le T}
        +g(N_i)\indc{N_i'>T}},
\end{equation}
where $T=c_2\log k$, $g(j)=\phi(j/n)+\frac{1}{2n}$, and we apply the estimator from polynomial approximation if $ N_i'\le T $ or the bias-corrected plug-in estimator otherwise (cf. the asymptotic expansion \prettyref{eq:hatH-asymp-bias} of the bias under the original sampling model). In view of the fact that $0 \leq H(P) \leq \log k $ for any distribution $ P $ with alphabet size $ k $, we define our final estimator by:
\begin{equation*}
    \hat{H} =(\tilde{H}\vee 0)  \wedge \log k.
\end{equation*}

The next result gives an upper bound on the above estimator under the Poisson sampling model, which, in view of the right inequality in \prettyref{eq:RRt} and \prettyref{eq:Rplug-rate}, implies the upper bound on the minimax risk $R^*(n,k)$ in \prettyref{eq:h-main}. 
\begin{proposition}
    \label{prop:err-rate}
    Assume that $\log n \leq  C \log k$ for some constant $C>0$. Then there exists $c_0,c_1,c_2$ depending on $C$ only, such that
    \[
        \sup_{P \in \calM_k} \Expect[(H(P)-\hat{H}(N))^2]
        \lesssim \pth{\frac{k}{n\log k}}^2+\frac{\log^2 k}{n},
    \]
    where $N=(N_1,\ldots,N_k) \inddistr \Poi(n p_i)$.
\end{proposition}

The above statistical guarantee is proved in \cite[Appendix C]{WY14}. Here we make a few comments on the optimal estimator.

\paragraph{Computation complexity. }  The estimate $\tilde H$ in \prettyref{eq:tH} can be expressed in terms of a linear combination of the fingerprints (see \prettyref{eq:hatH-linear}) of the second half of the sample. The coefficients $\{a_m\}$ can be \emph{pre-computed} using fast best polynomial approximation algorithms (\eg, \prettyref{algo:remez} due to Remez at the speed of linear convergence \cite[Theorem 8.2]{DL93}). It is clear that, after the coefficients are pre-computed, the estimator $\hat{H}$ can be computed in time that is linear in $n$, which is sublinear in the alphabet size.

\paragraph{Contribution from small probability masses. } The estimator in this section uses the polynomial approximation of $x\mapsto x\log \frac{1}{x}$ for those masses below $\frac{\log k}{n}$ and the bias-reduced plug-in estimator otherwise. This suggests that the major difficulty of entropy estimation lies in those probabilities in the interval $[0,\frac{\log k}{n}]$, which are individually small but collectively contribute significantly to the entropy. In the next section, to prove a minimax lower bound, the pair of unfavorable priors consists of randomized distributions whose masses are below $\frac{\log k}{n}$.

\begin{figure}[ht]
    \centering
    \includegraphics[width=0.7\linewidth]{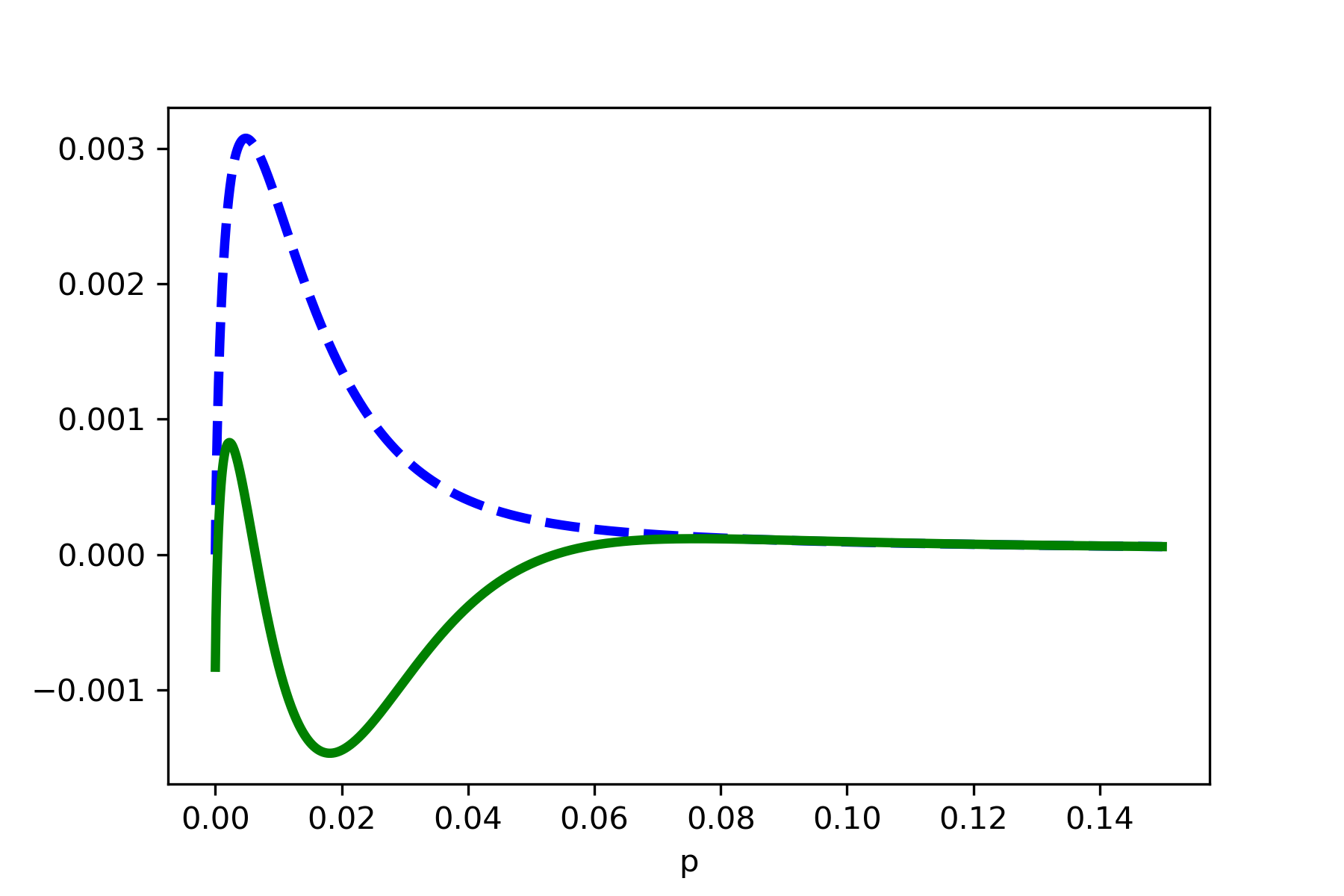}
    \caption{Comparison of the bias of estimators for $\phi(p)$ using $N\sim \Binom(n,p)$. The green solid line is the bias of the polynomial estimator $\tilde g(N)$ as a function of $p$; the blue dashed line shows the bias of the Miller-Madow estimator $\tilde g'(N)$. }
    \label{fig:hatH-compare}
\end{figure}

\paragraph{Bias reduction from polynomial approximation. } To show the effect of bias reduction using the best polynomial approximation, we illustrate $\phi(p)-\Expect[\tilde g(N)]$ as a function of $p$, where $N\sim \Binom(n,p)$ and 
\[
    \tilde g(j)=
    \begin{cases}
        g_L(j),& j\le T,\\
        \phi(j/n)+\frac{1-(j/n)}{2n}, & j>T.
    \end{cases}
\]
Here $g_L$ is obtained by \prettyref{eq:gL} using the best polynomial approximation. We also compare with that of the Miller-Madow estimate where $\tilde g'(j)=\phi(j/n)+\frac{1-(j/n)}{2n}$ for every $j$. In \prettyref{fig:hatH-compare}, we take a sample size $n=100$; $g_L(j)$ is obtained using the best polynomial of degree four to approximation $\phi$ on $[0,0.06]$, and is applied with $T=3$. We can clearly see the improvement on the bias as compared to the Miller-Madow estimate when $p$ is small.

\paragraph{Sample splitting. } The benefit of sample splitting is that we can first condition on the realization of $N'$ and treat the indicators in \prettyref{eq:tH} as deterministic, which has also been used in the entropy estimator in \cite{JVHW15}. Although not ideal operationally or aesthetically, this is a frequently-used idea in statistics and learning to simplify the analysis (also known as sample cloning in the Gaussian model \cite{Nemirovski03,CL11}) at the price of losing half of the sample thereby inflating the risk by a constant factor. It remains to be shown whether the optimality result in \prettyref{prop:err-rate} continues to hold if we can use the same sample in \prettyref{eq:tH} for both selection and estimation.

Note that the estimator \prettyref{eq:tH} is \emph{linear} in the fingerprint of the second half of the sample. We also note that for estimating other distribution functionals, \eg, support size, it is possible to circumvent sample splitting by directly using a linear estimator obtained from best polynomial approximation. See \prettyref{sec:supp} for details.

\begin{figure}[ht]
    \centering
    \includegraphics[width=\linewidth]{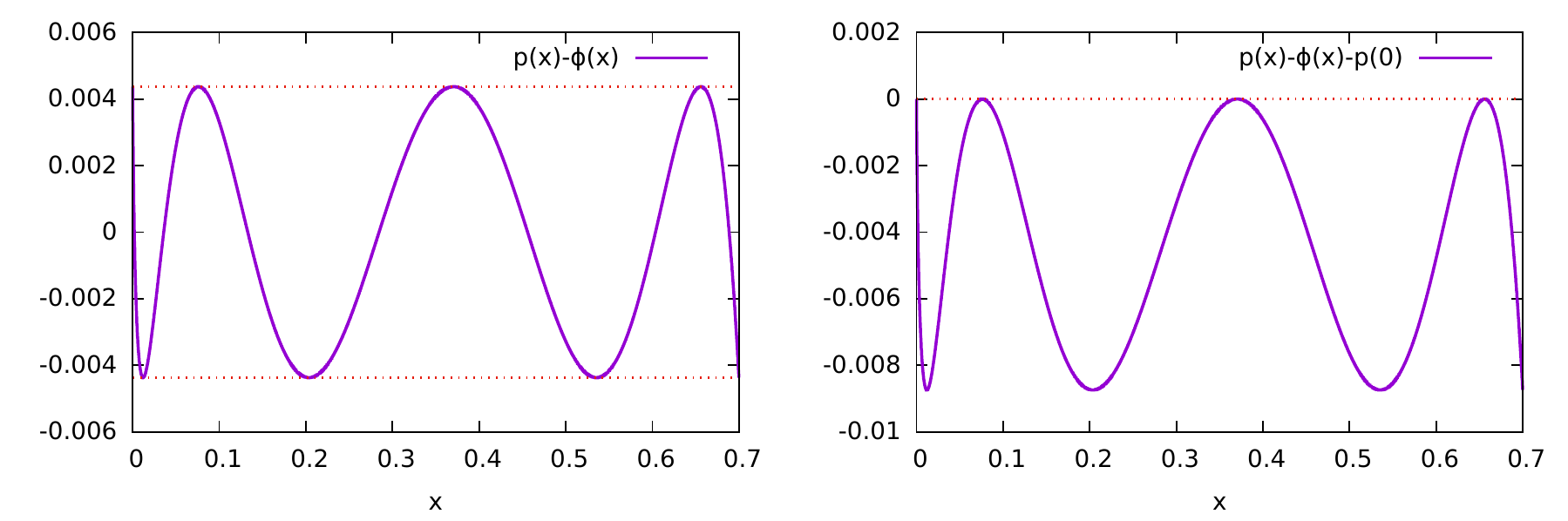}
    \caption{Bias of the degree-$ 6 $ polynomial estimator with and without the constant term.\label{fig:underbiased} }
\end{figure}

\paragraph{Adaptivity. } The estimator in \prettyref{eq:tH} depends on the alphabet size $k$ only through its logarithm; therefore the dependence on the alphabet size is rather insensitive. In many applications such as neuroscience the discrete data are obtained from quantizing an analog source and $k$ is naturally determined by the quantization level \cite{SLSKB97}. Nevertheless, it is also desirable to obtain an optimal estimator that is adaptive to $k$. To this end, we can replace all $\log k$ by $\log n$ and define the final estimator by $ \tilde{H}\vee 0 $. Moreover, we need to set $ g_L(0)=0 $ since the number of unseen symbols is unknown. Following \cite{JVHW15}, we can simply let the constant term $ a_0$ of the approximating polynomial \prettyref{eq:pL} go to zero and obtain the corresponding unbiased estimator \prettyref{eq:gL} through the falling factorials, which satisfies $ g_L(0)=0 $ by construction.\footnote{
    Alternatively, we can directly set $ g_L(0)=0 $ and use the original $g_L(j)$ in \prettyref{eq:gL} when $j \geq 1$. Then the bias becomes $ \sum_i(P_L(p_i)-\phi(p_i)-\prob{N_i=0}P_L(0)) $. In sublinear regime that $ n=o(k) $, we have $ \sum_i\prob{N_i=0}=\Theta(k) $; therefore this modified estimator also achieves the minimax rate.
}
The bias upper bound becomes $ \sum_i(P_L(p_i)-\phi(p_i)-P_L(0)) $ which is at most twice the original upper bound since $ P_L(0)\le \linf{P_L-\phi} $. The minimax rate in \prettyref{prop:err-rate} continues to hold in the regime of $\frac{k}{\log k} \lesssim n \lesssim \frac{k^2}{\log^2 k}$, where the plug-in estimator fails to attain the minimax rate. In fact, $ P_L(0) $ is always strictly positive and coincides with the uniform approximation error (see \cite[Appendix G]{WY14} for a short proof). Therefore, removing the constant term leads to $ g_L(N_i) $ which is always underbiased as shown in \prettyref{fig:underbiased}. A better choice for adaptive estimation is to find the best polynomial satisfying $ p_L(0)=0 $ that uniformly approximates $ \phi $.


\subsection{Numerical experiments}
In this subsection, we compare the performance of our entropy estimator to other estimators using synthetic data.\footnote{
    The C++ and Python implementation of our estimator is available at \url{https://github.com/Albuso0/entropy}.
}
Note that the coefficients of best polynomial to approximate $ \phi $ on $ [0,1] $ are independent of data so they can be pre-computed and tabulated to facilitate the computation of the estimator. The coefficients of the polynomial approximant can be efficiently obtained by applying the Remez algorithm, which provably converges linearly for all continuous functions (see, \eg, \cite[Theorem 1.10]{petrushev2011rational}). Since the choice of the polynomial degree is only logarithmic in the alphabet size, we pre-compute the coefficients up to degree $ 400 $ which suffices for all practical purposes. In the implementation of our estimator we replace $ N_i' $ by $ N_i $ in \prettyref{eq:tH} without conducting sample splitting. Though in proving the theoretical guarantees we are conservative about choosing the constants $ c_0,c_1,c_2 $, in experiments we observe that the performance of our estimator is in fact not sensitive to their values within the reasonable range. In the subsequent experiments the parameters are fixed to be $ c_0=c_2=1.6, c_1=3.5 $.

\begin{figure}[ht]
    \centering
    \includegraphics[width=0.9\linewidth]{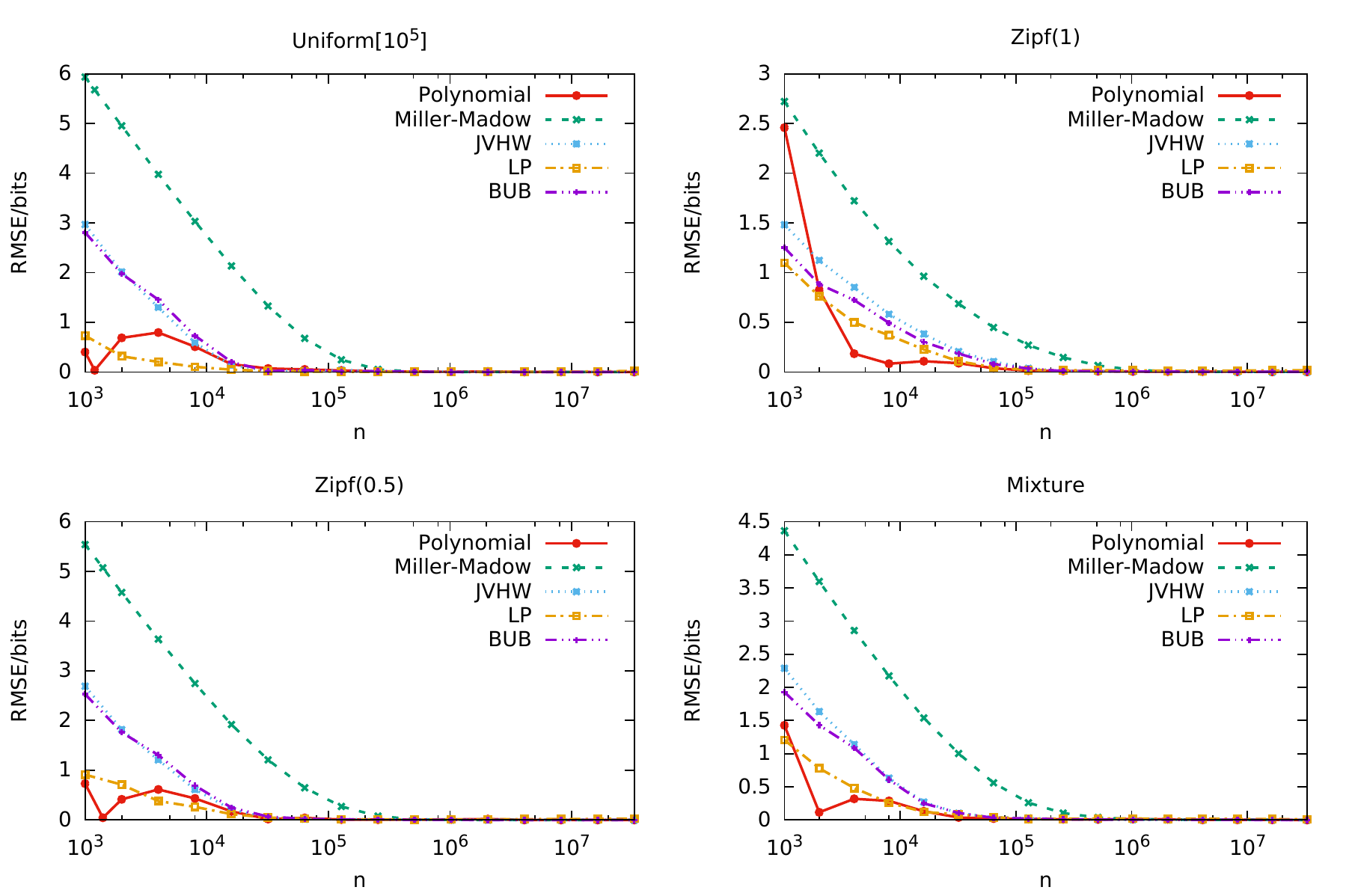}
    \caption{Performance comparison with sample size $n$ ranging from $ 10^3 $ to $3 \times 10^7 $. \label{fig:full} }
\end{figure}

\begin{figure}[ht]
    \centering
    \includegraphics[width=0.9\linewidth]{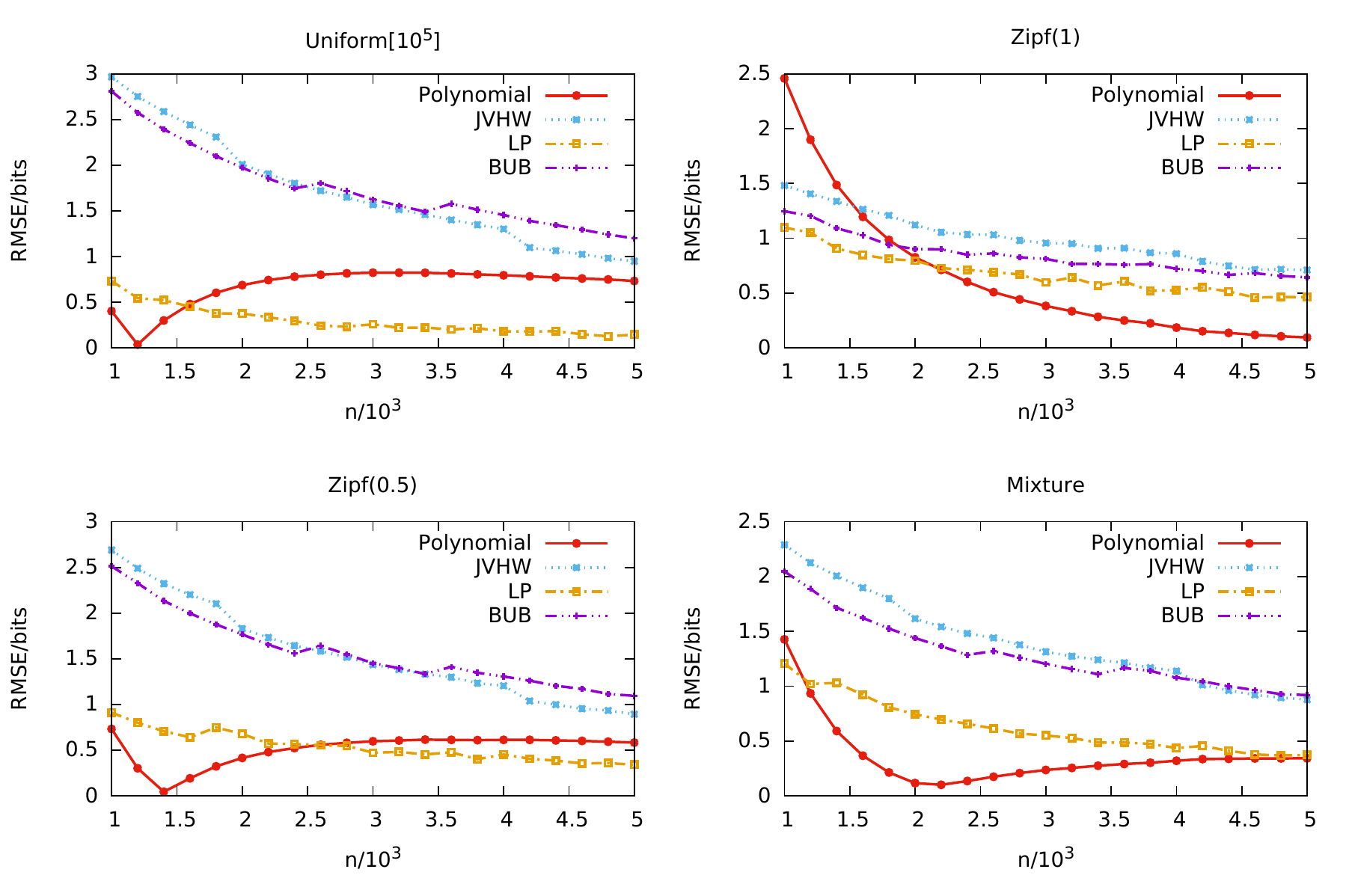}
    \caption{Performance comparison when sample size $n$ ranges from $ 1000 $ to $ 5000 $. \label{fig:scarce} }
\end{figure}

We generate data from four types of distributions over an alphabet of $ k=10^5 $ elements, namely, 
the uniform distribution with $ p_i=\frac{1}{k} $, Zipf distributions with $ p_i\propto i^{-\alpha} $ and $ \alpha$ being either $1$ or $0.5 $, and an ``even mixture'' of geometric distribution and Zipf distribution where  for the first half of the alphabet $ p_i \propto 1/i$ and  for the second half $ p_{i+k/2} \propto (1-\frac{2}{k})^{i-1} $, $ 1\le i\le \frac{k}{2} $. Using parameters mentioned above, the approximating polynomial has degree $ 18 $, the parameter determining the approximation interval is $ c_1\log k=40 $, and the threshold to decide which estimator to use in \prettyref{eq:tH} is $ 18 $; namely, we apply the polynomial estimator $ g_L $ if a symbol appeared at most 18 times and the  bias-corrected plug-in estimator otherwise. After obtaining the preliminary estimate $ \tilde{H} $ in \prettyref{eq:tH}, our final output is $ \tilde{H}\vee 0 $.\footnote{
    We can, as in \prettyref{prop:err-rate}, output $ (\tilde{H}\vee 0)\wedge \log k $, which yields a better performance. We elect not to do so for a stricter comparison.
} 
Since the plug-in estimator suffers from severe bias when the sample size is small, we forgo the comparison with it to save space in the figures and instead compare with its bias-corrected version, \ie, the Miller-Madow estimator \prettyref{eq:hatH-MM}. We also compare the performance with the LP-based estimator in \cite{VV13}, the best upper bound (BUB) estimator \cite{Paninski03}, and the estimator based on similar polynomial approximation techniques\footnote{
    The estimator in \cite{JVHW15} uses a smooth cutoff function in lieu of the indicator function in \prettyref{eq:tH}; this seems to improve neither the theoretical error bound nor the empirical performance.
} 
proposed by \cite{JVHW15} using their implementations with default parameters. Our estimator is implemented in C++ which is much faster than those from \cite{VV13,JVHW15,Paninski03} implemented in MATLAB so the running time comparison is ignored. We notice that the LP-based estimator in \cite{VV13} is much slower than the polynomial estimator in \cite{JVHW15}, especially when the sample size is large.

We compute the root mean squared error (RMSE) for each estimator over $ 50 $ trials. The full performance comparison is shown in \prettyref{fig:full} where the sample size ranges from one percent to $ 300 $ folds of the alphabet size. In \prettyref{fig:scarce} we further zoom into the more interesting regime of smaller sample size ranging from one to five percent of the alphabet size. In this regime our estimator, as well as those from \cite{VV13,JVHW15,Paninski03}, outperforms the classical Miller-Madow estimator significantly; furthermore, our estimator performs better than those in \cite{JVHW15,Paninski03} in most cases tested and comparably with that in \cite{VV13}. When the observations are abundant, all estimators achieve very small error; however, it has been empirically observed in \cite{JVHW15} that the performance of LP estimator starts to deteriorate when the sample size is very large, which is also observed in our experiments (see \cite{yang-msthesis}). The specific figures of that regime are ignored since the absolute errors are very small and even the plug-in estimator without bias correction is accurate. By \prettyref{eq:tH}, for large sample size our estimator tends to the Miller-Madow estimator when every symbol is observed many times.

\section{Fundamental limits of entropy estimation}
\label{sec:h-lb}
Thus far, we have described the empirical entropy and the construction of an estimator using the idea of polynomial approximation with significantly reduced bias. The worst-case MSE of both estimators are analyzed. To establish the minimax rate of 
entropy estimation in \prettyref{thm:entropy}, we need a matching minimax lower bound. This is the goal of the present section.

To obtain the lower bound part of \prettyref{eq:h-main}, it suffices to show that the minimax risk is lower bounded by the two terms in \prettyref{eq:h-main} separately. This follows from combining Propositions \ref{prop:lb1} and \ref{prop:lb2} below.
\begin{proposition}
    \label{prop:lb1}
    For all $k,n\in \naturals$,
    \begin{equation}
        \label{eq:lb1}
        R_H^*(k,n) \gtrsim  \frac{\log^2 k}{n}.   
    \end{equation}
\end{proposition}

\begin{proposition}
    \label{prop:lb2}
    For all $k,n\in \naturals$,
    \begin{equation}
        \label{eq:lb2}
        R_H^*(k,n) \gtrsim \pth{\frac{k}{n \log k}}^2\vee 1.
    \end{equation}
\end{proposition}

\prettyref{prop:lb1} follows from a simple application of Le Cam's \emph{two-point method} in \prettyref{thm:lb-two-pts}: If two input distributions $P$ and $Q$ are sufficiently close such that it is impossible to reliably distinguish between them using $n$ \iid~observations with error probability less than, say, $\frac{1}{2}$, then any estimator suffers a quadratic risk proportional to the separation of the functional values $|H(P)-H(Q)|^2$. 
\begin{proof}
    For any pair of distributions $P$ and $Q$, applying \prettyref{thm:lb-two-pts} yields that
    \begin{equation}
        \label{eq:lb-KL}
        R_H^*(k,n) \geq \frac{1}{4} (H(P)-H(Q))^2 \exp(-nD(P\|Q)).
    \end{equation}
    Therefore it boils down to solving the optimization problem: 
    \begin{equation}
        \label{eq:HD}
        \sup\{H(P)-H(Q): D(P\|Q) \leq 1/n\},
    \end{equation}
		where $P$ and $Q$ are distributions on $k$ elements. 
    Without loss of generality, assume that $k\geq 2$. Fix an $\epsilon \in (0,1)$ to be specified. Let 
    \begin{equation}
        \label{eq:PQ}
        P=\pth{\frac{1}{3k'},\ldots,\frac{1}{3k'},\frac{2}{3}}, \quad
        Q=\pth{\frac{1+\epsilon}{3k'},\ldots,\frac{1+\epsilon}{3k'},\frac{2-\epsilon}{3}},
    \end{equation}
    where $k'=k-1$. Direct computation yields
    $D(P\|Q) = \frac{2}{3} \log \frac{2}{2-\epsilon} + \frac{1}{3} \log \frac{1}{\epsilon+1}  \leq \epsilon^2$ and $H(Q)-H(P) = \frac{1}{3} (\epsilon \log k' + \log 4 + (2-\epsilon) \log \frac{1}{2-\epsilon} + (1+\epsilon) \log \frac{1}{\epsilon+1}) \geq \frac{1}{3} \log(2k') \epsilon- \epsilon^2$. Choosing $\epsilon = \frac{1}{\sqrt{n}}$ and applying \prettyref{eq:lb-KL}, we obtain the desired \prettyref{eq:lb1}.
\end{proof}
\begin{remark}
    \label{rmk:HD}
		It turns out that \prettyref{eq:lb1} is the best lower bound offered by the two-point method. To see this, recall the Pinsker inequality $D(P\|Q) \geq 2 \TV^2(P,Q)$ \cite[p. 58]{ckbook} and the continuity property of entropy with respect to the total variation distance: for any $P$ and $Q$ on $k$ elements, $|H(P)-H(Q)| \leq \TV(P,Q) \log \frac{k}{\TV(P,Q)}$ for $\TV(P,Q) \leq \frac{1}{4}$ \cite[Lemma 2.7]{ckbook}. Thus we conclude that the best two-point lower bound, \ie, the supremum in \prettyref{eq:HD}, is on the order of $\frac{\log k}{\sqrt{n}}$, achieved by the pair in \prettyref{eq:PQ}.
		Consequently, to prove the lower bound in \prettyref{prop:lb2}, one need to go beyond the simple two-point method.
\end{remark}

The remainder of this section is devoted to explaining  the main idea in proving \prettyref{prop:lb2};  for details see \cite{WY14}. Since the best lower bound provided by the two-point method is $\frac{\log^2 k}{n}$ (see \prettyref{rmk:HD}), proving \prettyref{eq:lb2} requires more powerful techniques. To this end, we use a generalized version of Le Cam's method involving two \emph{composite} hypotheses introduced in \prettyref{thm:lb-two-priors}:
\begin{equation}
    \label{eq:compHT}
    H_0: H(P) \leq t \quad \text{versus} \quad H_1: H(P) \geq t+ d,
\end{equation}
which is more general than the two-point argument using only simple hypothesis testing. Similarly, if we can establish that no test can distinguish \prettyref{eq:compHT} reliably, then we obtain a lower bound for the quadratic risk on the order of $d^2$. By the minimax theorem, the optimal probability of error for the composite hypotheses test is given by the Bayesian version with respect to the least favorable priors. For  \prettyref{eq:compHT} we need to choose a pair of priors, which, in this case, are distributions on the probability simplex $\calM_k$, to ensure that the entropy values are separated.

\subsection{Construction of the priors}
\label{sec:prior}
The main idea for constructing the priors is as follows: First of all, the symmetry of the entropy functional implies that the least favorable prior must be permutation-invariant. This inspires us to use the following \emph{\iid construction}. For conciseness, we focus on the special case of $n \asymp \frac{k}{\log k}$ for now with the goal of obtaining an $\Omega(1)$ lower bound. Let $U$ be a $\reals_+$-valued random variable with unit mean. Consider the random vector 
\begin{equation}
\sfP=  \frac{1}{k} (U_1,\ldots,U_k),
\label{eq:PU}
\end{equation}
consisting of \iid copies of $U$. Note that $\sfP$ itself is \emph{not} a probability distribution; however, the key observation is that, since $\Expect[U]=1$, as long as the variance of $ U $ is not too large, the weak law of large numbers ensures that $\sfP$ is \emph{approximately} a probability vector. 
\begin{example}
    A deterministic $U = 1$ generates $\sfP=\pth{\frac{1}{k},\dots,\frac{1}{k}}$, 
    which is the uniform distribution over $k$ elements. For the binary-valued $U \sim \frac{1}{2} (\delta_0+\delta_2)$, in \prettyref{eq:PU}
		approximately half of $U_i$ is two and others are zero. In this case $\sfP$ roughly corresponds to an uniform distribution over $k/2$ elements whose support set is uniformly chosen at random.
\end{example}
From this viewpoint, the CDF of the random variable $\frac{U}{k}$ plays the role of the \emph{histogram of the distribution} $\sfP$, which is the central object in the Valiant-Valiant lower bound construction (see \cite[Definition 3]{VV10}). 
Although $\sfP$ is not a probability vector per se, using a conditioning argument we can show that the distribution of $\sfP$ can effectively serve as a prior; cf.~\cite[Lemma 2]{WY14}.

Next we apply Le Cam's method outlined in \prettyref{thm:lb-two-priors} 
by choosing two random approximate probability vectors $\sfP=  \frac{1}{k} (U_1,\ldots,U_k)$ and $\sfP'=  \frac{1}{k} (U_1',\ldots,U_k')$,
where $U_i,U_i'$ are \iid copies of $U$ and $U'$ whose distributions need to be carefully chosen to achieve the following goals:
\begin{enumerate}
    \item \emph{Functional value separation}: 
		Recall $\phi(x) =x \log \frac{1}{x}$. 
		Note that
    \[
        H(\sfP)= \sum_{i=1}^k\phi\pth{\frac{U_i}{k}} =\frac{1}{k}\sum_{i=1}^k\phi(U_i)+\frac{\log k}{k}\sum_{i=1}^kU_i,
    \]
    which, by law of large numbers, concentrates near its mean $ \Expect[H(\sfP)]=\Expect[\phi(U)]+\Expect[U]\log k$. Therefore, given another random variable $U'$ with unit mean, we can obtain $\sfP'$ similarly using \iid copies of $U'$. Then with high probability, $H(\sfP)$ and $H(\sfP')$ are separated by the difference of their mean values, namely,
    \[
        \Expect[H(\sfP)]-\Expect[H(\sfP')] = \Expect[\phi(U)] - \Expect[\phi(U')],
    \]
    which we aim to maximize.

    \item \emph{Indistinguishability}: 
    Recall that given $P$, the sufficient statistics are distributed as $N_i\inddistr \Poi(n p_i)$. Therefore, if $P$ is drawn from the distribution of $\sfP$, then $N=(N_1,\ldots,N_k)$ are \iid distributed according the \emph{Poisson mixture}, which we denote by $\Expect[\Poi(\frac{n}{k} U)]$. Similarly, if $P$ is drawn from the prior of $\sfP'$, then $N$ is distributed according to $(\Expect[\Poi(\frac{n}{k} U')])^{\otimes k}$. To establish the impossibility of testing, we need the total variation distance between the two $k$-fold product distributions to be strictly bounded away from one, for which a sufficient condition is\footnote{This follows from the simple fact that $\TV(P^{\otimes k}, Q^{\otimes k}) \leq k \TV(P,Q)$.}
    \begin{equation}
        \label{eq:hu2}
        \TV(\Expect[\Poi(n U/k)],\Expect[\Poi(n U'/k)]) \leq c/k
    \end{equation}
    for some $c<1$.
\end{enumerate}
To conclude, we see that the \iid construction \prettyref{eq:PU} fully exploits the independence blessed by the Poisson sampling, thereby reducing the problem to \emph{one dimension}. This allows us to sidestep the difficulty encountered in \cite{VV10} when dealing with fingerprints which are high-dimensional random vectors with dependent entries.

What remains is the following optimization problem in one dimension: choose $U,U'$ to maximize $|\Expect[\phi(U)] - \Expect[\phi(U')]|$ subject to the constraint \prettyref{eq:hu2}. As previously discussed in \prettyref{sec:d-mm}, a commonly used proxy for bounding the total variation distance is \emph{moment matching}, \ie, $\Expect[U^j] = \Expect[U'^j]$ for all $j = 1,\ldots,L$. Together with the $L_\infty$-norm constraints that $|U|,|U'|\lesssim \frac{kL}{n}$, a sufficiently large degree $L$ ensures the total variation bound \prettyref{eq:hu2}; cf.~\prettyref{thm:tv-bound}. Combining the above steps, our lower bound is proportional to the value of the following convex optimization problem (in fact, an infinite-dimensional LP over probability measures):
\begin{equation}
    \label{eq:H-FL}
    \begin{aligned}
     \calF_L(\lambda) \triangleq  \sup & ~ \Expect[\phi(U)] - \Expect[\phi(U')]  \\
        \text{s.t.}     
        & ~ \Expect[U] = \Expect[U']=1 \\
        & ~ \Expect[U^j] = \Expect[U'^j], \quad j = 1,\ldots,L, \\
        & ~ U,U' \in [0, \lambda]
    \end{aligned}
\end{equation}
for some appropriately chosen $L \in \naturals$ and $\lambda > 1$ depending on $n$ and $k$.

Finally, the optimization problem \prettyref{eq:H-FL} to the machinery of \emph{best polynomial approximation} (see \cite[Appendix B.2]{WY14}):
\begin{equation}
    \label{eq:FLEL}
    \calF_L(\lambda) \geq 2 E_L(\log, [1/\lambda,1]).
\end{equation}
Due to the singularity of the logarithm at zero, the approximation error can be made bounded away from zero if $\lambda$ grows \emph{quadratically} with the degree $L$ (see \prettyref{lmm:sep} below). Choosing $L \asymp \log k$ and $\lambda \asymp \log^2 k$ leads to the impossibility of consistent estimation for $n \asymp \frac{k}{\log k}$. For $n \gg \frac{k}{\log k}$, the lower bound for the quadratic risk follows from relaxing the unit-mean constraint in \prettyref{eq:H-FL} to $\Expect[U] = \Expect[U'] \leq 1$ and a simple scaling argument; we refer the reader to \cite[Appendix B.2]{WY14} for details. Analogous construction of priors and proof techniques have subsequently been used in \cite{JVHW15} to obtain sharp minimax lower bound for estimating the power sum in which case the $\log p$ function is replaced by $p^\alpha$.


\begin{lemma}
    \label{lmm:sep}
    There exist universal positive constants $c , c', L_0$ such that for any $L\geq L_0$,
    \begin{equation}
        E_{\floor{c L}}(\log, [L^{-2},1]) \geq c'.
        \label{eq:sep}
    \end{equation}
\end{lemma}
The proof of \prettyref{lmm:sep} uses the converse result of best uniform approximation in \prettyref{thm:Ivanov}.
See \cite[Appendix F]{WY14} for details.

\section{Bibliographic notes}

The problem of entropy estimation has a long and rich history in the fields of information theory and statistics, with early work dating back to \cite{Dobrushin58,Basharin59,Miller55,Harris75,Stein1986}. Notably, the inadequacy of plug-in estimator on large domain was already recognized in \cite{Dobrushin58}. In Charles Stein's 1976 lectures in Leningrad \cite[Lecture 5]{Stein1986}, entropy estimation was used as an example to study models with many parameters and Taylor expansion was applied to tackle the difficulty arising from large alphabets.

It is well-known that to estimate the distribution $P$ itself, say, with a total variation loss at most 0.1, 
we need $\Theta(k)$ \iid~observations (see, \eg, \cite{BFSS02}). 
Using non-constructive arguments, Paninski first proved that it is possible to consistently estimate the entropy with \emph{sublinear} sample size, \ie, there exists $n_k=o(k)$, such that the minimax risk $R_H^*(k,n_k)\to 0$ as $k\to \infty$ \cite{Paninski04}. 
Valiant proved that no consistent estimator exists, \ie, $R^*(k,n_k) \gtrsim 1$ if $n \lesssim \frac{k}{\exp(\sqrt{\log k})}$ \cite{Valiant08}.
The sharp scaling of the minimal sample size of consistent estimation is shown to be $\frac{k}{\log k}$ in the breakthrough results of Valiant and Valiant \cite{VV10,VV11}. 
Later an estimator based on linear programming is shown to achieve an additive error of $\epsilon$ using $\frac{k}{\epsilon^2 \log k}$ observations \cite[Theorem 1]{VV13}, while $\frac{k}{\epsilon \log k}$ observations are shown to be necessary \cite[Corollary 10]{VV10}. 
This gap is partially amended in \cite{VV11-focs} by a different estimator, which requires $ \frac{k}{\epsilon \log k} $ observations but only valid when $ \epsilon>k^{-0.03} $.
\prettyref{thm:entropy} improve upon these results by characterizing the full minimax rate and the sharp sample complexity.

\prettyref{thm:entropy} has been independently developed by \cite{WY14} and \cite{JVHW15} using the idea of best polynomial approximation, and is sharpened in \cite{HJW15adaptive} by showing that the minimax rate is $\log^2(\frac{k}{n\log n})$ in the regime of $ n\lesssim \frac{k}{\log k} $ using similar techniques. 
For more recent results on estimating Shannon entropy, R\'enyi entropy and other distributional functionals on large alphabets, see \cite{JVHW2017,JOST15,WY15,HJW15,HJW15adaptive}.
Unified methodologies for symmetric functional estimation have also been developed including local moment matching \cite{HJW18lmm} and profile maximum likelihood \cite{ADOS17}.

Finally, we mention that entropy estimation based on polynomial approximation techniques has been recently extended to continuous settings and non-\iid data. 
Estimating the differential entropy of a smooth density was studied in \cite{HJWW17}, where optimal estimator is obtained based on polynomial approximation in conjunction with kernel methods. Estimation of the entropy rate of Markov chains was considered \cite{HJLWWY17}, where the optimal sample complexity was determined in certain regime in terms of both the alphabet size and the spectral gap of the Markov chain; in particular, the minimax lower bound uses on more sophisticated moment matching techniques to construct two random probability transition matrices with prescribed entropy rates and spectral gaps.



\chapter{Estimating the unseen}
\label{chap:unseen}

Estimating the support size of a distribution from data is a classical problem in statistics with widespread applications. For example, a major task for ecologists is to estimate the number of species \cite{FCW43} from field experiments; a classical task in linguistics is to estimate the vocabulary size of Shakespeare based on his complete works \cite{McNeil73,ET76,TE87}; in population genetics it is of great interest to estimate the number of different alleles in a population \cite{HW01}. Estimating the support size is equivalent to estimating the number of unseen symbols. This is particularly challenging when the sample size is relatively small compared to the total population size, since a significant portion of the population are never observed in the data. This chapter discusses the support size estimation problem in \prettyref{sec:supp} and the closely related distinct elements problem in \prettyref{sec:distinct}, which refers to the problem of estimating the number of distinct colors by sampling from an urn of colored balls.

\section{Estimating the support size}
\label{sec:supp}
We adopt the following statistical model \cite{BO79,RRSS09}. Let $ P $ be a discrete distribution over some countable alphabet. Without loss of generality, we assume the alphabet is $ \naturals $ and denote $ P=(p_1,p_2,\dots) $. Given $ n $ independent observations $ X\triangleq(X_1,\dots,X_n) $ drawn from $P$, the goal is to estimate the support size
\begin{equation}
    \label{eq:supp}
    S=S(P)\triangleq \sum_{i}\indc{p_i>0}.
\end{equation}
As discussed in \prettyref{sec:poi}, since the support size is a symmetric function of the distribution, both the histogram \prettyref{eq:histogram} and the fingerprint \prettyref{eq:fp} are sufficient statistics for estimating $S(P)$. We refer to the estimation of $S(P)$ as the \SupportSize problem.

It is clear that unless we impose further assumptions on the distribution $P$, it is impossible to estimate $S(P)$ within a given accuracy; otherwise there can be arbitrarily many masses in the support of $P$ that, with high probability, are never sampled and the worst-case risk for estimating $S(P)$ is obviously infinite. To prevent the triviality, a conventional assumption \cite{RRSS09} is to impose a lower bound on the non-zero probabilities. Therefore we restrict our attention to the parameter space $ \calD_k $, which consists of all probability distributions on $ \naturals $ whose minimum non-zero mass is at least $ \frac{1}{k} $; consequently $ S(P)\le k $ for any $P \in \calD_k$.
We refer to this as the \SupportSize problem.

\subsection{Fundamental limits of the \SupportSize problem}
\label{sec:supp-rate}
To quantify the decision-theoretic fundamental limit, we again consider the minimax MSE:
\begin{equation}
    \label{eq:S-risk}
    R_{\sf S}^*(k,n)\triangleq \inf_{\hat S}\sup_{P\in\calD_k}\Expect(\hat S- S)^2,
\end{equation}
where $\hat S$ is an integer-valued estimator measurable with respect to $X_1,\dots,X_n \iiddistr P$.

\begin{theorem}
    \label{thm:main}
    For all $k,n \geq 2$, 
    \begin{equation}
        \label{eq:S-main}
        R_{\sf S}^*(k,n) = k^2\exp\pth{-\Theta\pth{\sqrt{\frac{n\log k}{k}}\vee \frac{n}{k} \vee 1}}.
    \end{equation}
    Furthermore, if $ \frac{k}{\log k}\ll n\ll k\log k $, as $ k\rightarrow \infty $,
    \begin{equation}
        \label{eq:S-main-asymp}
        k^2\exp\pth{-c_1\sqrt{\frac{n\log k}{k}}}\le R_{\sf S}^*(k,n)\le k^2\exp\pth{-c_2\sqrt{\frac{n\log k}{k}}},
    \end{equation}
    where $c_1=\sqrt{2}e+o(1)$ and $c_2=1.579+o(1)$.
 \end{theorem}
To interpret the rate of convergence in \prettyref{eq:S-main}, we consider three cases:
\begin{description}
    \item[Simple regime] $ n \gtrsim k\log k $: we have $ R_{\sf S}^*(k,n) = k^2\exp(-\Theta(\frac{n}{k})) $ which can be achieved by the simple plug-in estimator
    \begin{equation}
        \label{eq:splug}
        \Splug\triangleq \sum_{i} \indc{N_i>0},
    \end{equation}
    that is, the number of observed symbols or the support size of the empirical distribution. Furthermore, if $\frac{n}{k \log k}$ exceeds a sufficiently large constant, all symbols are present in the data and $\Splug$ is in fact exact with high probability, namely, $\Prob[\Splug \neq S]\le \Expect(\Splug - S)^2 \to 0$. This can be understood as the classical coupon collector's problem (cf.~\eg, \cite{MU06}).    

    \item[Non-trivial regime] $\frac{k}{\log k} \ll n\ll k\log k $: In this case the observations are relatively scarce and the naive plug-in estimator grossly underestimate the true support size as many symbols are simply not observed. Nevertheless, accurate estimation is still possible and the optimal risk is given by $ R_{\sf S}^*(k,n) = k^2\exp(-\Theta(\sqrt{\frac{n\log k}{k}})) $. This can be achieved by a linear estimator based on the Chebyshev polynomial and its approximation-theoretic properties. Although more sophisticated than the plug-in estimator, this procedure can be evaluated in $O(n+\log^2 k)$ time.

    \item[Impossible regime] $n \lesssim \frac{k}{\log k}$: any estimator suffers an error proportional to $ k $ in the worst case.
\end{description}


Following \prettyref{def:sample-complexity}, we define the sample complexity in the \SupportSize problem by
\begin{equation}
    \label{eq:sample-support}
    n_{\sf S}^*(k,\Delta) \triangleq \min\{n \geq 0\colon \exists \hat S, \text{ s.t. } \Prob[|\hat S - S(P)| \ge \Delta] \leq 0.1,  \forall P\in \calD_k \},
\end{equation}
where $\hat S$ is an integer-valued estimator measurable with respect to the $X_1,\ldots,X_n\iiddistr P$.
The next result characterizes the sample complexity within universal constant factors (in fact, within a factor of six asymptotically). 
\begin{theorem}
    \label{thm:sample}
    Fix a constant $c_0 < \frac{1}{2}$. For all $1\le \Delta \leq c_0 k$,
    \begin{equation}
        \label{eq:sample-complexity}
        n_{\sf S}^*(k,\Delta) \asymp \frac{k}{\log k}\log^2\frac{k}{\Delta}.
    \end{equation}
    Furthermore, if $ \Delta=o(k) $ and $\Delta = k^{1-o(1)}$, as $ k\rightarrow \infty $,
    \begin{equation}
        \label{eq:sample-complexity-asymp}
        \frac{\tilde c_1k}{\log k}\log^2\frac{k}{\Delta}\le n_{\sf S}^*(k,\Delta)\le \frac{\tilde c_2k}{\log k}\log^2\frac{k}{\Delta},
    \end{equation}
    where $\tilde  c_1=\frac{1}{2e^2}+o(1)$ and $\tilde c_2=\frac{1}{2.494}+o(1)$.
\end{theorem}
Compared to \prettyref{thm:main}, the only difference is that here we are dealing with the zero-one loss $\indc{|S-\hat S| \ge \Delta}$ instead of the quadratic loss $(S-\hat S)^2$.
In this section, we shall show the upper bound for the quadratic risk (\prettyref{sec:optimal}) and the lower bound for the zero-one loss (\prettyref{sec:minimax}), thereby proving both \prettyref{thm:main} and \ref{thm:sample} simultaneously.
Furthermore, the choice of 0.1 as the probability of error in the definition of the sample complexity in \prettyref{eq:sample-support} 
is entirely arbitrary; replacing it by $1-\delta$ for any constant $\delta \in (0,1)$ only affect $n_{\sf S}^*(k,\Delta)$ up to constant factors.\footnote{Specifically, upgrading the confidence to $1-\delta$ can be achieved by oversampling by merely a factor of $\log \frac{1}{\delta}$: Let $T = \log \frac{1}{\delta}$. With $nT$ observations, divide them into $T$ batches, apply the $n$-sample estimator to each batch and aggregate by taking the median. Then Hoeffding's inequality implies the desired confidence.}

\subsection{Optimal estimator via Chebyshev polynomials}
\label{sec:optimal}
In this subsection we show the upper bound part of \prettyref{thm:main} and describe the rate-optimal support size estimator in the non-trivial regime. 
By \prettyref{thm:poisson-sampling}, since $0\leq S\leq k$, for the purpose of proving \prettyref{thm:main} it suffices to consider the Poisson sampling model, where the sample size is $ \Poi(n) $ instead of a fixed number $ n $ and hence the sufficient statistics $ N=(N_1,\dots,N_k)\inddistr\Poi(np_i) $. 
In the next proposition, we first analyze the risk of the plug-in estimator $ \Splug $ in \prettyref{eq:splug}, which yields the optimal upper bound of \prettyref{thm:main} in the regime of $ n\gtrsim k\log k $. This is consistent with the intuition of coupon collection. 
\begin{proposition}
    \label{prop:main-ub-plug}
    For all $ n, k\ge 1 $,
    \begin{equation}
        \label{eq:S-main-ub-plug}
        \sup_{P\in\calD_k}\Expect(S(P)-\Splug(N))^2\le k^2e^{-2n/k}+ke^{-n/k},
    \end{equation}
    where $ N=(N_1,N_2,\dots) $ and $ N_i\inddistr\Poi(np_i) $.
    
    Conversely, for $P$ that is uniform over $[k]$, for any fixed $\delta \in (0,1)$, if $n \leq (1-\delta) k \log \frac{k}{\Delta}$, then as $k\diverge$,
    \begin{equation}
        \label{eq:plug-converse}
        \Prob[|S(P)-\Splug(N)|\leq \Delta]\leq e^{-\Omega(k^\delta)}.
    \end{equation}    
\end{proposition}

In order to remedy the inaccuracy of the plug-in estimate $ \Splug $ in the regime of $ n\lesssim k\log k $, our proposed estimator adds a linear correction term:
\begin{equation}
    \label{eq:uj}
    \hat{S}=\Splug+\sum_{j\ge 1}u_j\Phi_j,
\end{equation}
where the coefficients $ u_j$'s are to be specified. Equivalently, the estimator can be expressed in terms of the histogram as
\begin{equation}
    \hat{S}=\sum_i g(N_i) 
    \label{eq:hatS}
\end{equation}
where $g:\integers_+\to\reals$ is given by $ g(j)=u_j+1 $ for $ j\ge 1 $ and $ g(0)=0 $. Then the bias of $ \hat{S} $ is
\begin{equation}
    \label{eq:bias-term}
    \Expect[\hat{S}-S]
    =\sum_{i:p_i>0}e^{-np_i}\pth{\sum_{j\ge 1}u_j\frac{(np_i)^j}{j!}-1}
    \triangleq \sum_{i:p_i>0}e^{-np_i}P(p_i),
\end{equation}
where $ P(0)=-1 $ by design. Therefore the bias of $\hat S$ is at most
\[
    S\max_{x\in[p_{\min},1]}|e^{-nx}P(x)|,
\]
and the variance can be upper bounded by $ 2S\Norm{g}_\infty^2 $ using the Efron-Stein inequality \cite{steele86}. Next we choose the coefficients in order to balance the bias and variance.

The construction is done using Chebyshev polynomials $T_L$ \prettyref{eq:cheby}.
Note that $T_L$ is bounded in magnitude by one over the interval $[-1,1]$. The shifted and scaled Chebyshev polynomial over the interval $ [l,r] $ is given by
\begin{equation}
    P_L(x)
    =-\frac{T_L(\frac{2x-r-l}{r-l})}{T_L(\frac{-r-l}{r-l})}
    \triangleq \sum_{m=1}^{L}a_mx^m-1,
    \label{eq:PL}
\end{equation}
the coefficients $a_1,\dots,a_L $ can be obtained from those of the Chebyshev polynomial \cite[2.9.12]{timan63} and the binomial expansion, or more directly,
\begin{equation}
    \label{eq:aj}
    a_j=\frac{P_L^{(j)}(0)}{j!}=-\pth{\frac{2}{r-l}}^j\frac{1}{j!}\frac{T_L^{(j)}(-\frac{r+l}{r-l})}{T_L(-\frac{r+l}{r-l})}.
\end{equation}
We now let
\begin{equation}
    \label{eq:constants-ref}
    L \triangleq \floor{c_0\log k}, \quad r\triangleq\frac{c_1\log k}{n}, \quad l\triangleq\frac{1}{k},
\end{equation} 
where $ c_0<c_1 $ are constants to be specified, and choose the coefficients of the estimator as 
\begin{equation}
    \label{eq:uj-choice}
    u_j=
    \begin{cases}
        \frac{a_jj!}{n^j} & j=1,\dots,L  \\
        0 & \text{otherwise},\\
    \end{cases}
\end{equation}
The estimator $\hat S$ is defined according to \prettyref{eq:uj}.

By choosing the Chebyshev polynomial of degree $L$, we only use the first $L$ fingerprints to estimate the unseen $\Phi_0$ by setting $u_j = 0$ in \eqref{eq:uj} for all $j > L$, which while possibly incurring bias, reduces the variance. 
A further reason for only using the first few fingerprints is that higher-order fingerprints are \emph{almost uncorrelated} with the number of unseens $\Phi_0$. For instance, if $n/2$ balls drawn from an urn are red, the only information this reveals is that approximately half of the urn are red. In fact, the correlation between $\Phi_0$ and $\Phi_j$ decays exponentially (see \cite[Appendix C]{WY2016sample}). Therefore for $L=\Theta(\log k)$, $\{\Phi_j\}_{j > L}$ offer little predictive power about $\Phi_0$.

We proceed to explain the reasoning behind the choice \prettyref{eq:uj-choice} and the role of the Chebyshev polynomial. The main intuition is that since $c_0<c_1$, then with high probability, whenever $ N_i\le L=\floor{c_0\log k} $ the corresponding mass must satisfy $ p_i\le \frac{c_1\log k}{n} $. That is, if $ p_i>0 $ and $ N_i\le L $ then $ p_i\in [l,r] $ with high probability, and hence $ P_L(p_i) $ is bounded by the sup-norm of $ P_L $ over the interval $ [l,r] $, which controls the bias in view of \prettyref{eq:bias-term}. In view of the extremal property of Chebyshev polynomials \cite[Ex.~2.13.14]{timan63}, \prettyref{eq:PL} is the unique degree-$L$ polynomial that passes through the point $(0,-1)$ and deviates the least from zero over the interval $ [l,r] $. This explains the coefficients \prettyref{eq:hatS} which are chosen to minimize the bias. The degree of the polynomial is only logarithmic so that the variance is small.

The next proposition gives an upper bound of the quadratic risk of our estimator \prettyref{eq:hatS}:
\begin{proposition}
    \label{prop:main-ub-poly}
    Assume the Poissonized sampling model where the histograms are distributed as  $ N=(N_1,N_2,\dots)\inddistr\Poi(np_i) $. Let $ c_0 = 0.558$ and $c_1=0.5 $. As $ \delta\triangleq \frac{n}{k\log k}\rightarrow 0 $ and $ k\diverge $, the bias and variance of $ \hat{S} $ are upper bounded by
    \begin{align*}
        & |\Expect(\hat S-S)|\le 2S(1+o_{k}(1))\exp\pth{-(1+o_\delta(1))\sqrt{\kappa\frac{n\log k}{k}}},\\
        & \var[\hat{S}]\le O\pth{S k^c},
    \end{align*}
    for some absolute constant $ c<1 $, and consequently,
    \begin{equation}
        \label{eq:S-main-ub-poly}
        \sup_{P\in\calD_k}\Expect(\hat S(N)-S(P))^2\le  4k^2(1+o_{k}(1))\exp\pth{-(2+o_\delta(1))\sqrt{\kappa\frac{n\log k}{k}}},
    \end{equation}
    where $\kappa = 2.494$.
\end{proposition}

The minimax upper bounds in Theorems \ref{thm:main} and \ref{thm:sample} follow from combining Propositions \ref{prop:main-ub-plug} and \ref{prop:main-ub-poly}.

Note that the optimal estimator \prettyref{eq:hatS} relies on the choice of parameters in \prettyref{eq:constants-ref}, which, in turn, depends on the knowledge of $1/k$, the lower bound on the minimum non-zero probability $p_{\min}$. While $k$ is readily obtainable in certain applications where the sample is uniformly drawn from a database or corpus of known size (see \cite{BJKST02,ET76} as well as the empirical results in \prettyref{sec:S-exp}), it is desirable to construct estimators that are agnostic to $ p_{\min} $ and retains the same optimality guarantee. To this end, we provide the following alternative choice of parameters. Let $\tilde{S}$ be the linear estimator defined using the same coefficients in \prettyref{eq:uj-choice}, with the approximation interval $ [l,r] $ and the degree $L$ in \prettyref{eq:constants-ref} replaced by 
\begin{equation}
    \label{eq:constants-ref-adapt}
    l=\frac{c_1}{c_0^2}\frac{\log^2(1/\epsilon)}{n\log n}, \quad r=\frac{c_1\log n}{n}, \quad L=\Floor{c_0\log n}.
\end{equation}
Here $\epsilon$ is the desired relative accuracy and the constants $c_0,c_1$ are the same as in \prettyref{prop:main-ub-poly}. Following the same analysis as in the proof of \prettyref{prop:main-ub-poly}, the above choice of parameters leads to the following upper bound of the quadratic risk:
\begin{proposition}
    \label{prop:alternative}
    Let $ c_0,c_1,c $ be the same constants as \prettyref{prop:main-ub-poly}. There exist constants $C,C'$ such that, if $ \epsilon>n^{-C} $, then
    \[
        \Expect(\tilde{S}-S)^2 \leq 
        C' (S^2 \epsilon^{2(1-\sqrt{\alpha})} + Sn^{c}),
    \]
    where $\alpha = \max\pth{1-\frac{c_0^2}{c_1}\frac{n\log n}{k \log^2(1/\epsilon) }, 0}$.
\end{proposition}
Therefore, whenever the sample size satisfies $n \geq (\frac{c_1}{c_0^2}+o_k(1))\frac{k}{\log k} \log^2 \frac{1}{\epsilon}$ and $ n\le (\epsilon^2 k)^{\frac{1}{c}} $, the upper bound is at most $O((\epsilon k)^2)$, recovering the optimal risk bound in \prettyref{prop:main-ub-poly}. The new result in \prettyref{prop:alternative} is that even when $n$ is not that large the risk degrades gracefully.

For proofs of Propositions \ref{prop:main-ub-plug}--\ref{prop:alternative}, see \cite{WY15}.
We finish this subsection with a few remarks.
\begin{remark}
    Combined with standard concentration inequalities, the mean-square error bound in \prettyref{prop:main-ub-poly} can be easily converted to a high-probability bound. In the regime of $ n\lesssim k\log k $, for any distribution $P\in\calD_k $, the bias of our estimate $ \hat{S} $ is at most the uniform approximation error:
    \begin{equation*}
        |\Expect[\hat{S}]-S|\le S\exp\pth{-\Theta\pth{\sqrt{\frac{n\log k}{k}}}}.
    \end{equation*}
    The standard deviation is significantly smaller than the bias. Indeed, the coefficients of the linear estimator \prettyref{eq:hatS} is uniformly bounded by $ \Norm{g}_\infty^2\le k^{c} $ for some absolute constant $ c<1 $ (see \prettyref{fig:coeffs} for numerical results). Therefore, by Hoeffding's inequality, we have the following concentration bound:
    \begin{equation*}
        \Prob[|\hat{S}-\Expect[\hat{S}]|\ge t k]\le 2\exp\pth{-\frac{t^2k}{2\Norm{g}_\infty^2}}=\exp\pth{-t^2k^{\Omega(1)}}.
    \end{equation*}
		This type of concentration result was crucially used in  \cite{ADOS17} to 
		show, via an indirect method, that profile maximum likelihood also achieves the optimal rate of estimating the support size (and other functionals).
\end{remark}

\begin{remark}
    Estimators of the form \prettyref{eq:hatS} is also known as \emph{linear estimators}:
    \begin{equation}
        \label{eq:linear}
        \hat S = \sum_i g(N_i) =\sum_{j \geq 1} g(j) \Phi_j ,       
    \end{equation} 
    which is a linear combination of fingerprints $\Phi_j$'s defined in \prettyref{eq:fp}.

    Other notable examples of linear estimators include:
    \begin{itemize}
        \item Plug-in estimator \prettyref{eq:splug}: $\Splug = \Phi_1+\Phi_2+\dots$.
        \item Good-Toulmin estimator \cite{GT56}: for some $t>0$,
        \begin{equation}
            \label{eq:GT}
            \hat S_{\rm GT} = \Splug + 
            t \Phi_1 - t^2 \Phi_2 + t^3 \Phi_3 - t^4 \Phi_4 + \ldots
        \end{equation}
        \item Efron-Thisted estimator \cite{ET76}: for some $t>0$ and $J\in\naturals$,
        \begin{equation}
            \label{eq:ET}
            \hat S_{\rm ET} = \Splug + \sum_{j=1}^J (-1)^{j+1} t^j b_j \Phi_j,
        \end{equation}
        where $b_j = \Prob[\Binom(J,1/(t+1)) \geq j]$.
    \end{itemize}
    By definition, our estimator \prettyref{eq:hatS} can be written as
    \begin{equation}
        \label{eq:Shat1}
        \hat S =  \sum_{j=1}^L g(j) \Phi_j + \sum_{j > L} \Phi_j.
    \end{equation}
    By \prettyref{eq:PL}, $ P_L $ is also a polynomial of degree $ L $, which is oscillating and results in coefficients with alternating signs (see \prettyref{fig:coeffs}).
    \begin{figure}[hbt]
        \centering
        \includegraphics[width=.6\linewidth]{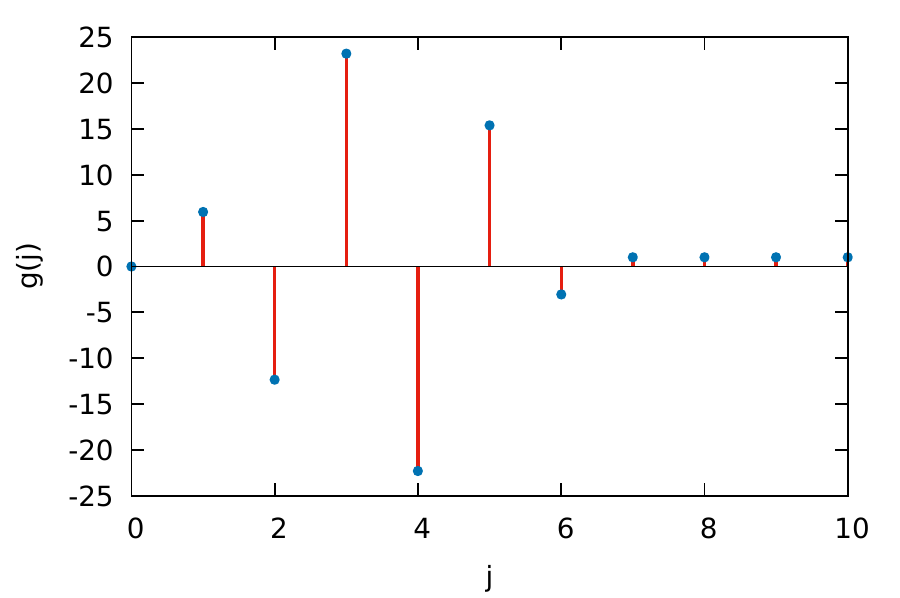}
        \caption{Coefficients of estimator $ g(j) $ in \prettyref{eq:hatS} with $ c_0=0.45,c_1=0.5 $, $ k=10^6$ and $ n=2\times 10^5$.\label{fig:coeffs} }
    \end{figure}
    Interestingly, this behavior, although counterintuitive, coincides with many classical estimators, such as \prettyref{eq:GT} and \prettyref{eq:ET}. The occurrence of negative coefficients can be explained as follows. Note that the rationale of the linear estimator is to form a linear prediction of the number of unseen $\Phi_0$ using the observed fingerprints $\Phi_1,\Phi_2,\ldots$; this is possible because the fingerprints are correlated. Indeed, since the sum of all fingerprints coincides with the support size, i.e., $\sum_{j\geq 0} \Phi_j = S$, for each $j\geq 1$, the random variable $ \Phi_j $ is negatively correlated with $ \Phi_0 $ and hence some of the coefficients in the linear estimator are negative. 
\end{remark}

\begin{remark}[Time complexity]
    \label{rmk:complexity}
    The evaluation of the estimator \prettyref{eq:linear} consists of three parts:
    \begin{enumerate}
        \item Construction of the estimator: $ O(L^2)=O(\log^2k) $, which amounts to computing the coefficients $ g(j) $ per \prettyref{eq:aj};
        \item Computing the histograms $ N_i $ and fingerprints $ \Phi_j $: $ O(n) $;
        \item Evaluating the linear combination: $ O(n\wedge k) $, since the number of non-zero terms in the second summation of \prettyref{eq:linear} is at most $ n\wedge k $. 
    \end{enumerate}
    Therefore the total time complexity is $ O(n+\log^2k) $.
\end{remark}

\subsection{Experiments}
\label{sec:S-exp}
We evaluate the performance of our estimator on both synthetic and real datasets in comparison with popular existing procedures.\footnote{
    The implementation of our estimator is available at \url{https://github.com/Albuso0/support}.
}
In the experiments we choose the constants $ c_0=0.45, c_1=0.5 $ in \prettyref{eq:constants-ref}, instead of $c_0=0.558$ which is optimized to yield the best rate of convergence in \prettyref{prop:main-ub-poly} under the \iid~sample model. The reason for such a choice is that in the real-data experiments the observations are not necessarily generated independently and dependency leads to a higher variance. By choosing a smaller $c_0$, the Chebyshev polynomials have a slightly smaller degree, which results in smaller variance and more robustness to model mismatch. Each experiment is averaged over $ 50 $ independent trials and the standard deviations are shown as error bars.

\paragraph{Synthetic data}
We consider data independently sampled from the following distributions:
\begin{itemize}
    \item the uniform distribution with $ p_i=\frac{1}{k} $;
    \item Zipf distributions with $ p_i\propto i^{-\alpha} $ and $ \alpha$ being either $1$ or $0.5 $;
    \item an even mixture of geometric distribution and Zipf distribution where for the first half of the alphabet $ p_i \propto 1/i$ and  for the second half $ p_{i+k/2} \propto (1-\frac{2}{k})^{i-1} $, $ 1\le i\le \frac{k}{2} $.
\end{itemize}
The alphabet size $ k $ varies in each distribution so that the minimum non-zero mass is roughly $ 10^{-6} $. Accordingly, a degree-6 Chebyshev polynomial is applied. Therefore, according to \prettyref{eq:Shat1}, we apply the polynomial estimator $ g $ to symbols appearing at most six times and the plug-in estimator otherwise. 
\begin{figure}[ht]
    \centering
    \includegraphics[width=1.0\linewidth]{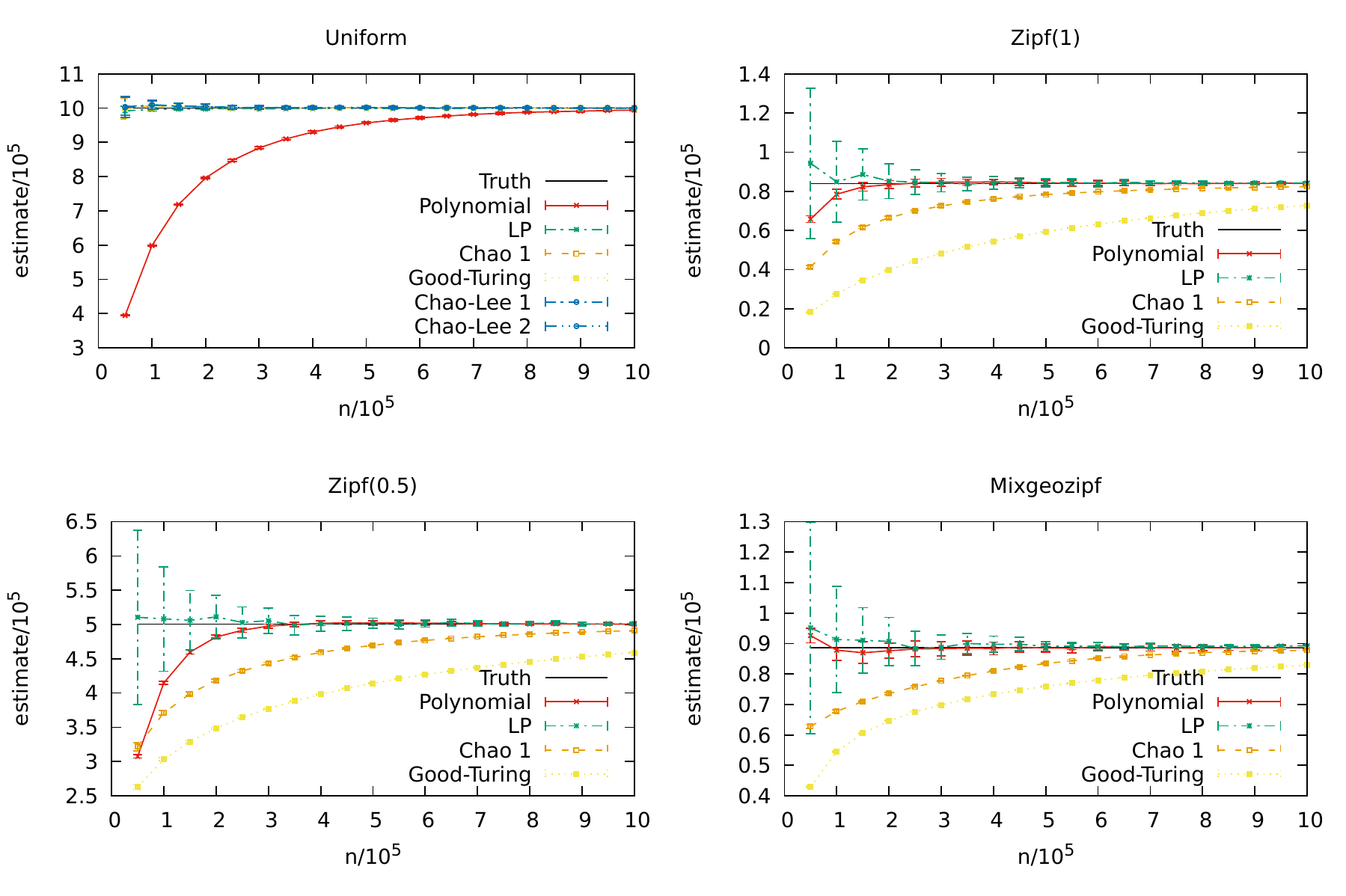}
    \caption{Performance comparison under four data-generating distributions.\label{fig:synthetic} }
\end{figure}
We compare our results with the Good-Turing estimator \cite{Good1953}, the Chao 1 estimator \cite{Chao84,GC11}, the two estimators proposed by Chao and Lee \cite{Chao92}, and the linear programming approach proposed by Valiant and Valiant \cite{VV13}. The results are shown in \prettyref{fig:synthetic}. Here the Good-Turing estimator refers to first estimate the total probability of seen symbols (sample coverage) by $\hat C=1- \frac{\Phi_1}{n}$ then estimate the support size by $\SGT = {\Splug}/{\hat C}$; the Chao 1 estimator refers to the bias-corrected form $ \SChao = \Splug+\frac{\Phi_1(\Phi_1-1)}{2(\Phi_2+1)} $. The plug-in estimator simply counts the number of distinct elements observed, which is always outperformed by the Good-Turing estimator in our experiments and hence omitted in the comparison.

Good-Turing's estimate on sample coverage performs remarkably well in the special case of uniform distributions. This has been noticed and analyzed in \cite{Chao92,DR80}. Chao-Lee's estimators are based on Good-Turing's estimate with further correction terms for non-uniform distributions. However, with limited sample size, if no symbol appears more than once, the sample coverage estimate $\hat C$ is zero and consequently the Good-Turing estimator and Chao-Lee estimators are not even well-defined. For Zipf and mixture distributions, the output of Chao-Lee's estimators is highly unstable and thus is omitted from the plots; the convergence rates of Good-Turing estimator and Chao 1 estimator are much slower than our estimator and the LP estimator, partly because they only use the information of how many symbols occurred exactly once and twice, namely the first two fingerprints $\Phi_1$ and $ \Phi_2 $, as opposed to the full spectrum of fingerprints $\{\Phi_j\}_{j\geq 1}$\footnote{In fact, both  Good-Turing and Chao 1 estimator fail to achieve the optimal rate in \prettyref{thm:main} and they suffer provably large bias under non-uniform distributions as simple as mixtures of two uniform distributions; cf.~\cite[Appendix C]{WY15}.}.
The linear programming approach has a similar convergence rate to ours but suffers from large variance when the sample size is small.
Finally, we mention that the bias of our estimators can be reduced by using a higher degree approximation at the cost of an increased variance, and we choose a slightly smaller degree such that our estimator is more robust to model mismatch.

\begin{figure}[ht]
    \centering
    \includegraphics[width=0.7\linewidth]{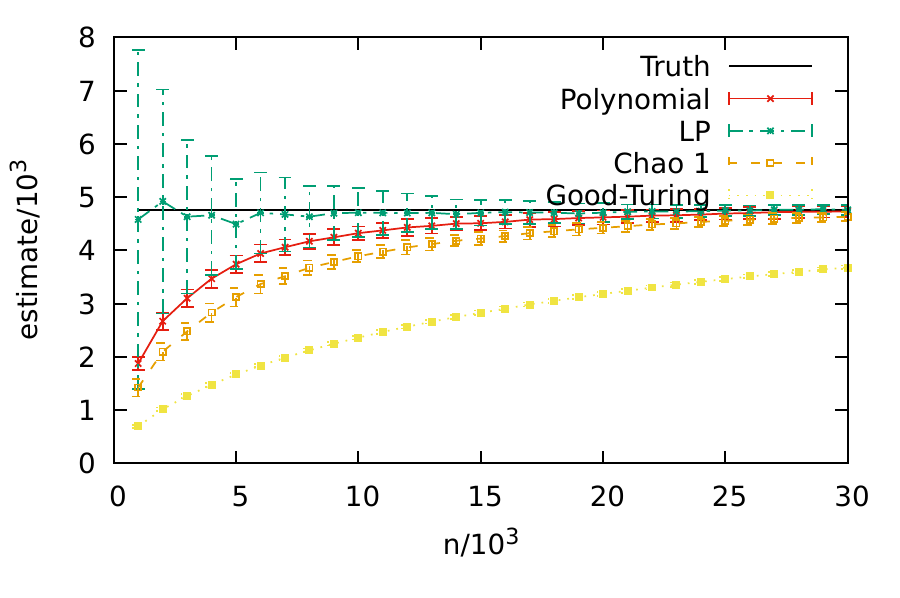}
    \caption{Comparison of various estimates of the total number of distinct words in \emph{Hamlet}. \label{fig:estimation} }
\end{figure}

\paragraph{Real data}
Next we evaluate our estimator by a real data experiment based on the text of \emph{Hamlet}, which contains about $ 32,000 $ words in total consisting of about $ 4,800 $ distinct words. Here and below the definition of ``distinct word'' is any distinguishable arrangement of letters that are delimited by spaces, insensitive to cases, with punctuations removed. We randomly sample the text with replacement and generate the fingerprints for estimation. The minimum non-zero mass is naturally the reciprocal of the total number of words, $ \frac{1}{32,000} $. In this experiment we use the degree-$ 4 $ Chebyshev polynomial. We also compare our estimator with the one in \cite{VV13}. The results are plotted in \prettyref{fig:estimation},
which shows that the estimator in \cite{VV13} has similar convergence rate to ours; however, the variance is again much larger and the computational cost of linear programming is significantly higher than linear estimators, which amounts to computing linear combinations with pre-determined coefficients.

\begin{figure}[ht]
    \centering
    \includegraphics[width=0.7\linewidth]{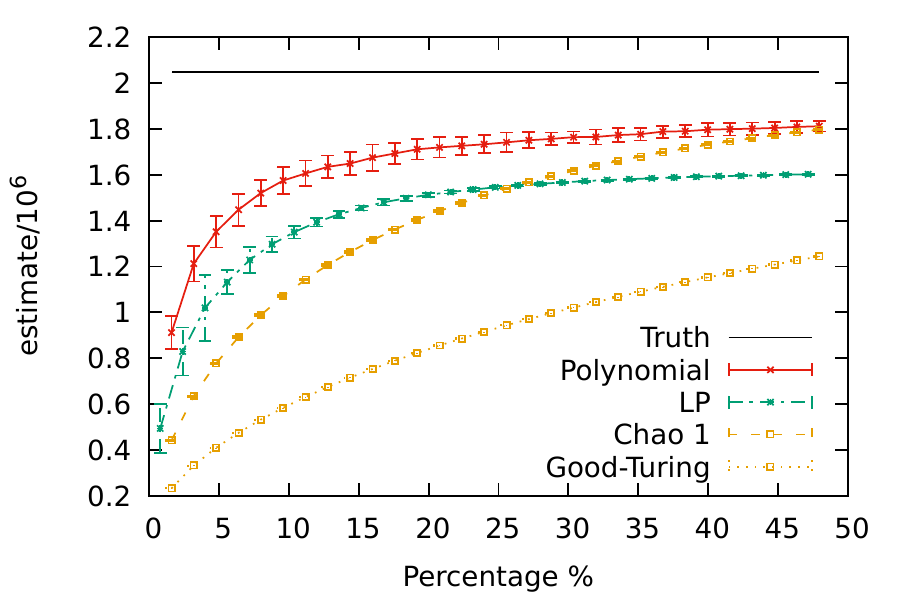}
    \caption{Performance comparison using \emph{New York Times Corpus}. \label{fig:NYT} }
\end{figure}

Next we conduct a larger-scale experiment using the \emph{New York Times Corpus} from the years 1987 -- 2007.\footnote{
    Dataset available at \url{https://catalog.ldc.upenn.edu/LDC2008T19}.
}
This corpus has a total of 25,020,626 paragraphs consisting of 996,640,544 words with 2,047,985 distinct words. We randomly sample 1\% -- 50\% out of the all paragraphs with replacements and feed the fingerprint to our estimator. The minimum non-zero mass is also the reciprocal of the total number of words, $ 1/10^9 $, and thus the degree-9 Chebyshev polynomial is applied. Using only 20\% of the sample our estimator achieves a relative error of about 10\%, which is a systematic error due to the model mismatch: the sampling here is paragraph by paragraph rather than word by word, which induces dependence across the observations as opposed to the \iid~sampling model for which the estimator is designed; in comparison, the LP estimator\footnote{
    In this large-scale experiment, the original MATLAB code of the linear programming estimator given in \cite{VV13} is extremely slow; the result in \prettyref{fig:NYT} is obtained using an optimized version provided by the author \cite{Valiant-email}. 
} 
suffers a larger bias from this model mismatch. Furthermore, the proposed linear estimator is significantly faster than linear programming based methods: given the sampled data, the curve in \prettyref{fig:NYT} corresponding to the LP estimator takes over 5 hours to compute; in contrast, the proposed linear estimator takes only 2 seconds on the same computer, which clearly demonstrate its computational advantage even if one takes into account the fact that our implementation is based on C++ while the LP estimator is in MATLAB.

Finally, we perform the classical experiment of ``how many words did Shakespeare know''. We feed the fingerprint of the entire Shakespearean canon (see \cite[Table 1]{ET76}), which contains 31,534 word types, to our estimator. We choose the minimum non-zero mass to be the reciprocal of the total number of English words, which, according to known estimates, is between 600,000 \cite{OED} to 1,000,000 \cite{GLM}, and obtain an estimate of 63,148 to 73,460 for Shakespeare's vocabulary size, as compared to 66,534 obtained by Efron-Thisted \cite{ET76}. Using the alternative choice of parameters that are agnostic to $k$ in \prettyref{prop:alternative}, by setting the desired accuracy to be $0.05$ and $0.1$, we obtain an estimate of 62,355 to 72,454.

\subsection{Minimax lower bound}
\label{sec:minimax}
The lower bound argument follows the idea in \cite{LNS99,CL11,WY14} and relies on the generalized Le Cam's method involving two composite hypothesis introduced in \prettyref{thm:lb-two-priors}. The main idea is similar to \prettyref{sec:h-lb}. Specifically, suppose the following (composite) hypothesis testing problem
\[
    H_0: S(P) \leq s, P \in \calD_k  \quad \text{versus} \quad H_1: S(P) \geq s+\delta, P \in\calD_k
\]
cannot be tested with vanishing probability of error on the basis of $n$ \iid~observations, then the sample complexity of estimating $S(P)$ within $\delta$ with high probability must exceed $n$. In particular, the impossibility to test the above composite hypotheses is shown by constructing two priors (i.e., two random probability vectors) so that the induced distributions of the sample are close in total variation. Next we elaborate the main ingredients of \prettyref{thm:lb-two-priors}:
\begin{itemize}
    \item construction of the two priors;
    \item separation between functional values;
    \item bound on the total variation.
\end{itemize}

Let $\lambda > 1$. Given unit-mean random variables $ U$ and $U'$ that take values in $\{0\} \cup [1,\lambda]$, define the following random vectors
\begin{equation}
    \label{eq:iid-PP}
    \sfP=\frac{1}{k}(U_1,\dots,U_k),\quad \sfP'=\frac{1}{k}(U_1',\dots,U_k'),
\end{equation}
where $ U_i $ and $U_i' $ are \iid~copies of $ U$ and $U' $, respectively. Although $ \sfP $ and $ \sfP' $ need not be probability distributions, as long as the standard deviations of $U$ and $U'$ are not too big, the law of large numbers ensures that with high probability $ \sfP$ and $\sfP' $ lie in a small neighborhood near the probability simplex, which we refer to as the set of \emph{approximate} probability distributions. Furthermore, the minimum non-zeros in $\sfP$ and $\sfP'$ are at least $\frac{1}{k}$. It can be shown that the minimax risk over approximate probability distributions is close to that over the original parameter space $\calD_k$ of probability distributions. This allows us to use $ \sfP$ and $\sfP' $ as priors and apply \prettyref{thm:lb-two-priors}. Note that both $ S(\sfP) $ and $ S(\sfP') $ are binomially distributed, which, with high probability, differ by the difference in their mean values:
\begin{align*}
    \Expect[S(\sfP)]-\Expect[S(\sfP')]
    =k(\Prob[U>0]-\Prob[U'>0])
    =k(\Prob[U'=0]-\Prob[U=0]).
\end{align*}
If we can establish the impossibility of testing whether data are generated from $\sfP$ or $\sfP'$, the resulting lower bound is proportional to $k(\Prob[U'=0]-\Prob[U=0])$.

To simplify the argument we apply the Poissonization technique where the sample size is a $ \Poi(n) $ random variable instead of a fixed number $n$. This provably does not change the statistical nature of the problem due to the concentration of $\Poi(n)$ around its mean $ n $. Under Poisson sampling, the histograms \prettyref{eq:histogram} still constitute a sufficient statistic, which are distributed as $ N_i\inddistr \Poi(np_i) $, as opposed to multinomial distribution in the fixed-sample-size model. Therefore through the \iid~construction in \prettyref{eq:iid-PP}, $ N_i\iiddistr \Expect[\Poi(\frac{n}{k}U)] $ or $ \Expect[\Poi(\frac{n}{k}U')] $. Then \prettyref{thm:lb-two-priors} is applicable if $   \TV(\Expect[\Poi(\frac{n}{k} U)]^{\otimes k},\Expect[\Poi(\frac{n}{k} U')]^{\otimes k})$ is strictly bounded away from one, for which it suffices to show
\begin{equation}
    \label{eq:tv-bd}
    \TV(\Expect[\Poi(nU/k)],\Expect[\Poi(nU'/k)])\le \frac{c}{k},
\end{equation}
for some constant $ c<1 $.

The above construction provides a recipe for the lower bound. To optimize the ingredients it boils down to the following optimization problem (over one-dimensional probability distributions): Construct two priors $ U,U' $ with unit mean that maximize the difference $ \prob{U'=0}-\prob{U=0} $ subject to the total variation distance constraint \prettyref{eq:tv-bd}, which, by \prettyref{thm:tv-bound}, can be guaranteed by \emph{moment matching}, \ie, ensuring $ U$ and $U' $ have identical first $L$ moments for some large $L$, and the $L_\infty$-norms of $ U,U' $ are not too large such that $|\frac{nU}{k}|,|\frac{nU'}{k}|\lesssim L$. To summarize, our lower bound entails solving the following optimization problem:
\begin{equation}
    \label{eq:FL}
    \begin{aligned}
        \sup & ~ \Prob[U'=0]-\Prob[U=0]  \\
        \text{s.t.}     
        & ~ \Expect[U] = \Expect[U']=1 \\
        & ~ \Expect[U^j] = \Expect[U'^j], \quad j = 1,\ldots,L \\
        & ~ U,U' \in \sth{0}\cup[1, \lambda].
    \end{aligned}
\end{equation}
The final lower bound is obtained from \prettyref{eq:FL} by choosing $L \asymp \log k$ and $\lambda \asymp \frac{k\log k}{n}$.

In order to evaluate the infinite-dimensional linear programming problem \prettyref{eq:FL}, we consider its dual program.
It is shown in \eqref{eq:poly-primal} that the problem of best polynomial and moment matching are dual to each other; however, unlike the standard moment matching problem which impose the equality of moments, the extra constraint in \prettyref{eq:FL} is that the values of the first moment must equal to one. Therefore its dual is no longer the best polynomial approximation problem. Nevertheless, for the specific problem \prettyref{eq:FL} which deals with the indicator function $x \mapsto \indc{x=0}$, via a change of variable we show in \prettyref{sec:opt-UX} that \prettyref{eq:FL} coincides exactly with the best uniform approximation error of the function $ x\mapsto\frac{1}{x} $ over the interval $ [1,\lambda] $ by degree-$ (L-1) $ polynomials:
\begin{equation*}
    \inf_{p\in \calP_{L-1}}\sup_{x\in [1,\lambda]}\abs{\frac{1}{x}-p(x)},
\end{equation*}
where $ \calP_{L-1} $ denotes the set of polynomials of degree at most $ L-1 $. This best polynomial approximation problem has been well-studied, cf.~\cite{timan63,DS08}; in particular, the exact formula for the best polynomial that approximates $ x\mapsto\frac{1}{x} $ and the optimal approximation error have been obtained in \cite[Sec. 2.11.1]{timan63}.

Applying the procedure described above, we obtain the following sample complexity lower bound:
\begin{proposition}
    \label{prop:main-lb2}
    Let $ \delta\triangleq \frac{\log(k/\Delta)}{\log k}$ and $ \tau\triangleq \frac{\sqrt{\log k}/k^{1/4}}{1-(2\Delta/k)} $. As $ k\rightarrow\infty $, $\delta \to 0$ and $\tau\to 0$, 
    \begin{equation}
        n_{\sf S}^*(k,\Delta) \ge (1-o_{\delta}(1)-o_{k}(1)-o_\tau(1))\frac{k}{2e^2\log k}\log^2\frac{k}{2\Delta}.
        \label{eq:S-main-lb2}
    \end{equation}
    Consequently,    if $ k^{1-c}\le \Delta \le \frac{k}{2}-c'k^{3/4}\sqrt{\log k} $ for some constants $ c,c' $, then $ n_{\sf S}^*(k,\Delta)\gtrsim \frac{k}{\log k}\log^2\frac{k}{2\Delta} $.
\end{proposition}
The lower bounds announced in Theorems \ref{thm:main} and \ref{thm:sample} follow from \prettyref{prop:main-lb2} combined with a simple two-point argument. See \cite{WY15}.

\subsection{Dual program of \eqref{eq:FL}}
\label{sec:opt-UX}
Define the following infinite-dimensional linear program:
\begin{equation}
    \begin{aligned}
        \calE_1^*\triangleq
        \sup & ~ \prob{U'=0}-\prob{U=0}  \\
        \text{s.t.}     
        & ~ \expect{U} = \expect{U'}=1 \\
        & ~ \expect{U^j} = \expect{U'^j}, \quad j = 1,\ldots,L+1, \\
        & ~ U,U' \in \sth{0}\cup I,
    \end{aligned}
    \label{eq:FL1}
\end{equation}
where $ I=[a,b] $ with $ b>a\ge 1 $ and the variables are probability measures on $I$ (distributions of the random variables $U,U'$). Then \eqref{eq:FL} is a special case of \prettyref{eq:FL1} with $ I=[1,\lambda] $. 
\begin{lemma}
    \label{lmm:FL}
    $\calE_1^* =2\inf_{p\in\calP_L}\sup_{x\in I}\abs{\frac{1}{x}-p(x)}$.
\end{lemma}
\begin{proof}
    We first show that \prettyref{eq:FL} coincides with the following optimization problem:
    \begin{equation}
        \begin{aligned}
            \calE_2^*\triangleq
            \sup & ~ \expect{\frac{1}{X}}-\expect{\frac{1}{X'}}  \\
            \text{s.t.}     
            & ~ \expect{X^j} = \expect{X'^j}, \quad j = 1,\ldots,L, \\
            & ~ X,X' \in I.
        \end{aligned}
        \label{eq:FL2}
    \end{equation}
    Given any feasible solution $ U,U' $ to \prettyref{eq:FL}, construct $ X,X' $ with the following distributions:
    \begin{equation}
        \begin{aligned}
            &P_X(\diff x)=xP_{U}(\diff x),\\
            &P_{X'}(\diff x)=xP_{U'}(\diff x),
        \end{aligned}
        \label{eq:U-X}
    \end{equation}
    It is straightforward to verify that $ X,X' $ are feasible for \prettyref{eq:FL2} and
    \begin{equation*}
        \calE_2^*
        \ge \expect{\frac{1}{X}}-\expect{\frac{1}{X'}}
        =\prob{U'=0}-\prob{U=0}.
    \end{equation*}
    Therefore $ \calE_2^*\ge \calE_1^* $.

    On the other hand, given any feasible $ X,X' $ for \prettyref{eq:FL2}, construct $ U,U' $ with the distributions:
    \begin{equation}
        \begin{aligned}
            &P_U(\diff u)=\pth{1-\expect{\frac{1}{X}}}\delta_0(\diff u)+\frac{1}{u}P_{X}(\diff u),\\
            &P_{U'}(\diff u)=\pth{1-\expect{\frac{1}{X'}}}\delta_0(\diff u)+\frac{1}{u}P_{X'}(\diff u),
        \end{aligned}
        \label{eq:X-U}
    \end{equation}
    which are well-defined since $ X,X'\ge 1 $ and hence $ \expect{\frac{1}{X}}\le 1, \expect{\frac{1}{X'}}\le 1 $.
    Then $ U,U' $ are feasible for \prettyref{eq:FL} and hence 
    \begin{equation*}
        \calE_1^*
        \ge \prob{U'=0}-\prob{U=0}
        =\expect{\frac{1}{X}}-\expect{\frac{1}{X'}}.
    \end{equation*}
    Therefore $ \calE_1^*\ge \calE_2^* $. Finally, the dual of \prettyref{eq:FL2} is precisely the best polynomial approximation problem (see \prettyref{eq:dual-mm}).
\end{proof}

\section{Distinct elements problem}
\label{sec:distinct}

The \DistinctElements problem \cite{CCMN00} refers to the following:
\begin{quote}
    \emph{Given $ n $ balls randomly drawn from an urn containing $ k $ colored balls, how to estimate the total number of distinct colors in the urn?}
\end{quote}
This is an instance of 
the \SupportSize problem considered in \prettyref{sec:supp}, in that each probability mass $p_i$ has the special form $p_i=\frac{k_i}{k}$ with $k_i\in\integers_+$.
Originating from ecology, numismatics, and linguistics, the \DistinctElements problem is also known as the \emph{species problem} in the statistics literature \cite{Lo92,BF93}. Apart from the theoretical interests, it has a wide array of applications in various fields, such as estimating the number of species in a population of animals \cite{FCW43,Good1953}, the number of dies used to mint an ancient coinage \cite{Esty86}, and the vocabulary size of an author \cite{ET76}. In computer science, this problem frequently arises in large-scale databases, network monitoring, and data mining \cite{RRSS09,BJKST02,CCMN00}, where the objective is to estimate the types of database entries or IP addresses from limited observations, since it is typically impossible to have full access to the entire database or keep track of all the network traffic. 

Similar to the \SupportSize problem, the key challenge in the \DistinctElements problem is: given a small sample where most of the colors are not observed, how to accurately extrapolate the number of unseens? 
Analogous to \prettyref{eq:S-risk} and \prettyref{eq:sample-support}, we define the minimax risk and the sample complexity of the \DistinctElements problem as $R_{\sf DS}^*(k,n)$ and  $n^*_{\sf DS}(k,\Delta)$, respectively, where the supremum is taken over all distributions $P$ that arise from $k$-ball urns instead\footnote{
In the presentation of our main results we focus on sampling with replacement. The main results hold for sampling without replacements as well. See \cite[Appendix A]{WY2016sample} for details.
}.
By definition, we have $R_{\sf DS}^* \leq R_{\sf S}^*$ and  $n^*_{\sf DS} \leq n^*_{\sf S}$. 
In particular, by \prettyref{eq:sample-complexity}, we have
\begin{equation}
    \label{eq:sample-supp}
    n^*_{\sf DS}(k,\Delta)\le O\pth{\frac{k}{\log k}\log^2\frac{k}{\Delta}}.
\end{equation}

Using the polynomial approximation methods in \prettyref{sec:functional}, the key is to approximate the function to be estimated by a polynomial, whose degree is chosen to balance the approximation error (bias) and the estimation error (variance). 
In \prettyref{eq:sample-complexity}, the worst-case performance guarantee for the \SupportSize problem is governed by the uniform approximation error over an interval where the probabilities may reside. In contrast, in the \DistinctElements problem, the observations are generated from a distribution supported on a \emph{discrete} set of values. Uniform approximation over a discrete subset leads to smaller approximation error and, in turn, improved sample complexity. It turns out that $ O(\frac{k}{\log k}\log\frac{k}{\Delta}) $ observations are sufficient to achieve an additive error of $ \Delta $ that satisfies $ k^{0.5+O(1)} \le \Delta \le O(k) $, which strictly improves \prettyref{eq:sample-supp} thanks to the discrete structure of the \DistinctElements problem.

\subsection{A summary of the sample complexity}

The main results of this section provide bounds and constant-factor approximations of the sample complexity in various regimes (see \prettyref{tab:main} for a summary), as well as computationally efficient algorithms. Below we highlight a few important conclusions drawn from \prettyref{tab:main}:
\begin{description}
    \item[From linear to sublinear: ] From the result for $ k^{0.5+\delta}\le \Delta \le ck $ in \prettyref{tab:main}, we conclude that the sample complexity is sublinear in $k$ if and only if $ \Delta=k^{1-o(1)} $, which also holds for sampling without replacement. To estimate within a constant fraction of balls $ \Delta=ck $ for any small constant $ c $, the sample complexity is $ \Theta(\frac{k}{\log k}) $, which coincides with the general support size estimation problem. However, in other regimes we can achieve better performance by exploiting the discrete nature of the \DistinctElements problem.
    
    \item[From linear to superlinear: ] The transition from linear to superlinear sample complexity occurs near $ \Delta=\sqrt{k} $. Although the exact sample complexity near $ \Delta=\sqrt{k} $ is not completely resolved in the current chapter, the lower bound and upper bound in \prettyref{tab:main} differ by a factor of at most $ \log\log k $. In particular, the estimator via interpolation can achieve $ \Delta=\sqrt{k} $ with sample size $ n=O(k\log\log k) $, and achieving a precision of $ \Delta\le k^{0.5-o(1)} $ requires strictly superlinear sample size.
\end{description}


\begin{table}[ht]
    \centering
    \begin{tabular}{|c|c|c|}  
      \toprule
      $ \Delta $                               & Lower bound                        & Upper bound                              \\
      \midrule
      $ \le 1 $                               & \multicolumn{2}{c|}{$ \Theta(k\log k) $}                                     \\
      \hline
                                               & \multicolumn{2}{c|}{\multirow{3}{*}{$ \Theta\pth{k\log \frac{k}{\Delta^2}} $}}   \\
      $ \qth{1,\sqrt{k}(\log k)^{-\delta}} $         & \multicolumn{2}{c|}{}                                                         \\
                                               & \multicolumn{2}{c|}{}     \\
      \cline{1-3}
      $ \qth{\sqrt{k}(\log k)^{-\delta},k^{0.5+\delta}} $ & $ \Omega\pth{k \pth{1\vee \log \frac{k}{\Delta^2}}} $ & $ O\pth{k\log\frac{\log k}{1 \vee \log\frac{\Delta^2}{k}}} $  \\
      \hline
                                                & \multicolumn{2}{c|}{\multirow{3}{*}{$ \Theta\pth{\frac{k}{\log k}\log\frac{k}{\Delta}} $}}  \\
      $[k^{^{0.5+\delta}},ck]$                     & \multicolumn{2}{c|}{} \\
                                                & \multicolumn{2}{c|}{}   \\
      \cline{1-3}
      $ [ck,(0.5-\delta)k] $                    & $k e^{-O(\sqrt{\log k \log \log k})} $\cite{RRSS09}  & $O\pth{\frac{k}{\log k}}$   \\
      \bottomrule
    \end{tabular}
    \caption[Summary]{Summary of the sample complexity $ n_{\sf DS}^*(k,\Delta) $, where $ \delta $ is any sufficiently small constant, $c$ is an absolute positive constant less than 0.5 (same over the table).
		The estimators are linear with coefficients obtained from either interpolation or $L_2$-approximation.
    \label{tab:main}}
\end{table}

In addition to the sample complexity, we also obtain the following characterization of the minimax mean squared error (MSE) of the \DistinctElements problem:
\begin{align}
  \frac{1}{k^2} R_{\sf DS}^*(k,n)
  & =\exp\sth{-\Theta\pth{\pth{1\vee \frac{n\log k}{k}} \wedge \pth{\log k \vee \frac{n}{k}} }}\\
  & =\begin{cases}
      \Theta(1) ,& n\le \frac{k}{\log k},\\
      \exp(-\Theta(\frac{n\log k}{k})) ,& \frac{k}{\log k}\le n\le k,\\
      \exp(-\Theta(\log k)) ,& k\le n \le k\log k,\\
      \exp(-\Theta(\frac{n}{k})) , &n \ge k\log k.
  \end{cases}
	\label{eq:distinct-rates}
\end{align}

\subsection{Linear estimators via discrete polynomial approximation}
\label{sec:distinct-linear}
%

Recall that $S$ denotes the number of distinct colors in a urn containing $k$ colored balls. 
Let $k_i$ denote the number of balls of the $ i\Th$ color in the urn. Then $\sum_i k_i = k$ and $S = \sum_i \indc{k_i > 0}$. 
Let $X_1,X_2,\ldots$ be independently drawn with replacement from the urn. Equivalently, the $X_i$'s are \iid according to a distribution $P = (p_i)_{i\geq 1}$, where $p_i=k_i/k$ is the fraction of balls of the $i\Th$ color. 
We assume the Poisson sampling model in \prettyref{sec:poi} where the observed data are $X_1,\ldots,X_N$, and the sample size $N\sim \Poi(n)$ is independent from $(X_i)_{i\geq 1}$.


We again consider the linear estimator of the form \eqref{eq:uj} (or equivalently \prettyref{eq:hatS}); however, the choice of the coefficients requires new consideration.
Following the analysis in \eqref{eq:bias-term}, we can bound the bias in the \DistinctElements problem as follows: 
\begin{equation}
    \label{eq:lp-full}
    \Expect[\hat S-S]
    \le \sum_i |\Expect[g(N_i)-1]|
    \le S e^{-n/k}\max_{a\in[k]}\abs{\phi(a)-1},
\end{equation}
where $\phi(a) \triangleq \sum_{j\ge 1}a^j\frac{u_j(n/k)^j}{j!}$ is a (formal) power series with $\phi(0)=0$. The right-hand side of \prettyref{eq:lp-full} can be made zero by choosing $\phi$ to be, \eg, the Lagrange interpolating polynomial that satisfies $ \phi(0)=-1 $ and $ \phi(i)=0 $ for $ i\in[k] $, namely, $\phi(a) = \frac{(-1)^{k+1}}{k!}\prod_{i=1}^k (a-i)$; however, this strategy results in a high-degree polynomial $\phi$ with large coefficients, which, in turn, leads to a large variance of the estimator.

As previously explained in \prettyref{sec:optimal}, due to the correlation decay between $\Phi_0$ and $\Phi_j$, we choose the threshold $L$ to be proportional to $\log k$.
Moreover, if a color is observed at most $L$ times, say, $ N_i\le L $, this implies that, with high probability, $ k_i\le M $, where $ M = O(kL/n)$, thanks to the concentration of Poisson random variables. Therefore, effectively we only need to consider those colors that appear in the urn for at most $M$ times, i.e., $ k_i\in [M] $, for which the bias is at most
\begin{align}
    |\Expect[g(N_i)-1]|
    &\le e^{-n/k}\max_{a\in[M]} \abs{\phi(a)-1} = e^{-n/k}\max_{x\in[M]/M} \abs{p(x)-1} \nonumber\\
    &= e^{-n/k} \norm{Bw-\ones}_{\infty}, \label{eq:lp}
\end{align}
where $p(x) \triangleq \phi(Mx) = \sum_{j=1}^L w_j x^j$, $ w=(w_1,\dots,w_L)^\top $, and
\begin{equation}
    w_j 
    \triangleq \frac{u_j(Mn/k)^j}{j!},\quad
    B \triangleq 
    \begin{pmatrix}
        1/M & (1/M)^2 & \cdots & (1/M)^L  \\
        2/M & (2/M)^2 & \cdots & (2/M)^L  \\
        \vdots  & \vdots  & \ddots & \vdots  \\
        1 & 1 & \cdots & 1
    \end{pmatrix}
    \label{eq:Bw}
\end{equation}
is a (partial) Vandermonde matrix. Lastly, since $ \Splug\le S\le k $, we define the final estimator to be $\hat S$ projected to the interval $[\Splug, k]$. We have the following error bound:
\begin{proposition}
    \label{prop:rate-w}
    Assume the Poisson sampling model. Let
    \begin{equation}
        L=\alpha \log k, \quad M=\frac{\beta  k\log k}{n},
        \label{eq:LM}
    \end{equation}
    for any $ \beta > \alpha $ such that $ L $ and $ M $ are integers. Let $ w\in\reals^{L} $. Let $\hat S$ be defined in \prettyref{eq:uj} with 
		\begin{equation}
		u_j = 
		\begin{cases}
w_j j! (\frac{k}{nM})^j  			&  j\in [L]\\
0		 & \text{otherwise}. \\
		\end{cases}
		\label{eq:ujgj}
		\end{equation}		
		Define $ \tilde S\triangleq(\hat S \vee \Splug)\wedge k $. Then
    \begin{gather}
        \Expect{(\tilde S-S)^2}
        \le k^2 e^{-2n/k}\norm{Bw-\ones}_{\infty}^2+ke^{-n/k}+k\max_{m\in[M]}\Expect_{N\sim\Poi(nm/k)}[u_{N}^2]\nonumber\\+k^{-(\beta -\alpha \log\frac{e\beta }{\alpha }-3)}.\label{eq:rate-w}
    \end{gather}
\end{proposition}
\begin{proof}
    Since $ \Splug\le S\le k $, $ \tilde{S} $ is always an improvement of $ \hat{S} $. Define the event $ E\triangleq \cap_{i=1}^k\{N_i\le L\Rightarrow kp_i\le M\} $, which means that whenever $N_i\le L$ we have $p_i\le M/k$. Since $ \beta>\alpha $, applying the Chernoff bound and the union bound yields $ \Prob[E^c]\le k^{1-\beta +\alpha \log\frac{e\beta }{\alpha }} $, and thus
    \begin{equation}
        \Expect{(\tilde S-S)^2}\le \Expect((\tilde S-S)\Indc_E)^2+k^2\Prob[E^c]
        \le \Expect((\hat S-S)\Indc_E)^2+k^{3-\beta +\alpha \log\frac{e\beta }{\alpha }}.
        \label{eq:hatC-risk1}
    \end{equation}
    By the definition of $ \hat{S} $ in \prettyref{eq:hatS},
    \begin{equation*}
        (\hat S-S)\Indc_E
        =\sum_{i:k_i\in [M]}(g(N_i)-1)\indc{N_i\le L}\triangleq\calE.
    \end{equation*}
    In view of the bias analysis in \prettyref{eq:lp}, we have
    \begin{equation}
        |\Expect[{\cal E}]|
        \le \sum_{i:k_i\in [M]}e^{-nk_i/k}\norm{Bw-\ones}_{\infty}
        \le ke^{-n/k}\norm{Bw-\ones}_{\infty}.
        \label{eq:calE-bias}
    \end{equation}
    Recall that $g(0)=0$ and $g(j)=u_j+1$ for $j\in[L]$. Since $ N_i$ is independently distributed as $ \Poi(nk_i/k)$, we have
    \begin{align}
      \var[{\cal E}]
      & = \sum_{i:k_i\in [M]}\var\qth{(g(N_i)-1)\indc{N_i\le L}}
        \nonumber\\&\le \sum_{i:k_i\in [M]}\Expect\qth{(g(N_i)-1)^2\indc{N_i\le L}}\nonumber\\
      & = \sum_{i:k_i\in [M]}\pth{e^{-nk_i/k}+\Expect [u_{N_i}^2]}
        \nonumber\\&\le ke^{-n/k}+k\max_{m\in[M]}\Expect_{N\sim\Poi(nm/k)}[u_N^2]. \label{eq:calE-var}
    \end{align}

    Combining the upper bound on the bias in \prettyref{eq:calE-bias} and the variance in \prettyref{eq:calE-var} yields an upper bound on $ \Expect[\calE^2] $. Then the MSE in \prettyref{eq:rate-w} follows from \prettyref{eq:hatC-risk1}.
\end{proof}

\prettyref{prop:rate-w} suggests that the coefficients of the linear estimator can be chosen by solving the following linear programming (LP):
\begin{equation}
\min_{w \in \reals^L} \norm{Bw-\ones}_{\infty}
\label{eq:LP}
\end{equation}
and, furthermore, proving that the coefficients of the minimizer is bounded in magnitude. Instead of the $L_\infty$-approximation problem \prettyref{eq:LP}, whose optimal value is difficult to analyze, we solve the $L_2$-approximation problem as a relaxation:
\begin{equation}
\min_{w \in \reals^L} \Norm{Bw-\ones}_2,
\label{eq:L2}
\end{equation}
which is an upper bound of \prettyref{eq:LP}, and is in fact within an $O(\log k)$ factor since $M=O(k\log k/n)$ and $n=\Omega(k/\log k)$. 
We consider two separate cases:
\begin{itemize}
    \item $M>L$ ($n \lesssim k$): In this case, the linear system in \prettyref{eq:L2} is overdetermined and the minimum is non-zero. Surprisingly, the exact optimal value can be found in closed form \cite[Lemma 1]{WY2016sample}:
    \[
    \min_{w \in \reals^L} \Norm{Bw-\ones}_2 
    =\qth{{\frac{\binom{M+L+1}{L+1}}{\binom{M}{L+1}}-1}}^{-1/2}, 
    \] 
		using discrete Chebyshev polynomials previously introduced in \eqref{eq:tn-full}  and their orthogonality properties.
		The coefficients of the minimizer can be bounded using the minimum singular value of the matrix $ B $.

    \item $M\leq L$ ($n \gtrsim k$): In this case, the linear system is underdetermined and the minimum in \prettyref{eq:L2} is zero. To bound the variance, it turns out that the coefficients bound obtained from the minimum singular value is too crude for this regime. Instead, we can express the coefficients in terms of Lagrange interpolating polynomials (see \prettyref{sec:interpol}) and use Stirling numbers to obtain sharp variance bounds. 
\end{itemize}


For details of determining the values of \prettyref{eq:L2} and deriving the MSE using \prettyref{eq:rate-w} in the above two cases, as well as statistical lower bound justifying the optimality of the sample complexity in \prettyref{tab:main}, see \cite{WY2016sample}. 
Since the size of the least square problem in \prettyref{eq:lp} is only logarithmic in $k$, similar to \prettyref{rmk:complexity}, the overall time complexity is $O(n)$.

\begin{figure}[ht]
        \centering
        \begin{subfigure}[ht]{.31\linewidth}
            \includegraphics[width=\linewidth]{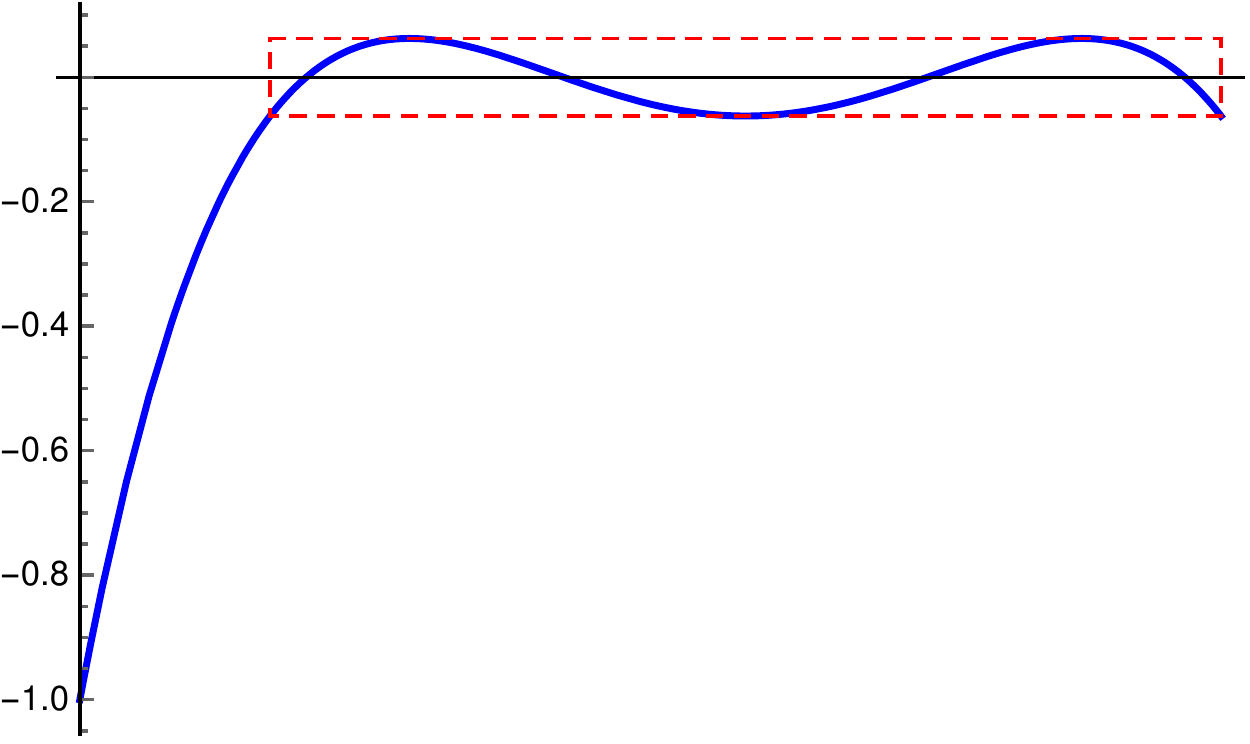}
            \caption{Continuous approximation}
            \label{fig:cont-approx} 
        \end{subfigure}
        \begin{subfigure}[ht]{.31\linewidth}
            \includegraphics[width=\linewidth]{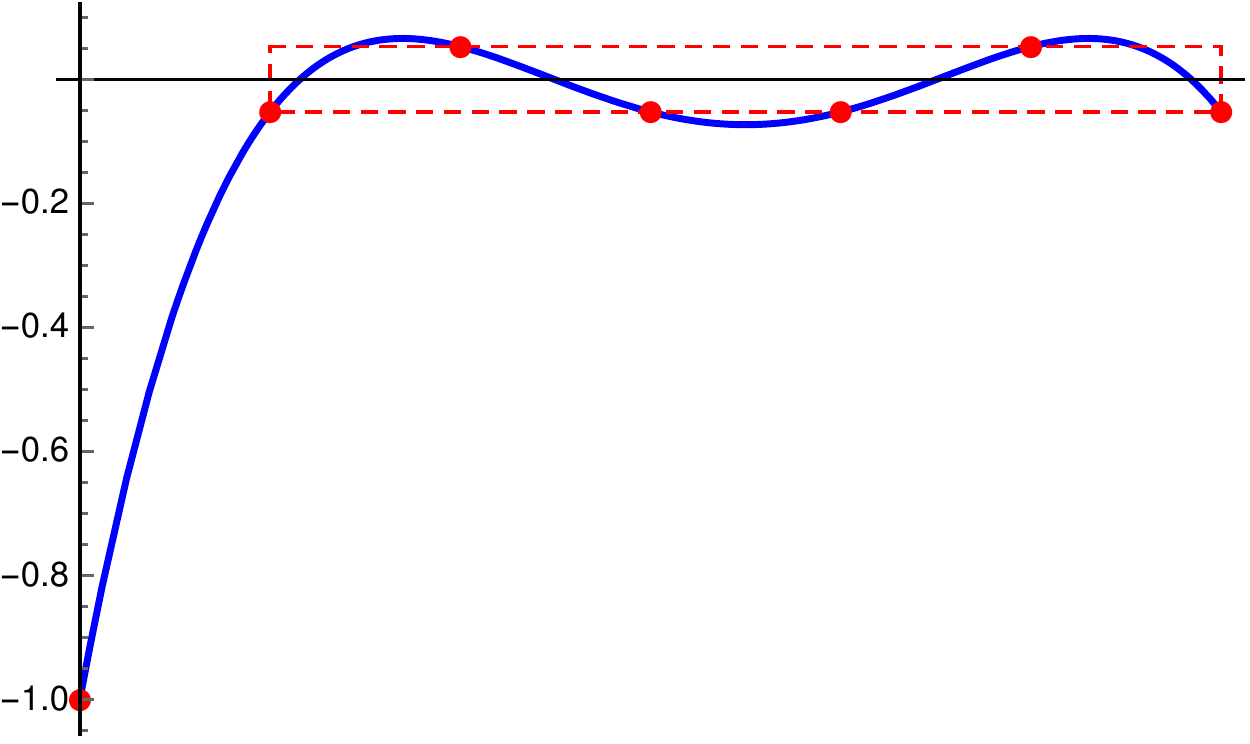}
            \caption{Discrete approximation}
            \label{fig:discrete-approx} 
        \end{subfigure}
        \begin{subfigure}[ht]{.31\linewidth}
            \includegraphics[width=\linewidth]{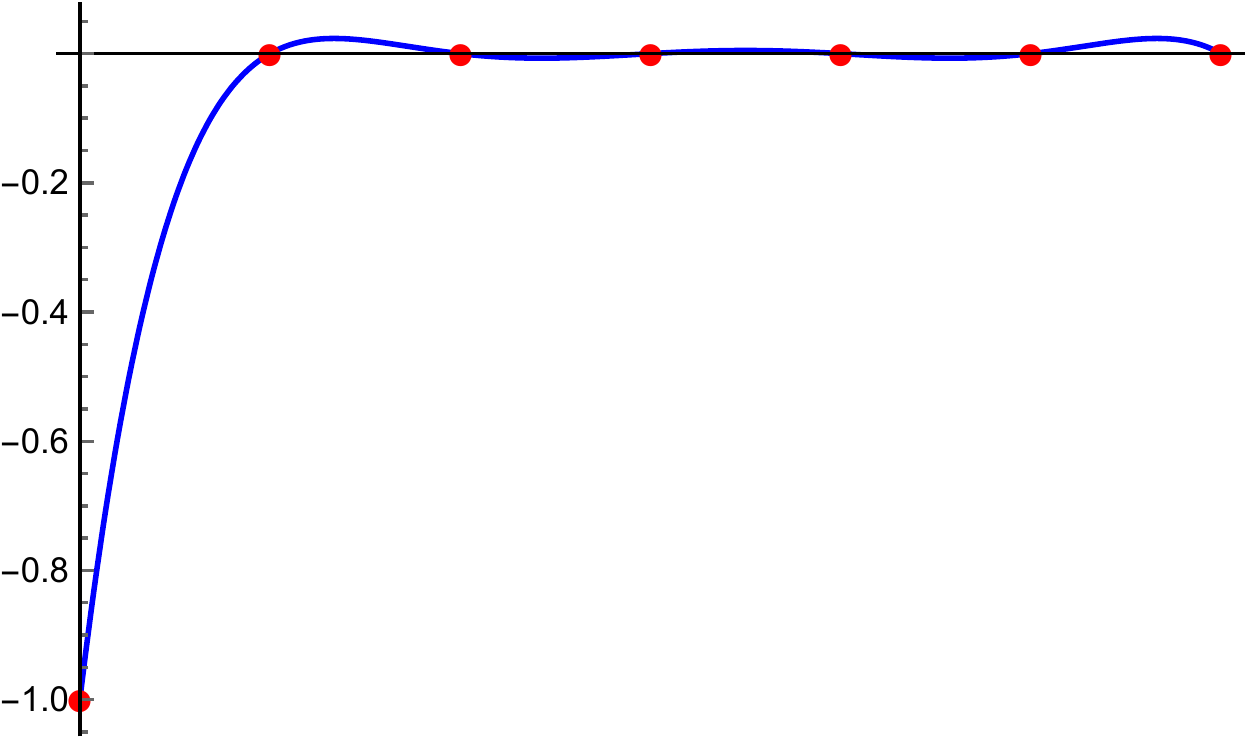}
            \caption{Interpolation}
            \label{fig:interpolation} 
        \end{subfigure}
        \caption{Continuous and discrete polynomial approximations for $M= 6 $ and degree $L=4$, where (\subref{fig:cont-approx}) and (\subref{fig:discrete-approx}) plot the optimal solution to \prettyref{eq:cheby-cont} and \prettyref{eq:cheby-discrete} respectively. The interpolating polynomial in (\subref{fig:interpolation}) requires a higher degree $L=6$. \label{fig:approx}}
    \end{figure}

We conclude this section by comparing the solutions to the \SupportSize and the \DistinctElements problems, both of which are based on the polynomial approximation method in \prettyref{chap:approx}; however, due to the definition of the parameter space, the approximation problem is of a continuous and discrete character, respectively.
    The optimal estimator for the both problems are of the linear form \prettyref{eq:uj}, but with different choice of the coefficients.
    In the \SupportSize problem, since the probabilities can take any values in an interval, the coefficients are found to be the solution of the continuous polynomial approximation problem
    \begin{equation}
        \inf_{p} \max_{x \in [\frac{1}{M},1]}  |p(x) - 1| = \exp\Big(-\Theta\Big(\frac{L}{\sqrt{M}}\Big)\Big).
        \label{eq:cheby-cont}
    \end{equation}
    where the infimum is taken over all degree-$L$ polynomials such that $p(0)=0$, achieved by the (appropriately shifted and scaled) Chebyshev polynomial \cite{timan63}. 
    In contrast, in the \DistinctElements problem, the discrete version of \prettyref{eq:cheby-cont}, which is equivalent to the LP \prettyref{eq:LP}, satisfies
    \begin{equation}
        \inf_{p} \max_{x \in \{\frac{1}{M},\frac{2}{M},\ldots,1\}}  |p(x) - 1| = \poly(M) \exp\Big(-\Theta\Big(\frac{L^2}{M}\Big)\Big),
        \label{eq:cheby-discrete}
    \end{equation}
    provided $L < M$ (see \cite{WY2016sample}). The difference between \prettyref{eq:cheby-cont} and \prettyref{eq:cheby-discrete} explains why the sample complexity \prettyref{eq:sample-supp} for the \SupportSize problem has an extra log factor compared to that of the \DistinctElements problem in \prettyref{tab:main}. When the sample size $n$ is large enough, interpolation is used in lieu of approximation. See \prettyref{fig:approx} for an illustration.

\section{Bibliographic notes}

For the \SupportSize problem, in \prettyref{sec:supp}, to estimate $S(P)$ within $\pm \epsilon k$, the naive plug-in estimator requires $\Theta(k \log \frac{1}{\epsilon})$ independent observations, which scales logarithmically in $\frac{1}{\epsilon}$ but linearly in $k$, the same scaling for estimating the distribution $P$ itself. 
Valiant and Valiant \cite{VV11} showed that the sample complexity is in fact sub-linear in $k$; however, the performance guarantee of the proposed estimators are still far from being optimal.
Specifically, an estimator based on an LP that is a modification of \cite[Program 2]{ET76} is proposed and shown to achieve $n_{\sf S}^*(k,\epsilon k) \lesssim \frac{k}{\epsilon^{2+\delta}\log k} $ for any arbitrary $ \delta>0 $ \cite[Corollary 11]{VV11}, which has subsequently been improved to $ \frac{k}{\epsilon^2\log k} $ in \cite[Theorem 2, Fact 9]{VV13}. 
The lower bound $ n_{\sf S}^*(k,\epsilon k) \gtrsim \frac{k}{\log k} $ in \cite[Corollary 9]{VV10} is optimal in $k$ but provides no dependence on $ \epsilon $.
These results show that the optimal scaling in terms of $k$ is $\frac{k}{\log k}$ but the achieved dependence on the accuracy $\epsilon$ is $\frac{1}{\epsilon^2}$, which is even worse than the plug-in estimator. 
From \prettyref{thm:sample} we see that the dependence on $\epsilon$ can be improved from polynomial to polylogarithmic $\log^2\frac{1}{\epsilon}$, and this is sharp.
Furthermore, this can be attained by a linear estimator which is far more scalable than linear programming on massive datasets.
A general framework of designing and analyzing linear estimators is given in \cite{VV11-focs} based on linear programming (as opposed to the approximation-theoretic approach in this survey).

The interests in the \DistinctElements problem also arise in the database literature, where various heuristic estimators \cite{HOT88,NS90} have been proposed under simplifying assumptions such as uniformity, and few performance guarantees are available.
More recent work in \cite{CCMN00,BKS01} obtained the optimal sample complexity under the \emph{multiplicative} error criterion\footnote{Here the multiplicative error means $\max\{\frac{\hat S}{S},\frac{S}{\hat S}\}$. This is not to be confused with the relative error $|\frac{\hat S}{S}-1|$. }, where the minimum sample size to estimate the number of distinct elements within a factor of $\alpha $ is shown to be $ \Theta(k/\alpha^2) $.
For this task, it turns out the least favorable scenario is to distinguish an urn with unitary color from one with \emph{almost} unitary color, the impossibility of which implies large multiplicative error.
However, the optimal estimator performs poorly compared with others on an urn with many distinct colors \cite{CCMN00}, the case where most estimators enjoy small multiplicative error.
In view of the limitation of multiplicative error, additive error is later considered by \cite{RRSS09,Valiant11}.
To achieve an additive error of $ ck $ for a constant $ c\in(0,\frac{1}{2}) $, the result in \cite{CCMN00} only implies an $ \Omega(1/c) $ sample complexity lower bound. A much stronger lower bound $ k^{1 - O(\sqrt{\frac{\log\log k}{\log k}})} $ was obtained in \cite{RRSS09}, which scales almost linearly but still suboptimal compared to \prettyref{tab:main}.
Recently, the optimal rate has been sharpened in \cite[Theorem 9]{PW18} in the linear sampling regime,\footnote{To be more precise, \cite{PW18} considered the Bernoulli sampling model (each ball is sampled with probability $p$) as opposed to the multinomial sampling model. Nevertheless, the result continuous to hold in view of the approximate equivalence between the two sampling models (see \cite[Corollary 1]{WY2016sample}).} where the sample size is $n=p\cdot k$ for some constant $p \in (0,1)$, 
by showing that the minimax rate is within logarithmic factors of $k^{-\min\{\frac{p}{1-p},1\}}$, improving over \prettyref{eq:distinct-rates} which only shows $k^{-\Theta(1)}$.

Finally we mention the closely related problem of \emph{species extrapolation}, which is a classical question dating back to Fisher \cite{FCW43} and arises ubiquitously in many scientific disciplines such as ecology \cite{FCW43,Chao92}, computational linguistics \cite{ET76}, genomics \cite{ICL09}. In the \iid sampling framework of \prettyref{sec:supp}, we observe $X_1,\ldots,X_n$ drawn independently from an unknown discrete distribution $P$, and the goal is to estimate the number of hitherto unobserved symbols that would be observed if $m$ new observations $X_1', \ldots, X_{m}'$ were collected from $P$, i.e.,
 \[
U = U_{n,m}  \triangleq |\{X_1', \ldots, X_m'\} \backslash \{X_1,\ldots,X_n\}|.
\]
In particular, the sequence $m\mapsto U_{n,m}$ is called the species discovery curve, which provides a guideline on how many new species would be observed were an additional sample of size $m$ to be collected. The extreme cases of $m=1$ and $m=\infty$ corresponds to the Good-Turing problem (estimating the probability that the next sample is new)~\cite{Good1953,OS15} and the \SupportSize problem, respectively.
Clearly, the bigger the sample to be extrapolated, the more difficult it is to obtain a reliable prediction of $U$. 
Of central importance is the largest \emph{extrapolation ratio} $t = \frac{m}{n}$ between
the number of new and existing observations for which $U$ can be accurately
predicted.  It was shown in \cite{OSW15} the minimax risk of estimating $U$ (mean squared error normalized by $n^2$) is given by $n^{-\Theta(\frac{1}{t})}$. 
Thus, it is possible to extrapolate with a vanishing risk provided that $t = o(\log n)$, and this condition is the best possible (see also \cite{zou2016quantifying}).
In particular, the upper bound is achieved by linear estimators with coefficients obtained by a random truncation of the unbiased estimator of Good-Toulmin \cite{GT56}, while the lower bound is shown by a reduction to the \SupportSize problem and applying \prettyref{thm:main}. 
For constant $t$, the minimax rate is recently shown to be within polylogarithmic factors of $n^{-\min\{\frac{2}{t+1},1\}}$ \cite{PW18}.



\chapter{Mixture models and moment comparison theorems}
\label{chap:framework}
In both this chapter and the next, we consider learning mixture models, another important task in statistical inference. 
Compared to the property estimation problems in Chapters~\ref{chap:approx}--\ref{chap:unseen} that can be successfully solved by polynomial approximation methods, here the key to efficiently learning mixture models is the geometry of the moment space that builds on the theory of positive polynomials. 
In this chapter, we first propose a framework of learning mixture models in \prettyref{sec:why-wass} and introduce the Wasserstein distance as the loss function in \prettyref{sec:wass}.
Then we present a number of moment comparison theorems in Sections~\ref{sec:compare-wass} and \ref{sec:compare-higher-moments} which will be used for deriving statistical guarantees of moment-based procedures in \prettyref{chap:gm} as well as certifying their optimality.
The main techniques for establishing the moment comparison theorems are polynomial interpolation introduced in Sections~\ref{sec:interpol} and \ref{sec:hermite}.

Learning mixture models has a long history in statistics \cite{Pearson1894} and enjoyed recent renewed interest in latent variable models in machine learning. In a $k$-component mixture model from a parametric family of distributions $\calP=\{P_\theta:\theta\in\Theta\}$, each observation is distributed as
\begin{equation*}
    X\sim \sum_{i=1}^k w_i P_{\theta_i}.
\end{equation*}
Here $w_i$ is the mixing weight such that $w_i\ge 0$ and $\sum_i w_i = 1$, $\theta_i\in \Theta$ is the parameter of the $i\Th$ component. Equivalently, we can write the distribution of an observation $ X $ as 
\begin{equation}
    \label{eq:model-mixing}
    X\sim \pi_\nu = \int P_{\theta} \diff\nu(\theta),
\end{equation}
where $\nu = \sum_{i=1}^k w_i \delta_{\theta_i}$ denotes the \emph{mixing distribution} and $U\sim\nu$ is referred to as the latent variable.

Generally speaking, there are three common formulations of learning mixture models:
\begin{itemize}
    \item \textbf{Parameter estimation}:
    estimate the parameter $\theta_i$'s and the weights $w_i$'s up to a global permutation.
    \item \textbf{Density estimation}:
    estimate the probability density function of the mixture model under certain loss such as $L_2$ or Hellinger distance. This task is further divided into the cases of \emph{proper} and \emph{improper} learning, depending on whether the estimate is required to be a mixture of distributions in $\calP$ or not; in the latter case, there are more flexibility in designing the estimator but less interpretability.
    \item \textbf{Clustering}: estimate the latent variable of each observation (i.e.~$U_i$, such that $i$th observation is distributed as $P_{\theta}$ conditioned on $U_i=\theta$) with a small misclassification rate.
\end{itemize}
It is clear that to ensure the possibility of clustering it is necessary to impose certain separation conditions between the clusters; however, as far as estimation is concerned, both parametric and non-parametric, no separation condition should be needed and one can obtain accurate estimates of the parameters even when clustering is impossible. Furthermore, one should be able to learn from the data the order of the mixture model, that is, the number of components. However, in the present literature, most of the estimation procedures with finite sample guarantees are either clustering-based, or rely on separation conditions in the analysis (e.g.~\cite{balakrishnan2017statistical,lu2016statistical,hopkins2018mixture}). Bridging this conceptual divide is one of the main motivations of the present chapter.

\section{Estimating the mixing distribution}
\label{sec:why-wass}
Following the framework proposed in \cite{Chen95,HK2015}, in this chapter we consider the estimation of the mixing distribution $\nu$, rather than estimating the parameters $(w_i,\theta_i)$ of each component. The main benefit of this formulation include the following:
\begin{itemize}
    \item Assumption-free: to recover individual components it is necessary to impose certain assumptions to ensure identifiability, such as lower bounds on the mixing weights and separations between components, none of which is needed for estimating the mixing distribution.
    Furthermore, under the usual assumption such as separation conditions, statistical guarantees on estimating the mixing distribution can be naturally translated to those for estimating the individual parameters.
    \item Inference on the number of components: this formulation allows us to deal with misspecified models and estimate the order of the mixture model.
\end{itemize}

In this framework, a meaningful and flexible loss function for estimating the mixing distribution is the \emph{$1$-Wasserstein distance} defined by
\begin{equation}
    \label{eq:w1}
    W_1(\nu,\nu')\triangleq \inf \{\Expect[\Norm{X-Y}]: X\sim \nu, Y\sim\nu'\},
\end{equation}
where the infimum is taken over all couplings, i.e., joint distributions of $X$ and $Y$ which are marginally distributed as $\nu$ and $\nu'$ respectively.
In one dimension, the $W_1$ distance coincides with the $L_1$-distance between the CDFs \cite{villani.topics}. This is a natural criterion, which is not too stringent to yield trivial result (e.g.~the Kolmogorov-Smirnov (KS) distance\footnote{
    Consider two mixing distributions $\delta_0$ and $\delta_{\epsilon}$ with arbitrarily small $\epsilon$, whose KS distance is always one.
    })
and, at the same time, strong enough to provide meaningful guarantees on the means and weights. In fact, the commonly used criterion $\min_\Pi\sum_i \Norm{\theta_i-\hat\theta_{\Pi(i)}}$ over all permutations $\Pi$ is precisely ($k$ times) the Wasserstein distance between two equally weighted distributions \cite{villani.topics}.

Furthermore, we can obtain statistical guarantees on the support sets and weights of the estimated mixing distribution under the usual assumptions in literature \cite{Dasgupta1999,KMV2010,HP15} that include separation between the means and lower bound on the weights. See \prettyref{sec:wass} for a detailed discussion. We highlight the following result, phrased in terms of the parameter estimation error up to a permutation:
\begin{lemma}
    \label{lmm:w1-parameters}
    Let
    \[
        \nu=\sum_{i=1}^k w_i \delta_{\theta_i}, \quad \hat \nu=\sum_{i=1}^k \hat w_i \delta_{\hat \theta_i}.
    \]
 Let
    \begin{align*}
        \epsilon&=W_1(\nu,\hat\nu),\\
        \epsilon_1&=\min\{\Norm{\theta_i-\theta_j},\Norm{\hat\theta_i-\hat\theta_j}:1\le i< j\le k\},\\
        \epsilon_2&=\min\{w_i,\hat w_i:i\in [k]\}.
    \end{align*}
    If $\epsilon<\epsilon_1\epsilon_2/4$, then, there exists a permutation $\Pi$ on $[k]$ such that
    \[
        \Norm{\theta_i-\hat \theta_{\Pi(i)}} \le \epsilon/\epsilon_2,\quad |w_i-\hat w_{\Pi(i)}| \le 2\epsilon/\epsilon_1,\quad \forall~i.
    \]
\end{lemma}

\section{Wasserstein distance}
\label{sec:wass}
A central quantity in the theory of optimal transportation, the Wasserstein distance is the minimum cost of mapping one distribution to another.
In this part, we will be mainly concerned with the 1-Wasserstein distance defined in \prettyref{eq:w1}, which can be equivalently expressed, through the Kantorovich duality \cite{villani.topics}, as
\begin{equation}
    \label{eq:W1-dual}
    W_1(\nu,\nu')=\sup \{\Expect_\nu[\varphi]-\Expect_{\nu'}[\varphi]:\varphi\text{ is 1-Lipschitz}\}.
\end{equation}
The optimal coupling in \prettyref{eq:w1} has many equivalent characterization \cite{optimal.transport.old.new} but is often difficult to compute analytically in general. Nevertheless, the situation is especially simple for distributions on the real line, where the quantile coupling is known to be optimal and hence
\begin{equation}
    \label{eq:W1-CDF}
    W_1(\nu,\nu')=\int |F_\nu(t)-F_{\nu'}(t)|\diff t.
\end{equation}
Both \prettyref{eq:W1-dual} and \prettyref{eq:W1-CDF} provide convenient characterizations to bound the Wasserstein distance in \prettyref{sec:compare-wass}.

As previously mentioned in \prettyref{sec:why-wass}, two discrete distributions close in the Wasserstein distance have similar support sets and weights. This is made precise by the next two lemmas.

\begin{lemma}
    \label{lmm:w1-hausdorff}
    Suppose $\nu$ and $\nu'$ are discrete distributions supported on $S$ and $S'$, respectively. 
		Define the Hausdorff distance between $S$ and $S'$ as
		\[
    d_H(S,S') \triangleq \max\sth{\sup_{x\in S}\inf_{x'\in S'}\Norm{x-x'},\sup_{x'\in S'}\inf_{x\in S}\Norm{x-x'}}.
		\]
		Then
    \[
        d_H(S,S')\le W_1(\nu,\nu')/\epsilon,
    \]
		where $\epsilon=\min\{\nu(x):x\in S\}\wedge\min\{\nu'(x):x\in S'\}$ denotes the minimum probability mass of $\nu$ and $\nu'$. 
\end{lemma}
\begin{proof}
    For any coupling $P_{XY}$ such that $X\sim\nu$ be $Y\sim\nu'$,
    \begin{align*}
        \Expect \Norm{X-Y} &= \sum_x \Prob[X=x] \Expect[\Norm{X-Y} | X=x]\\
        &\ge \sum_x \epsilon\cdot \inf_{x'\in S'} \Norm{x-x'} \ge \epsilon \cdot\sup_{x\in S} \inf_{x'\in S'} \Norm{x-x'}.
    \end{align*}
    Interchanging $X$ and $Y$ completes the proof.
\end{proof}

\begin{lemma}
    \label{lmm:w1-prokhorov}
    For any $\delta>0$,
    \begin{gather*}
        \nu(x)- \nu'(B(x,\delta))\le W_1(\nu,\nu')/\delta,\\
        \nu'(x)- \nu(B(x,\delta))\le W_1(\nu,\nu')/\delta. 
    \end{gather*}
\end{lemma}
\begin{proof}
    Using the optimal coupling $P_{XY}^*$ such that $X\sim\nu$ be $Y\sim\nu'$, applying Markov inequality yields that
    \[
        \Prob[|X-Y|>\delta]\le \Expect|X-Y|/\delta=W_1(\nu,\nu')/\delta.
    \]
    By Strassen's theorem (see \cite[Corollary 1.28]{villani.topics}), for any Borel set $B$, we have $\nu(B)\le \nu'(B^\delta)+W_1(\nu,\nu')/\delta$ and $\nu'(B)\le \nu(B^\delta)+W_1(\nu,\nu')/\delta$, where $B^\delta\triangleq \{x:\inf_{y\in B}\Norm{x-y}\le \delta\}$ denotes the $\delta$-fattening of $B$. The conclusion follows by considering a singleton $B=\{x\}$.
\end{proof}

\prettyref{lmm:w1-hausdorff} and \ref{lmm:w1-prokhorov} together yield a bound on the parameter estimation error (up to a permutation) in terms of the Wasserstein distance, which was previously given in \prettyref{lmm:w1-parameters}: 
\begin{proof}[Proof of \prettyref{lmm:w1-parameters}]
    Denote the support sets of $\nu$ and $\hat\nu$ by $S=\{\theta_1,\ldots,\theta_k\}$ and $S'=\{\hat\theta_1,\ldots,\hat\theta_k\}$, respectively. Applying \prettyref{lmm:w1-hausdorff} yields that $d_H(S,S')\le \epsilon/\epsilon_2$, which is less than $\epsilon_1/4$ by the assumption $\epsilon<\epsilon_1\epsilon_2/4$. Since $\Norm{\theta_i-\theta_j}\ge \epsilon_1$ for every $i\ne j$, then there exists a permutation $\Pi$ such that 
    \[
        \Norm{\theta_i-\hat \theta_{\Pi(i)}} \le \epsilon/\epsilon_2,\quad \forall~i.
    \]
    Applying \prettyref{lmm:w1-prokhorov} twice with $\delta=\epsilon_1/2$, $x=\theta_i$ and $x=\hat \theta_{\Pi(i)}$, respectively, we obtain that
    \[
        w_i-\hat w_{\Pi(i)}\le 2\epsilon/\epsilon_1,\quad \hat w_{\Pi(i)}-w_i\le 2\epsilon/\epsilon_1.\tag*{\qedhere}
    \]
\end{proof}

\section{Moment comparison in Wasserstein distance for discrete distributions }
\label{sec:compare-wass}
Under the framework of learning finite-order mixture models in \prettyref{sec:why-wass}, in this section we prove comparison theorems in Wasserstein distance for discrete distributions in terms of moments. They form the foundation of proving statistical guarantees for moment-based learning algorithms of mixture models in \prettyref{chap:gm}.

Moment comparison is a classical topic in the probability theory. Classical moments comparison theorems aim to show convergence of distributions by comparing a \emph{growing} number of moments. For example, Chebyshev's theorem states that if the first $r$ moments of a distribution agrees with those of the normal distribution, then their KS distance is at most $O(\frac{1}{\sqrt{r}})$.
More precisely, if ${\bf m}_r(\pi)={\bf m}_r(N(0,1))$, then (see \cite[Theorem 2]{Diaconis1987})
\[
    \sup_{x\in\reals}|F_{\pi}(x)-\Phi(x)|\le \sqrt{\frac{\pi}{2r}},
\]
where $F_{\pi}$ and $\Phi$ denote the CDFs of $\pi$ and $N(0,1)$, respectively. For two compactly supported distributions with identical first $r$ moments, the estimate for their KS distance can be sharpened to $O(\frac{\log r}{r})$ \cite{Krawtchouk1932}. In contrast, in the context of estimating finite mixtures we are dealing with finitely-supported mixing distributions, which can be identified by a \emph{fixed} number of moments. However, for finite sample size, it is impossible to exactly determine the moments, and measuring the error in the KS distance is too much to ask (see the counterexample in \prettyref{sec:why-wass}). It turns out that $W_1$-distance is a suitable metric for this purpose, and the closeness of moments does imply the closeness of distribution in the $W_1$ distance, which is the integrated difference ($L_1$-distance) between two CDFs as opposed the uniform error ($L_\infty$-distance).

We first present the identifiability of discrete distributions from a finite set of moments in order to familiarize the readers with the techniques of polynomial interpolation in \prettyref{sec:interpol}.
Note that a discrete distribution with $k$ atoms has $2k-1$ free parameters. Therefore it is reasonable to expect that it can be uniquely determined by its first $2k-1$ moments. Indeed, we have the following simple identifiability results for discrete distributions: 
\begin{lemma}
    \label{lmm:identify} 
	Let $\nu$ and $\nu'$ be distributions on the real line.
    \begin{enumerate}
        \item If $\nu$ and $\nu'$ are both $k$-atomic, then $\nu=\nu'$ if and only if $\bfm_{2k-1}(\nu)=\bfm_{2k-1}(\nu')$.
        \item If $\nu$ is $k$-atomic, then $\nu=\nu'$ if and only if $\bfm_{2k}(\nu)=\bfm_{2k}(\nu')$.
    \end{enumerate}    
\end{lemma}
\begin{proof}
    We only need to prove the ``if'' part. We prove this lemma using the apparatus of interpolating polynomials previously introduced in \prettyref{sec:interpol}.
    \begin{enumerate}
        \item Denote the union of the support sets of $\nu$ and $\nu'$ by $S$. Here $S$ is of size at most $2k$. For any $t\in\reals$, there exists a polynomial $P$ of degree at most $2k-1$ that interpolates the step function $x\mapsto \indc{x\le t}$ on $S$. Since $m_i(\nu)=m_i(\nu')$ for $i=1,...,2k-1$, we have 
        \[
            F_{\nu}(t)=\Expect_\nu[\indc{X\le t}]=\Expect_{\nu}[P(X)]=\Expect_{\nu'}[P(X)]=\Expect_{\nu'}[\indc{X\le t}]=F_{\nu'}(t).
        \]
        \item Denote the support set of $\nu$ by $S'=\{x_1,\dots,x_k\}$. Let $Q(x)=\prod_i(x-x_i)^2$. Then $Q$ is a non-negative polynomial of degree $2k$. Since $m_i(\nu)=m_i(\nu')$ for $i=1,...,2k$, we have 
        \[
            \Expect_{\nu'}[|Q(X)|]=\Expect_{\nu'}[Q(X)]=\Expect_{\nu}[Q(X)]=0.
        \]
        Since $Q(x)=0$ if and only $x\in S'$, $\nu'$ is also supported on $S'$ and thus is $k$-atomic. The conclusion follows from the first case of \prettyref{lmm:identify}. \qedhere
    \end{enumerate}
\end{proof}

In the context of statistical estimation, we only have access to samples and noisy estimates of moments. To solve the inverse problems from moments to distributions, our theory relies on the following quantitative version of the identifiability result in \prettyref{lmm:identify}:
\begin{proposition}
    \label{prop:stable1}
    Let $\nu$ and $\nu'$ be $k$-atomic distributions supported on $[-1,1]$. If $|m_i(\nu)-m_i(\nu')|\le \delta$ for $i=1,\dots,2k-1$, then
    \[
        W_1(\nu,\nu')\le O\pth{k\delta^{\frac{1}{2k-1}}}.
    \]
\end{proposition}
\begin{proposition}
    \label{prop:stable2}        
    Let $\nu$ be a $k$-atomic distribution supported on $[-1,1]$. If $|m_i(\nu)-m_i(\nu')|\le \delta$ for $i=1,\dots,2k$, then
    \[
        W_1(\nu,\nu')\le O\pth{k\delta^{\frac{1}{2k}}}.
    \]
\end{proposition}
Analogous to the proof of \prettyref{lmm:identify}, the proofs of Propositions \ref{prop:stable1} and \ref{prop:stable2} also rely on techniques of polynomial interpolation and will be presented in \prettyref{sec:compare-proof}.
The main concern of those results is the optimal dependence on $\delta$, and the linear factor of $k$ might be suboptimal. 
In the sequel we refer to Propositions \ref{prop:stable1} and \ref{prop:stable2} as moment comparison theorems, which show that closeness of in the moment space implies closeness of distributions in Wasserstein distance.
An upper bound on the Wasserstein distance is obtained in \cite{KV17} involving the differences of the first $k$ moments and a $\Theta(\frac{1}{k})$ term that does not vanish for fixed $k$. 

\begin{remark}
    \label{rmk:stable}
    The exponents in \prettyref{prop:stable1} and \ref{prop:stable2} are optimal. To see this, we first note that the number of moments needed for identifiability in \prettyref{lmm:identify} cannot be reduced:
    \begin{enumerate}
        \item Given any $2k$ distinct points, there exist two $k$-atomic distributions with disjoint support sets but identical first $2k-2$ moments (see \prettyref{lmm:match2k-2} in \prettyref{sec:aux-lemmas}).
        \item Given any continuous distribution, its $k$-point Gaussian quadrature is $k$-atomic and has identical first $2k-1$ moments (see \prettyref{sec:GQ}).
    \end{enumerate}
    By the first observation, there exists two $k$-atomic distributions $\nu$ and $\nu'$ such that 
    \begin{gather*}
        m_i(\nu)=m_i(\nu'),\quad i=1,\ldots,2k-2,\\ |m_{2k-1}(\nu)-m_{2k-1}(\nu')|= c_k,\quad
        W_1(\nu,\nu')=d_k,
    \end{gather*}
    where $c_k$ and $d_k$ are strictly positive constants that depend on $k$. Let $\tilde\nu$ and $\tilde\nu'$ denote the distributions of $\epsilon X$ and $\epsilon X'$ such that $X\sim\nu$ and $X'\sim\nu'$, respectively. Then, we have
    \[
        \max_{i\in [2k-1]}|m_i(\tilde\nu)-m_i(\tilde\nu')|= \epsilon^{2k-1}c_k,
        \quad W_1(\tilde\nu,\tilde\nu')=\epsilon d_k.
    \]
    This concludes the tightness of the exponent in \prettyref{prop:stable1}. Similarly, the exponent in \prettyref{prop:stable2} is also tight using the second observation.
\end{remark}

When the atoms of the discrete distributions are separated, we have the following adaptive version of the moment comparison theorems (cf.~Propositions \ref{prop:stable1} and \ref{prop:stable2}), which naturally lead to the adaptive rates of our moment-based algorithms in \prettyref{chap:gm}. 
\begin{proposition}
    \label{prop:stable1-separation}
    Suppose both $\nu$ and $\nu'$ are supported on a set of $\ell$ atoms in $[-1,1]$, and each atom is at least $\gamma$ away from all but at most $\ell'$ other atoms. Let $\delta=\max_{i\in [\ell-1]}|m_i(\nu)-m_i(\nu')|$. Then,
    \[
        W_1(\nu,\nu')\le \ell \pth{\frac{\ell 4^{\ell-1}\delta}{\gamma^{\ell-\ell'-1}}}^{\frac{1}{\ell'}}.
    \]
\end{proposition}
\begin{proposition}
    \label{prop:stable2-separation}
    Suppose $\nu$ is supported on $k$ atoms in $[-1,1]$ and any $t\in\reals$ is at least $\gamma$ away from all but $k'$ atoms. Let $\delta=\max_{i\in [2k]}|m_i(\nu)-m_i(\nu')|$. Then,
    \[
        W_1(\nu,\nu')\le 8k\pth{\frac{k 4^{2k} \delta}{\gamma^{2(k-k')}}}^{\frac{1}{2k'}}.
    \]
\end{proposition}

\subsection{Proofs}
\label{sec:compare-proof}
First we prove \prettyref{prop:stable1-union}, which is slightly stronger than \prettyref{prop:stable1}. We provide three proofs: the first two are based on the primal (coupling) formulation of $W_1$ distance \prettyref{eq:W1-CDF}, and the third proof uses the dual formulation \prettyref{eq:W1-dual}. Specifically,
\begin{itemize}
    \item The first proof uses polynomials to interpolate step functions, whose expected values are the CDFs. The closeness of moments imply the closeness of distribution functions and thus, by \prettyref{eq:W1-CDF}, a small Wasserstein distance. Similar idea applies to the proof of \prettyref{prop:stable2} later.
    \item The second proof finds a polynomial that preserves the sign of the difference between two CDFs, and then relate the Wasserstein distance to the integral of that polynomial. Similar idea is used in \cite[Lemma 20]{MV2010} which uses a polynomial that preserves the sign of the difference between two density functions. 
    
		\item The third proof uses the dual formulation \prettyref{eq:W1-dual} of the Wasserstein distance and 
		polynomial interpolation of 1-Lipschitz functions.
\end{itemize}
In this subsection we present the main proofs, and the technical lemmas are collected in \prettyref{sec:aux-lemmas}.
\begin{proposition}
    \label{prop:stable1-union}    
    Let $\nu$ and $\nu'$ be discrete distributions supported on $\ell$ atoms in $[-1,1]$. If
    \begin{equation}
        \label{eq:momentdiff}
        |m_i(\nu)-m_i(\nu')|\le \delta,\quad i=1,\ldots,\ell-1,
    \end{equation}
    then
    \[
        W_1(\nu,\nu')\le O\pth{\ell \delta^{\frac{1}{\ell-1}}}.
    \]
\end{proposition}
\begin{proof}[First proof of \prettyref{prop:stable1-union}]
    Suppose $\nu$ and $\nu'$ are supported on
    \begin{equation}
        \label{eq:supports-1}
        S=\{t_1,\dots,t_\ell\},\quad t_1<t_2<\dots<t_{\ell}.    
    \end{equation}
    Then, using the integral representation \prettyref{eq:W1-CDF}, the $W_1$ distance reduces to 
    \begin{equation}
        \label{eq:W1-CDF-discrete}
        W_1(\nu,\nu')=\sum_{r=1}^{\ell-1}|F_{\nu}(t_r)-F_{\nu'}(t_r)|\cdot|t_{r+1}-t_r|.
    \end{equation}
    For each $r$, let $f_r(x)=\indc{x\le t_r}$, and $P_r$ be the unique polynomial of degree $\ell-1$ to interpolate $f_r$ on $S$. In this way we have $f_r=P_r$ almost surely under both $\nu$ and $\nu'$, and thus
    \begin{equation}
        \label{eq:CDF-poly}    
        |F_{\nu}(t_r)-F_{\nu'}(t_r)|=|\Expect_{\nu}f_r-\Expect_{\nu'}f_r|
        =|\Expect_{\nu}P_r-\Expect_{\nu'}P_r|.
    \end{equation}
    $P_r$ can expressed using Newton formula \prettyref{eq:interpolation-Newton} as
    \begin{equation}
        \label{eq:Pr-1}
        P_r(x)=1+\sum_{i=r+1}^{\ell}f_r[t_1,\dots,t_i]g_{i-1}(x),
    \end{equation}
    where $g_{r}(x)=\prod_{j=1}^{r}(x-t_j)$ and we used $f_r[t_1,\dots,t_i]=0$ for $i=1,\ldots,r$. In \prettyref{eq:Pr-1}, the absolute values of divided differences are obtained in \prettyref{lmm:divided} in \prettyref{sec:aux-lemmas}:
    \begin{equation}
        \label{eq:fr-div}
        |f_r[t_1,\dots,t_i]|\le \frac{\binom{i-2}{r-1}}{(t_{r+1}-t_{r})^{i-1}}.
    \end{equation}
    In the summation of \prettyref{eq:Pr-1}, let $g_{i-1}(x)=\sum_{j=0}^{i-1} a_jx^j$. Since $|t_j|\le 1$ for every $j$, we have $\sum_{j=0}^{i-1}|a_j|\le 2^{i-1}$ (see \prettyref{lmm:Px-expand} in \prettyref{sec:aux-lemmas}). Applying \prettyref{eq:momentdiff} yields that
    \begin{equation}
        \label{eq:polydiff}
        |\Expect_\nu[g_{i-1}]-\Expect_{\nu'}[g_{i-1}]|\le \sum_{j=1}^{i-1} |a_j|\delta\le 2^{i-1}\delta.
    \end{equation}
    Then we obtain from \prettyref{eq:CDF-poly} and \prettyref{eq:Pr-1} that
    \begin{equation}
        \label{eq:CDFdiff}
        |F_{\nu}(t_r)-F_{\nu'}(t_r)|\le \sum_{i=r+1}^{\ell}\frac{\binom{i-2}{r-1}2^{i-1}\delta}{(t_{r+1}-t_{r})^{i-1}}
        \le \frac{\ell 4^{\ell-1}\delta}{(t_{r+1}-t_r)^{\ell-1}}.
    \end{equation}
    Also, $|F_{\nu}(t_r)-F_{\nu'}(t_r)|\le 1$ trivially. Therefore,
    \begin{equation}
        \label{eq:W1-example}
        W_1(\nu,\nu')\le \sum_{r=1}^{\ell-1}\pth{\frac{\ell 4^{\ell-1}\delta}{(t_{r+1}-t_r)^{\ell-1}}\wedge 1}\cdot|t_{r+1}-t_r|\le 4e(\ell-1)\delta^{\frac{1}{\ell-1}},
    \end{equation}
    where we used $\max\{\frac{\alpha}{x^{\ell-2}}\wedge x:x>0\}=\alpha^{\frac{1}{\ell-1}}$ and $x^{\frac{1}{x-1}}\le e$ for $x\ge 1$.
\end{proof}

\begin{proof}[Second proof of \prettyref{prop:stable1-union}]
    Suppose on the contrary that
    \begin{equation}
        W_1(\nu,\nu')\ge C\ell\delta^{\frac{1}{\ell-1}},
    \end{equation}
    for some absolute constant $C$. We will show that $\max_{i\in[\ell-1]}|m_i(\nu)-m_i(\nu')|\ge \delta$. We continue to use $S$ in \prettyref{eq:supports-1} to denote the support of $\nu$ and $\nu'$. Let $\Delta F(t)=F_{\nu}(t)-F_{\nu'}(t)$ denote the difference between two CDFs. Using \prettyref{eq:W1-CDF-discrete}, there exists $r\in [\ell-1]$ such that
    \begin{equation}
        \label{eq:DeltaF-2}
        |\Delta F(t_r)|\cdot |t_{r+1}-t_r|\ge C\delta^{\frac{1}{\ell-1}}.
    \end{equation}
    We first construct a polynomial $L$ that preserves the sign of $\Delta F$. To this end, let $S'=\{s_1,\dots,s_m\}\subseteq S$ such that $t_1=s_1<s_2<\dots<s_m=t_{\ell}$ be the set of points where $\Delta F$ changes sign, \ie, $\Delta F(x)\Delta F(y)\le 0$ for every $x\in(s_i,s_{i+1})$, $y\in(s_{i+1},s_{i+2})$, for every $i$. Let $L(x)\in \pm \prod_{i=2}^{m-1}(x-s_i)$ be a polynomial of degree at most $\ell-2$ that also changes sign on $S'$ such that
    \[
        \Delta F(x)L(x)\ge 0,\quad t_1\le x\le t_\ell.
    \]
    Consider the integral of the above positive function. Applying integral by parts and using the boundary condition $\Delta F(t_{\ell})=\Delta F(t_{1})=0$, we obtain
    \begin{equation}
        \label{eq:int-EP}
        \int_{t_1}^{t_\ell} \Delta F(x)L(x)\diff x= -\int_{t_1}^{t_\ell}P(x)\diff \Delta F(x)
        =\Expect_{\nu'}[P(X)]-\Expect_\nu [P(X)],
    \end{equation}
    where $P(x)$ is a polynomial of degree at most $\ell-1$ such that $P'(x)=L(x)$. If we write $L(x)=\sum_{j=0}^{\ell-2}a_jx^j$, then $P(x)=\sum_{j=0}^{\ell-2}\frac{a_j}{j+1}x^{j+1}$. Since $|s_j|\le 1$ for every $j$, we have $\sum_{j=0}^{\ell-2}|a_j|\le 2^{\ell-2}$ (see \prettyref{lmm:Px-expand} in \prettyref{sec:aux-lemmas}), and thus $\sum_{j=0}^{\ell-2}\frac{|a_j|}{j+1}\le 2^{\ell-2}$. Hence, 
    \begin{equation}
        \label{eq:EP-ub}
        |\Expect_{\nu'}[P(X)]-\Expect_\nu [P(X)]|\le 2^{\ell-2}\max_{i\in[\ell-1]}|m_i(\nu)-m_i(\nu')|.
    \end{equation}

    Since $\Delta F(x)L(x)$ is always non-negative, applying \prettyref{eq:DeltaF-2} to \prettyref{eq:int-EP} yields that
    \begin{align}
        |\Expect_{\nu'}[P(X)]-\Expect_\nu [P(X)]|&\ge \int_{t_r}^{t_{r+1}} |\Delta F(x)L(x)|\diff x\nonumber\\
        &\ge \frac{C\delta^{\frac{1}{\ell-1}}}{|t_{r+1}-t_r|}\int_{t_r}^{t_{r+1}}|L(x)|\diff x.\label{eq:diff-EP}
    \end{align}
    Recall that $|L(x)|=\prod_{i=2}^{m-1}|x-s_i|$. Then for $x\in (t_r,t_{r+1})$, we have $|x-s_i|\ge x-t_r$ if $s_i\le t_r$, and $|x-s_i|\ge t_{r+1}-x$ if $s_i\ge t_{r+1}$. Hence,
    \[
        |L(x)|\ge (t_{r+1}-x)^a(x-t_r)^b,
    \]
    for some $a,b\in\naturals$ such that $a,b\ge 1$ and $a+b\le \ell-2$. Using the density of the Beta distribution, the integral of the right-hand side of the above inequality can be
		found as
    \[
        \int_{t_r}^{t_{r+1}}(t_{r+1}-x)^a(x-t_{r})^b\diff x=\frac{(t_{r+1}-t_r)^{a+b+1}}{(a+1)\binom{a+b+1}{b}}.
    \]
    Since $|t_{r+1}-t_r|\ge |\Delta F(t_r)|\cdot |t_{r+1}-t_r|\ge C\delta^{\frac{1}{\ell-1}}$ and $\binom{a+b+1}{b}\le 2^{a+b+1}$, and $a+b+1\le \ell-1$, we obtain from \prettyref{eq:diff-EP} that
    \begin{equation}
        \label{eq:EP-lb}
        |\Expect_{\nu'}[P(X)]-\Expect_\nu [P(X)]|\ge \delta\frac{(C/2)^{\ell-1}}{\ell}.
    \end{equation}
    We obtain from \prettyref{eq:EP-ub} and \prettyref{eq:EP-lb} that 
    \[
        \max_{i\in[\ell-1]}|m_i(\nu)-m_i(\nu')|\ge \delta\frac{(C/4)^{\ell-1}}{\ell}.\tag*{\qedhere}
    \]
\end{proof}

\begin{proof}[Third proof of \prettyref{prop:stable1-union}]
    We continue to use $S$ in \prettyref{eq:supports-1} to denote the support of $\nu$ and $\nu'$. For any 1-Lipschitz function $f$, $\Expect_\nu f$ and $\Expect_{\nu'}f$ only pertain to function values $f(t_1),\dots,f(t_{\ell})$, which can be interpolated by a polynomial of degree $\ell-1$. 
    However, the coefficients of the interpolating polynomial can be arbitrarily large.\footnote{
        For example, the polynomial to interpolate $f(-\epsilon)=f(\epsilon)=\epsilon, f(0)=0$ is $P(x)=x^2/\epsilon$.
    }
    To fix this issue, we slightly modify the function $f$ on $S$ to $\tilde f$, and then interpolate $\tilde f$ with bounded coefficients. In this way we have 
    \begin{equation*}
        |\Expect_\nu f-\Expect_{\nu'}f|
        \le 2 \max_{x\in \{t_1,\dots,t_{\ell}\}}|\tilde f(x)-f(x)| + |\Expect_\nu P-\Expect_{\nu'}P|.
    \end{equation*}
    To this end, we define the values of $\tilde f$ recursively by
    \begin{equation}
        \label{eq:tilde-f-def}
        \tilde f(t_1)=f(t_1),\quad \tilde f(t_i)=\tilde f(t_{i-1})+(f(t_i)-f(t_{i-1}))\indc{t_i-t_{i-1}>\tau},
    \end{equation}
    where $\tau\le 2$ is a parameter we will optimize later. 
    From the above definition and the fact that $f$ is 1-Lipschitz $|\tilde f(t_i)-f(t_i)|\le |\tilde f(t_{i-1})-f(t_{i-1})|+\tau$, and thus $|\tilde f(x)-f(x)|\le \tau\ell$ for $x\in S$ by induction. The interpolating polynomial $P$ can be expressed using Newton formula \prettyref{eq:interpolation-Newton} as
    \[
        P(x)=\sum_{i=1}^\ell \tilde f[t_1,\dots,t_i] g_{i-1}(x),
    \]
    where $g_{r}(x)=\prod_{j=1}^{r}(x-t_j)$ such that $|\Expect_\nu[g_{r}]-\Expect_{\nu'}[g_{r}]|\le 2^{r}\delta$ by \prettyref{eq:polydiff} for $r\le \ell-1$. Since $f$ is 1-Lipschitz, we have $|\tilde f[t_i,t_{i+1}]|\le 1$ for every $i$. Higher order divided differences are recursively evaluated by \prettyref{eq:div-recursion}. We now prove
    \begin{equation}
        \label{eq:tilde-f-div}
        \tilde f[t_i,\ldots,t_{i+j}]\le (2/\tau)^{j-1},~\forall~i,j.
    \end{equation}
    by induction on $j$. Assume \prettyref{eq:tilde-f-div} holds for every $i$ and some fixed $j$. The recursion \prettyref{eq:div-recursion} gives 
    \[
        \tilde f[t_i,\ldots,t_{i+j+1}]=\frac{\tilde f[t_{i+1},\ldots,t_{i+j+1}]-\tilde f[t_{i},\ldots,t_{i+j}]}{t_{i+j+1}-t_{i}}.
    \]
    If $t_{i+j+1}-t_{i}<\tau$, then $\tilde f[t_i,\ldots,t_{i+j+1}]=0$ by \prettyref{eq:tilde-f-def}; otherwise, $\tilde f[t_i,\ldots,t_{i+j+1}]\le (\frac{2}{\tau})^{j}$ by triangle inequality. Using \prettyref{eq:tilde-f-div}, we obtain that
    \begin{equation*}
        |\Expect_\nu f-\Expect_{\nu'}f|
        \le 2\tau\ell+\sum_{i=2}^\ell \pth{\frac{2}{\tau}}^{i-2}2^{i-1}\delta
        \le 2\ell\pth{\tau+\frac{4^{\ell-2}}{\tau^{\ell-2}}\delta}.
    \end{equation*}
    The conclusion follows by letting $\tau=4\delta^{\frac{1}{\ell-1}}$.
\end{proof}

The proof of \prettyref{prop:stable2} uses a similar idea as the first proof of \prettyref{prop:stable1-union} to approximate step functions for all values of $\nu$ and $\nu'$; however, this is clearly impossible for non-discrete $\nu'$. For this reason, we turn from interpolation to majorization. A classical method to bound a distribution function by moments is to construct two polynomials that majorizes and minorizes a step function, respectively. Then the expectations of these two polynomials provide a sandwich bound for the distribution function. This idea is used, for example, in the proof of Chebyshev-Markov-Stieltjes inequality (cf.~\cite[Theorem 2.5.4]{Akhiezer1965}).
\begin{proof}[Proof of \prettyref{prop:stable2}]
    Suppose $\nu$ is supported on $x_1<x_2<\ldots<x_k$. Fix $t\in\reals$ and let $I_t(x)=\indc{x\le t}$. Suppose $x_m<t<x_{m+1}$. Similar to \prettyref{ex:Hermite}, we construct polynomial majorant and minorant using Hermite interpolation. To this end, let $P_t$ and $Q_t$ be the unique degree-$2k$ polynomials to interpolate $I_t$ with the following:
    \begin{center}
        \begin{tabular}{ c | c  c  c | c | c c c }
        \toprule
                &$x_1$  &\dots  &$x_m$  &$t$      &$x_{m+1}$  &\dots  &$x_k$ \\
        \hline
        $P$     &1     &\dots   &1      &1      &0          &\dots  &0 \\
        $P'$    &0     &\dots   &0      &any    &0          &\dots  &0 \\
        \hline
        $Q$     &1     &\dots   &1      &0      &0          &\dots  &0 \\
        $Q'$    &0     &\dots   &0      &any    &0          &\dots  &0 \\
        \bottomrule
        \end{tabular}
    \end{center}
    \vspace{0.5em}
    As a consequence of Rolle's theorem, $P_t\ge I_t\ge Q_t$ (cf.~\cite[p.~65]{Akhiezer1965}, and an illustration in \prettyref{fig:major-minor}).
    \begin{figure}[ht]
        \centering
        \includegraphics[width=.5\linewidth]{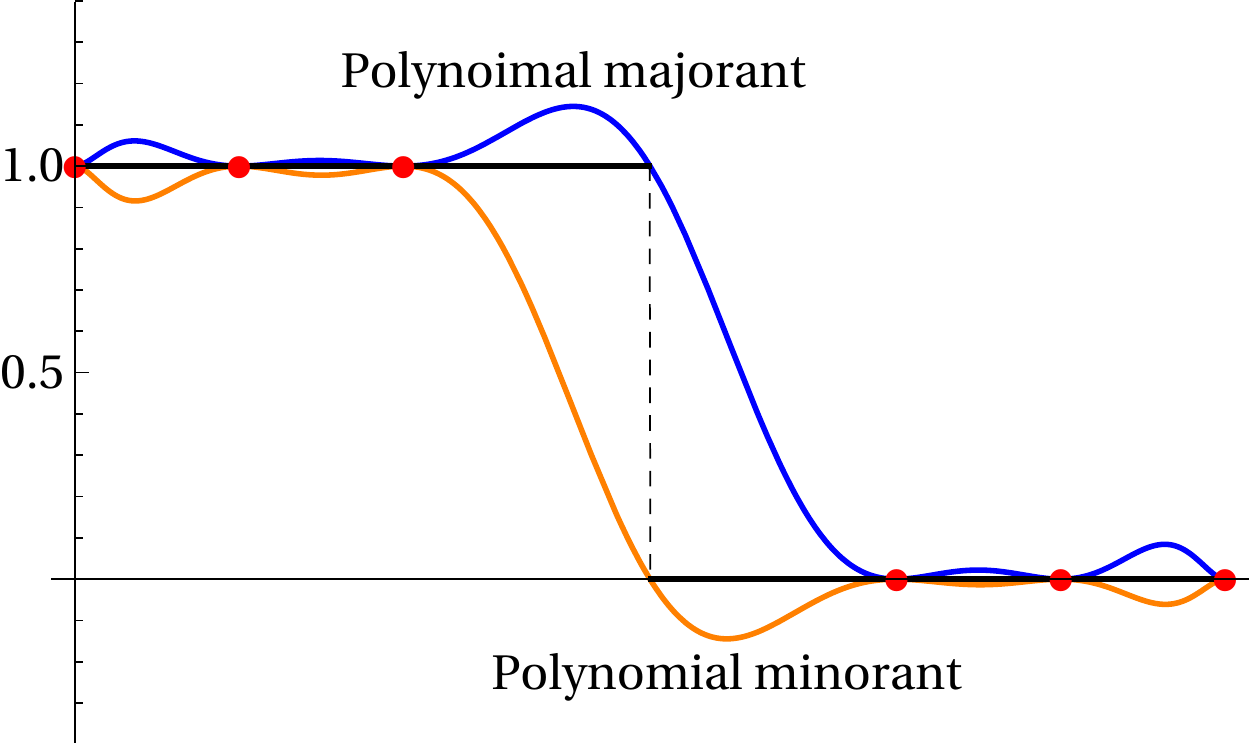}
        \caption{Polynomial majorant $P_t$ and minorant $Q_t$ that coincide with the step function on 6 red points.
            The polynomials are of degree 12, obtained by Hermite interpolation in \prettyref{sec:interpol}.
        }
        \label{fig:major-minor}
    \end{figure}
    Using Lagrange formula of Hermite interpolation \cite[pp.~52--53]{stoer.2002}, $P_t$ and $Q_t$ differ by
    \[
        P_t(x)-Q_t(x)=R_t(x)\triangleq\prod_i\pth{\frac{x-x_i}{t-x_i}}^2.
    \]
    The sandwich bound for $I_t$ yields a sandwich bound for the CDFs:
    \begin{align*}
        & \Expect_{\nu'}[Q_t]\le F_{\nu'}(t) \le \Expect_{\nu'}[P_t]=\Expect_{\nu'}[Q_t]+\Expect_{\nu'}[R_t],\\
        & \Expect_{\nu}[Q_t]\le F_{\nu}(t) \le \Expect_{\nu}[P_t]=\Expect_{\nu}[Q_t].
    \end{align*}
    Then the CDFs differ by 
    \begin{gather}
        |F_{\nu}(t)-F_{\nu'}(t)|\le (f(t)+g(t))\wedge 1\le f(t)\wedge 1+g(t)\wedge 1,\label{eq:CDF-diff-fg}\\
        f(t)\triangleq |\Expect_{\nu'}[Q_t]-\Expect_{\nu}[Q_t]|, \quad g(t)\triangleq\Expect_{\nu'}[R_t].\nonumber
    \end{gather}
    The conclusion will be obtained from the integral of CDF difference using \prettyref{eq:W1-CDF}. Since $R_t$ is almost surely zero under $\nu$, we also have $g(t)=|\Expect_{\nu'}[R_t]-\Expect_{\nu}[R_t]|$. Similar to \prettyref{eq:polydiff}, we obtain that 
    \[
        g(t)=|\Expect_{\nu'}[R_t]-\Expect_{\nu}[R_t]|\le \frac{2^{2k}\delta}{\prod_{i=1}^k (t-x_i)^2}.
    \]
    Hence, 
    \begin{equation}
        \label{eq:int-g}
        \int (g(t)\wedge 1) \diff t \le \int \pth{\frac{2^{2k}\delta}{\prod_{i=1}^k (t-x_i)^2}\wedge 1} \diff t \le 16k\delta^{\frac{1}{2k}},
    \end{equation}
    where the last inequality is proved in \prettyref{lmm:wedge1-integral} in \prettyref{sec:aux-lemmas}.

    Next we analyze $f(t)$. The polynomial $Q_t$ (and also $P_t$) can be expressed using Newton formula \prettyref{eq:interpolation-Newton} as
    \begin{equation}
        \label{eq:Qt}
        Q_t(x)=1+\sum_{i=2m+1}^{2k+1}I_t[t_1,\dots,t_i]g_{i-1}(x),
    \end{equation}
    where $t_1,\ldots,t_{2k+1}$ denotes the expanded sequence
    \[
        x_1,x_1,\ldots,x_m,x_m,t,x_{m+1},x_{m+1},\ldots,x_{k},x_{k}
    \]
    obtained by \prettyref{eq:new-sequence}, $g_{r}(x)=\prod_{j=1}^{r}(x-t_j)$, and we used $I_t[t_1,\dots,t_i]=0$ for $i=1,\ldots, 2m$. In \prettyref{eq:Qt}, the absolute values of divided differences are obtained in \prettyref{lmm:divided} in \prettyref{sec:aux-lemmas}: 
    \[
        I_t[t_1,\dots,t_i]\le \frac{\binom{i-2}{2m-1}}{(t-x_m)^{i-1}}.
    \]
    Using \prettyref{eq:Qt}, and applying the upper bound for $|\Expect_{\nu}[g_{i-1}]-\Expect_{\nu'}[g_{i-1}]|$ in \prettyref{eq:polydiff}, we obtain that, 
    \begin{gather*}
        f(t)=|\Expect_{\nu'}[Q_t]-\Expect_{\nu}[Q_t]|
        \le \sum_{i=2m+1}^{2k+1}\frac{\binom{i-2}{2m-1}2^{i-1}\delta}{(t-x_m)^{i-1}}
        \le \frac{k4^{2k}\delta}{(t-x_m)^{2k}},\\ \forall~x_m<t<x_{m+1},~m\ge 1.
    \end{gather*}
    If $t<x_1$, then $Q_t=0$ and thus $f(t)=0$. Then, analogous to \prettyref{eq:int-g}, we obtain that
    \begin{equation}
        \label{eq:int-f}
        \int (f(t)\wedge 1) \diff t \le 16k\delta^{\frac{1}{2k}}.
    \end{equation}
    Using \prettyref{eq:int-g} and \prettyref{eq:int-f}, the conclusion follows by applying \prettyref{eq:CDF-diff-fg} to the integral representation of Wasserstein distance \prettyref{eq:W1-CDF}.
\end{proof}

\begin{proof}[Proof of \prettyref{prop:stable1-separation}]
    The proof is analogous to the first proof of \prettyref{prop:stable1-union}, apart from a more careful analysis of polynomial coefficients. When each atom is at least $\gamma$ away from all but at most $\ell'$ other atoms, the left-hand side of \prettyref{eq:CDFdiff} is upper bounded by
    \[
        |F_{\nu}(t_r)-F_{\nu'}(t_r)|\le \frac{\ell 4^{\ell-1}\delta}{(t_{r+1}-t_r)^{\ell'}\gamma^{\ell-\ell'-1}},
    \]
    The remaining proof is similar. 
\end{proof}
\begin{proof}[Proof of \prettyref{prop:stable2-separation}]
    Similar to the proof of \prettyref{prop:stable1-separation}, this proof is analogous to \prettyref{prop:stable2} apart from a more careful analysis of polynomial coefficients. When every $t\in\reals$ is at least $\gamma$ away from all but $k'$ atoms, the left-hand sides of \prettyref{eq:int-g} and \prettyref{eq:int-f} are upper bounded by
    \begin{gather*}
        \int (g(t)\wedge 1) \diff t \le 4k\pth{\frac{2^{2k} \delta}{\gamma^{2(k-k')}}}^{1/(2k')},\\
        \int (f(t)\wedge 1) \diff t \le 4k\pth{\frac{k 4^{2k} \delta}{\gamma^{2(k-k')}}}^{1/(2k')}.
    \end{gather*}
    The remaining proof is similar. 
\end{proof}

\subsection{Auxiliary lemmas}
\label{sec:aux-lemmas}
In this subsection, we collect technical lemmas that are used to prove moment comparison theorems in \prettyref{sec:compare-proof}.
In Lemmas~\ref{lmm:divided} and \ref{lmm:Px-expand}, we provide upper bounds for the coefficients of the interpolating polynomials. 
\begin{lemma}
    \label{lmm:divided}
    Let $t_1\le t_2\le \dots$ be an ordered sequence (not necessarily distinct) and $t_r<t<t_{r+1}$. Let $I_t(x)=\indc{x\le t}$. Then 
    \begin{equation}
        \label{eq:divided}
        I_t[t_i,\dots,t_j]=(-1)^{i-r}\sum_{L\in\calL(i,j)}\prod_{(x,y)\in L}\frac{1}{t_x-t_y}, \quad i\le r < r+1 \le j,
    \end{equation}
    where $ \calL(i,j) $ is the set of lattice paths from $ (r,r+1) $ to $ (i,j) $ using steps $ (0,1) $ and $ (-1,0) $\footnote{
        Formally, for $a,b\in \naturals^2$, a lattice path from $ a $ to $ b $ using a set of steps $ S $ is a sequence $a=x_1,x_2,\dots,x_n=b$ with all increments $x_{j+1}-x_j\in S$. In the matrix representation shown in the proof, this corresponds to a path from $a_{r,r+1}$ to $a_{i,j}$ going up and right. This path consists of entries $(i,j)$ such that $i\le r<r+1\le j$, and thus in \prettyref{eq:divided} we always have $t_x\le t_r<t_{r+1}\le t_y$.
    }. 
    Furthermore,
    \begin{equation}
        \label{eq:divided-ub}
        |I_t[t_1,\dots,t_i]|\le \frac{\binom{i-2}{r-1}}{(t_{r+1}-t_{r})^{i-1}},\quad i\ge r+1.
    \end{equation}
\end{lemma}
\begin{proof}
    Denote by $a_{i,j}=I_t[t_i,\dots,t_j]$ when $i\le j$. It is obvious that $a_{i,i}=1$ for $i\le r$; $a_{i,i}=0$ for $i\ge r+1$; $a_{i,j}=0$ for both $i<j\le r$ and $j>i\ge r+1$. For $i\le r<r+1\le j$, the values can be obtained recursively by 
    \begin{equation}
        \label{eq:recursion}
        a_{i,j}=\frac{a_{i,j-1}-a_{i+1,j}}{t_i-t_j}.
    \end{equation}
    The above recursion can be represented in Neville's diagram previously introduced in \prettyref{sec:interpol}. In this proof, it is equivalently represented as an upper triangular matrix as follows:
    \begin{equation*}
        \begin{bmatrix}
            1 & 0 & \cdots      & 0      & a_{1,r+1}    &  \cdots  &      &      \\
            & 1   & \ddots      & \vdots & \vdots      &          &      &      \\
            &     & 1           & 0      & a_{r-1,r+1}  &  \cdots  &      &       \\
            &     &             & 1      & a_{r,r+1}    &  \cdots  &      &      \\
            &     &             &        & 0 & \cdots   &  0   &      \\
            &     &\text{\huge0}&        &   & \ddots   &\vdots&      \\
            &     &             &        &   &          &   0  &      
        \end{bmatrix}.
    \end{equation*}
    In the matrix, every $a_{i,j}$ is calculated using the value on its left and the value below it. The values on any path from $a_{r,r+1}$ to $a_{i,j}$ going up and right will contribute to the formula of $a_{i,j}$ in \prettyref{eq:divided}. The paths consist of two types: first go to $a_{i,j-1}$ and then go right; first go to $a_{i+1,j}$ and then go up. Formally, $\{L,(i,j):L\in\calL_{i,j-1}\}\cup \{L,(i,j):L\in\calL_{i+1,j}\}=\calL_{i,j}$. This will be used in the proof of \prettyref{eq:divided} by induction present next. The base cases ($r\Th$ row and $(r+1)\Th$ column) can be directly computed:
    \[
        a_{r,j}=\prod_{v=r+1}^j\frac{1}{t_{r}-t_v},\quad a_{i,r+1}=(-1)^{i-r}\prod_{v=i}^{r}\frac{1}{t_v-t_{r+1}}.
    \]
    Suppose \prettyref{eq:divided} holds for both $a_{i,j-1}$ and $a_{i+1,j}$. Then $a_{i,j}$ can be evaluated by
    \begin{align*}
        a_{i,j}&=\frac{(-1)^{i-r}}{t_i-t_j}\pth{\sum_{L\in\calL(i,j-1)}\prod_{(x,y)\in L}\frac{1}{t_x-t_y}+\sum_{L\in\calL(i+1,j)}\prod_{(x,y)\in L}\frac{1}{t_x-t_y}}\\
        &=(-1)^{i-r}\pth{\sum_{L\in\calL(i,j)}\prod_{(x,y)\in L}\frac{1}{t_x-t_y}}.
    \end{align*}
    For the upper bound in \prettyref{eq:divided-ub}, we note that $|\calL(i,j)|\le\binom{(r-1)+(i-(r+1))}{r-1}$ in \prettyref{eq:divided}, and each summand is at most $\frac{1}{(t_{r+1}-t_{r})^{i-1}}$ in magnitude.
\end{proof}

\begin{lemma}
    \label{lmm:Px-expand}
    Let
    \[
        P(x)=\prod_{i=1}^\ell(x-x_i)=\sum_{j=0}^{\ell}a_jx^j.
    \]
    If $|x_i|\le \beta$ for every $i$, then
    \[
        |a_j|\le \binom{\ell}{j}\beta^{\ell-j}.
    \]
\end{lemma}
\begin{proof}
    $P$ can be explicitly expanded and we obtain that
    \[
        a_{\ell-j}=(-1)^j\sum_{\{i_1, i_2,\dots, i_j\}\subseteq [\ell]}x_{i_1}\cdot x_{i_2}\cdot \ldots\cdot x_{i_j}.
    \]
    The summation consists of $\binom{\ell}{j}$ terms, and each term is at most $\beta^j$ in magnitude. 
\end{proof}

In the proofs of moment comparison theorems, the coefficients of the interpolating polynomials for $I_t(x)=\indc{x\le t}$ diverge when $t$ approaches atoms in the support set due to the discontinuity of the step function. 
Nevertheless, the difference between two CDFs is at most one.
Then, we obtain the upper bound of the Wasserstein distance by an integration in the following \prettyref{lmm:wedge1-integral}.
\begin{lemma}
    \label{lmm:wedge1-integral}
    Let $r\ge 2$.
    Then,
    \begin{equation*}
        \int \pth{\frac{\delta}{\prod_{i=1}^r|t-x_i|}\wedge 1}\diff t \le 4r\delta^{\frac{1}{r}}.
    \end{equation*}
\end{lemma}
\begin{proof}
    Without loss of generality, let $x_1\le x_2\le \dots \le x_r$.
    Note that
    \begin{equation*}
        \int \pth{\frac{\delta}{\prod_{i=1}^r|t-x_i|}\wedge 1} \diff t
        =\int_{-\infty}^{x_1}+\int_{x_1}^{\frac{x_1+x_2}{2}}+\int_{\frac{x_1+x_2}{2}}^{x_2}+\dots + \int_{x_r}^{\infty}.
    \end{equation*}
    There are $2r$ terms in the summation and each term can be upper bounded by
    \begin{equation*}
        \int_{x_i}^{\infty}\pth{\frac{\delta}{|t-x_i|^{r}}\wedge 1} \diff t 
        =\int_{0}^{\infty}\pth{\frac{\delta}{t^{r}}\wedge 1} \diff t 
        =\frac{r}{r-1}\delta^{\frac{1}{r}}.
    \end{equation*}
    The conclusion follows.
\end{proof}

To demonstrate the tightness of the exponents in \prettyref{prop:stable1}, we need two $k$-atomic distributions with identical first $2k-2$ moments in \prettyref{rmk:stable}, the existence of which is guaranteed by \prettyref{lmm:match2k-2}.
\begin{lemma}
    \label{lmm:match2k-2}
    Given any $2k$ distinct points $x_1<x_2<\dots<x_{2k}$, there exist two distributions $\nu$ and $\nu'$ supported on $\{x_1,x_3,\dots,x_{2k-1}\}$ and $\{x_2,x_4,\dots,x_{2k}\}$, respectively, such that ${\bf m}_{2k-2}(\nu)={\bf m}_{2k-2}(\nu')$.
\end{lemma}
\begin{proof}
    Consider the following linear equation
    \begin{equation*}
        \begin{pmatrix}
            1 & 1 & \cdots & 1 \\
            x_1 & x_2 & \cdots & x_{2k} \\
            \vdots  & \vdots  & \ddots & \vdots  \\
            x_1^{2k-2} & x_{2}^{2k-2} & \cdots & x_{2k}^{2k-2} 
        \end{pmatrix}
        \begin{pmatrix}
            w_1 \\
            w_2 \\
            \vdots  \\
            w_{2k}
        \end{pmatrix}
        =0,
    \end{equation*}
    This underdetermined system has a non-zero solution. Let $w$ be a solution with $\|w\|_1=2$.
    Since all weights sum up to zero, then positive weights in $w$ sum up to $1$ and negative weights sum up to $-1$.
    Let one distribution be supported on $x_i$ with weight $w_i$ for $w_i>0$, and the other one be supported on the remaining $x_i$'s with the corresponding weights $|w_i|$.
    Then these two distribution match the first $2k-2$ moments.

    It remains to show that the weights in any non-zero solution have alternating signs.
    Note that all weights are non-zero: if one $w_i$ is zero, then the solution must be all zero since the Vandermonde matrix is of full row rank. To verify the signs of the solution, without loss generality, assume that $w_{2k}=-1$ and then
    \begin{equation*}
        \begin{pmatrix}
            1 & \cdots & 1 \\
            x_1 & \cdots & x_{2k-1} \\
            \vdots & \ddots & \vdots  \\
            x_1^{2k-2} & \cdots & x_{2k-1}^{2k-2} 
        \end{pmatrix}
        \begin{pmatrix}
            w_1 \\
            w_2 \\
            \vdots  \\
            w_{2k-1}
        \end{pmatrix}
        =
        \begin{pmatrix}
            1 \\
            x_{2k} \\
            \vdots  \\
            x_{2k}^{2k-2} 
        \end{pmatrix}.
    \end{equation*}
    The solution has an explicit formula that $w_i=P_i(x_{2k})$ where $P_i$ is an interpolating polynomial of degree $2k-2$ satisfying $P_i(x_j)=1$ for $j=i$ and $P_i(x_j)=0$ for all other $j\le 2k-1$.
    Specifically, we have $w_i=\frac{\prod_{j\ne i,j\le 2k-1}(x_{2k}-x_j)}{\prod_{j\ne i,j\le 2k-1}(x_{i}-x_j)}$, which satisfies $w_i>0$ for odd $i$ and $w_i<0$ for even $i$.
    The proof is complete.
\end{proof}

\section{Moment comparison for Gaussian mixture densities}
\label{sec:compare-higher-moments}
In addition to comparing discrete distributions in the Wasserstein distance by means of moments, in this section we show that the information geometry of finite Gaussian mixtures can be characterized by moments of the mixing distributions. We start with a comparison result of higher moments by lower ones. 
\begin{lemma}
    \label{lmm:high-moments}
    If $U$ and $U'$ each takes at most $k$ values in $[-1,1]$, and $|\Expect[U^j]-\Expect[U'^j]|\le \epsilon$ for $j=1,\dots,2k-1$, then, for any $\ell\ge 2k$,
    \[
        |\Expect[U^\ell]-\Expect[U'^\ell]|\le 3^\ell \epsilon.
    \]
\end{lemma}
\begin{proof}
    Let $f(x)=x^\ell$ and denote the atoms of $U$ and $U'$ by $x_1<\dots<x_{m}$ for $m\le 2k$. The function $f$ can be interpolated on $x_1,\dots,x_{m}$ using a polynomial $P$ of degree at most $2k-1$, which, in the Newton form \prettyref{eq:interpolation-Newton}, is
    \[
        P(x)=\sum_{i=1}^{m}f[x_1,\dots,x_i]g_{i-1}(x)
        =\sum_{i=1}^{m}\frac{f^{(i-1)}(\xi_i)}{(i-1)!}g_{i-1}(x),
    \]
    for some $\xi_i\in[x_1,x_i]$, where $g_{r}(x)=\prod_{j=1}^{r}(x-x_j)$ and we used the intermediate value theorem for the divided differences (see \cite[(2.1.4.3)]{stoer.2002}).  Note that for $\xi_i\in[-1,1]$, $|f^{(i-1)}(\xi_i)|\le \frac{\ell!}{(\ell-1+i)!}$. Similar to \prettyref{eq:polydiff}, we obtain that
    \[
        |\Expect[U^\ell]-\Expect[U'^\ell]|
        =|\Expect[P(U)]-\Expect[P(U')]|
        \le \sum_{i=1}^{m}\binom{\ell}{i-1}2^{i-1}\epsilon
        \le 3^\ell \epsilon.\tag*{\qedhere}
    \]
\end{proof}

In the context of learning Gaussian mixtures, we can bound the distance between two mixture densities by comparing the moments of the mixing distributions.
\begin{lemma}[Bound $\chi^2$-divergence using moments]
    \label{lmm:chi2-moments}
    Suppose all moments of $\nu$ and $\nu'$ exist, and $\nu'$ is centered with variance $\sigma^2$.
    Then,
    \[
        \chi^2(\nu*N(0,1)\|\nu'*N(0,1))\le e^{\frac{\sigma^2}{2}}\sum_{\ell\ge 1}\frac{(m_{\ell}(\nu)-m_{\ell}(\nu'))^2}{\ell!}.
    \]
\end{lemma}
\begin{proof}
    The proof is similar to \prettyref{thm:chi2-gm}. Let $\Delta m_{\ell} =  m_{\ell}(\nu)-m_{\ell}(\nu')$, $f$ and $f'$ denote the density of $\nu*N(0,1)$ and $\nu'*N(0,1)$, respectively.
    By an analogous argument to \eqref{eq:chi2-matching}, we obtain that
    \begin{align*}
        &\phantom{{}={}}\chi^2(f\|f') = \int \frac{(f(x)-f'(x))^2}{f'(x)}\diff x\\
        &\leq e^{\frac{\sigma^2}{2}} \expect{\pth{\sum_{\ell\ge 1}H_{\ell}(Z)\frac{\Delta m_{\ell}}{\ell!}}^2}       
        =  e^{\frac{\sigma^2}{2}} \sum_{\ell\ge 1}\frac{(\Delta m_{\ell})^2}{\ell!},
    \end{align*}
    where $Z\sim N(0,1)$ and the last step follows from the orthogonality property of Hermite polynomials in \prettyref{eq:Hermite-ortho}.
\end{proof}


Next we show that the distance (be it Hellinger, KL divergence, or $\chi^2$-divergence; cf.~\prettyref{sec:notation}) between two $k$-component Gaussian mixture models is characterized up to constant factors by the differences of their first $2k-1$ moments. This result forms the basis for proving density estimation guarantees of moment-based learning algorithms in \prettyref{chap:gm}.
\begin{proposition}
\label{prop:compare-density}
Let $\nu$ and $\nu'$ be $k$-atomic distributions supported on $[-1,1]$. 
Let $f$ and $f'$ denote the density of $\nu*N(0,1)$ and $\nu'*N(0,1)$, respectively. 
Then, for any $\rho \in \{H^2,D,\chi^2\}$,
\[
c\max_{\ell \le 2k-1} |m_\ell(\nu)-m_\ell(\nu')|^2
\le
\rho(f, f') \le  C\max_{\ell \le 2k-1} |m_\ell(\nu)-m_\ell(\nu')|^2,
\]
for constants $c,C$ only depending on $k$. 
\end{proposition}
\begin{proof}
Since
\[
H^2(P, Q) \le D(P \| Q) \le \chi^2(P \| Q ),
\]
(see, e.g., \cite[Section 2.4.1]{Tsybakov09}), 
it suffices to prove the lower bound for $H^2$ and the upper bound for $\chi^2$. 
We first prove the upper bound. 
Let $U'\sim\nu'$. 
If $\Expect[U']=0$, then the upper bound follows from Lemmas~\ref{lmm:high-moments} and \ref{lmm:chi2-moments}. 
If $\mu\triangleq\Expect[U']\ne 0$ (which satisfies $|\mu|\le 1$), by the shift-invariance of $\chi^2$-divergence, applying the following simple binomial expansion yields the desired upper bound:
\[
|m_\ell(U-\mu)-m_\ell(U'-\mu)|
\le \sum_{k=0}^\ell \binom{\ell}{k}|m_k(U)-m_k(U')| |\mu|^{\ell-k}.
\]

Next we show the lower bound that, for any $(2k-1)$-times differentiable function $h$,
\[
H(f, f') \ge c |\Expect[h(U)]-\Expect[h(U')]|,
\]
for a constant $c$ depending on $k$ and $h$, where $U\sim \nu$ and $U'\sim \nu'$.
The proof relies on Newton interpolating polynomials and orthogonal polynomials in Sections~\ref{sec:interpol} and \ref{sec:orthogo}, respectively.

Note that $f=\nu*\phi$ and $f'=\nu'*\phi$, where $\phi$ denotes the standard normal density. 
Equivalently, we have $f(y) = \Expect[\phi(y-U)]$ and $f(y) = \Expect[\phi(y-U')]$. 
Since $f$ and $f'$ are bounded above by $\frac{1}{\sqrt{2\pi}}$, it suffices to prove that
\begin{equation}
\label{eq:l2_moments}
\|f-f'\|_2 \ge c |\Expect[h(U)]-\Expect[h(U')]|.
\end{equation}
Define $\alpha_j(y) \triangleq \sqrt{\frac{\sqrt{2} \phi(\sqrt{2}y)}{j!}} H_j(\sqrt{2} y)$ which forms an orthonormal basis on $L^2(\reals, \diff y)$ in view of \prettyref{eq:Hermite-ortho}.
Then $f$ admits an orthogonal expansion\footnote{Note that unlike \prettyref{eq:GM-orthogonal} which is with respect to the Gaussian weight, this orthogonal expansion is unweighted.} $f(y)=\sum_{j\geq 0} a_j \alpha_j(y)$ with coefficients
\[
a_j = \int f(y)\alpha_j(y) \diff y = \Expect[(\alpha_j * \phi)(U)] 
= \frac{1}{2^{\frac{j+1}{2}}     \pi ^{\frac{1}{4}} \sqrt{j!}}\Expect[U^j e^{-\frac{U^2}{4}}],
\]
where the last equality follows from the fact that \cite[7.374.6, p.~803]{GR}
\[
    (\alpha_j*\phi) (y) = \frac{1}{2^{\frac{j+1}{2}}  \pi ^{\frac{1}{4}} \sqrt{j!}} y^j e^{-\frac{y^2}{4}}.
\]
By a similar orthogonal expansion for $f'$, we obtain that
\[
\|f- f'\|_2^2 
= \sum_{j \ge 0} \frac{1}{j! 2^{j+1} \sqrt{\pi} } \pth{ \Expect[U^j e^{-U^2/4}] - \Expect[U'^j e^{-U'^2/4}] }^2. 
\]
In particular, for each $j\geq 0$,
\begin{align}
 \|f- f'\|_2 \geq \frac{1}{\sqrt{j! 2^{j+1} \sqrt{\pi}}} |\Expect[U^j e^{-U^2/4}] - \Expect[U'^j e^{-U'^2/4}]|.
 \label{eq:l2_expmoments}
\end{align}
Thus, our strategy of proving \eqref{eq:l2_moments} is to interpolate $h(y)$ by linear combinations of
$\{y^j e^{-y^2/4}: j=0,\ldots,2k-1\}$ over the atoms of $\nu$ and $\nu'$, a total of at most $2k$ points.
Clearly, this is equivalent to the standard polynomial interpolation of $\tilde h(y) \triangleq  h(y) e^{y^2/4}$ by a degree-($2k-1$) polynomial.
Specifically, denote by $T$ the union of the support of $\nu$ and $\nu'$. Let $P(y) = \sum_{j=0}^{2k-1}  b_j y^j$ be the interpolating polynomial of $\tilde h$ on the set $T$.
Then
\begin{align*}
|\Expect[h(U)] - \Expect[h(U')] | 
&\stepa{=} |\Expect[P(U)e^{-U^2/4} - \Expect[P(U')e^{-U'^2/4}  | \\
&\le \sum_{j=0}^{2k-1} |b_j| | \Expect[U^j e^{-U^2/4}] - \Expect[U'^j e^{-U'^2/4}] |  \\
&\stepb{\le} \|f- f'\|_2 \sum_{j=0}^{2k-1} |b_j| \sqrt{j! 2^{j+1} \sqrt{\pi}},
\end{align*}
where 
(a) is by virtue of interpolation so that $P(U)=\tilde h(U) = e^{U^2/4} h(U)$ and similarly for $U'$; 
(b) follows from \prettyref{eq:l2_expmoments}.
It remains to bound the coefficient $b_j$ independently of the set $T$, 
 which is accomplished by the next lemma.
\end{proof}

\begin{lemma}
\label{lmm:interp}  
    Let $h$ be an $m$-times differentiable function on the interval $[-R,R]$, whose derivatives are bounded by
    $ |h^{(i)}(x)| \leq M$ for all $0\leq i\leq m$ and all $x\in[-R,R]$. 
    Then for any $m\geq 1$ and $R>0$, there exists a positive constant $C=C(m,R,M)$, such that the following holds.
    For any set of distinct nodes $T=\{x_0,\ldots,x_m\} \subset [-R,R]$, denote by $P(x)=\sum_{i=0}^{m} b_j x^j$ 
    the unique interpolating polynomial of degree at most $m$ of $h$ on $T$.
    Then $\max_{0\leq j \leq m} |b_j| \leq C$.  
\end{lemma}
\begin{proof}
Express $P$ in the Newton form~\prettyref{eq:interpolation-Newton}:
\begin{align*}
P(y) &= \sum_{i = 0}^{m} h[x_0, \hdots, x_i] \prod_{j = 0}^{i-1} (y - x_j)
\end{align*}
By the intermediate value theorem, finite differences can be bounded by derivatives as follows: (c.f.~\cite[(2.1.4.3)]{stoer.2002})
\[
|h[x_0, \hdots, x_i]| \leq \frac{1}{i!} \sup_{|\xi|\leq R}|h^{(i)}(\xi)|.
\]
Let $\prod_{j = 0}^{i-1} (y - x_j)  = \sum_{j = 0}^{i} c_{ij} y^j$. Since $|x_j|\leq R$, 
$|c_{ij}| \leq C_1=C_1(R,m)$ all $i,j$.
This completes the proof.
\end{proof}

\prettyref{prop:compare-density} has a natural extension to multiple dimensions in terms of \emph{moment tensors}. The moment tensor of order $\ell$ of a distribution $\nu$ on $\reals^d$ is a symmetric tensor defined as
\begin{equation}
M_{\ell}(\nu) \triangleq  \Expect_{U\sim\nu}[\underbrace{U\otimes \cdots \otimes U}_{\text{$\ell$ times}}].
\label{eq:moment-tensor}
\end{equation}
In particular, $M_1$ and $M_2-M_1 M_1^\top$ are the mean and the covariance matrix, respectively. 
In the special case where $\nu = \sum_{j=1}^k w_j \mu_j $ is $k$-atomic, $M_\ell(\nu)=\sum_{j=1}^k w_j \mu_j^{\otimes \ell}$ is a rank-$k$ tensor. 
Generalizing \prettyref{prop:compare-density}, the following result provides a moment comparison for multivariate Gaussian mixtures with \emph{dimension-independent} factors. 
\begin{proposition}[{\cite[Theorem 4.1]{DWYZ20}}]
\label{prop:compare-density-d-dim}
Let $\nu$ and $\nu'$ be $k$-atomic distributions supported on the ball $B(0,R)\subseteq \reals^d$.
Let $f$ and $f'$ denote the density of $\nu*N(0,I_d)$ and $\nu'*N(0,I_d)$, respectively. 
Then, for any $\rho \in \{H^2,D,\chi^2\}$,
\[
c\max_{\ell \le 2k-1} \|M_\ell(\nu)-M_\ell(\nu')\|_F^2
\le
\rho(f,f') 
\le  C\max_{\ell \le 2k-1} \|M_\ell(\nu)-M_\ell(\nu')\|_F^2,
\]
where the constants $c,C$ may depend on $k$ and $R$ but \emph{not} $d$, and $\|\cdot\|_F$ denotes the Frobenius norm (the $\ell_2$-norm of the vectorization) of a tensor.
\end{proposition}

\begin{remark}
	Related to \prettyref{lmm:chi2-moments}, for general (not necessarily finite) Gaussian mixtures, the following inequality is shown in 
	\cite[Theorem 9 and Appendix A]{bandeira2017optimal}:
	\[
\sum_{\ell = 1}^\infty \frac{ c^\ell}{\ell!} \|M_\ell(\nu)-M_\ell(\nu')\|_F^2
\le D(f\|f')  \le \sum_{\ell = 1}^\infty \frac{ C^\ell}{\ell!} \|M_\ell(\nu)-M_\ell(\nu')\|_F^2,
\]
where $c$ and $C$ depends only on the support of $\nu$ and $\nu'$. 
Compared with this result which involves all moments, the point of Propositions \ref{prop:compare-density} and \ref{prop:compare-density-d-dim} is that for $k$-component mixtures it suffices to compare moments up to degree $2k-1$. 
\end{remark}





\chapter{Learning Gaussian mixtures}
\label{chap:gm}
In this chapter, we discuss the problem of learning Gaussian mixtures using the framework and tools developed in \prettyref{chap:framework}.
Consider a $k$-component Gaussian location mixture model, where each observation is distributed as
\begin{equation}
    \label{eq:model}
    X\sim \sum_{i=1}^k w_i N(\mu_i, \sigma^2).
\end{equation}
Here $w_i$ is the mixing weight such that $w_i\ge 0$ and $\sum_i w_i = 1$, $\mu_i$ is the mean (center) of the $i\Th$ component, and $\sigma$ is the common standard deviation. Equivalently, we can write the distribution of an observation $ X $ as a convolution
\begin{equation}
    \label{eq:model-conv}
    X\sim\nu*N(0,\sigma^2),
\end{equation}
where $\nu = \sum_{i=1}^k w_i \delta_{\mu_i}$ denotes the \emph{mixing distribution}. Thus, we can write $X=U+\sigma Z$, where $U\sim\nu$ is referred to as the latent variable, and $Z$ is standard normal and independent of $U$.
We adopt the framework in \prettyref{chap:framework} with the goal of estimating the mixing distribution $\nu$ with respect to the Wasserstein distance \prettyref{eq:w1}. 
Equivalently, estimating the mixing distribution can be viewed as a deconvolution problem, where the goal is to recover the distribution $\nu$ using sample drawn from the convolution \prettyref{eq:model-conv}.

The focus of the present chapter is moment-based learning algorithms.
We first discuss the problems with the \emph{method of moments} and the \emph{generalized method of moments} in \prettyref{sec:gm-moment}, and the key idea to overcome those issues.
We study the mixture model \prettyref{eq:model} for known and unknown $\sigma$ in Sections~\ref{sec:known} and \ref{sec:unknown}, respectively. 
The optimality of the proposed methods is demonstrated in \prettyref{sec:gm-lb}, and open problems are discussed in \prettyref{sec:discuss}.
Finally, \prettyref{sec:gm-bib} contains bibliographical notes.

\section{Moment-based learning algorithms}
\label{sec:gm-moment}
\subsection{Failure of the ordinary method of moments}
\label{sec:mmfail}
The method of moments, commonly attributed to Pearson \cite{Pearson1894}, produces an estimator by equating the population moments to the sample moments.
While conceptually simple, this method suffers from the following problems, especially in the context of mixture models:
\begin{itemize}
    \item \emph{Solubility}: 
    the method of moments entails solving a multivariate polynomial system, in which one frequently encounters non-existence or non-uniqueness of statistically meaningful solutions.
    \item \emph{Computation}: solving moment equations can be computationally intensive. For instance, for $k$-component Gaussian mixture models, the system of moment equations consist of $2k-1$ polynomial equations with $2k-1$ variables.
    \item \emph{Accuracy}: 
    existing statistical literature on the method of moments \cite{VdV00,Hansen1982} either shows mere consistency under weak assumptions, or proves asymptotic normality assuming very strong regularity conditions (so that delta method works), which generally do not hold in mixture models since the convergence rates can be slower than parametric.
    Some results on nonparametric rates are known (cf.~\cite[Theorem 5.52]{VdV00} and \cite[Theorem 14.4]{Kosorok}) but the conditions are extremely hard to verify.
\end{itemize}

To explain the failure of the vanilla method of moments in Gaussian mixture models, we analyze the following simple two-component example:
\begin{example}
    \label{ex:mmfail}
    Consider a Gaussian mixture model with two unit variance components: $X\sim w_1 N(\mu_1, 1)+w_2 N(\mu_2, 1)$. Since there are three parameters $\mu_1,\mu_2$ and $w_1=1-w_2$, we use the first three moments and solve the following system of equations:
    \begin{equation}
        \begin{aligned}
            \Expect_n[X]&=\Expect[X]=w_1\mu_1+w_2\mu_2,\\
            \Expect_n[X^2]&=\Expect[X^2]=w_1\mu_1^2+w_2\mu_2^2+1,\\
            \Expect_n[X^3]&=\Expect[X^3]=w_1\mu_1^3+w_2\mu_2^3+3(w_1\mu_1+w_2\mu_2),
        \end{aligned}    
        \label{eq:MM-X}
    \end{equation}
    where $\Expect_n[X^i] \triangleq \frac{1}{n} \sum_{j=1}^n X_j^i$ denotes the $i^{\rm th}$ moment of the empirical distribution from $n$ \iid~observations. The right-hand sides of \prettyref{eq:MM-X} are related to the moments of the mixing distribution by a linear transformation, which allow us to equivalently rewrite the moment equations \prettyref{eq:MM-X} as:
    \begin{equation}
        \begin{aligned}
            \Expect_n[X]&=\Expect[U]=w_1\mu_1+w_2\mu_2,\\
            \Expect_n[X^2-1]&=\Expect[U^2]=w_1\mu_1^2+w_2\mu_2^2,\\
            \Expect_n[X^3-3X]&=\Expect[U^3]=w_1\mu_1^3+w_2\mu_2^3,
        \end{aligned}
        \label{eq:MM-U}
    \end{equation}
    where $U \sim w_1\delta_{\mu_1}+w_1\delta_{\mu_2}$. 
        It turns out that with finitely many observations, there is always a non-zero chance that \prettyref{eq:MM-U} has no solution; even when sample size tends to infinity, it is still possible that the solution remains non-existent with constant probability. To see this, note that, from the first two equations of \prettyref{eq:MM-U}, the solution does not exist whenever
    \begin{equation}
        \label{eq:CS-fail}
        \Expect_n[X^2]-1<\Expect_n^2[X],
    \end{equation}
    that is, the Cauchy-Schwarz inequality fails. Consider the extreme case $\mu_1=\mu_2=0$, i.e., $X\sim N(0,1)$. Then \prettyref{eq:CS-fail} is equivalent to
    \[
        n(\Expect_n[X^2]-\Expect_n^2[X])< n,
    \]
    where the left-hand side follows the $\chi^2$-distribution with $n-1$ degrees of freedom. Thus, \prettyref{eq:CS-fail} occurs with probability  approaching $\frac{1}{2}$ as $n$ diverges, according to the central limit theorem.
\end{example}

\subsection{Generalized method of moments}
\label{sec:gmm}
The Generalized Method of Moments (GMM) is
introduced by Hansen \cite{Hansen1982}. GMM is a widely used methodology for analyzing economic and financial data (cf.~\cite{Hall2005} for a thorough review).
 Instead of solving the moment equations, GMM aims to minimize the difference between estimated and fitted moments:
\begin{equation}
    \label{eq:GMM-opt}
    Q(\theta)=(\hat m- m(\theta))^\top W(\hat m - m(\theta)),
\end{equation}
where $\hat m$ is the estimated moment, $\theta$ is the model parameter, and $W$ is a positive semidefinite weighting matrix. The minimizer of $Q(\theta)$ serves as the GMM estimate for the unknown model parameter $\theta_0$. In general the objective function $Q$ is nonconvex in $\theta$, notably under the Gaussian mixture model with $\theta$ corresponding to the unknown means and weights, which is hard to optimize. 

In theory, the optimal weighting matrix $W^*$ that minimizes the asymptotic variance is the inverse of $\lim_{n\diverge}\cov[\sqrt{n}(\hat m- m(\theta_0))]$, which depends the unknown model parameters $\theta_0$. Thus, a popular approach is a two-step estimator \cite{Hall2005}:  
\begin{enumerate}
    \item A suboptimal weighting matrix, \eg, identify matrix, is used in the GMM to obtain a consistent estimate of $\theta_0$ and hence a consistent estimate $\hat W$ for $W^*$;
    \item Re-estimate $\theta_0$ using the weighting matrix $\hat W$. 
\end{enumerate}
The above two-step approach minimizes the asymptotic variance, but again both steps involve optimizing nonconvex objectives that are computationally challenging. 
Empirical comparison of GMM with other methods are implemented in \cite{WY18}.

\subsection{Efficiently solving GMM via moment space optimization}
In view of \prettyref{ex:mmfail}, we note that the main issue with the ordinary method of moments is the following: although individually each moment estimate is accurate ($\sqrt{n}$-consistent), jointly they do not correspond to the moments of any distribution. 
Moment vectors satisfy many geometric constraints as introduced in \prettyref{sec:moment} and lie in the moment space described by \prettyref{thm:moment-psd}.
Thus for any model parameters, with finite sample size the method of moments fails with nonzero probability whenever the noisy estimates escape the moment space. Furthermore, even with infinite sample size, it can also provably happen with constant probability when the order of the true mixture model is strictly less than the postulated order $k$, or equivalently, the population moments lie on the boundary of the moment space.

The GMM resolves the issue of solubility of the ordinary method of moments, but the prevailing approach requires solving nonconvex optimization that is computationally expensive.
Next we propose a procedure, named the \emph{denoised method of moments} (DMM), that simultaneously resolves all three issues of the method of moments.
The key idea is to jointly denoise the moment estimates by projecting them onto the moment space, which is a convex set characterized by PSD conditions of moment matrices in \prettyref{thm:moment-psd}.
The details of DMM will be presented in \prettyref{sec:known}. 
In particular, the program to be presented in \prettyref{eq:project} with the Euclidean norm is \emph{equivalent} to GMM with the identity weighting matrix. Therefore DMM can be used as an exact solver for GMM in the Gaussian location mixture model. 
The two-step approach of GMM can be similarly implemented by incorporating the corresponding weighting matrices in \prettyref{eq:project}.

\section{Known variance}
\label{sec:known}
\subsection{Denoised method of moments}
The \emph{denoised method of moments} consists of
three main steps: (1) compute noisy estimates of moments, e.g., the unbiased estimates; (2) jointly denoise the moment estimates by projecting them onto the moment space as described in \prettyref{sec:moment} by semidefinite programming (SDP); (3) execute the usual method of moments. 
It turns out that the extra step of projection resolves the three issues of the vanilla version of the method of moments identified in \prettyref{sec:mmfail} simultaneously:
\begin{itemize}
    \item \emph{Solubility}: a unique statistically meaningful solution is guaranteed to exist by the classical theory of moments;
    \item \emph{Computation}: the solution can be found through an efficient algorithm (Gaussian quadrature) instead of invoking generic solvers of polynomial systems;
    \item \emph{Accuracy}: the solution provably achieves the optimal rate of convergence, and automatically adaptive to the clustering structure of the population.
\end{itemize}
The complete algorithm is summarized in \prettyref{algo:known}.

\begin{algorithm}[ht]
    \caption{Denoised method of moments (DMM) with known variance.}
    \label{algo:known}
    \begin{algorithmic}[1]
        \REQUIRE $n$ independent observations $X_1,\dots,X_n$, order $k$, variance $\sigma^2$, interval $I=[a,b]$. 
        \ENSURE estimated mixing distribution.
        \FOR{$r=1$ \TO $2k-1$}\label{line:for}
        \STATE{$\hat{\gamma}_r=\frac{1}{n}\sum_i X_i^r$}
        \STATE{$\tilde{m}_r=r!\sum_{i=0}^{\Floor{r/2}}\frac{(-1/2)^i}{i!(r-2i)!}\hat{\gamma}_{r-2i}\sigma^{2i}$}\label{line:tildem}
        \ENDFOR\label{line:endfor}
        \STATE{Let $\hat m$ be the optimal solution of the following:}
        \begin{equation}
            \label{eq:project}
            \min\{\Norm{\tilde m - \hat m}: \hat m\in\calM_{2k-1}([a,b])\},
        \end{equation} 
        where $\tilde m=(\tilde m_1,\dots,\tilde m_{2k-1})$ and $\calM_{2k-1}([a,b])$ denotes the moment space characterized by \prettyref{eq:moment-psd}.
        \STATE{Report the outcome of Gaussian quadrature (\prettyref{algo:quadrature}) with input $\hat m$.}
    \end{algorithmic}
\end{algorithm}

We estimate the moments of the mixing distribution in lines~\ref{line:for} to \ref{line:endfor}. The unique unbiased estimators for the polynomials of the mean parameter in a Gaussian location model are Hermite polynomials \prettyref{eq:Hermite1}
such that $\Expect H_r(X)=\mu^r$ when $X\sim N(\mu,1)$. Thus, if we define
\begin{equation}
    \label{eq:Hermite}
    h_r(x,\sigma)=\sigma^rH_r(x/\sigma)=r!\sum_{j=0}^{\floor{r/2}} \frac{(-1/2)^j }{j!(r-2j)!}\sigma^{2j}x^{r-2j},
\end{equation}
then $\Expect h_r(X,\sigma)=\mu^r$ when $X\sim N(\mu,\sigma^2)$. Hence, by linearity, $\tilde m_r$ is an unbiased estimate of $m_r(\nu)$. The variance of $\tilde m_r$ is analyzed in the following lemma by an upper bound for the variance of $h_r$. See a proof in \cite{WY18}.
\begin{lemma}
    \label{lmm:var-tildem}
    If $X_1,\dots,X_n\iiddistr \nu*N(0,\sigma^2)$ and $\nu$ is supported on $[-M,M]$, then
    \[
        \var[\tilde m_r]\le \frac{1}{n}(O(M+\sigma\sqrt{r}))^{2r}.
    \]
\end{lemma}

As observed in \prettyref{sec:mmfail}, the major cause for the failure of the usual method of moments is that the unbiased estimate $\tilde m$ needs not constitute a legitimate moment sequence, despite the consistency of each individual $\tilde m_i$. To resolve this issue, we project $\tilde m$ to the moment space using \prettyref{eq:project}. As explained in \prettyref{thm:moment-psd}, the description of the moment space consists of positive semidefinite constraints, and thus the optimal solution of \prettyref{eq:project} can be obtained by solving an SDP.
 In fact, to achieve the optimal theoretical guarantee,
it suffices to solve a \emph{feasibility} program and find any valid moment vector $\hat m$ that is within the desired $\frac{1}{\sqrt{n}}$ statistical accuracy. 
Finally, we note that \prettyref{eq:project} can be replaced by optimizing weighted quadratic objectives instead of the $L_2$ norm, and the two-step procedure of GMM in \prettyref{sec:gmm} can also be similarly implemented.

Now that $\hat m$ is indeed a valid moment sequence, we use the Gaussian quadrature introduced in \prettyref{sec:GQ} (see \prettyref{algo:quadrature}) to find the unique $k$-atomic distribution $\hat\nu$ such that ${\bf m}_{2k-1}(\hat \nu)=\hat m$. Using \prettyref{algo:known}, $\tilde m$ is computed in $O(kn)$ time, the semidefinite programming is solvable in $O(k^{6.5})$ time using the interior-point method (see \cite{WSV2012}), and the Gaussian quadrature can be evaluated in $O(k^3)$ time \cite{GW1969}. 
In view of the global assumption \prettyref{eq:k}, \prettyref{algo:known} can be executed in $O(kn)$ time.

\subsection{Optimal rates of the DMM estimator}
We now present the statistical guarantee for the DMM estimator.
Throughout this chapter, we assume that the number of components satisfies
\begin{equation}
k = O\pth{\frac{\log n}{\log\log n}}.
\label{eq:k}
\end{equation}
If the order of mixture is large, namely, $k\ge \Omega(\frac{\log n}{\log\log n})$, including continuous mixtures, then one can approximate it by a finite mixture with $O(\frac{\log n}{\log\log n})$ components and estimate the mixing distribution using the DMM estimator. Furthermore, this method turns out to be optimal (see \prettyref{thm:large-k} at the end of this subsection).
\begin{theorem}[Optimal rate]
    \label{thm:main-W1}
    Suppose that $\nu$ is supported on $[-M,M]$ for $M\ge 1$ and $\sigma$ is bounded from above 
    by a constant, and both $k$ and $M$ are given.
    Then, \prettyref{algo:known} produces $\hat\nu$ such that, with probability at least $1-\delta$,
    \begin{equation}
        \label{eq:known-variance-W1}
        W_1(\nu,\hat \nu)\le O\pth{Mk^{1.5}\pth{\frac{n}{\log(1/\delta)}}^{-\frac{1}{4k-2}}};
    \end{equation}
\end{theorem}
The above convergence rates are minimax rate-optimal for constant $k$, which will be shown in \prettyref{sec:gm-lb}. 
Note that \prettyref{thm:main-W1} is proved under the worst-case scenario where the centers can be arbitrarily close, \eg, components completely overlap. It is reasonable to expect a faster convergence rate when the components are better separated, and, in fact, a parametric rate in the best-case scenario
where the components are fully separated and weights are bounded away from zero. 
To capture the clustering structure of the mixture model, we introduce the following definition:
\begin{definition}
    \label{def:sep}
    The Gaussian mixture \prettyref{eq:model} has $k_0$ $(\gamma,\omega)$-separated clusters if there exists a partition $S_1,\dots,S_{k_0}$ of $[k]$ such that
    \begin{itemize}
        \item $|\mu_i - \mu_{i'}| \ge \gamma$ for any $i\in S_\ell$ and $i'\in S_{\ell'}$ such that $\ell\ne \ell'$;
        \item $\sum_{i\in S_\ell}w_i\ge \omega$ for each $\ell$.
    \end{itemize}
    In the absence of the minimal weight condition (i.e.~$\omega=0$), we say the Gaussian mixture has $k_0$ $\gamma$-separated clusters.    
\end{definition}

The next result shows that the DMM estimators attain the following adaptive rates: 
\begin{theorem}[Adaptive rate]
    \label{thm:main-W1-adaptive}
    Under the conditions of \prettyref{thm:main-W1}, suppose there are $k_0$ $(\gamma,\omega)$-separated clusters such that $\gamma\omega\ge C\epsilon$ for some absolute constant $C>2$, 
        where $\epsilon$ denotes the right-hand side of \prettyref{eq:known-variance-W1}.
    Then, with probability at least $1-\delta$,\footnote{Here $O_k(\cdot)$ denotes a constant factor that depends on $k$ only.}
    \begin{equation}
        \label{eq:known-variance-W1-adaptive}
        W_1(\nu,\hat \nu)\le  O_k\pth{M\gamma^{-\frac{2k_0-2}{2(k-k_0)+1}}\pth{\frac{n}{\log(k/\delta)}}^{-\frac{1}{4(k-k_0)+2}}}.
    \end{equation}
\end{theorem}
The result \prettyref{eq:known-variance-W1-adaptive} is also minimax rate-optimal when $k,k_0$ and $\gamma$ are constants, in view of the lower bounds in \cite{HK2015}. We also provide a simple proof in \prettyref{rmk:lb-adaptive-known} by extending the lower bound argument of \prettyref{prop:lb-known-sigma} in \prettyref{sec:gm-lb} (see \cite[Proposition 8]{GM-suppl} for a complete proof). 

Next we discuss the implication on density estimation (\emph{proper} learning), where the goal is to estimate the density function of the Gaussian mixture by another $k$-Gaussian mixture density. 
Given that the estimated mixing distribution $\hat \nu$ from \prettyref{thm:main-W1}, a natural density estimate is the convolution $\hat f =\hat \nu * N(0,\sigma^2)$. \prettyref{thm:main-density} below shows that the density estimate $\hat f$ is  $O(\frac{1}{n})$-close to the true density $f$ in $\chi^2$-divergence, which bounds other common distance measures such as the KL divergence, total variation, and Hellinger distance.
\begin{theorem}[Density estimation]
    \label{thm:main-density}
    Under the conditions of \prettyref{thm:main-W1}, denote the density of the underlying model by $f = \nu * N(0,\sigma^2)$. If $\sigma$ is given, then there exists an estimate $\hat f$ such that
    \[
        \chi^2(\hat f\|f) + \chi^2(f\|\hat f)\le O_k(\log(1/\delta)/n),
    \]
    with probability $1-\delta$.
\end{theorem}

So far we have been focusing on well-specified models.
In the case of misspecified models, the data need not be generated from a $k$-Gaussian mixture. In this case, the DMM procedure still reports a meaningful estimate that is close to the best $k$-Gaussian mixture fit of the unknown distribution. This is 
made precise by the next result of oracle inequality type. We state the version for the total variation. Analogous results hold for $\chi^2$-divergence, KL divergence, and Hellinger distance as well.

\begin{theorem}[Misspecified model]
    \label{thm:oracle}
    Assume that $X_1,\ldots,X_n$ is independently drawn from a density $f$ which is 1-subgaussian.
        Suppose there exists a $k$-component Gaussian location mixture $g$ with variance $\sigma^2$ such that $\TV(f,g) \le \epsilon$. Then, there exists an estimate $\hat f$ such that
    \[
        \TV(\hat f,f)\le O_k\pth{\epsilon\sqrt{\log(1/\epsilon)}+\sqrt{\log(1/\delta)/n}},
    \]
    with probability $1-\delta$.
\end{theorem}

Finally, we present a result for estimating mixtures of an arbitrarily large order, including continuous mixtures.
In this situation we apply the DMM method to produce a mixture of order $\min\{k,O(\frac{\log n}{\log\log n})\}$. The convergence rate is minimax optimal in view of the matching lower bound in \cite{WY18}.
\begin{theorem}[Higher-order mixture]
\label{thm:large-k}
Suppose $|\mu_i|\le M$ for $M\ge 1$ and $\sigma$ is a bounded constant, where $M,\sigma$ are given. Then there exists an estimate $\hat\nu$ such that, with probability at least $1-\delta$,
\[
W_1(\nu,\hat\nu)\le O\pth{M\pth{\frac{\log\log n}{\log n}+\sqrt{\frac{\log(1/\delta)}{n^{1-c}}}}},
\]
for some constant $c<1$.
\end{theorem}

To conclude this subsection, we prove \prettyref{thm:main-W1} using Propositions \ref{prop:stable1}. Theorems \ref{thm:main-W1-adaptive} -- \ref{thm:large-k} can be proved similarly using other moment comparison theorems in \prettyref{chap:framework}; see \cite{WY18}.
\begin{proof}[Proof of \prettyref{thm:main-W1}]
    By scaling it suffices consider $M=1$. We use \prettyref{algo:known} with Euclidean norm in \prettyref{eq:project}. Using the variance of $\tilde m$ in \prettyref{lmm:var-tildem} and Chebyshev inequality yield that, for each $r=1,\ldots,2k-1$, with probability $1-\frac{1}{8k}$, 
    \begin{equation}
        \label{eq:tilde-mr-rate}
        |\tilde m_r-m_r(\nu)| \le \sqrt{k/n}(c\sqrt{r})^r,
    \end{equation}
    for some absolute constant $c$. By the union bound, with probability $3/4$, \prettyref{eq:tilde-mr-rate} holds simultaneously for every $r=1,\dots,2k-1$, and thus
    \[
        \Norm{\tilde m-{\bf m}_{2k-1}(\nu)}_2\le \epsilon,\quad \epsilon\triangleq \frac{(\sqrt{ck})^{2k+1}}{\sqrt{n}}.
    \]
    Since ${\bf m}_{2k-1}(\nu)$ satisfies \prettyref{eq:moment-psd} and thus is a feasible solution for \prettyref{eq:project}, we have $\Norm{\tilde m-\hat m}_2\le \epsilon$. Note that $\hat m={\bf m}_{2k-1}(\hat \nu)$. Hence, by triangle inequality, we obtain the following statistical accuracy:
    \begin{equation}
        \label{eq:bfm-example}
        \Norm{{\bf m}_{2k-1}(\hat \nu)-{\bf m}_{2k-1}(\nu)}_2\le \epsilon,
    \end{equation}
    Applying \prettyref{prop:stable1} yields that, with probability $3/4$,
    \[
        W_1(\hat\nu,\nu)\le O\pth{k^{1.5}n^{-\frac{1}{4k-2}}}.
    \]
    The confidence $1-\delta$ in \prettyref{eq:known-variance-W1} can be obtained by the usual ``median trick'' (which has previously appeared in \prettyref{sec:supp-rate}): divide the sample into $T=\log\frac{2k}{\delta}$ batches, apply \prettyref{algo:known} to each batch of size $n/T$, and take $\tilde m_r$ to be the median of these estimates. Then Hoeffding's inequality and the union bound imply that, with probability $1-\delta$,
    \begin{equation}
        \label{eq:Hoeffding-example}
        |\tilde m_r-m_r(\nu)| \le \sqrt{\frac{\log(2k/\delta)}{n}}(c\sqrt{r})^r,\quad \forall~r=1,\ldots,2k-1,
    \end{equation}
    and the conclusion follows. 
\end{proof}


\section{Unknown variance (Lindsay's algorithm)}
\label{sec:unknown}
When the variance parameter $\sigma^2$ is unknown, an unbiased estimator for the moments of the mixing distribution no longer exists (see \cite[Lemma 31]{WY18}). 
A natural idea is to find a consistent estimate of the scale parameter $\sigma$ (e.g.~$\hat\sigma=\max_{i\in[n]} X_i / \sqrt{2\log n}$ which satisfies $|\sigma-\hat \sigma|=O_P(1/\sqrt{\log n})$), then plug into the DMM estimator in \prettyref{sec:known} to obtain a consistent estimate of the mixing distribution $\nu$; however, the resulting convergence rate is far from optimal. In fact, to achieve the optimal rate, 
 it is crucial to \emph{simultaneously} estimate both the location and the scale parameters. To this end, again we adopt a moment-based approach. The following result provides a guarantee for any joint estimate of both the mixing distribution and the variance parameter in terms of the moments accuracy.
\begin{proposition}[{\cite[Proposition 3]{WY18}}]
    \label{prop:accuracy2}
    Let
    \[
        \pi=\nu*N(0,\sigma^2),\quad \hat\pi=\hat\nu*N(0,\hat\sigma^2),
    \]
    where $\nu,\hat\nu$ are $k$-atomic distributions supported on $[-M,M]$, and $\sigma,\hat\sigma$ are bounded. If $|m_r(\pi)-m_r(\hat\pi)|\le \epsilon$ for $r=1,\dots,2k$, then
    \[
        |\sigma^2-\hat\sigma^2|\le O(M^2\epsilon^{\frac{1}{k}}),\quad W_1(\nu,\hat\nu)\le O(Mk^{1.5}\epsilon^{\frac{1}{2k}}).
    \]
\end{proposition}
To apply \prettyref{prop:accuracy2}, we can solve the method of moments equations, namely, find a $k$-atomic distribution $\hat\nu$ and $\hat\sigma^2$ such that
\begin{equation}
    \label{eq:MM-sigma}
    \Expect_n[X^r] = \Expect_{\hat\pi}[X^r], \qquad r=1,\ldots,2k
\end{equation}
where $\hat\pi=\hat\mu * N(0,\hat\sigma^2)$ is the fitted Gaussian mixture. Here both the number of equations and the number of variables are equal to $2k$. Suppose \prettyref{eq:MM-sigma} has a solution $(\hat\mu,\hat\sigma)$. Then applying \prettyref{prop:accuracy2}  with $\delta = O_k(\frac{1}{\sqrt{n}})$ achieves the rate $O_k(n^{-1/(4k)})$, which is minimax optimal (see \prettyref{sec:gm-lb}). 
In stark contrast to the known $\sigma$ case, where we have shown in \prettyref{sec:mmfail} that the vanilla method of moments equation can have no solution unless we denoise by projection to the moment space, here with one extra scale parameter $\sigma$, one can show that \prettyref{eq:MM-sigma} has a solution \emph{with probability one}!\footnote{
    It is possible that the equation \prettyref{eq:MM-sigma} has no solution, for instance, when $k=2, n=7$ and the empirical distribution is $\pi_7=\frac{1}{7}\delta_{-\sqrt{7}}+\frac{1}{7}\delta_{\sqrt{7}}+\frac{5}{7}\delta_0$. The first four empirical moments are ${\bf m}_4(\pi_7)=(0,2,0,14)$, which cannot be realized by any two-component Gaussian mixture \prettyref{eq:model}. Indeed, suppose $\hat \pi=w_1N(\mu_1,\sigma^2)+(1-w_1)N(\mu_2,\sigma^2)$ is a solution to \prettyref{eq:MM-sigma}. Eliminating variables leads to the contradiction that $2\mu_1^4+2=0$. Assuringly, as we will show later in \prettyref{lmm:U-exist}, such instances occur with probability zero.
} 
Furthermore, an efficient method of finding \emph{a} solution to \prettyref{eq:MM-sigma} is due to Lindsay \cite{Lindsay1989} and summarized in \prettyref{algo:Lindsay}. Indeed, the sample moments are computable in $O(kn)$ time, and the smallest non-negative root of the polynomial of degree $k(k+1)$ can be found in $O(k^2)$ time using Newton's method (see \cite{Atkinson89}). So overall Lindsay's estimator can be evaluated in $O(kn)$ time.
\begin{algorithm}[ht]
    \caption{Lindsay's estimator for normal mixtures with an unknown common variance}
    \label{algo:Lindsay}
    \begin{algorithmic}[1]
        \REQUIRE $ n $ observations $ X_1,\dots,X_n $.
        \ENSURE estimated mixing distribution $\hat\nu$, and estimated variance $\hat\sigma^2$.
        \FOR{$ r=1 $ \TO $ 2k $}
        \STATE{$ \hat{\gamma}_r=\frac{1}{n}\sum_i X_i^r $}
        \STATE{$ \hat{m}_r(\sigma)=r!\sum_{i=0}^{\Floor{r/2}}\frac{(-1/2)^i}{i!(r-2i)!}\hat{\gamma}_{r-2i}\sigma^{2i} $}
        \ENDFOR
        \STATE{Let $ \hat{d}_k(\sigma) $ be the determinant of the matrix $\{\hat{m}_{i+j}(\sigma)\}_{i,j=0}^k$.}\label{line:hat-dk}
        \STATE{Let $ \hat{\sigma} $ be the smallest positive root of $\hat{d}_k(\sigma)=0$.}\label{line:hat-sigma}
        \FOR{$ r=1 $ \TO $ 2k $}
        \STATE{$ \hat m_r=\hat{m}_r(\hat{\sigma}) $} \label{line:hat-mr}
        \ENDFOR
        \STATE{Let $\hat\nu$ be the outcome of the Gaussian quadrature (\prettyref{algo:quadrature}) with input $\hat m_1,\dots,\hat m_{2k-1}$}
        \STATE{Report $\hat\nu$ and $\hat\sigma^2$.}
    \end{algorithmic}
\end{algorithm}

In \cite{Lindsay1989} the consistency of this estimator was proved under the extra condition that $\hat \sigma$ (which is a random variable) as a root of $\hat d_k$ 
 has multiplicity one.
It is unclear whether this condition is guaranteed to hold.
We will show that, unconditionally, Lindsay's estimator is not only consistent, but in fact achieves the minimax optimal rate that will be shown in \prettyref{thm:main-W1-unknown}.
In this section we show that Lindsay's algorithm produces an estimator $\hat\sigma$ so that the corresponding the moment estimates lie in the moment space with probability one. In this sense, although no explicit projection is involved, the noisy estimates are \emph{implicitly} denoised.

We first describe the intuition of the choice of $\hat\sigma$ in Lindsay's algorithm, \ie, line \ref{line:hat-sigma} of \prettyref{algo:Lindsay}. Let $X\sim \nu*N(0,\sigma^2)$. For any $\sigma' \leq \sigma$, we have 
\[
    \Expect[h_j(X,\sigma')]=m_j(\nu*N(0,\sigma^2-\sigma'^2)).
\]
Let $d_k(\sigma')$ denote the determinant of the moment matrix $\{\Expect[h_{i+j}(X,\sigma')]\}_{i,j=0}^k$, which is an even polynomial in $\sigma'$ of degree $k(k+1)$. According to \prettyref{thm:supp-detM}, $d_k(\sigma')>0$ when $0\le \sigma'<\sigma$ and becomes zero at $\sigma'=\sigma$, and thus $\sigma$ is characterized by the smallest positive zero of $d_k$. In lines \ref{line:hat-dk} -- \ref{line:hat-sigma}, $d_k$ is estimated by $\hat d_k$ using the empirical moments, and $\sigma$ is estimated by the smallest positive zero of $\hat d_k$. We first note that $\hat d_k$ indeed has a positive zero:
\begin{lemma}
    \label{lmm:sigma-exist}
    Assume $ n>k $ and the mixture distribution has a density. Then, almost surely, $\hat d_k$ has a positive root within $(0,s]$, where $s^2 \triangleq \frac{1}{n}\sum_{i=1}^n(X_i-\Expect_n[X])^2$ denotes the sample variance. 
\end{lemma}

The next result shows that, with the above choice of $\hat\sigma$, the moment estimates $\hat{m}_j=\Expect_n[h_j(X,\hat\sigma)]$ for $j=1,\dots,2k$ given in line \ref{line:hat-mr} are implicitly denoised and lie in the moment space with probability one. Thus \prettyref{eq:MM-sigma} has a solution, and the estimated mixing distribution $\hat\nu$ can be found by the Gaussian quadrature. This result was previously shown in \cite{Lindsay1989} assuming that $\hat \sigma$ is of multiplicity one.
\begin{lemma}
    \label{lmm:U-exist}
    Assume $ n\ge 2k-1 $ and the mixture distribution has a density. Then, almost surely, there exists a $k$-atomic distribution $\hat\nu$ such that $m_j(\hat\nu)=\hat{m}_j$ for $j\le 2k$, where $\hat{m}_j$ is from \prettyref{algo:Lindsay}.
\end{lemma}

For the proofs of Lemmas~\ref{lmm:sigma-exist} and \ref{lmm:U-exist}, see \cite{WY18}.
With the above analysis, we now present the statistical guarantee for Lindsay's algorithm including the minimax and the adaptive rates (cf. Theorems~\ref{thm:main-W1} and \ref{thm:main-W1-adaptive}):
\begin{theorem}[Optimal rate]
    \label{thm:main-W1-unknown}
    Suppose that $|\mu_i|\le M$ for $M\ge 1$ and $\sigma$ is bounded 
    from above
    by a constant, and both $k$ and $M$ are given.
    Then, \prettyref{algo:Lindsay} produces $(\hat \nu,\hat \sigma)$ such that, with probability at least $1-\delta$,
    \begin{equation}
        \label{eq:unknown-variance-W1}
        W_1(\nu,\hat \nu)\le O\pth{Mk^{2}\pth{\frac{n}{\log(1/\delta)}}^{-\frac{1}{4k}}},
    \end{equation}
    and 
    \begin{equation}
        \label{eq:unknown-variance-sigma}        
        |\sigma^2-\hat \sigma^2| \leq  O\pth{M^2k\pth{\frac{n}{\log(1/\delta)}}^{-\frac{1}{2k}}}.
    \end{equation}     
\end{theorem}

\begin{theorem}[Adaptive rate]
    \label{thm:main-W1-adaptive-unknown}
    Under the conditions of \prettyref{thm:main-W1-unknown}, suppose there are $k_0$ $(\gamma,\omega)$-separated clusters such that $\gamma\omega\ge C\epsilon$ for some absolute constant $C>2$, 
    where $\epsilon$ denotes the right-hand side of \prettyref{eq:unknown-variance-W1} when $\sigma$.
    Then, with probability at least $1-\delta$,\footnote{Note that the estimation rate for the mean part $\nu$ is the square root of the rate for estimating the variance parameter $\sigma^2$. Intuitively, this phenomenon is due to the infinite divisibility of the Gaussian distribution: for the location mixture model $\nu * N(0,\sigma^2)$ with $\nu\sim N(0,\epsilon^2)$ and $\sigma^2=1$ has the same distribution as that of $\nu\sim \delta_0$ and $\sigma^2=1+\epsilon^2$.}
    \begin{equation}
        \label{eq:unknown-variance-W1-adaptive}
        \sqrt{|\sigma^2-\hat \sigma^2|},~W_1(\nu,\hat \nu)\le  O_k\pth{M\gamma^{-\frac{k_0-1}{k-k_0+1}}\pth{\frac{n}{\log(k/\delta)}}^{-\frac{1}{4(k-k_0+1)}}}.
    \end{equation}
\end{theorem}

The proofs are similar to Theorems~\ref{thm:main-W1} and \ref{thm:main-W1-adaptive}, where we will apply the moment comparisons with one discrete distribution in Propositions~\ref{prop:stable2} and \ref{prop:stable2-separation} instead of Propositions~\ref{prop:stable1} and \ref{prop:stable1-separation}.
Next we prove \prettyref{thm:main-W1-unknown}; \prettyref{thm:main-W1-adaptive-unknown} can be proved similarly using \prettyref{prop:stable2-separation} (see \cite{WY18} for details).
Technical lemmas are deferred to \prettyref{sec:aux-lemma-unknown}.

\begin{proof}[Proof of \prettyref{thm:main-W1-unknown}]
    It suffices to consider $M=1$. Let $\hat\pi=\hat\nu*N(0,\hat\sigma^2)$ and $\pi=\nu*N(0,\sigma^2)$ denote the estimated mixture distribution and the ground truth, respectively. Let $\hat m_r=\Expect_n[X^r]$ and $m_r=m_r(\pi)$. The variance of $\hat m_r$ is upper bounded by
    \[
        \var[\hat m_r]=\frac{1}{n}\var[X_1^r]\le \frac{1}{n}\Expect[X^{2r}]
        \le \frac{(\sqrt{cr})^{2r}}{n},
    \]
    for some absolute constant $c$. Using Chebyshev's inequality, for each $r=1,\dots,2k$, with probability $1-\frac{1}{8k}$, we have, 
    \begin{equation}
        \label{eq:tildemr-rate2}    
        |\hat m_r-m_r|\le (\sqrt{cr})^r\sqrt{k/n}.
    \end{equation}
    By the union bound, with probability 3/4, the above holds holds simultaneously for every $r=1,\dots,2k$. It follows from \prettyref{lmm:sigma-exist} and \ref{lmm:U-exist} that \prettyref{eq:MM-sigma} holds with probability one. Therefore,
    \[
        |m_r(\hat\pi)-m_r(\pi)|\le (\sqrt{cr})^r\sqrt{k/n},\quad r=1,\dots,2k.
    \]
    for some absolute constant $c$. In the following, the error of variance estimate is denoted by $\tau^2=|\sigma^2-\hat\sigma^2|$.

    \begin{itemize}
        \item If $\sigma\le \hat\sigma$, let $\nu'=\hat\nu*N(0,\tau^2)$. Using $\Expect_\pi[h_r(X,\sigma)]=m_r(\nu)$ and $\Expect_{\hat\pi}[h_r(X,\sigma)]=m_r(\nu')$, where $h_r$ is the Hermite polynomial \prettyref{eq:Hermite}, we obtain that (see \prettyref{lmm:Hermite-moments} in \prettyref{sec:aux-lemma-unknown})
        \begin{equation}
            \label{eq:mr-deconv-example}
            |m_r(\nu')-m_r(\nu)|\le (\sqrt{c'k})^{2k}\sqrt{k/n},\quad r=1,\dots,2k,
        \end{equation}
        for an absolute constant $c'$. Applying \prettyref{prop:accuracy2} yields that,
        \[
            |\sigma^2-\hat\sigma^2|\le O(kn^{-\frac{1}{2k}}),\quad W_1(\nu,\hat\nu)\le O(k^2n^{-\frac{1}{4k}}).
        \]
        \item If $\sigma\ge \hat\sigma$, let $\nu'=\nu*N(0,\tau^2)$. Similar to \prettyref{eq:mr-deconv-example}, we have 
        \[
            |m_r(\hat\nu)-m_r(\nu')|\le (\sqrt{c'k})^{2k}\sqrt{k/n}\triangleq \epsilon,\quad r=1,\dots,2k.
        \]
        To apply \prettyref{prop:accuracy2}, we also need to ensure that $\hat\nu$ has a bounded support, which is not obvious. To circumvent this issue, we apply a truncation argument thanks to the following tail probability bound for $\hat\nu$ (see \prettyref{lmm:tail-hatU} in \prettyref{sec:aux-lemma-unknown}):
        \begin{equation}
            \label{eq:hatnu-tail}
            \Prob[|\hat U|\ge \sqrt{c_0k}]\le \epsilon(\sqrt{c_1k}/t)^{2k},\quad \hat U\sim \hat \nu,
        \end{equation}
        for absolute constants $c$ and $c'$. To this end, consider $\tilde U=\hat U\indc{|\hat U|\le \sqrt{c_0k}}$ with distribution $\tilde\nu$. Note that $\tilde U$ is $k$-atomic supported on the interval $[-\sqrt{c_0k},\sqrt{c_0k}]$, we have $W_1(\nu,\hat\nu)\le \epsilon e^{O(k)}$ and $|m_r(\tilde\nu)-m_r(\hat\nu)|\le k\epsilon(c_1k)^k$ for $r=1,\dots,2k$. Using the triangle inequality yields that 
        \[
            |m_r(\tilde \nu)-m_r(\nu')|\le \epsilon+k\epsilon(c_1k)^k.
        \]
        Now we apply \prettyref{prop:accuracy2} with $\tilde\nu$ and $\nu*N(0,\tau^2)$ where both $\tilde \nu$ and $\nu$ are $k$-atomic supported on $[-\sqrt{c_0k},\sqrt{c_0k}]$. In the case $\tilde\nu$ is discrete, the dependence on $k$ in \prettyref{prop:accuracy2} can be improved 
        and we obtain that
        \[
            |\sigma^2-\hat\sigma^2|\le O(kn^{-\frac{1}{2k}}),\quad W_1(\nu,\tilde\nu)\le O(k^2n^{-\frac{1}{4k}}).
        \]
        Using $k\le O(\frac{\log n}{\log\log n})$, we also obtain $W_1(\nu,\hat\nu)\le O(k^2n^{-\frac{1}{2k}})$ by the triangle inequality.
    \end{itemize}
    To obtain a confidence $1-\delta$ in \prettyref{eq:unknown-variance-W1} and \prettyref{eq:unknown-variance-sigma}, we can replace the empirical moments $\hat m_r$ by the median of $T=\log\frac{1}{\delta}$ independent estimates similar to \prettyref{eq:Hoeffding-example}. 
\end{proof}

\subsection{Auxiliary lemmas}
\label{sec:aux-lemma-unknown}

In the subsection we present auxiliary lemmas used to prove \prettyref{thm:main-W1-unknown}.
\begin{lemma}
    \label{lmm:Hermite-moments}
    If $|\Expect[X^\ell]-\Expect[X'^\ell]|\le (C\sqrt{\ell})^\ell\epsilon$ for $\ell=1,\dots,r$, then, for $h_r$ in \prettyref{eq:Hermite},
    \[
        |\Expect[h_r(X,\sigma)]-\Expect[h_r(X',\sigma)]|\le \epsilon\pth{(2\sigma \sqrt{r/e})^r+(2C\sqrt{r})^r}.
    \]
\end{lemma}
\begin{proof}
    Note that $|\Expect[X^\ell]-\Expect[X'^\ell]|\le \Expect|C\sqrt{e}Z'|^r\epsilon$ by \prettyref{lmm:normal-moments} below, where $Z'\sim N(0,1)$. Then, 
    \begin{gather*}
        |\Expect[h_r(X,\sigma)]-\Expect[h_r(X',\sigma)]|
        \le \sum_{i=0}^{\Floor{r/2}}\frac{r!\sigma^{2i}}{i!(r-2i)!2^i}\Expect[|C\sqrt{e}Z'|^r]\epsilon\\
        =\epsilon\cdot\Expect[(\sigma Z+|C\sqrt{e}Z'|)^r],
    \end{gather*}
    where $Z\sim N(0,1)$ independent of $Z'$. Applying $(a+b)^r\le 2^{r-1}(|a|^r+|b|^r)$ and \prettyref{lmm:normal-moments} completes the proof.
\end{proof}

\begin{lemma}
    \label{lmm:normal-moments}
    \[
        (p/e)^{p/2}\le \Expect|Z|^p\le \sqrt{2}(p/e)^{p/2},\quad p\ge 0.
    \]
\end{lemma}
\begin{proof}
    Note that
    \[
        \frac{\Expect|Z|^p}{(p/e)^{p/2}}=\frac{2^{p/2}\Gamma(\frac{p+1}{2})}{\sqrt{\pi}(p/e)^{p/2}}\triangleq f(p),\quad \forall~p\ge 0.
    \]
    Since $f(0)=1$ and $f(\infty)=\sqrt{2}$, it suffices to show that $f$ is increasing in $[0,\infty)$.
    Equivalently, $\frac{x}{2}\log\frac{2e}{x}+\log\Gamma(\frac{x+1}{2})$ is increasing, which is equivalent to $\psi(\frac{x+1}{2})\ge \log\frac{x}{2}$ by the derivative, where $\psi(x)\triangleq \frac{\diff}{\diff x}\log\Gamma(x)$.
    The last inequality holds for any $x>0$ (see, \eg, \cite[(3)]{DS2016}).
\end{proof}

\begin{lemma}
    \label{lmm:tail-hatU}
    Let
    \[
        \pi=\nu*N(0,\tau^2),\quad \hat\pi= \hat\nu,
    \]
    where $\nu$ and $\nu$ are both $k$-atomic, $\nu$ is supported on $[-1,1]$, and $\tau\le 1$. If $|m_i(\pi)-m_i(\hat{\pi})|\le \epsilon$ for $i\le 2k$, then, for any $t\ge \sqrt{18k}$,
    \begin{equation*}
        \Prob[|\hat U|\ge t]\le \frac{2^{2k+1}\epsilon}{(\frac{t}{\sqrt{18k}}-1)^{2k}},\quad \hat U\sim\hat \nu.
    \end{equation*}
\end{lemma}
\begin{proof}
    Let $g$ be the $(k+1)$-point Gaussian quadrature of the standard normal distribution. Furthermore, $g$ is supported on $[-\sqrt{4k+6},\sqrt{4k+6}]$ for some absolute constant $c$ (see the bound on the zeros of Hermite polynomials in \cite[p.~129]{orthogonal.poly}). Let $G\sim g$, $U\sim \nu$, and $\hat U\sim\hat\nu$. Denote the maximum absolute value of $U+\tau G$ by $M$ which is at most $1+\sqrt{4k+6}\le \sqrt{18k}$ for $k\ge 1$. Applying \prettyref{lmm:tail} below with the distributions of $\frac{U+\tau G}{\sqrt{18k}}$ and $\frac{\hat U}{\sqrt{18k}}$ yields the conclusion.
\end{proof}

\begin{lemma}
    \label{lmm:tail}
    Let $\epsilon=\max_{i\in [2k]}|m_i(\nu)-m_i(\nu')|$. If either $\nu$ or $\nu'$ is $k$-atomic, and $\nu$ is supported on $[-1,1]$, then, for any $t>1$,
    \[
        \Prob[|Y|\ge t]\le 2^{2k+1}\epsilon/(t-1)^{2k},\quad Y\sim\nu'.
    \]
\end{lemma}
\begin{proof}
    We only show the upper tail bound $ \Prob[Y\ge t] $. The lower tail bound of $ Y $ is equal to the upper tail bound of $ -Y $. 
    \begin{itemize}
        \item Suppose $\nu$ is $k$-atomic supported on $\{x_1,\ldots,x_k\}$. Consider a polynomial $ P(x)=\prod_i (x-x_i)^2 $ of degree $2k$ that is almost surely zero under $\nu$. Since every $|x_i|\le 1$, similar to \prettyref{eq:polydiff}, we obtain that 
        \[
            \Expect_{\nu'}[P]=|\Expect_\nu[P]-\Expect_{\nu'}[P]|\le 2^{2k} \epsilon.
        \]
        Using Markov inequality, for any $t>1$, we have
        \[
            \Prob[Y\ge t]\le \Prob[P(Y)\ge P(t)] \le \frac{\Expect[P(Y)]}{P(t)} \le \frac{2^{2k} \epsilon}{(t-1)^{2k}}.
        \]
        \item Suppose $\nu'$ is $k$-atomic supported on $\{x_1,\ldots,x_k\}$. If those values are all within $[-1,1]$, then we are done. If there are at most $ k-1 $ values, denoted by $ \{x_1,\dots,x_{k-1}\} $, are within $ [-1,1] $, then we consider a polynomial $ P(x)=(x^2-1)\prod_i (x-x_i)^2 $ of degree $2k$ that is almost surely non-positive under $\nu$. Similar to \prettyref{eq:polydiff}, we obtain that
        \[
            \Expect_{\nu'}[P]\le \Expect_{\nu'}[P]-\Expect_{\nu}[P]\le 2^{2k} \epsilon.
        \]
        Since $P\ge 0$ almost surely under $\nu'$, the conclusion follows follows analogously using Markov inequality. \qedhere
    \end{itemize}
\end{proof}

\section{Lower bounds}
\label{sec:gm-lb}

In this section we discuss minimax lower bounds for estimating Gaussian location mixture models which certify the optimality of our estimators. We will apply Le Cam's two-point method in \prettyref{thm:lb-two-pts}, which entails finding two Gaussian mixtures that are statistically close but have different parameters. 

As mentioned in \prettyref{sec:d-mm}, moment matching is useful to ensure a small statistical distance between two mixture models. In the context of Gaussian mixtures, $\nu*N(0,1)$ and $\nu'*N(0,1)$ are statistically close if $\nu$ and $\nu'$ share many moments.
This phenomenon is demonstrated in \prettyref{fig:tv-latent-mixture} and quantified in \prettyref{thm:chi2-gm}. The best lower bound follows from two distinct mixing distributions $\nu$ and $\nu'$ that have as many identical moments as possible ($2k-2$ when both are $k$-atomic and $2k-1$ when only one of them is $k$-atomic; see \prettyref{lmm:identify} and \prettyref{rmk:stable}). 
Next we provide the precise minimax lower bounds for the case of known and unknown variance separately.

\paragraph{Known variance. } For simplicity, we shall assume that $\sigma^2=1$. First, we define the space of all $k$ Gaussian location mixtures as 
\[
    \calP_k=\{\nu*N(0,1): \nu \text{ is $k$-atomic supported on }[-1,1]\},
\]
and consider the worst-case risk over all mixtures in $\calP_k$. From the identifiability of discrete distributions in \prettyref{lmm:identify}, two different $k$-atomic distributions can match up to $2k-2$ moments. Therefore, using \prettyref{thm:chi2-gm}, the best minimax lower bound using \prettyref{thm:lb-two-pts} is obtained from the optimal pair of distributions for the following:
\begin{equation}
    \label{eq:optimal-lb}
    \begin{aligned}
        \max & ~W_1(\nu,\nu')\\
        \mathrm{s.t.}&~{\bf m}_{2k-2}(\nu)={\bf m}_{2k-2}(\nu'),\\
        &~ \nu,\nu' \text{ are $k$-atomic on }[-\epsilon,\epsilon].
    \end{aligned}
\end{equation}
The value of the above optimization problem is $\Omega(\epsilon/k)$ by the next lemma:
\begin{lemma}
    \label{lmm:max-W1}
    \[
        \sup\{W_1(\nu,\nu'):{\bf m}_\ell(\nu) = {\bf m}_\ell(\nu'),~\nu,\nu'\text{ on }[-\epsilon,\epsilon]\}=\Theta(\epsilon/\ell).
    \]
    Furthermore, the supremum is $\frac{\epsilon(\pi-o(1))}{\ell}$ as $\ell\rightarrow\infty$, achieved by two distributions whose support sizes differ by at most one and sum up to $\ell+2$. 
\end{lemma}
\begin{proof}
    It suffices to prove for $\epsilon=1$. Using the dual characterization of the $W_1$ distance in \prettyref{sec:wass}, the supremum is equal to
    \[
        \sup_{f:1-\mathrm{Lipschitz}}\sup\sth{\Expect_\nu f-\Expect_{\nu'}f:{\bf m}_\ell(\nu) = {\bf m}_\ell(\nu'),~\nu,\nu'\text{ on }[-1,1]}.
    \]
    Using the duality between moment matching and best polynomial approximation (see \prettyref{sec:dual-best}), the optimal value is further equal to
    \begin{equation}
            2\sup_{f:1-\mathrm{Lipschitz}}E_\ell(f,[-1,1]).
        \label{eq:Elip}
        \end{equation}
    The supremum in \prettyref{eq:Elip} is the worst-case uniform approximation error of 1-Lipschitz functions, a well-studied quantity in the approximation theory (see, \eg, \cite[Section 4.1]{Jorge2011}). It is attained by a certain 1-Lipschitz function $f^*$ and the value satisfies $\frac{\pi-o(1)}{2\ell}$ as the degree $\ell\diverge$. A pair of optimal distributions are supported on the maxima and the minima of $P^*-f^*$, respectively, where $P^*$ is the best degree-$\ell$ polynomial approximation of $f^*$. The numbers of maxima and minima differ by at most one by Chebyshev's alternating theorem (see \prettyref{thm:cheby}).
\end{proof}

Using $\epsilon=\sqrt{k}n^{-\frac{1}{4k-2}}$ in \eqref{eq:optimal-lb}, we obtain the following minimax lower bound:
\begin{proposition}
    \label{prop:lb-known-sigma}
    \[
        \inf_{\hat \nu}\sup_{P\in\calP_k}\Expect_P W_1(\nu,\hat\nu) \ge \Omega\pth{\frac{1}{\sqrt{k}}n^{-\frac{1}{4k-2}}}
    \]
    where $\hat\nu$ is an estimator measurable with respect to $X_1,\ldots,X_n\iiddistr P=\nu*N(0,1)$.
\end{proposition}

\begin{remark}
    \label{rmk:lb-adaptive-known}
    The construction of a pair of mixing distributions in \prettyref{eq:optimal-lb}
		can be easily extended to prove the optimality of the adaptive rate \prettyref{eq:known-variance-W1-adaptive} in \prettyref{thm:main-W1-adaptive}, where the mixture satisfies further separation conditions in the sense of \prettyref{def:sep}. In this case, the main difficulty is to estimate parameters in the biggest cluster. When there are $k_0$ $\gamma$-separated clusters, the biggest cluster is of order at most $k'=k-k_0+1$. Similar to \prettyref{eq:optimal-lb}, let $\tilde\nu$ and $\tilde\nu'$ be two $k'$-atomic distributions on $[-\epsilon,\epsilon]$. Consider the following mixing distributions
    \[
        \nu=\frac{k_0-1}{k_0}\nu_0+\frac{1}{k_0}\tilde\nu,\quad \nu'=\frac{k_0-1}{k_0}\nu_0+\frac{1}{k_0}\tilde\nu',
    \]
    where $\nu_0$ is the uniform distribution over $\{\pm 2\gamma,\pm 3\gamma,\ldots\}$ of cardinality $k_0-1$. Then both mixture models have $k_0$ $(\gamma,\frac{1}{k_0})$-separated clusters. Thus the minimax lower bound $\Omega(\frac{1}{\sqrt{k'}}n^{-\frac{1}{4k'-2}})$ analogously follows from \prettyref{thm:lb-two-pts}.
\end{remark}

\paragraph{Unknown variance. } In this case the collection of mixture models is defined as 
\[
    \calP_k'=\{\nu*N(0,\sigma^2): \nu \text{ is $k$-atomic supported on }[-1,1],~\sigma\le 1\}.
\]
In \prettyref{thm:chi2-gm}, mixing distributions are not restricted to be $k$-atomic but can be Gaussian location mixtures themselves, thanks to the infinite divisibility of the Gaussian distributions, \eg, $N(0,\epsilon^2)*N(0,0.5)=N(0,0.5+\epsilon^2)$. Let $g_k$ be the $k$-point Gaussian quadrature of $N(0,\epsilon^2)$. Then $g_k$ has the same first $2k-1$ moments as $N(0,\epsilon^2)$, and $g_k*N(0,0.5)$ is an order-$k$ Gaussian mixture. Applying \prettyref{eq:chi2-subg} yields that
\[
    \chi^2(g_k*N(0,0.5)\|N(0,0.5+\epsilon^2))\le O(\epsilon^{4k}).
\]
Using $W_1(g_k,\delta_0)\ge \Omega(\epsilon/\sqrt{k})$ (see \prettyref{lmm:quadrature-l1}), and choosing $\epsilon=n^{-\frac{1}{4k}}$, we obtain the following minimax lower bound:
\begin{proposition}
    \label{prop:lb-unknown-sigma}
    For $k\ge 2$,
    \begin{gather*}
        \inf_{\hat \nu}\sup_{P\in\calP_k}\Expect_P W_1(\nu,\hat\nu) \ge \Omega\pth{\frac{1}{\sqrt{k}}n^{-\frac{1}{4k}}},\\
        \inf_{\hat \nu}\sup_{P\in\calP_k}\Expect_P |\sigma^2-\hat\sigma^2| \ge \Omega\pth{n^{-\frac{1}{2k}}},
    \end{gather*}
    where the infimum is taken over estimators $\hat\nu, \hat\sigma^2$ measurable with respect to $X_1,\ldots,X_n\iiddistr P=\nu*N(0,\sigma^2)$.
\end{proposition}

\section{Discussion and open problems}
\label{sec:discuss}


\paragraph{Gaussian location-scale mixtures}
We have been focusing on the Gaussian location mixture model \prettyref{eq:model}, where all components share the same (possibly unknown) variance.
One immediate extension is the Gaussian location-scale mixture model with heteroscedastic components:
\begin{equation}
\sum_{i=1}^k w_i N(\mu_i, \sigma_i^2).
\label{eq:model-locationscale}
\end{equation}
Parameter estimation for this model turns out to be significantly more difficult than the location mixture model, in particular
\begin{itemize}
        \item The likelihood function is unbounded. In fact, it is well-known that the maximum likelihood estimator is ill-defined \cite[p.~905]{KW56}.
        For instance, 
consider $k=2$, for any sample size $n$, we have
\[
\sup_{p_1,p_2,\theta_1,\theta_2,\sigma} \prod_{i=1}^n \qth{\frac{p_1 }{\sigma_1}\varphi\pth{\frac{X_i-\theta_1}{\sigma_1} }+
\frac{p_2 }{\sigma_2}\varphi\pth{\frac{X_i-\theta_2}{\sigma_2} }} = \infty,
\]
achieved by, e.g., $\theta_1=X_1, p_1=1/2$, $\sigma_2=1$, and $\sigma_1\to0$.
        
        \item 
In this model, the identifiability result based on moments is not completely settled and we do not have a counterpart of \prettyref{lmm:identify}. Note that the model \prettyref{eq:model-locationscale} comprises $3k-1$ free parameters ($k$ means, $k$ variances, and $k$ weights normalized to one), so it is expected to be identified through its first $3k-1$ moments.
However, the intuition of equating the number of parameters and the number of equations is already known to be wrong as pointed out by Pearson \cite{Pearson1894}, who showed that for $k=2$, five moments are insufficient and six moments are enough.
The recent result \cite{ARS2016} showed that, if the parameters are in general positions, then $3k-1$ moments can identify the Gaussian mixture distribution up to finitely many solutions (known as algebraic identifiability). Whether $3k$ moments can uniquely identify the model (known as rational identifiability) in general positions remains open, except for $k=2$.
In the worst case, we need at least $4k-2$ moments for identifiability since for scale-only Gaussian mixtures all odd moments are zero. 

\end{itemize}

Besides the issue of identifiability, the optimal estimation rate under the Gaussian location-scale mixture model is resolved only in special cases. 
The sharp rate is only known in the case of two components to be $\Theta(n^{-1/12})$ for estimating means and $\Theta(n^{-1/6})$ for estimating variances \cite{HP15}, achieved by a robust variation of Pearson's method of moment equations \cite{Pearson1894}; however, this approach is difficult to generalize to higher order mixtures.
For $k$ components, the optimal rate is known to be $n^{-\Theta(1/k)}$ \cite{MV2010,KMV2010}, achieved by an exhaustive search over the discretized parameter space so that the fitted moments is close to the empirical moments. In addition to being computationally expensive, this method achieves the estimation accuracy $n^{-C/k}$ for some constant $C$, which is suboptimal. 
In addition, the above results all aim to recover parameters of all components (up to a global permutation), which necessarily requires many assumptions such as lower bounds on mixing weights and separation between components; recovering the mixing distribution with respect to, say, Wasserstein distance, remains open.

\paragraph{Multiple dimensions}
This chapter focuses on learning Gaussian mixtures in one dimension. The multivariate version of this problem has been studied in the context of clustering, or classification, which requires separation between components \cite{Dasgupta1999,VW2004}. 
One commonly used approach is dimensionality reduction: projecting data onto some lower dimensional subspace, clustering samples in that subspace, and mapping back to the original space. Common choices of the subspace include random subspaces and subspaces obtained from the singular value decomposition.
The approach using random subspace is analyzed in \cite{Dasgupta1999,AK2001}, and requires a pairwise separation polynomial in the dimensions;
the subspace from singular value decomposition is analyzed in \cite{VW2004,KSV2005,AM2005,BV2008}, and requires a pairwise separation that grows polynomially in the number of components.
Tensor decomposition for spherical Gaussian mixtures has been studied in \cite{HK2013}, which requires the stronger assumption that
the means are linear independent and is inapplicable in lower dimensions, say, two or three dimensions.

Combined with dimension reduction, the DMM methodology is extended to multivariate Gaussian location models in \cite{DWYZ20} achieving the optimal Wasserstein rate in high dimensions. The main idea is as follows: First estimate the subspace spanned by the centers via principal component analysis, then project the sample to the learned subspace and estimate the mixing distribution in the lower-dimensional subspace through a sliced version of the DMM. The analysis relies on the extension of the moment comparison result in \prettyref{prop:compare-density-d-dim}.
Although the procedure runs in time that is polynomial in the dimension and the sample size, it involves some exhaustive search and is still far from being practical.



\paragraph{General finite mixtures}
Though we focused on Gaussian location mixture models, the moments comparison theorems in \prettyref{chap:framework} are independent of properties of Gaussian. As long as moments of the mixing distribution are estimated accurately, similar theory and algorithms can be obtained. Unbiased estimate of moments exists in many useful mixture models, including exponential mixtures \cite{Jewell1982}, Poisson mixtures \cite{KX2005}, and more generally the quadratic variance exponential family (QVEF) whose variance is at most a quadratic function of the mean \cite[(8.8)]{morris1982natural}.


As a closely related topic of this chapter, we discuss the Gaussian scale mixture model \eqref{eq:Gauss-scale} in details, which has been extensively studied in the statistics literature \cite{AM1974} and is widely used in image and video processing \cite{WS2000,PSWS2003}. 
In a Gaussian scale mixture, each observation can be represented as $X=\sqrt{V}Z$, where $V\sim \nu=\sum_{i=1}^kw_i\delta_{\sigma_i^2}$ is a $k$-atomic mixing distribution and $Z$ is standard normal independent of $V$. In this model, samples from different components significantly overlap, so clustering-based algorithms will fail. Nevertheless, moments of $\nu$ can be easily estimated, for instance, using $\Expect_n[X^{2r}]/\Expect[Z^{2r}]$ for $m_r(\nu)$ with accuracy $O_r(1/\sqrt{n})$. Applying a similar algorithm to DMM in \prettyref{sec:known}, we obtain an estimate $\hat \nu$ such that
\[
    W_1(\nu,\hat \nu)\le O_k(n^{-\frac{1}{4k-2}}),
\]
with high probability.

Moreover, using similar recipe to \prettyref{sec:gm-lb}, a minimax lower bound can be established. Analogous to \prettyref{eq:optimal-lb}, let $\nu$ and $\nu'$ be a pair of $k$-atomic distributions supported on $[0,\epsilon]$ such that they match the first $2k-2$ moments, and let
\begin{equation}
\pi=\int N(0,\sigma^2)\diff \nu(\sigma^2),\quad 
    \pi'=\int N(0,\sigma^2)\diff \nu'(\sigma^2),
\label{eq:gaussian-scale-lb}
\end{equation}
which match their first $4k-3$ moments and are $\sqrt{\epsilon}$-subgaussian. Applying \prettyref{thm:chi2-gm} with $\pi*N(0,0.5)$, $\pi'*N(0,0.5)$, and $\epsilon=O_k(n^{-\frac{1}{4k-2}})$ yields a minimax lower bound 
\[
    \inf_{\hat\nu}\sup_{P\in\calG_k}\Expect_P W_1(\nu,\hat \nu)\ge \Omega_k\pth{n^{-\frac{1}{4k-2}}},
\]
where the estimator $\hat \nu$  is measurable with respect to $X_1,\dots,X_n\sim P$, and the space of $k$ Gaussian scale mixtures is defined as
\[
    \calG_k=\sth{\int N(0,\sigma^2)\diff \nu(\sigma^2): \nu \text{ is $k$-atomic supported on }[0,1]}.
\]

\section{Bibliographic notes}
\label{sec:gm-bib}
There exist a vast literature on mixture models, in particular Gaussian mixtures, and the method of moments. For a comprehensive review we refer the reader to the monographs \cite{Lindsay1995,Fruhwirth2006}. 

\paragraph{Likelihood-based methods.} Maximum likelihood estimation (MLE) is one of the most useful method for parameter estimation. Under strong separation assumptions, MLE is consistent and asymptotically normal \cite{RW1984}; however, those assumptions are difficult to verify, and it is computationally hard to obtain the global maximizer due to the non-convexity of the likelihood function in the location parameters.

Expectation-Maximization (EM) \cite{DLR1977} is an iterative algorithm that aims to approximate the MLE. It has been widely applied in Gaussian mixture models \cite{RW1984,XJ1996} and more recently in high-dimensional settings \cite{balakrishnan2017statistical}. In general, this method is only guaranteed to converge to a local maximizer of the likelihood function rather than the global MLE. In practice we need to employ heuristic choices of the initialization \cite{KX2003} and stopping criteria \cite{SMA2000}, as well as possibly data augmentation techniques \cite{meng1997algorithm,PL2001}. Furthermore, its slow convergence rate is widely observed in practice \cite{RW1984,KX2003}. Global convergence of the EM algorithm is recently analyzed by \cite{XHM16,DTZ17} but only in the special case of two equally weighted components. Additionally, the EM algorithm accesses the entire data set in each iteration, which is particularly expensive for large sample size and high dimensions.

Lastly, we mention the nonparametric maximum likelihood estimation (NPMLE) for mixture models proposed by \cite{KW56}, where the maximization is taken over all mixing distributions which need not be $k$-atomic. This is an infinite-dimensional convex optimization problem, which has been studied in \cite{Laird1978,Lindsay1981,Lindsay1995} and more recently in \cite{koenker2014convex} for its computational aspect. One of the drawbacks of NPMLE is its lack of interpretability since the solution is a discrete distribution with at most $n$ atoms cf.~\cite[Theorem 2]{koenker2014convex}. Furthermore, existing statistical guarantees are suboptimal compared to \prettyref{thm:main-density}. For example, the squared Hellinger rate of $O(\frac{\log^2 n}{n})$ is shown in \cite{Zhang_2009} for compactly supported mixing distributions in one dimension and extended to multiple dimensions in \cite{Saha_2017}. While the logarithmic factors are unavoidable for continuous mixtures \cite{kim2014minimax}, it is not the case for finite mixtures.

\paragraph{Moment-based methods.} The simplest moment-based method is the method of moments (MM) introduced by Pearson \cite{Pearson1894}. The failure of the vanilla MM described in \prettyref{sec:mmfail} has motivated various modifications including, notably, 
the \emph{Generalized Method of Moments} (GMM) introduced by Hansen \cite{Hansen1982}. GMM is a widely used methodology for analyzing economic and financial data (cf.~\cite{Hall2005} for a thorough review). Instead of exactly solving the MM equations, GMM aims to minimize the sum of squared differences between the sample moments and the fitted moments. Despite its nice asymptotic properties \cite{Hansen1982}, GMM involves a non-convex optimization problem which is computationally challenging to solve. 
In practice, heuristics such as gradient descent are used \cite{Chausse2010} which converge slowly and lack theoretical guarantees.

For Gaussian mixture models (and more generally finite mixture models), our results can be viewed as a solver for GMM which is provably exact and computationally efficient, improving over existing heuristic methods in terms of both speed and accuracy significantly.
The key is to switch the view from optimizing over $k$-atomic mixing distributions (which is non-convex) to moment space (which is convex and efficiently optimizable via SDP).
We also note that minimizing the sum of squares in GMM is not crucial and minimizing any distance yields the same theoretical guarantee. 

There are a number of recent work in the theoretical computer science literature on provable results for moment-based estimators in Gaussian location-scale mixture models, see, e.g., \cite{MV2010,KMV2010,BS10,HP15,li2017robust}.
For instance, the algorithm \cite{MV2010} is based on exhaustive search over the discretized parameter space such that the population moments is close to the empirical moments. In addition to being computationally expensive, this method achieves the estimation accuracy $n^{-C/k}$ for some constant $C$, which is suboptimal in view of \prettyref{thm:main-W1}.
By carefully analyzing Pearson's method of moments equations \cite{Pearson1894},  \cite{HP15} showed that the optimal rate for two-component location-scale mixtures is $\Theta(n^{-1/12})$; however, this approach is difficult to generalize to higher order mixtures.
Finally, for moment-based methods in multiple dimensions, such as spectral and tensor decomposition, see the discussion in \prettyref{sec:discuss}.

\paragraph{Minimum distance estimators.} In the case of known variance, the minimum distance estimator is studied by \cite{DK68,Chen95,HK2015}. Specifically, the estimator is a $k$-atomic distribution $\hat\nu$ such that $\hat\nu*N(0,\sigma^2)$ is closest to the empirical distribution of the sample in a certain distance. The minimax optimal rate $O(n^{-\frac{1}{4k-2}})$ for estimating the mixing distribution under the Wasserstein distance is shown in \cite{HK2015} (which corrects the previous result in \cite{Chen95}), by bounding the $W_1$ distance between the mixing distributions in terms of the KS distance of the Gaussian mixtures \cite[Lemma 4.5]{HK2015}. 
However, the minimum distance estimator is in general computationally expensive and suffers from the same non-convexity issue of the MLE. 
In contrast, denoised method of moments is efficiently computable and adaptively achieves the optimal rate of accuracy as given in \prettyref{thm:main-W1-adaptive}.
For arbitrary Gaussian location mixtures in one dimension, the minimum distance estimator was considered in \cite{edelman1988estimation} in the context of empirical Bayes. Under the assumptions of bounded first moment, it is shown in \cite[Corollary 2]{edelman1988estimation} that the mixing distribution can be estimated at rate $O((\log n)^{-1/4})$ under the $L_2$-distance between the CDFs; this loss is, however, weaker than the $W_1$-distance (i.e.~$L_1$ distance between the CDFs).

\paragraph{Density estimation}

If the estimator is allowed to be any density (\emph{improper} learning), it is known that as long as the mixing distribution has a bounded support, the rate of convergence is close to parametric regardless of the number of components. Specifically, the optimal squared $L_2$-risk is found to be $\Theta(\frac{\sqrt{\log n}}{n})$ \cite{kim2014minimax}, achieved by the kernel density estimator designed for analytic densities \cite{ibragimov2001estimation}. 
As mentioned before, \emph{proper} density estimate (which is required to be a $k$-Gaussian mixture) is more desirable for the sake of interpretability; however, 
finding the $k$-Gaussian mixture that best approximates a given function such as a kernel density estimate can be computationally challenging due to, again, the non-convexity in the location parameters. 
In this regard, another contribution of Theorems \ref{thm:main-density} and \ref{thm:oracle} is the observation that proper and near optimal estimates/approximates can be found efficiently via the method of moments.
Finally, we note that MLE for estimating the density of general Gaussian mixtures has been studied in \cite{genovese.wasserman,GV2001}.

%

 \begin{acknowledgements}
     This work is supported in part by the NSF Grant CCF-1527105, CCF-1900507, an NSF CAREER award CCF-1651588, and an Alfred Sloan fellowship.
		The authors are grateful to the anonymous referees for helpful comments and suggestions and to Prof.~Alexandre Tsybakov for the reference \cite{Stein1986}.
 \end{acknowledgements}






\backmatter
\printbibliography

@inproceedings{OS15,
  title={Competitive Distribution Estimation: Why is {G}ood-{T}uring Good},
  author={Orlitsky, Alon and Suresh, Ananda Theertha},
  booktitle={Proceedings of the Twenty-ninth Conference on Neural Information Processing Systems},
  pages={2143--2151},
  publisher={Curran Associates, Inc.},
  year={2015}
}

@article{ICL09,
  title={Estimating the number of unseen variants in the human genome},
  author={Ionita-Laza, Iuliana and Lange, Christoph and Laird, Nan M},
  journal={Proceedings of the National Academy of Sciences},
  volume={106},
  number={13},
  pages={5008--5013},
  year={2009},
  publisher={National Acad Sciences}
}

@article{ET76,
  title={Estimating the number of unseen species: How many words did {S}hakespeare know?},
  author={Efron, B. and Thisted, R.},
  journal={Biometrika},
  volume={63},
  number={3},
  pages={435--447},
  year={1976},
  publisher={Biometrika Trust}
}

@article{CL11,
	Author = {Cai, T.T. and Low, M. G.},
	Journal = {The Annals of Statistics},
	Number = {2},
	Pages = {1012--1041},
	Title = {Testing composite hypotheses, {Hermite} polynomials and optimal estimation of a nonsmooth functional},
	Volume = {39},
	Year = {2011}
}

@ARTICLE{RRSS09,
  author = {Raskhodnikova, Sofya and Ron, Dana and Shpilka, Amir and Smith, Adam},
  title = {Strong lower bounds for approximating distribution support size and
	the distinct elements problem},
  journal = {SIAM Journal on Computing},
  year = {2009},
  volume = {39},
  pages = {813--842},
  number = {3},
  publisher = {SIAM}
}

@article{LNS99,
  title={On estimation of the {$L_r$} norm of a regression function},
  author={Lepski, Oleg and Nemirovski, Arkady and Spokoiny, Vladimir},
  journal={Probability Theory and Related Fields},
  volume={113},
  number={2},
  pages={221--253},
  year={1999}
}

@article{WY14,
  author = {Yihong Wu and Pengkun Yang},
  title = {Minimax rates of entropy estimation on large alphabets via best polynomial approximation},
  journal = {to appear in IEEE Transactions on Information Theory, arxiv:1407.0381},
  month = {Jul},
  year = {2015}
}

@article{JVHW15,
  title={Minimax estimation of functionals of discrete distributions},
  author={Jiao, Jiantao and Venkat, Kartik and Han, Yanjun and Weissman, Tsachy},
  journal={IEEE Transactions on Information Theory},
  volume={61},
  number={5},
  pages={2835--2885},
  year={2015},
  publisher={IEEE}
}

@ARTICLE{BF93,
  author = {Bunge, John and Fitzpatrick, M},
  title = {Estimating the number of species: a review},
  journal = {Journal of the American Statistical Association},
  year = {1993},
  volume = {88},
  pages = {364--373},
  number = {421},
  publisher = {Taylor \& Francis Group}
}

@inproceedings{VV11,
  author = {Valiant, Gregory and Valiant, Paul},
  title = {Estimating the unseen: an $n/\log (n)$-sample estimator for entropy
	and support size, shown optimal via new {CLT}s},
  booktitle = {Proceedings of the 43rd annual ACM symposium on Theory of computing},
  year = {2011},
  pages = {685--694}
}

@book{AS64,
	Address = {New York, NY},
	Author = {Abramowitz, M. and Stegun, I. A.},
	Publisher = {Wiley-Interscience},
	Title = {{Handbook of mathematical functions with formulas, graphs, and mathematical tables}},
	Year = {1964}
}

@article{FCW43,
  author = {Fisher, Ronald Aylmer and Corbet, A Steven and Williams, Carrington B},
  title = {The relation between the number of species and the number of individuals in a random sample of an animal population},
  journal = {The Journal of Animal Ecology},
  year = {1943},
  volume={12}, number={1},
  pages = {42--58}
}

@article{GT56,
  title={The number of new species, and the increase in population coverage, when a sample is increased},
  author={Good, I.J. and Toulmin, G.H.},
  journal={Biometrika},
  volume={43},
  number={1-2},
  pages={45--63},
  year={1956}
}

@article{TE87,
  title={Did {S}hakespeare write a newly-discovered poem?},
  author={Thisted, Ronald and Efron, Bradley},
  journal={Biometrika},
  volume={74},
  number={3},
  pages={445--455},
  year={1987},
  publisher={Biometrika Trust}
}

@INPROCEEDINGS{VV13,
  author = {Valiant, Paul and Valiant, Gregory},
  title = {Estimating the Unseen: Improved Estimators for Entropy and other
	Properties},
  booktitle = {Advances in Neural Information Processing Systems},
  year = {2013},
  pages = {2157--2165}
}

@article{zou2016quantifying,
  title={Quantifying unobserved protein-coding variants in human populations provides a roadmap for large-scale sequencing projects},
  author={Zou, James and Valiant, Gregory and Valiant, Paul and Karczewski, Konrad and Chan, Siu On and Samocha, Kaitlin and Lek, Monkol and Sunyaev, Shamil and Daly, Mark and MacArthur, Daniel G},
  journal={Nature communications},
  volume={7},
  pages={13293},
  year={2016},
  publisher={Nature Publishing Group}
}

@book{Ismail2005,
  title={Classical and quantum orthogonal polynomials in one variable},
  author={Ismail, Mourad E. H.},
  volume={13},
  year={2005},
  publisher={Cambridge university press}
}

@book{KN1977,
  title={The Markov moment problem and extremal problems},
  author={Krein, Mark Grigorevich and Adol'f Abramovich Nudel'man},
  year={1977},
  publisher={American Mathematical Society}
}

@article{bandeira2017optimal,
  title={Optimal rates of estimation for multi-reference alignment},
  author={Bandeira, Afonso S and Rigollet, Philippe and Weed, Jonathan},
  journal={arXiv:1702.08546},
  year={2017}
}

@article{DWYZ20,
  title={Optimal estimation of high-dimensional Gaussian mixtures},
  author={Doss, Natalie and Wu, Yihong and Yang, Pengkun and Zhou, Harrison H},
  journal={arXiv:2002.05818},
  year={2020}
}

@book{Reimer2012,
  title={Multivariate polynomial approximation},
  author={Reimer, Manfred},
  volume={144},
  year={2012},
  publisher={Birkh{\"a}user}
}

@article{Stein1986,
  title={Lectures on the theory of estimation of many parameters},
  author={Stein, Charles},
  journal={Journal of Soviet Mathematics},
  volume={34},
  number={1},
  pages={1373--1403},
  year={1986},
  note = {[Zap. Nauchn. Sem. LOMI, 1977, Volume 74, Pages 4–65]}
}

@article{nikolsky1946mean,
  title={Approximation of functions in the mean by trigonometrical polynomials},
  author={Nikolsky, S},
  journal={Izvestiya Rossiiskoi Akademii Nauk. Seriya Matematicheskaya},
  volume={10},
  number={3},
  pages={207--256},
  year={1946}
}

@inproceedings{tsai1999analysis,
  title={Analysis of functional {MRI} data using mutual information},
  author={Tsai, Andy and Fisher, John W and Wible, Cindy and Wells, William M and Kim, Junmo and Willsky, Alan S},
  booktitle={International Conference on Medical Image Computing and Computer-Assisted Intervention},
  pages={473--480},
  year={1999},
  organization={Springer}
}

@article{tsai2004mutual,
  title={Mutual information in coupled multi-shape model for medical image segmentation},
  author={Tsai, Andy and Wells, William and Tempany, Clare and Grimson, Eric and Willsky, Alan},
  journal={Medical image analysis},
  volume={8},
  number={4},
  pages={429--445},
  year={2004},
  publisher={Elsevier}
}

@article{pluim2003mutual,
  title={Mutual-information-based registration of medical images: a survey},
  author={Pluim, Josien PW and Maintz, JB Antoine and Viergever, Max A},
  journal={IEEE transactions on medical imaging},
  volume={22},
  number={8},
  pages={986--1004},
  year={2003},
  publisher={Citeseer}
}

@inproceedings{kybic2004high,
  title={High-dimensional mutual information estimation for image registration},
  author={Kybic, Jan},
  booktitle={2004 International Conference on Image Processing, 2004. ICIP'04.},
  volume={3},
  pages={1779--1782},
  year={2004},
  organization={IEEE}
}

@article{PW18,
  title={Dualizing {L}e {C}am's method, with applications to estimating the unseens},
  author={Polyanskiy, Yury and Wu, Yihong},
  journal={arXiv:1902.05616},
  year={2019}
}

@article{Aktulga07,
  title={Identifying statistical dependence in genomic sequences via mutual information estimates},
  author={Aktulga, Hasan Metin and Kontoyiannis, Ioannis and Lyznik, L Alex and Szpankowski, Lukasz and Grama, Ananth Y and Szpankowski, Wojciech},
  journal={EURASIP Journal on Bioinformatics and Systems Biology},
  volume={2007},
  pages={3},
  year={2007},
  publisher={Hindawi Publishing Corp.}
}

@article{GM-suppl,
  title={Supplement to `optimal estimation of {G}aussian mixtures via denoised method of moments'},
  Author={Yihong Wu and Pengkun Yang},
  Year={2020},
  doi={10.1214/19-AOS1873SUPP}
}

@book{DK2012,
  title={Moduli of smoothness},
  author={Ditzian, Zeev and Totik, Vilmos},
  volume={9},
  year={2012},
  publisher={Springer Science \& Business Media}
}

@article{Ivanov1983,
  title={On a new characteristic of functions. II. {D}irect and converse theorems for the best algebraic approximation in {$C[-1,1]$} and {$L_p[-1,1]$}},
  author={Ivanov, Kamen G},
  journal={Pliska},
  volume={5},
  pages={151--163},
  year={1983}
}

@inproceedings{ADOS17,
  title={A unified maximum likelihood approach for estimating symmetric properties of discrete distributions},
  author={Acharya, Jayadev and Das, Hirakendu and Orlitsky, Alon and Suresh, Ananda Theertha},
  booktitle={Proceedings of the 34th International Conference on Machine Learning},
  pages={11--21},
  publisher={PLMR},
  year={2017}
}

@inproceedings{HJW18lmm,
  title={Local moment matching: A unified methodology for symmetric functional estimation and distribution estimation under {W}asserstein distance},
  author={Han, Yanjun and Jiao, Jiantao and Weissman, Tsachy},
  booktitle={Proceedings of the 31st  Conference On Learning Theory},
  pages={3189--3221},
  publisher={PMLR},
  year={2018}
}

@book{MS1999,
  title={Foundations of statistical natural language processing},
  author={Manning, Christopher D and Sch{\"u}tze, Hinrich},
  year={1999},
  publisher={MIT press}
}

@book{TV2005,
  title={Fundamentals of wireless communication},
  author={Tse, David and Viswanath, Pramod},
  year={2005},
  publisher={Cambridge university press}
}

@book{Oppenheim1999,
  title={Discrete-time signal processing},
  author={Oppenheim, Alan V},
  year={1999},
  publisher={Pearson Education India}
}

@book{Mitchell1997,
  title={Machine Learning},
  author={Mitchell, Tom M},
  year={1997},
  publisher={McGraw Hill}
}

@inproceedings{KSV2005,
  title={The spectral method for general mixture models},
  author={Kannan, Ravindran and Salmasian, Hadi and Vempala, Santosh},
  booktitle={International Conference on Computational Learning Theory},
  pages={444--457},
  year={2005},
  organization={Springer}
}

@inproceedings{Dasgupta1999,
  title={Learning mixtures of Gaussians},
  author={Dasgupta, Sanjoy},
  booktitle={Foundations of computer science, 1999. 40th annual symposium on},
  pages={634--644},
  year={1999},
  organization={IEEE}
}

@book{Lasserre2009,
  title={Moments, positive polynomials and their applications},
  author={Lasserre, Jean Bernard},
  volume={1},
  year={2009},
  publisher={World Scientific}
}

@Article{Chao92,
  Title                    = {Estimating the number of classes via sample coverage},
  Author                   = {Chao, Anne and Lee, Shen-Ming},
  Journal                  = {Journal of the American statistical Association},
  Year                     = {1992},
  Number                   = {417},
  Pages                    = {210--217},
  Volume                   = {87},

  Publisher                = {Taylor \& Francis}
}

@Article{HW01,
  Title                    = {Estimating the total number of alleles using a sample coverage method},
  Author                   = {Huang, Shu-Pang and Weir, BS},
  Journal                  = {Genetics},
  Year                     = {2001},
  Number                   = {3},
  Pages                    = {1365--1373},
  Volume                   = {159},

  Publisher                = {Genetics Soc America}
}

@Article{McNeil73,
  Title                    = {Estimating an author's vocabulary},
  Author                   = {McNeil, Donald R},
  Journal                  = {Journal of the American Statistical Association},
  Year                     = {1973},
  Number                   = {341},
  Pages                    = {92--96},
  Volume                   = {68},

  Publisher                = {Taylor \& Francis}
}

@InProceedings{Valiant08,
  Title                    = {Testing Symmetric Properties of Distributions},
  Author                   = {Valiant, Paul},
  Booktitle                = {Proceedings of the Fortieth Annual ACM Symposium on Theory of Computing},
  Year                     = {2008},
  Pages                    = {383--392},
  Series                   = {STOC '08},

  Numpages                 = {10}
}

@InProceedings{BFSS02,
  Title                    = {How to achieve minimax expected {K}ullback-{L}eibler distance from an unknown finite distribution},
  Author                   = {Braess, Dietrich and Forster, J{\"u}rgen and Sauer, Tomas and Simon, Hans U},
  Booktitle                = {Algorithmic Learning Theory},
  Year                     = {2002},
  Organization             = {Springer},
  Pages                    = {380--394}
}

@Article{Miller55,
  Title                    = {Note on the bias of information estimates},
  Author                   = {Miller, George A.},
  Journal                  = {Information theory in psychology: Problems and methods},
  Year                     = {1955},
  Pages                    = {95--100},
  Volume                   = {2}
}

@InCollection{Harris75,
  Title                    = {The statistical estimation of entropy in the non-parametric case},
  Author                   = {B. Harris},
  Booktitle                = {Topics in Information Theory},
  Publisher                = {Springer Netherlands},
  Year                     = {1975},
  Editor                   = {I Csisz\'ar and P. Elias},
  Pages                    = {323-355},
  Volume                   = {16}
}

@article{Lindsay1989,
  title={Moment matrices: applications in mixtures},
  author={Lindsay, Bruce G},
  journal={The Annals of Statistics},
  pages={722--740},
  year={1989},
  publisher={JSTOR}
}

@article{Hansen1982,
  title={Large sample properties of generalized method of moments estimators},
  author={Hansen, Lars Peter},
  journal={Econometrica: Journal of the Econometric Society},
  pages={1029--1054},
  year={1982},
  publisher={JSTOR}
}

@book{Hall2005,
  title={Generalized method of moments},
  author={Hall, Alastair R},
  year={2005},
  publisher={Oxford University Press}
}

@inproceedings{HP15,
  title={Tight bounds for learning a mixture of two gaussians},
  author={Hardt, Moritz and Price, Eric},
  booktitle={Proceedings of the Forty-Seventh Annual ACM on Symposium on Theory of Computing},
  pages={753--760},
  year={2015},
  organization={ACM}
}

@inproceedings{KMV2010,
  title={Efficiently learning mixtures of two {G}aussians},
  author={Kalai, Adam Tauman and Moitra, Ankur and Valiant, Gregory},
  booktitle={Proceedings of the forty-second ACM symposium on Theory of computing},
  pages={553--562},
  year={2010},
  organization={ACM}
}

@inproceedings{MV2010,
  title={Settling the polynomial learnability of mixtures of {G}aussians},
  author={Moitra, Ankur and Valiant, Gregory},
  booktitle={Foundations of Computer Science (FOCS), 2010 51st Annual IEEE Symposium on},
  pages={93--102},
  year={2010},
  organization={IEEE}
}

@article{meng2016mllib,
  title={Mllib: Machine learning in apache spark},
  author={Meng, Xiangrui and Bradley, Joseph and Yavuz, Burak and Sparks, Evan and Venkataraman, Shivaram and Liu, Davies and Freeman, Jeremy and Tsai, DB and Amde, Manish and Owen, Sean and others},
  journal={The Journal of Machine Learning Research},
  volume={17},
  number={1},
  pages={1235--1241},
  year={2016},
  publisher={JMLR. org}
}

@inproceedings{abadi2016tensorflow,
  title={Tensorflow: A system for large-scale machine learning.},
  author={Abadi, M. and Barham, P. and Chen, J. and Chen, Z. and Davis, A. and Dean, J. and Devin, M. and Ghemawat, S. and Irving, G. and Isard, M. and Kudlur, M. and Levenberg, J. and Monga, R. and Moore, S. and Murray, D. G. and Steiner, B. and Tucker, P. and Vasudevan, V. and Warden, P. and Wicke, M. and Yu, Y. and X. Zheng},
  booktitle={12th USENIX Symposium on Operating Systems Design and Implementation (OSDI 16)},
  volume={16},
  pages={265--283},
  publisher={USENIX Association},
  year={2016}
}

@article{scikit-learn,
 title={Scikit-learn: Machine Learning in {P}ython},
 author={Pedregosa, F. and Varoquaux, G. and Gramfort, A. and Michel, V.
         and Thirion, B. and Grisel, O. and Blondel, M. and Prettenhofer, P.
         and Weiss, R. and Dubourg, V. and Vanderplas, J. and Passos, A. and
         Cournapeau, D. and Brucher, M. and Perrot, M. and Duchesnay, E.},
 journal={Journal of Machine Learning Research},
 volume={12},
 pages={2825--2830},
 year={2011}
}

@article{Pearson1894,
  title={Contributions to the mathematical theory of evolution},
  author={Pearson, Karl},
  journal={Philosophical Transactions of the Royal Society of London. A},
  volume={185},
  pages={71--110},
  year={1894},
  publisher={JSTOR}
}

@Article{WVK11,
  Title                    = {Probability estimation in the rare-events regime},
  Author                   = {Wagner, Aaron B and Viswanath, Pramod and Kulkarni, Sanjeev R},
  Journal                  = IEEE_J_IT,
  Year                     = {2011},
  Number                   = {6},
  Pages                    = {3207--3229},
  Volume                   = {57},

  Publisher                = {IEEE}
}

@Article{KWTV13,
  Title                    = {Classification of Homogeneous Data With Large Alphabets},
  Author                   = {Kelly, B.G. and Wagner, A.B. and Tularak, T. and Viswanath, P.},
  Journal                  = {IEEE Transactions on Information Theory},
  Year                     = {2013},
  Number                   = {2},
  Pages                    = {782-795},
  Volume                   = {59}
}

@InProceedings{BS09,
  Title                    = {Knowing the unseen: estimating vocabulary size over unseen samples},
  Author                   = {Bhat, Suma and Sproat, Richard},
  Booktitle                = {Proceedings of the Joint Conference of the 47th Annual Meeting of the ACL and the 4th International Joint Conference on Natural Language Processing of the AFNLP: Volume 1},
  Year                     = {2009},
  Pages                    = {109--117}
}

@Article{OSZ04,
  Title                    = {Universal compression of memoryless sources over unknown alphabets},
  Author                   = {Orlitsky, Alon and Santhanam, Narayana P and Zhang, Junan},
  Journal                  = {IEEE Transactions on Information Theory},
  Year                     = {2004},
  Number                   = {7},
  Pages                    = {1469--1481},
  Volume                   = {50},

  Publisher                = {IEEE}
}

@Article{SLSKB97,
  Title                    = {Reproducibility and variability in neural spike trains},
  Author                   = {de Ruyter van Steveninck, Rob R. and Lewen, Geoffrey D. and Strong, Steven P. and Koberle, Roland and Bialek, William},
  Journal                  = {Science},
  Year                     = {1997},
  Number                   = {5307},
  Pages                    = {1805--1808},
  Volume                   = {275}
}

@Article{mainen1995reliability,
  Title                    = {Reliability of spike timing in neocortical neurons},
  Author                   = {Mainen, Zachary F and Sejnowski, Terrence J},
  Journal                  = {Science},
  Year                     = {1995},
  Number                   = {5216},
  Pages                    = {1503--1506},
  Volume                   = {268},

  Publisher                = {American Association for the Advancement of Science}
}

@Article{Berry13051997,
  Title                    = {The structure and precision of retinal spike trains},
  Author                   = {Berry, Michael J. and Warland, David K. and Meister, Markus},
  Journal                  = {Proceedings of the National Academy of Sciences},
  Year                     = {1997},
  Number                   = {10},
  Pages                    = {5411-5416},
  Volume                   = {94}
}

@InProceedings{benevenuto2009characterizing,
  Title                    = {Characterizing user behavior in online social networks},
  Author                   = {Benevenuto, Fabr{\'\i}cio and Rodrigues, Tiago and Cha, Meeyoung and Almeida, Virg{\'\i}lio},
  Booktitle                = {Proceedings of the 9th ACM SIGCOMM conference on Internet measurement conference},
  Year                     = {2009},
  Pages                    = {49--62}
}

@Book{VdV00,
  Title                    = {Asymptotic statistics},
  Author                   = {Van der Vaart, Aad W.},
  Publisher                = {Cambridge university press},
  Year                     = {2000},

  Address                  = {Cambridge, United Kingdom}
}

@Article{Berkson80,
  Title                    = {Minimum chi-square, not maximum likelihood! (with discussion)},
  Author                   = {Berkson, Joseph},
  Journal                  = A_S,
  Year                     = {1980},
  Pages                    = {457--487}
}

@Article{Efron82,
  Title                    = {Maximum Likelihood and Decision Theory},
  Author                   = {Efron, Bradley},
  Journal                  = A_S,
  Year                     = {1982},
  Number                   = {2},
  Pages                    = {pp. 340-356},
  Volume                   = {10},

  ISSN                     = {00905364},
  Language                 = {English}
}

@Article{DL91,
  Title                    = {Geometrizing rates of convergence, {II}},
  Author                   = {Donoho, David L. and Liu, Richard C.},
  Journal                  = A_S,
  Year                     = {1991},
  Pages                    = {668--701},
  Volume                   = {19}
}

@Article{CL05,
  Title                    = {Nonquadratic estimators of a quadratic functional},
  Author                   = {Cai, T. T. and Low, M. G.},
  Journal                  = A_S,
  Year                     = {2005},
  Number                   = {6},
  Pages                    = {2930--2956},
  Volume                   = {33}
}

@Article{Stone80,
  Title                    = {Optimal rates of convergence for nonparametric estimators},
  Author                   = {Stone, Charles J.},
  Journal                  = A_S,
  Year                     = {1980},
  Number                   = {6},
  Pages                    = {1348--1360},
  Volume                   = {8}
}

@book{Rao2014,
  title={Nonparametric functional estimation},
  author={Rao, BLS Prakasa},
  year={2014},
  publisher={Academic press}
}

@ARTICLE{Basharin59,
  author = {Basharin, G.P.},
  title = {On a statistical estimate for the entropy of a sequence of independent
	random variables},
  journal = {Theory of Probability \& Its Applications},
  year = {1959},
  volume = {4},
  pages = {333--336},
  number = {3},
  __markedentry = {[pyang14:]},
  publisher = {SIAM}
}

@ARTICLE{BO79,
  author = {Burnham, Kenneth P and Overton, W Scott},
  title = {Robust estimation of population size when capture probabilities vary
	among animals},
  journal = {Ecology},
  year = {1979},
  volume = {60},
  pages = {927--936},
  number = {5},
  publisher = {Eco Soc America}
}

@ARTICLE{Chao84,
  author = {Chao, Anne},
  title = {Nonparametric estimation of the number of classes in a population},
  journal = {Scandinavian Journal of statistics},
  year = {1984},
  pages = {265--270},
  publisher = {JSTOR}
}

@ARTICLE{DR80,
  author = {Darroch, JN and Ratcliff, D},
  title = {A note on capture-recapture estimation},
  journal = {Biometrics},
  year = {1980},
  pages = {149--153},
  publisher = {JSTOR}
}

@ARTICLE{Good1953,
  author = {Good, Irving J},
  title = {The population frequencies of species and the estimation of population
	parameters},
  journal = {Biometrika},
  year = {1953},
  volume = {40},
  pages = {237--264},
  number = {3-4},
  publisher = {Biometrika Trust}
}

@ARTICLE{INK87,
  author = {Ibragimov, IA and Nemirovskii, AS and Khas' minskii, RZ},
  title = {Some problems on nonparametric estimation in Gaussian white noise},
  journal = {Theory of Probability \& Its Applications},
  year = {1987},
  volume = {31},
  pages = {391--406},
  number = {3},
  publisher = {SIAM}
}

@ARTICLE{Paninski04,
  author = {Paninski, Liam},
  title = {Estimating entropy on $m$ bins given fewer than $m$ samples},
  journal = {IEEE Transactions on Information Theory},
  year = {2004},
  volume = {50},
  pages = {2200--2203},
  number = {9},
  __markedentry = {[pyang14:]}
}

@ARTICLE{Paninski03,
  author = {Paninski, Liam},
  title = {Estimation of entropy and mutual information},
  journal = {Neural Computation},
  year = {2003},
  volume = {15},
  pages = {1191--1253},
  number = {6},
  __markedentry = {[pyang14:]}
}

@BOOK{timan63,
  title = {Theory of approximation of functions of a real variable},
  publisher = {Courier Corporation},
  year = {1963},
  author = {Timan, Aleksandr Filippovich},
  volume = {34}
}

@BOOK{Tsybakov09,
  title = {Introduction to Nonparametric Estimation},
  publisher = {Springer Verlag},
  year = {2009},
  author = {Tsybakov, A.B.},
  address = {New York, NY}
}

@ARTICLE{VV10,
  author = {Gregory Valiant and Paul Valiant},
  title = {A {CLT} and tight lower bounds for estimating entropy},
  journal = {Electronic Colloquium on Computational Complexity (ECCC)},
  year = {2010},
  __markedentry = {[pyang14:]}
}

@ARTICLE{Vinck12,
  author = {Vinck, Martin and Battaglia, Francesco P. and Balakirsky, Vladimir
	B. and Vinck, A.J. Han and Pennartz, Cyriel M.A.},
  title = {Estimation of the entropy based on its polynomial representation},
  journal = {Physical Review E},
  year = {2012},
  volume = {85},
  pages = {051139},
  number = {5}
}

@book{hardle2012wavelets,
  title={Wavelets, approximation, and statistical applications},
  author={H{\"a}rdle, Wolfgang and Kerkyacharian, Gerard and Picard, Dominique and Tsybakov, Alexander},
  volume={129},
  year={2012},
  publisher={Springer Science \& Business Media}
}

@inproceedings{jiao2016beyond,
  title={Beyond maximum likelihood: {B}oosting the {C}how-{L}iu algorithm for large alphabets},
  author={Jiao, Jiantao and Han, Yanjun and Weissman, Tsachy},
  booktitle={2016 50th Asilomar Conference on Signals, Systems and Computers},
  pages={321--325},
  year={2016},
  organization={IEEE}
}

@article{laurent1996efficient,
  title={Efficient estimation of integral functionals of a density},
  author={Laurent, B{\'e}atrice},
  journal={The Annals of Statistics},
  volume={24},
  number={2},
  pages={659--681},
  year={1996},
  publisher={Institute of Mathematical Statistics}
}

@InProceedings{PSW17-colt,
  author = {Y. Polyanskiy and A. T. Suresh and Y. Wu},
  title = {Sample complexity of population recovery},
 booktitle = {Proceedings of Conference on Learning Theory (COLT)},
  month = {Jul},
  year = {2017},
  address={Amsterdam, Netherland},
note={arXiv:1702.05574}}

@article{JN09,
  title={Nonparametric estimation by convex programming},
  author={Juditsky, Anatoli B and Nemirovski, Arkadi S},
  journal={The Annals of Statistics},
  pages={2278--2300},
  year={2009},
  volume={37},
  number={5A}
}

@article{KV17,
  title={Spectrum estimation from samples},
  author={Kong, Weihao and Valiant, Gregory},
  journal={The Annals of Statistics},
  volume={45},
  number={5},
  pages={2218--2247},
  year={2017},
  publisher={Institute of Mathematical Statistics}
}

@inproceedings{XHM16,
  title={Global analysis of expectation maximization for mixtures of two {G}aussians},
  author={Xu, Ji and Hsu, Daniel J and Maleki, Arian},
  booktitle={Proceedings of the Thirtieth Conference on Neural Information Processing Systems},
  pages={2676--2684},
  year={2016},
  publisher={Curran Associates, Inc.}
}

@inproceedings{DTZ17,
  title={Ten Steps of {EM} Suffice for Mixtures of Two {G}aussians},
  author={Daskalakis, Constantinos and Tzamos, Christos and Zampetakis, Manolis},
  booktitle={Proceedings of the 30th Annual Conference on Learning Theory (COLT 2017)},
  pages={704--710},
  publisher={PMLR},
  year={2017}
}

@book{Rudin2006,
  title={Real and complex analysis},
  author={Rudin, Walter},
  year={2006},
  publisher={Tata McGraw-Hill Education}
}

@book{Fruhwirth2006,
  title={Finite mixture and {M}arkov switching models},
  author={Fr{\"u}hwirth-Schnatter, Sylvia},
  year={2006},
  publisher={Springer Science \& Business Media}
}

@Book{horn-2nd,
  Title                    = {Matrix Analysis},
  Author                   = {Roger A. Horn and Charles R. Johnson},
  Publisher                = {Cambridge University Press},
  Year                     = {2012},
  edition                  = 2
}

@article{BS2004,
  title={Bernstein polynomials and learning theory},
  author={Braess, Dietrich and Sauer, Thomas},
  journal={Journal of Approximation Theory},
  volume={128},
  number={2},
  pages={187--206},
  year={2004},
  publisher={Elsevier}
}

@book{Attneave1959,
  title={Applications of information theory to psychology: A summary of basic concepts, methods, and results},
  author={Attneave, Fred},
  year={1959},
  publisher={Holt, Rinehart and Winston}
}

@book{Gautschi2004,
  title={Orthogonal polynomials: computation and approximation},
  author={Gautschi, Walter},
  year={2004},
  publisher={Oxford University Press on Demand}
}

@book{WSV2012,
  title={Handbook of semidefinite programming: theory, algorithms, and applications},
  author={Wolkowicz, Henry and Saigal, Romesh and Vandenberghe, Lieven},
  volume={27},
  year={2012},
  publisher={Springer Science \& Business Media}
}

@book{Rockafellar1974,
  title={Conjugate duality and optimization},
  author={Rockafellar, R Tyrrell},
  volume={16},
  year={1974},
  publisher={Siam}
}

@article{deBoor2005,
  title={Divided differences},
  author={{de Boor}, Carl},
  journal={Surveys in Approximation Theory},
  volume={1},
  pages={46--49},
  year={2005}
}

@book{Davis1975,
  title={Interpolation and approximation},
  author={Davis, Philip J},
  year={1975},
  publisher={Courier Corporation}
}

@book{Prasolov2009,
  title={Polynomials},
  author={Prasolov, Victor V},
  volume={11},
  year={2009},
  publisher={Springer Science \& Business Media}
}

@article{DS2016,
  title={Bounds for the logarithm of the {E}uler gamma function and its derivatives},
  author={Diamond, Harold G and Straub, Armin},
  journal={Journal of Mathematical Analysis and Applications},
  volume={433},
  number={2},
  pages={1072--1083},
  year={2016},
  publisher={Elsevier}
}

@article{KW56,
  title={Consistency of the maximum likelihood estimator in the presence of infinitely many incidental parameters},
  author={Kiefer, Jack and Wolfowitz, Jacob},
  journal={The Annals of Mathematical Statistics},
  pages={887--906},
  year={1956},
  publisher={JSTOR}
}

@inproceedings{HK2013,
  title={Learning mixtures of spherical {G}aussians: moment methods and spectral decompositions},
  author={Hsu, Daniel and Kakade, Sham M},
  booktitle={Proceedings of the 4th conference on Innovations in Theoretical Computer Science},
  pages={11--20},
  year={2013},
  organization={ACM}
}

@article{morris1982natural,
  title={Natural exponential families with quadratic variance functions},
  author={Morris, Carl N},
  journal={Annals of Statistics},
  pages={65--80},
  year={1982},
  volume={10},
  issue={1}
}

@article{KX2005,
  title={Mixed {P}oisson distributions},
  author={Karlis, Dimitris and Xekalaki, Evdokia},
  journal={International Statistical Review},
  volume={73},
  number={1},
  pages={35--58},
  year={2005},
  publisher={Wiley Online Library}
}

@article{Jewell1982,
  title={Mixtures of exponential distributions},
  author={Jewell, Nicholas P},
  journal={Annals of Statistics},
  pages={479--484},
  year={1982},
  volume={10},
  issue={2},
  publisher={JSTOR}
}

@inproceedings{WS2000,
  title={Scale mixtures of {G}aussians and the statistics of natural images},
  author={Wainwright, Martin J and Simoncelli, Eero P},
  booktitle={Proceedings of the Thirteenth Conference on Neural Information Processing Systems (NIPS 1999)},
  pages={855--861},
  year={2000},
  publisher={MIT Press}
}

@article{AM1974,
  title={Scale mixtures of normal distributions},
  author={Andrews, David F and Mallows, Colin L},
  journal={Journal of the Royal Statistical Society. Series B (Methodological)},
  pages={99--102},
  volume={36},
  issue={1},
  year={1974},
  publisher={JSTOR}
}

@article{PSWS2003,
  title={Image denoising using scale mixtures of {G}aussians in the wavelet domain},
  author={Portilla, Javier and Strela, Vasily and Wainwright, Martin J and Simoncelli, Eero P},
  journal={IEEE Transactions on Image processing},
  volume={12},
  number={11},
  pages={1338--1351},
  year={2003},
  publisher={IEEE}
}

@inproceedings{Krawtchouk1932,
  title={Sur le probl{\`e}me de moments},
  author={Krawtchouk, Michel},
  booktitle={ICM Proceedings},
  pages={127--128},
  year={1932},
  Note={Available at \url{https://www.mathunion.org/fileadmin/ICM/Proceedings/ICM1932.2/ICM1932.2.ocr.pdf}},
}

@incollection{Diaconis1987,
  title={Application of the method of moments in probability and statistics},
  author={Diaconis, Persi},
  booktitle={Moments in mathematics},
  volume={37},
  pages={125--139},
  year={1987},
  publisher={Amer. Math. Soc.: Providence, RI}
}

@inproceedings{BV2008,
  title={Isotropic {PCA} and affine-invariant clustering},
  author={Brubaker, Spencer Charles and Vempala, Santosh},
  booktitle={IEEE 49th Annual IEEE Symposium on Foundations of Computer Science, 2008.},
  pages={551--560},
  year={2008},
  organization={IEEE}
}

@inproceedings{AM2005,
  title={On spectral learning of mixtures of distributions},
  author={Achlioptas, Dimitris and McSherry, Frank},
  booktitle={Proceedings of the 18th Annual Conference on Learning Theory (COLT 2005)},
  pages={458--469},
  year={2005},
  organization={Berlin, Heidelberg: Springer}
}

@article{VW2004,
  title={A spectral algorithm for learning mixture models},
  author={Vempala, Santosh and Wang, Grant},
  journal={Journal of Computer and System Sciences},
  volume={68},
  number={4},
  pages={841--860},
  year={2004},
  publisher={Elsevier}
}

@inproceedings{AK2001,
  title={Learning mixtures of arbitrary {G}aussians},
  author={Arora, Sanjeev and Kannan, Ravi},
  booktitle={Proceedings of the thirty-third annual ACM symposium on Theory of computing},
  pages={247--257},
  year={2001},
  organization={ACM}
}

@article{ARS2016,
  title={Algebraic Identifiability of {G}aussian Mixtures},
  author={Am{\'e}ndola, Carlos and Ranestad, Kristian and Sturmfels, Bernd},
  journal={International Mathematics Research Notices},
  volume={2018},
  issue={21},
  pages={6556–-6580},
  year={2018}
}

@inproceedings{Lindsay1995,
  title={Mixture models: Theory, geometry and applications},
  author={Lindsay, Bruce G},
  booktitle={NSF-CBMS regional conference series in probability and statistics},
  volume={5},
  pages={I--163},
  year={1995},
  organization={American Statistical Association}
}

@article{Laird1978,
  title={Nonparametric maximum likelihood estimation of a mixing distribution},
  author={Laird, Nan},
  journal={Journal of the American Statistical Association},
  volume={73},
  number={364},
  pages={805--811},
  year={1978},
  publisher={Taylor \& Francis}
}

@incollection{Lindsay1981,
  title={Properties of the maximum likelihood estimator of a mixing distribution},
  author={Lindsay, Bruce G},
  booktitle={Statistical Distributions in Scientific Work},
  pages={95--109},
  year={1981},
  publisher={Springer}
}

@article{GV2001,
  title={Entropies and rates of convergence for maximum likelihood and {B}ayes estimation for mixtures of normal densities},
  author={Ghosal, Subhashis and van der Vaart, Aad W},
  journal={Annals of Statistics},
  pages={1233--1263},
  year={2001},
  volume={29},
  issue={5},
  publisher={JSTOR}
}

@article{DLR1977,
  title={Maximum likelihood from incomplete data via the {EM} algorithm},
  author={Dempster, Arthur P and Laird, Nan M and Rubin, Donald B},
  journal={Journal of the Royal Statistical Society: Series B (Methodological)},
  pages={1--38},
  year={1977},
  volume={39},
  issue={1},
  publisher={JSTOR}
}

@article{HK2015,
  title={Strong identifiability and optimal minimax rates for finite mixture estimation},
  author={Heinrich, Philippe and Kahn, Jonas},
  journal={The Annals of Statistics},
  volume={46},
  number={6A},
  pages={2844--2870},
  year={2018},
  publisher={Institute of Mathematical Statistics}
}

@article{PL2001,
  title={Alternative {EM} methods for nonparametric finite mixture models},
  author={Pilla, Ramani S and Lindsay, Bruce G},
  journal={Biometrika},
  volume={88},
  number={2},
  pages={535--550},
  year={2001},
  publisher={Oxford University Press}
}

@article{SMA2000,
  title={A cautionary note on likelihood ratio tests in mixture models},
  author={Seidel, Wilfried and Mosler, Karl and Alker, Manfred},
  journal={Annals of the Institute of Statistical Mathematics},
  volume={52},
  number={3},
  pages={481--487},
  year={2000},
  publisher={Springer}
}

@book{Schmudgen17,
  title={The moment problem},
  author={Konrad Schm\"{u}dgen},
  year={2017},
  publisher={Springer}
}

@book{ST1943,
  title={The problem of moments},
  author={Shohat, James Alexander and Tamarkin, Jacob David},
  number={1},
  year={1943},
  publisher={American Mathematical Soc.}
}

@book{KS1953,
  title={Geometry of moment spaces},
  author={Karlin, Samuel and Shapley, Lloyd S},
  number={12},
  year={1953},
  publisher={American Mathematical Soc.}
}

@Article{Chausse2010,
    title = {Computing Generalized Method of Moments and Generalized Empirical Likelihood with {R}},
    author = {Pierre Chauss{\'e}},
    journal = {Journal of Statistical Software},
    year = {2010},
    volume = {34},
    number = {11},
    pages = {1--35},
    url = {http://www.jstatsoft.org/v34/i11/},
}

@article{KX2003,
  title={Choosing initial values for the {EM} algorithm for finite mixtures},
  author={Karlis, Dimitris and Xekalaki, Evdokia},
  journal={Computational Statistics \& Data Analysis},
  volume={41},
  number={3},
  pages={577--590},
  year={2003},
  publisher={Elsevier}
}

@article{XJ1996,
  title={On convergence properties of the {EM} algorithm for {G}aussian mixtures},
  author={Xu, Lei and Jordan, Michael I},
  journal={Neural computation},
  volume={8},
  number={1},
  pages={129--151},
  year={1996},
  publisher={MIT Press}
}

@article{RW1984,
  title={Mixture densities, maximum likelihood and the {EM} algorithm},
  author={Redner, Richard A and Walker, Homer F},
  journal={SIAM review},
  volume={26},
  number={2},
  pages={195--239},
  year={1984},
  publisher={SIAM}
}

@book{Uspensky37,
  title={Introduction to mathematical probability},
  author={Uspensky, James Victor},
  year={1937},
  publisher={McGraw-Hill}
}

@inproceedings{BS10,
  title={Polynomial learning of distribution families},
  author={Belkin, Mikhail and Sinha, Kaushik},
  booktitle={Foundations of Computer Science (FOCS), 2010 51st Annual IEEE Symposium on},
  pages={103--112},
  year={2010},
  organization={IEEE}
}

@article{Chen95,
  title={Optimal rate of convergence for finite mixture models},
  author={Chen, Jiahua},
  journal={The Annals of Statistics},
  pages={221--233},
  year={1995},
  volume={23},
  issue={1},
  publisher={JSTOR}
}

@article{OSW15,
  title={Optimal Prediction of the Number of Unseen Species},
  author={Orlitsky, Alon and Suresh, Ananda Theertha and Wu, Yihong},
  journal={Proceedings of the National Academy of Sciences (PNAS)},
  volume={113},number={47},pages={13283-13288},
  year={2016}
}

@InProceedings{JOST15,
  Title                    = {The complexity of estimating {R}{\'e}nyi entropy},
  Author                   = {Acharya, Jayadev and Orlitsky, Alon and Suresh, Ananda Theertha and Tyagi, Himanshu},
  Booktitle                = {Proceedings of the Twenty-Sixth Annual ACM-SIAM Symposium on Discrete Algorithms},
  Year                     = {2015},
  Organization             = {SIAM},
  Pages                    = {1855--1869}
}

@Article{amigo2004estimating,
  Title                    = {Estimating the entropy rate of spike trains via Lempel-Ziv complexity},
  Author                   = {Amig{\'o}, Jos{\'e} M and Szczepa{\'n}ski, Janusz and Wajnryb, Elek and Sanchez-Vives, Maria V},
  Journal                  = {Neural Computation},
  Year                     = {2004},
  Number                   = {4},
  Pages                    = {717--736},
  Volume                   = {16}
}

@Book{Atkinson89,
  Title                    = {An introduction to numerical analysis},
  Author                   = {Atkinson, Kendall E},
  Publisher                = {John Wiley \& Sons},
  Year                     = {1989}
}

@InCollection{BJKST02,
  Title                    = {Counting distinct elements in a data stream},
  Author                   = {Bar-Yossef, Ziv and Jayram, TS and Kumar, Ravi and Sivakumar, D and Trevisan, Luca},
  Booktitle                = { Proceedings of the 6th Randomization and Approximation Techniques in Computer Science},
  Year                     = {2002},
  Pages                    = {1--10},
  publisher = {Springer-Verlag}
}

@InProceedings{BKS01,
  Title                    = {Sampling algorithms: lower bounds and applications},
  Author                   = {Bar-Yossef, Ziv and Kumar, Ravi and Sivakumar, D},
  Booktitle                = {Proceedings of the thirty-third annual ACM symposium on Theory of computing},
  Year                     = {2001},
  Organization             = {ACM},
  Pages                    = {266--275}
}

@Article{Birge83,
  Title                    = {Approximation dans les espaces m{\'e}triques et th{\'e}orie de l'estimation},
  Author                   = {Birg{\'e}, L.},
  Journal                  = {Z. f{\"u}r Wahrscheinlichkeitstheorie und Verw. Geb.},
  Year                     = {1983},
  Number                   = {2},
  Pages                    = {181--237},
  Volume                   = {65}
}

@inproceedings{Bresler15,
  title={Efficiently learning Ising models on arbitrary graphs},
  author={Bresler, Guy},
  booktitle={Proceedings of the forty-seventh annual ACM symposium on Theory of computing},
  pages={771--782},
  year={2015},
  organization={ACM}
}

@Book{combinatorics,
  Title                    = {Introductory Combinatorics},
  Author                   = {Richard A. Brualdi},
  Publisher                = {Prentice Hall},
  Year                     = {1999},

  Address                  = {Upper Saddle River, NJ},
  Edition                  = {3rd}
}

@Article{subgaussian,
  Title                    = {{Sub-Gaussian random variables}},
  Author                   = {Buldygin, V. V and Kozachenko, Y. V.},
  Journal                  = {Ukrainian Mathematical Journal},
  Year                     = {1980},
  Number                   = {6},
  Pages                    = {483--489},
  Volume                   = {32},

  Publisher                = {Springer}
}

@InProceedings{CCMN00,
  Title                    = {Towards estimation error guarantees for distinct values},
  Author                   = {Charikar, Moses and Chaudhuri, Surajit and Motwani, Rajeev and Narasayya, Vivek},
  Booktitle                = {Proceedings of the nineteenth ACM SIGMOD-SIGACT-SIGART Symposium on Principles of Database Systems (PODS)},
  Year                     = {2000},
  Organization             = {ACM},
  Pages                    = {268--279}
}

@Article{CL68,
  Title                    = {Approximating discrete probability distributions with dependence trees},
  Author                   = {Chow, C.K. and Liu, C.N.},
  Journal                  = IEEE_J_IT,
  Year                     = {1968},
  Number                   = {3},
  Pages                    = {462--467},
  Volume                   = {14}
}

@Book{ckbook,
  Title                    = {Information Theory: Coding Theorems for Discrete Memoryless Systems},
  Author                   = {Imre Csisz{\'a}r and J{\'a}nos K\"{o}rner},
  Publisher                = {Academic Press, Inc.},
  Year                     = {1982}
}

@Book{DL93,
  Title                    = {Constructive approximation},
  Author                   = {DeVore, Ronald A. and Lorentz, George G.},
  Publisher                = {Springer},
  Year                     = {1993}
}

@Misc{OED,
  title                   = {Oxford {E}nglish Dictionary},
  HowPublished             = {\url{http://public.oed.com/about/}},
  Note                     = {Accessed: 2016-02-16}
}

@Article{Dobrushin58,
  Title                    = {A statistical problem arising in the theory of detection of signals in the presence of noise in a multi-channel system and leading to stable distribution laws},
  Author                   = {Dobrushin, R.L.},
  Journal                  = TPIA,
  Year                     = {1958},
  Number                   = {2},
  Pages                    = {161--173},
  Volume                   = {3}
}

@Book{DS08,
  Title                    = {Theory of uniform approximation of functions by polynomials},
  Author                   = {Dzyadyk, Vladislav K and Shevchuk, Igor A},
  Publisher                = {Walter de Gruyter},
  Year                     = {2008}
}

@InProceedings{GKB06,
  Title                    = {From the Entropy to the Statistical Structure of Spike Trains},
  Author                   = {Yun Gao and Kontoyiannis, I. and Bienenstock, E.},
  Booktitle                = {2006 IEEE International Symposium on Information Theory},
  Year                     = {2006},
  Month                    = {July},
  Pages                    = {645-649}
}

@Article{genovese.wasserman,
  Title                    = {{Rates of convergence for the Gaussian mixture sieve}},
  Author                   = {Genovese, C. R. and Wasserman, L.},
  Journal                  = {Annals of Statistics},
  Year                     = {2000},
  Number                   = {4},
  Pages                    = {1105--1127},
  Volume                   = {28}
}

@Article{probability.metrics,
  Title                    = {{On choosing and bounding probability metrics}},
  Author                   = {Gibbs, A. L. and Su, F. E.},
  Journal                  = {International Statistical Review},
  Year                     = {2002},
  Number                   = {3},
  Pages                    = {419--435},
  Volume                   = {70}
}

@Book{GR,
  Title                    = {{Table of Integrals Series and Products}},
  Author                   = {I. S. Gradshteyn and I. M. Ryzhik},
  Publisher                = {Academic},
  Year                     = {2007},

  Address                  = {New York, NY},
  Edition                  = {Seventh}
}

@inproceedings{HJW15,
  title={Does {D}irichlet prior smoothing solve the {S}hannon entropy estimation problem?},
  author={Han, Yanjun and Jiao, Jiantao and Weissman, Tsachy},
  booktitle={2015 IEEE International Symposium on Information Theory (ISIT)},
  pages={1367--1371},
  year={2015},
  organization={IEEE}
}

@inproceedings{HJW15adaptive,
  title={Adaptive estimation of {S}hannon entropy},
  author={Han, Yanjun and Jiao, Jiantao and Weissman, Tsachy},
  booktitle={2015 IEEE International Symposium on Information Theory (ISIT)},
  pages={1372--1376},
  year={2015},
  organization={IEEE}
}

@Book{IKbook,
  Title                    = {{Statistical Estimation: Asymptotic Theory}},
  Author                   = {Ibragimov, I.A. and Has'minskii, R.Z.},
  Publisher                = {Springer},
  Year                     = {1981}
}

@Article{jiao2013di,
  Title                    = {Universal estimation of directed information},
  Author                   = {Jiao, Jiantao and Permuter, Haim H and Zhao, Lei and Kim, Young-Han and Weissman, Tsachy},
  Journal                  = IEEE_J_IT,
  Year                     = {2013},
  Number                   = {10},
  Pages                    = {6220--6242},
  Volume                   = {59}
}

@article{JVHW2017,
  title={Maximum likelihood estimation of functionals of discrete distributions},
  author={Jiao, Jiantao and Venkat, Kartik and Han, Yanjun and Weissman, Tsachy},
  journal={IEEE Transactions on Information Theory},
  volume={63},
  number={10},
  pages={6774--6798},
  year={2017},
  publisher={IEEE}
}

@InProceedings{knudson2013spike,
  Title                    = {Spike train entropy-rate estimation using hierarchical {D}irichlet process priors},
  Author                   = {Knudson, Karin C and Pillow, Jonathan W},
  Booktitle                = {Proceedings of the Twenty-seventh Conference on Neural Information Processing Systems},
  Year                     = {2013},
  Pages                    = {2076--2084},
  publisher={Curran Associates, Inc.}
}

@Article{LeCam73,
  Title                    = {Convergence of Estimates Under Dimensionality Restrictions},
  Author                   = {Le Cam, L.},
  Journal                  = A_S,
  Year                     = {1973},
  Number                   = {1},
  Pages                    = {38 -- 53},
  Volume                   = {1}
}

@Article{Lo92,
  Title                    = {From the species problem to a general coverage problem via a new interpretation},
  Author                   = {Lo, Shaw-Hwa},
  Journal                  = A_S,
  Year                     = {1992},
  Number                   = {2},
  Pages                    = {1094--1109},
  Volume                   = {20}
}

@Book{MU06,
  Title                    = {Probability and computing: Randomized algorithms and probabilistic analysis},
  Author                   = {Mitzenmacher, Michael and Upfal, Eli},
  Publisher                = {Cambridge University Press},
  Year                     = {2005}
}

@Misc{GLM,
  Title                    = {Global Language Monitor. Number of words in the English language},
  HowPublished             = {\url{https://www.languagemonitor.com/global-english/no-of-words/}},
  Note                     = {Accessed: 2016-02-16}
}

@Article{NBS04,
  Title                    = {Entropy and information in neural spike trains: Progress on the sampling problem},
  Author                   = {Nemenman, Ilya and Bialek, William and de Ruyter van Steveninck, Rob R.},
  Journal                  = {Physical Review E},
  Year                     = {2004},
  Number                   = {5},
  Pages                    = {056111},
  Volume                   = {69}
}

@InProceedings{Nemirovski03,
  Title                    = {On tractable approximations of randomly perturbed convext constaints},
  Author                   = {A. Nemirovski},
  Booktitle                = {Proceedings of the 42nd IEEE Conference on Decision and Control},
  Year                     = {2003},
  Pages                    = {2419--2422},
  publisher={IEEE}
}

@Book{petrushev2011rational,
  Title                    = {Rational approximation of real functions},
  Author                   = {Petrushev, Penco Petrov and Popov, Vasil Atanasov},
  Publisher                = {Cambridge University Press},
  Year                     = {2011}
}

@Article{pinsker.minimax,
  Title                    = {{Optimal filtering of square-integrable signals in Gaussian noise}},
  Author                   = {Pinsker, M. S.},
  Journal                  = {Problemy Peredachi Informatsii},
  Year                     = {1980},
  Number                   = {2},
  Pages                    = {52--68},
  Volume                   = {16},

  Publisher                = {Russian Academie of Sciences}
}

@InCollection{PW96,
  Title                    = {An Entropy Estimator Algorithm and Telecommunications Applications},
  Author                   = {Plotkin, Nina T. and Wyner, Abraham J.},
  Booktitle                = {Maximum Entropy and Bayesian Methods},
  Publisher                = {Springer Netherlands},
  Year                     = {1996},
  Pages                    = {351-363},
  Series                   = {Fundamental Theories of Physics},
  Volume                   = {62}
}

@Article{porta2001entropy,
  Title                    = {Entropy, entropy rate, and pattern classification as tools to typify complexity in short heart period variability series},
  Author                   = {Porta, A. and Guzzetti, S. and Montano, N. and Furlan, R. and Pagani, M. and Malliani, A. and Cerutti, S.},
  Journal                  = {IEEE Transactions on Biomedical Engineering},
  Year                     = {2001},
  Number                   = {11},
  Pages                    = {1282--1291},
  Volume                   = {48}
}

@Article{QKC13,
  Title                    = {Efficient methods to compute optimal tree approximations of directed information graphs},
  Author                   = {Quinn, Christopher J and Kiyavash, Negar and Coleman, Todd P},
  Journal                  = IEEE_J_SP,
  Year                     = {2013},
  Number                   = {12},
  Pages                    = {3173--3182},
  Volume                   = {61}
}

@Book{spikes-book,
  Title                    = {Spikes: Exploring the Neural Code},
  Author                   = {Rieke, Fred and Bialek, William and Warland, David and van Steveninck, Rob de Ruyter},
  Publisher                = {The MIT Press},
  Year                     = {1999}
}

@Article{shannon,
  Title                    = {A Mathematical Theory of Communication},
  Author                   = {Claude E. Shannon},
  Journal                  = {Bell System Technical Journal},
  Year                     = {1948},
  Pages                    = {379 -- 423, 623 -- 656},
  Volume                   = {27}
}

@Article{steele86,
  Title                    = {An {E}fron-{S}tein inequality for nonsymmetric statistics},
  Author                   = {Steele, J Michael},
  Journal                  = {The Annals of Statistics},
  Year                     = {1986},
  Pages                    = {753--758},
  volume={14},
  issue={2},
  Publisher                = {JSTOR}
}

@Book{stoer.2002,
  Title                    = {{Introduction to Numerical Analysis}},
  Author                   = {Stoer, J. and Bulirsch, R.},
  Publisher                = {Springer-Verlag},
  Year                     = {2002},
  Address                  = {New York, NY},
  Edition                  = {3rd}
}

@Book{Strasser85,
  Title                    = {Mathematical theory of statistics: {S}tatistical experiments and asymptotic decision theory},
  Author                   = {Strasser, Helmut},
  Publisher                = {Walter de Gruyter},
  Year                     = {1985},

  Address                  = {Berlin, Germany}
}

@Article{SKSB98,
  Title                    = {Entropy and Information in Neural Spike Trains},
  Author                   = {Strong, S. P. and Koberle, Roland and de Ruyter van Steveninck, Rob R. and Bialek, William},
  Journal                  = {Phys. Rev. Lett.},
  Year                     = {1998},

  Month                    = {Jan.},
  Pages                    = {197--200},
  Volume                   = {80},

  Issue                    = {1}
}

@Book{orthogonal.poly,
  Title                    = {{Orthogonal polynomials}},
  Author                   = {Szeg{\"o}, G.},
  Publisher                = {American Mathematical Society},
  Year                     = {1975},

  Address                  = {Providence, RI},
  Edition                  = {4th}
}

@InProceedings{VV11-focs,
  Title                    = {The power of linear estimators},
  Author                   = {Valiant, Gregory and Valiant, Paul},
  Booktitle                = {Foundations of Computer Science (FOCS), 2011 IEEE 52nd Annual Symposium on},
  Year                     = {2011},
  Organization             = {IEEE},
  Pages                    = {403--412}
}

@article{Valiant11,
  title={Testing symmetric properties of distributions},
  author={Valiant, Paul},
  journal={SIAM Journal on Computing},
  volume={40},
  number={6},
  pages={1927--1968},
  year={2011},
  publisher={SIAM}
}

@Book{optimal.transport.old.new,
  Title                    = {{Optimal Transport: Old and New}},
  Author                   = {Villani, C.},
  Publisher                = {Springer Verlag},
  Year                     = {2008},

  Address                  = {Berlin}
}

@Book{villani.topics,
  Title                    = {{Topics in optimal transportation}},
  Author                   = {Villani, C.},
  Publisher                = {American Mathematical Society},
  Year                     = {2003},

  Address                  = {Providence, RI}
}

@article{WY15,
  title={Chebyshev polynomials, moment matching, and optimal estimation of the unseen},
  author={Wu, Yihong and Yang, Pengkun},
  journal={The Annals of Statistics},
  volume={47},
  number={2},
  pages={857--883},
  year={2019},
  publisher={Institute of Mathematical Statistics}
}

@Article{WY18,
  Title                    = {Optimal estimation of {G}aussian mixtures via denoised method of moments},
  Author                   = {Wu, Yihong and Yang, Pengkun},
  journal = {Annals of Statistics},
  volume ={48},
  issue={4},
  pages={1981–-2007},
  Year                     = {2020}
}

@MastersThesis{yang-msthesis,
  Title                    = {Optimal property estimation on large alphabets:
fundamental limits and fast algorithms},
  Author                   = {Pengkun Yang},
  School                   = {University of Illinois at Urbana-Champaign},
  Year                     = {2016}
}

@Article{YB99,
  Title                    = {Information-theoretic determination of minimax rates of convergence},
  Author                   = {Yang, Y. and Barron, A. R.},
  Journal                  = A_S,
  Year                     = {1999},
  Number                   = {5},
  Pages                    = {1564--1599},
  Volume                   = {27}
}

@incollection{Yu97,
  Title                    = {Assouad, {F}ano, and {L}e {C}am},
  Author                   = {Bin Yu},
  booktitle                  = {Festschrift for Lucien Le Cam},
  Year                     = {1997},
  Pages                    = {423--435},

  Publisher                = {Springer, New York}
}

@inproceedings{HOT88,
  title={Statistical estimators for relational algebra expressions},
  author={Hou, Wen-Chi and Ozsoyoglu, Gultekin and Taneja, Baldeo K},
  booktitle={Proceedings of the seventh ACM SIGACT-SIGMOD-SIGART symposium on Principles of database systems},
  pages={276--287},
  year={1988},
  organization={ACM}
}

@inproceedings{NS90,
  title={On estimating the size of projections},
  author={Naughton, Jeffrey F and Seshadri, S},
  booktitle={International Conference on Database Theory},
  pages={499--513},
  year={1990},
  organization={Springer}
}

@article{Esty86,
  title={Estimation of the size of a coinage: A survey and comparison of methods},
  author={Esty, Warren W},
  journal={The Numismatic Chronicle (1966-)},
  volume={146},
  pages={185--215},
  year={1986},
  publisher={JSTOR}
}

@article{GC11,
  title={Estimating species richness},
  author={Gotelli, Nicholas J and Colwell, Robert K},
  journal={Biological diversity: frontiers in measurement and assessment},
  volume={12},
  pages={39--54},
  year={2011},
  publisher={Oxford University Press Oxford}
}

@Misc{Valiant-email,
  Author                   = {G. Valiant},
  HowPublished             = {Private communication},
  Month                    = {Mar.},
  Year                     = {2017}
}

@article{WY2016sample,
  title={Sample complexity of the distinct elements problem},
  author={Wu, Yihong and Yang, Pengkun},
  journal={Mathematical Statistics and Learning},
  volume={1},
  number={1},
  pages={37--72},
  year={2018}
}

@book{Akhiezer1965,
  title={The classical moment problem: and some related questions in analysis},
  author={Akhiezer, N. I.},
  volume={5},
  year={1965},
  publisher={Oliver \& Boyd}
}

@article{CF1991,
  title={Recursiveness, positivity, and truncated moment problems},
  author={Curto, Ra{\'u}l E and Fialkow, Lawrence A},
  journal={Houston Journal of Mathematics},
  volume={17},
  number={4},
  pages={603--635},
  year={1991},
  publisher={UNIV HOUSTON DEPT MATH, HOUSTON, TX 77204}
}

@article{GW1969,
  title={Calculation of {G}auss quadrature rules},
  author={Golub, Gene H and Welsch, John H},
  journal={Mathematics of computation},
  volume={23},
  number={106},
  pages={221--230},
  year={1969}
}

@inproceedings{WV2010,
  title={The impact of constellation cardinality on {G}aussian channel capacity},
  author={Wu, Yihong and Verd{\'u}, Sergio},
  booktitle={Communication, Control, and Computing (Allerton), 2010 48th Annual Allerton Conference on},
  pages={620--628},
  year={2010},
  organization={IEEE}
}

@book{Jorge2011,
  title={Algebraic Approximation: A Guide to Past and Current Solutions},
  author={Bustamante, Jorge},
  year={2011},
  publisher={Springer Science \& Business Media}
}

@book{Kosorok,
  title={Introduction to empirical processes and semiparametric inference},
  author={Kosorok, Michael R},
  year={2007},
  publisher={Springer Science \& Business Media}
}

@article{balakrishnan2017statistical,
  title={Statistical guarantees for the EM algorithm: From population to sample-based analysis},
  author={Balakrishnan, Sivaraman and Wainwright, Martin J and Yu, Bin},
  journal={The Annals of Statistics},
  volume={45},
  number={1},
  pages={77--120},
  year={2017},
  publisher={Institute of Mathematical Statistics}
}

@article{lu2016statistical,
  title={Statistical and Computational Guarantees of {L}loyd's Algorithm and its Variants},
  author={Lu, Yu and Zhou, Harrison H},
  journal={arXiv preprint arXiv:1612.02099},
  year={2016}
}

@article{DK68,
  title={Construction of sequences estimating the mixing distribution},
  author={Deely, JJ and Kruse, RL},
  journal={The Annals of Mathematical Statistics},
  volume={39},
  number={1},
  pages={286--288},
  year={1968},
  publisher={JSTOR}
}

@article{meng1997algorithm,
  title={The {EM} algorithm—an old folk-song sung to a fast new tune},
  author={Meng, Xiao-Li and Van Dyk, David},
  journal={Journal of the Royal Statistical Society: Series B (Statistical Methodology)},
  volume={59},
  number={3},
  pages={511--567},
  year={1997},
  publisher={Wiley Online Library}
}

@article{koenker2014convex,
  title={Convex optimization, shape constraints, compound decisions, and empirical {B}ayes rules},
  author={Koenker, Roger and Mizera, Ivan},
  journal={Journal of the American Statistical Association},
  volume={109},
  number={506},
  pages={674--685},
  year={2014},
  publisher={Taylor \& Francis}
}

@article{kim2014minimax,
  title={Minimax bounds for estimation of normal mixtures},
  author={Kim, Arlene KH},
  journal={Bernoulli},
  volume={20},
  number={4},
  pages={1802--1818},
  year={2014},
  publisher={Bernoulli Society for Mathematical Statistics and Probability}
}

@article{ibragimov2001estimation,
  title={Estimation of analytic functions},
  author={Ibragimov, I},
  journal={Lecture Notes-Monograph Series},
  volume={36},
  pages={359--383},
  year={2001},
  publisher={Institute of Mathematical Statistics}
}

@inproceedings{li2017robust,
  title={Robust and proper learning for mixtures of {G}aussians via systems of polynomial inequalities},
  author={Li, Jerry and Schmidt, Ludwig},
  booktitle={Proceedings of the 30th Annual Conference on Learning Theory (COLT 2017)},
  pages={1302--1382},
  year={2017},
  publisher={PMLR}
}

@inproceedings{hopkins2018mixture,
  title={Mixture models, robustness, and sum of squares proofs},
  author={Hopkins, Samuel B and Li, Jerry},
  booktitle={Proceedings of the 50th Annual ACM SIGACT Symposium on Theory of Computing},
  pages={1021--1034},
  year={2018},
  organization={ACM}
}

@article{edelman1988estimation,
  title={Estimation of the mixing distribution for a normal mean with applications to the compound decision problem},
  author={Edelman, David},
  journal={The Annals of Statistics},
  volume={16},
  number={4},
  pages={1609--1622},
  year={1988},
  publisher={Institute of Mathematical Statistics}
}

@book{rivlin2003introduction,
  title={An introduction to the approximation of functions},
  author={Rivlin, Theodore J},
  year={1981},
  publisher={Dover}
}

@article{HJWW17,
  Author = {Yanjun Han and Jiantao Jiao and Tsachy Weissman and Yihong Wu},
  Title = {Optimal rates of entropy estimation over {L}ipschitz balls},
  journal = {Annals of Statistics},
  volume={48},
  issue={6},
  pages={3228 -- 3250},
  Year = {2020}
}

@inproceedings{HJLWWY17,
  Author = {Yanjun Han and Jiantao Jiao and Chuan-Zheng Lee and Tsachy Weissman and Yihong Wu and Tiancheng Yu},
  Title = {Entropy Rate Estimation for {M}arkov Chains with Large State Space},
  Year = {2018},
  booktitle = {Proceedings of the Thirty-second Conference on Neural Information Processing Systems.},
  publisher={Curran Associates, Inc.},
pages = {9781--9792}
}

@book{Freud,
  title={Orthogonal polynomials},
  author={Freud, G{\'e}za},
  year={1971},
  publisher={Pergamon}
}

@article{Zhang_2009,
  Author = {Cun-Hui Zhang},
  Journal = {Statistica Sinica},
  Pages = {1297-1318},
  Title = {Generalized maximum likelihood estimation of normal mixture densities},
  Volume = {19},
  Year = {2009}}

@article{Saha_2017,
  title={On the nonparametric maximum likelihood estimator for gaussian location mixture densities with application to gaussian denoising},
  author={Saha, Sujayam and Guntuboyina, Adityanand and others},
  journal={Annals of Statistics},
  volume={48},
  number={2},
  pages={738--762},
  year={2020},
  publisher={Institute of Mathematical Statistics}
}

@String { A_S       = {The Annals of Statistics} }

@String { Allerton  = {Proc.\ Allerton Conf.\ Commun., Control, and Computing} }

@String { IEEE_J_IT = {{IEEE} Trans. Inf. Theory} }

@String { IEEE_J_SP = {{IEEE} Trans. Signal Process.} }

@String { ISIT      = {Proc.\ IEEE Int.\ Symp.\ Inform. Theory} }

@String { TPIA      = {Theory of Probability \& Its Applications} }

\end{document}